%% file: er_localization.tex
\documentclass[11pt]{aart}

\usepackage[letterpaper, hmargin=1in, top=1.2in, bottom=1.2in, footskip=0.6in]{geometry}

\usepackage{titlesec}
\titleformat{\section}[block]{\filcenter\normalfont\bfseries\large}{\thesection.}{.5em}{}\titlespacing*{\section}{0pt}{2\baselineskip}{1\baselineskip}
\titleformat{\subsection}[runin]{\normalfont\bfseries}{\thesubsection.}{.4em}{}[.]\titlespacing{\subsection}{0pt}{2ex plus .1ex minus .2ex}{.8em}
\titleformat{\subsubsection}[runin]{\normalfont\itshape}{\thesubsubsection.}{.3em}{}[.]\titlespacing{\subsubsection}{0pt}{1ex plus .1ex minus .2ex}{.5em}
\titleformat{\paragraph}[runin]{\normalfont\itshape}{\theparagraph.}{.3em}{}[.]\titlespacing{\paragraph}{0pt}{1ex plus .1ex minus .2ex}{.5em}

\usepackage[T1]{fontenc}
\usepackage[utf8]{inputenc}
\usepackage{lmodern}

\usepackage[labelfont=bf,font=small,labelsep=period]{caption}
\usepackage{microtype}

\let\originalleft\left
\let\originalright\right
\renewcommand{\left}{\mathopen{}\mathclose\bgroup\originalleft}
\renewcommand{\right}{\aftergroup\egroup\originalright}

\usepackage{amsmath}
\usepackage{amssymb}
\usepackage{amsfonts}
\usepackage{latexsym}
\usepackage{amsthm}
\usepackage{amsxtra}
\usepackage{amscd}
\usepackage{bbm}
\usepackage{mathrsfs}
\usepackage{bm}
\usepackage{mathtools}

\usepackage{graphicx, color}

\definecolor{darkred}{rgb}{0.9,0,0.3}
\definecolor{darkblue}{rgb}{0,0.3,0.9}
\definecolor{purple}{rgb}{0.6,0,0.7}

\usepackage{booktabs}
\usepackage[nottoc,notlof,notlot]{tocbibind}
\usepackage{cite} 

\numberwithin{equation}{section}
\numberwithin{figure}{section}

\usepackage{enumitem}

\definecolor{vdarkred}{rgb}{0.6,0,0.2}
\definecolor{vdarkblue}{rgb}{0,0.2,0.6}
\usepackage[pdftex, colorlinks, linkcolor=vdarkblue,citecolor=vdarkred]{hyperref}

\usepackage{cleveref} 

\theoremstyle{plain} 
\newtheorem{theorem}{Theorem}[section]
\newtheorem*{theorem*}{Theorem}
\newtheorem{lemma}[theorem]{Lemma}
\newtheorem{claim}[theorem]{Claim}
\newtheorem*{lemma*}{Lemma}
\newtheorem{corollary}[theorem]{Corollary}
\newtheorem*{corollary*}{Corollary}
\newtheorem{proposition}[theorem]{Proposition}
\newtheorem*{proposition*}{Proposition}

\newtheorem*{conjecture*}{Conjecture}

\theoremstyle{definition} 
\newtheorem{definition}[theorem]{Definition}
\newtheorem*{definition*}{Definition}

\newtheorem*{example*}{Example}
\newtheorem{remark}[theorem]{Remark}
\newtheorem*{remark*}{Remark}

\newtheorem*{assumption*}{Assumption}

\newtheorem*{convention*}{Convention}

\newcommand{\f}{\mathbf} 
\renewcommand{\r}{\mathrm}  
\newcommand{\bb}{\mathbb} 
\renewcommand{\cal}{\mathcal} 
 
\newcommand{\fra}{\mathfrak} 
\newcommand{\ol}[1]{\overline{#1} \!\,} 
\newcommand{\wh}{\widehat}
\newcommand{\wt}{\widetilde}
\newcommand{\op}{\operatorname}

\renewcommand{\P}{\mathbb{P}}
\newcommand{\E}{\mathbb{E}}
\newcommand{\R}{\mathbb{R}}
\newcommand{\C}{\mathbb{C}}
\newcommand{\N}{\mathbb{N}}
\newcommand{\Z}{\mathbb{Z}}

\newcommand{\ee}{\mathrm{e}}
\newcommand{\ii}{\mathrm{i}}
\newcommand{\dd}{\mathrm{d}}
\newcommand{\col}{\mathrel{\vcenter{\baselineskip0.75ex \lineskiplimit0pt \hbox{.}\hbox{.}}}}
\newcommand*{\deq}{\mathrel{\vcenter{\baselineskip0.5ex \lineskiplimit0pt\hbox{\scriptsize.}\hbox{\scriptsize.}}}=}
\newcommand*{\eqd}{=\mathrel{\vcenter{\baselineskip0.5ex \lineskiplimit0pt\hbox{\scriptsize.}\hbox{\scriptsize.}}}}
\newcommand{\eqdist}{\overset{\r d}{=}}

\renewcommand{\leq}{\leqslant}
\renewcommand{\geq}{\geqslant}
\renewcommand{\epsilon}{\varepsilon}

\newcommand{\floor}[1] {\lfloor #1 \rfloor}

\newcommand{\ind}[1]{\mathbbm 1_{#1}}

\newcommand{\pb}[1]{\bigl(#1\bigr)}
\newcommand{\pB}[1]{\Bigl(#1\Bigr)}
\newcommand{\pbb}[1]{\biggl(#1\biggr)}
\newcommand{\pBB}[1]{\Biggl(#1\Biggr)}

\newcommand{\qb}[1]{\bigl[#1\bigr]}
\newcommand{\qB}[1]{\Bigl[#1\Bigr]}
\newcommand{\qbb}[1]{\biggl[#1\biggr]}
\newcommand{\qBB}[1]{\Biggl[#1\Biggr]}

\newcommand{\h}[1]{\{#1\}}
\newcommand{\hb}[1]{\bigl\{#1\bigr\}}
\newcommand{\hB}[1]{\Bigl\{#1\Bigr\}}
\newcommand{\hbb}[1]{\biggl\{#1\biggr\}}

\newcommand{\abs}[1]{\lvert #1 \rvert}
\newcommand{\absb}[1]{\bigl\lvert #1 \bigr\rvert}
\newcommand{\absB}[1]{\Bigl\lvert #1 \Bigr\rvert}
\newcommand{\absbb}[1]{\biggl\lvert #1 \biggr\rvert}
\newcommand{\absBB}[1]{\Biggl\lvert #1 \Biggr\rvert}

\newcommand{\norm}[1]{\lVert #1 \rVert}
\newcommand{\normb}[1]{\bigl\lVert #1 \bigr\rVert}

\newcommand{\normbb}[1]{\biggl\lVert #1 \biggr\rVert}

\newcommand{\scalar}[2]{\langle#1 \mspace{2mu}, #2\rangle}

\newcommand{\condB}{\,\Big\vert\,}

\DeclareMathOperator{\diag}{diag}

\DeclareMathOperator{\supp}{supp}

\DeclareMathOperator{\im}{Im}

\DeclareMathOperator{\dist}{dist}
\DeclareMathOperator{\diam}{diam}
\DeclareMathOperator{\spec}{spec}
\DeclareMathOperator{\ran}{ran}
\newcommand{\Span}{\operatorname{span}}

\newcommand{\eps}{\varepsilon}

\newcommand*{\defeq}{\mathrel{\vcenter{\baselineskip0.5ex \lineskiplimit0pt\hbox{\scriptsize.}\hbox{\scriptsize.}}}=}

\usepackage{tocloft}
\renewcommand{\contentsname}{\hfill\bfseries\large Contents\hfill}

\title{Localized phase for the Erd\H{o}s-R\'enyi graph} 
\author{Johannes Alt \and Raphael Ducatez \and Antti Knowles}

\begin{document}

\maketitle

\begin{abstract}
We analyse the eigenvectors of the adjacency matrix of the Erd\H{o}s-R\'enyi graph $\bb G(N,d/N)$ for $\sqrt{\log N} \ll d \lesssim \log N$. We show the existence of a localized phase, where each eigenvector is exponentially localized around a single vertex of the graph. This complements the completely delocalized phase previously established in \cite{ADK20}. For large enough $d$, we establish a mobility edge by showing that the localized phase extends up to the boundary of the delocalized phase. We derive explicit asymptotics for the localization length up to the mobility edge and characterize its divergence near the phase boundary.

The proof is based on a rigorous verification of Mott's criterion for localization, comparing the tunnelling amplitude between localization centres with the eigenvalue spacing. The first main ingredient is a new family of global approximate eigenvectors, for which sharp enough estimates on the tunnelling amplitude can be established. The second main ingredient is a lower bound on the spacing of approximate eigenvalues. It follows from an anticoncentration result for these discrete random variables, obtained by a recursive application of a self-improving anticoncentration estimate due to Kesten.
\end{abstract} 

\bigskip
\bigskip

\tableofcontents

\bigskip

\section{Introduction}

\subsection{Overview} \label{sec:overview}
Let $A$ be the adjacency matrix of a graph with vertex set $[N] = \{1, \dots, N\}$. We are interested in the geometric structure of the eigenvectors of $A$, in particular their \emph{spatial localization}. An $\ell^2$-normalized eigenvector $\f w = (w_x)_{x \in [N]} \in \R^N$ gives rise to a probability measure $x \mapsto w_x^2$ on the set of vertices $[N]$. Informally, $\f w$ is \emph{delocalized} if its mass is approximately uniformly distributed throughout $[N]$, and \emph{localized} if its mass is essentially concentrated on a small number of vertices.

In this paper we study the spatial localization of eigenvectors for the Erd\H{o}s-R\'enyi graph $ \mathbb{G} \equiv \mathbb{G}(N,d/N)$. It is the simplest model of a random graph, where each edge of the complete graph on $N$ vertices is kept independently with probability $d/N$, with $0 \leq d \leq N$. Here, $d \equiv d_N$ is a parameter whose interpretation is the expected degree of a vertex.
It is well known that $\bb G$ undergoes a dramatic change in behaviour at the \emph{critical scale} $d \asymp \log N$, which is the scale at and below which the vertex degrees do not concentrate. For $d \gg \log N$, with high probability all degrees are approximately equal and the graph is \emph{homogeneous}. On the other hand, for $d \lesssim \log N$, the degrees do not concentrate and the graph becomes highly \emph{inhomogeneous}: it contains for instance hubs of large degree, leaves, and isolated vertices. As long as $d > 1$, the graph $\bb G$ has with high probability a unique giant component, and we shall always restrict our attention to it.

The Erd\H{o}s-R\'enyi graph $\bb G$ at and below criticality was proposed in \cite{ADK20} as a simple and natural model on which to address the question of spatial localization of eigenvectors. Its graph structure provides an intrinsic and nontrivial notion of distance, which allows for a study of the geometry of the eigenvectors. It can be interpreted as a model of quantum disorder, where the disorder arises from the random geometry of the graph. Moreover, its phase diagram turns out to be remarkably amenable to rigorous analysis.

In this paper we establish the existence of a fully \emph{localized phase} in a region of the phase diagram of $\bb G$ near the spectral edge.  This complements the completely \emph{delocalized phase} established in \cite{ADK20,ADK_delocalized, HeKnowlesMarcozzi2018, EKYY1}. Our results in both phases are quantitative with essentially optimal bounds.

As a consequence, for a range of critical densities $d \asymp \log N$, we establish a \emph{mobility edge} separating the localized and delocalized phases. We derive the explicit behaviour of the localization length on either side of the mobility edge. In particular, we show how the localization length diverges as one approaches the mobility edge from the localized phase (see Figure \ref{fig:ll} below). The Erd\H{o}s-R\'enyi graph at criticality is hence one of the very few models where a mobility edge can be rigorously established. Moreover, our proofs yield strong quantitative control of the localization length in the localized phase, as well as complete delocalization in the delocalized phase, all the way up to the mobility edge in both phases. To the best of our knowledge, this is the first time quantitative control is obtained in the vicinity of the mobility edge.

\begin{figure}[!ht]
\begin{center}
{\small 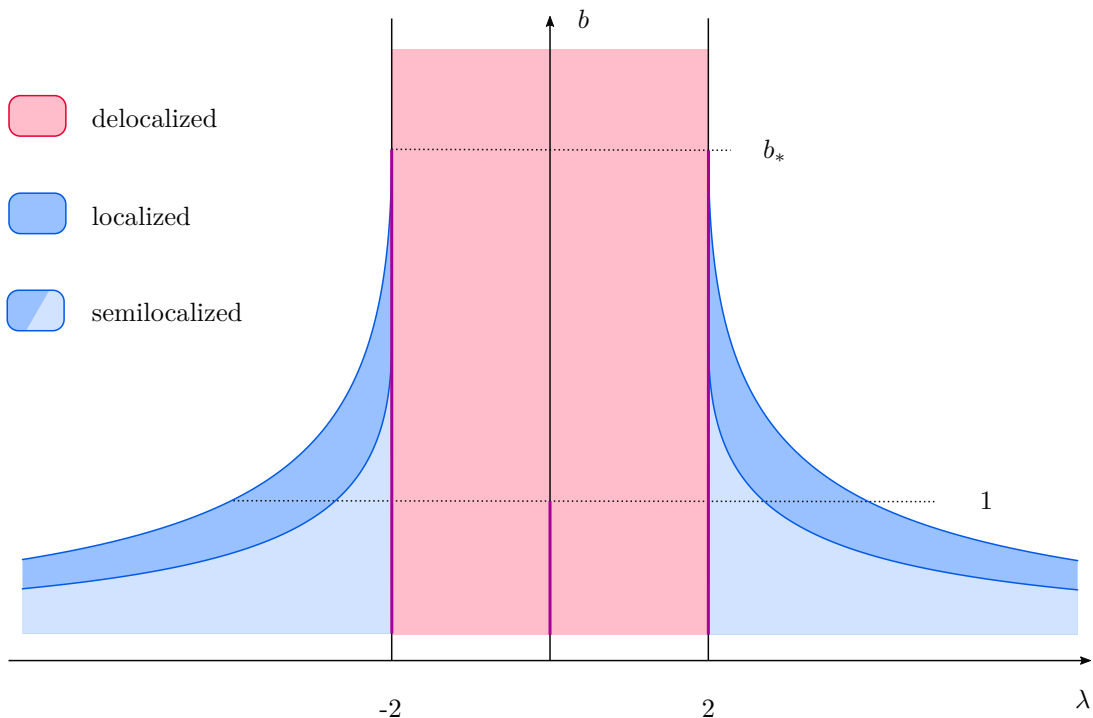}
\end{center}
\caption{The phase diagram of the rescaled adjacency matrix $A / \sqrt{d}$ of the (giant component of the) Erd\H{o}s-R\'enyi graph $\bb G(N,d/N)$ at criticality, where $d = b \log N$ with $b$ fixed. The horizontal axis records the location $\lambda$ in the spectrum and the vertical axis the sparseness parameter $b$. The spectrum is confined to the coloured region, which is split into the indicated phases. The thick purple lines correspond to phase boundaries, which are not covered by our results. (The phase boundary at energy $0$ for $b \leq 1$ is discussed in \cite{ADK20,ADK_delocalized}.) For large enough $b$, there is a mobility edge between the localized and the delocalized phases at energies $\pm 2$.
\label{fig:phases}}
\end{figure}

A graphical overview of the main result of this paper, and how it fits into the previous results of \cite{ADK20, ADK_delocalized}, is provided by the phase diagram of Figure \ref{fig:phases}. 
It depicts three phases, which are most conveniently characterized by  the $\ell^\infty$-norm $\norm{\f w}_\infty$ of an $\ell^2$-normalized eigenvector $\f w$. Clearly, $N^{-1} \leq \norm{\f w}_\infty^2 \leq 1$, where $\norm{\f w}_\infty^2 = 1$ corresponds to localization at a single vertex and $\norm{\f w}_\infty^2 = N^{-1}$ to perfect delocalization. We say that an eigenvalue $\lambda$ of the rescaled adjacency matrix $H \deq A / \sqrt{d}$ with eigenvector $\f w$ belongs to
\begin{enumerate}[label=(\roman*)]
\item \label{itm:phase1}
the \emph{localized phase} if $\norm{\f w}_\infty^2 \asymp 1$,
\item \label{itm:phase2}
the \emph{delocalized phase} if $\norm{\f w}_\infty^2 = N^{-1 + o(1)}$,
\item \label{itm:phase3}
the \emph{semilocalized phase} if $\norm{\f w}_\infty^2 \geq N^{-\gamma}$ for some constant $\gamma < 1$.
\end{enumerate}
In particular, the localized phase is a subphase of the semilocalized phase. The result of this paper is the existence of phase \ref{itm:phase1}, while phases \ref{itm:phase2} and \ref{itm:phase3} were previously established in \cite{ADK20, ADK_delocalized}.

We now briefly describe the structure of the phase diagram in Figure \ref{fig:phases}. It is well known \cite{Wigner1958} that, as long as $d \gg 1$, the global eigenvalue density of $H$ converges to the semicircle law supported in $[-2,2]$. We write $d = b \log N$ for some constant\footnote{Our results hold also for $d \ll \log N$, i.e.\ below the critical scale, but in this overview we suppose for simplicity that $b$ is a constant.} $b > 0$. The localized and semilocalized phases exist only if $b < b_*$, where
\begin{equation} \label{def_b_star}
b_* \deq \frac{1}{2\log 2 - 1} \approx 2.59\,.
\end{equation}
For fixed $b < b_*$, the spectrum splits into two disjoints parts, the delocalized phase $(-2,0) \cup (0,2)$ and the semilocalized phase $(-\lambda_{\max}, -2) \cup (2, \lambda_{\max})$, where $\lambda_{\max} > 2$ is an explicit function of $b$ (see \eqref{def_lambda_max} below). A region $(-\lambda_{\max}, -\lambda_{\r{loc}}) \cup (\lambda_{\r{loc}}, \lambda_{\max})$ near the spectral edges is in fact fully localized, where $2 \leq \lambda_{\r{loc}} < \lambda_{\max}$. In particular, for large enough $b$, the semilocalized phase consists entirely of the localized phase, i.e.\ $\lambda_{\r{loc}} = 2$. For smaller values of $b$, the diagram in Figure \ref{fig:phases} does not rule out the possibility of an eigenvector $\f w$ in the semilocalized phase satisfying $\norm{\f w}_\infty^2 = N^{-\gamma + o(1)}$ for some constant $\gamma \in (0,1)$. This latter scenario corresponds to eigenvectors that are neither fully localized nor fully delocalized\footnote{In the literature, this phenomenon is sometimes referred to as nonergodic delocalization.}, where $\gamma \in (0,1)$ plays the role of an anomalous fractal dimension. In fact, it is plausible that such fractal eigenvectors occur in the semilocalized phase for small enough $b$; for more details, we refer to \cite{tarzia2022} and the heuristic discussion on localization phenomena later in this subsection, as well as Appendix \ref{sec:optimal_loc} below.

The localization-delocalization transition for $\bb G$ described above is an example of an \emph{Anderson transition}, where a disordered quantum system exhibits localized or delocalized states depending on the disorder strength and the location in the spectrum, corresponding to an insulator or conductor, respectively. Originally proposed in the 1950s \cite{anderson1958absence} to model conduction in semiconductors with random impurities, this phenomenon is now recognized as a general feature of wave transport in disordered media, and is one of the most influential ideas in modern condensed matter physics \cite{lee1985disordered, evers2008anderson, lagendijk2009fifty, abrahams201050}. It is expected to occur in great generality whenever linear waves, such as quantum particles, propagate through a disordered medium. For weak enough disorder, the stationary states are expected to be delocalized, while a strong enough disorder can give rise to localized states.

The general heuristic behind localization is the following. A disordered quantum system is characterized by its Hamiltonian $H$, a large Hermitian random matrix. The disorder inherent in $H$ gives rise to spatial regions where the environment is in some sense exceptional, such as vertices of unusually large degree for the Erd\H{o}s-R\'enyi graph\footnote{For the Anderson model discussed below, these regions arise from an unusually large mean random potential. We note that, just like the large degree vertices for $\bb G$, these regions arise as a result of a collective phenomenon, requiring the conspiracy of a large number of independent random variables to yield an exceptional local configuration.}. These regions are possible localization centres, around which localized states may form. Whether they do so is captured by the following well-known rule of thumb, also known as Mott's criterion: localization occurs whenever \emph{the eigenvalue spacing is much larger than the tunnelling amplitude between localization centres}. The simplest illustration of this rule is for a two-state system whose Hamiltonian is the matrix
$\pb{
\begin{smallmatrix}
a & \tau
\\
\tau & b
\end{smallmatrix}
}.
$ Setting the tunnelling amplitude $\tau$ between the two sites to zero, we have an eigenvalue spacing $\abs{a - b}$. Denoting by $\f e_1, \f e_2$ the standard basis vectors of $\R^2$, we find that if $\abs{\tau} \ll \abs{a - b}$ then the eigenvectors are approximately $\f e_1, \f e_2$, corresponding to localization, and if $\abs{\tau} \gg \abs{a - b}$ then the eigenvectors are approximately $\frac{1}{\sqrt{2}} (\f e_1 \pm \f e_2)$, corresponding to delocalization.

\begin{figure}[!ht]
\begin{center}
{\small 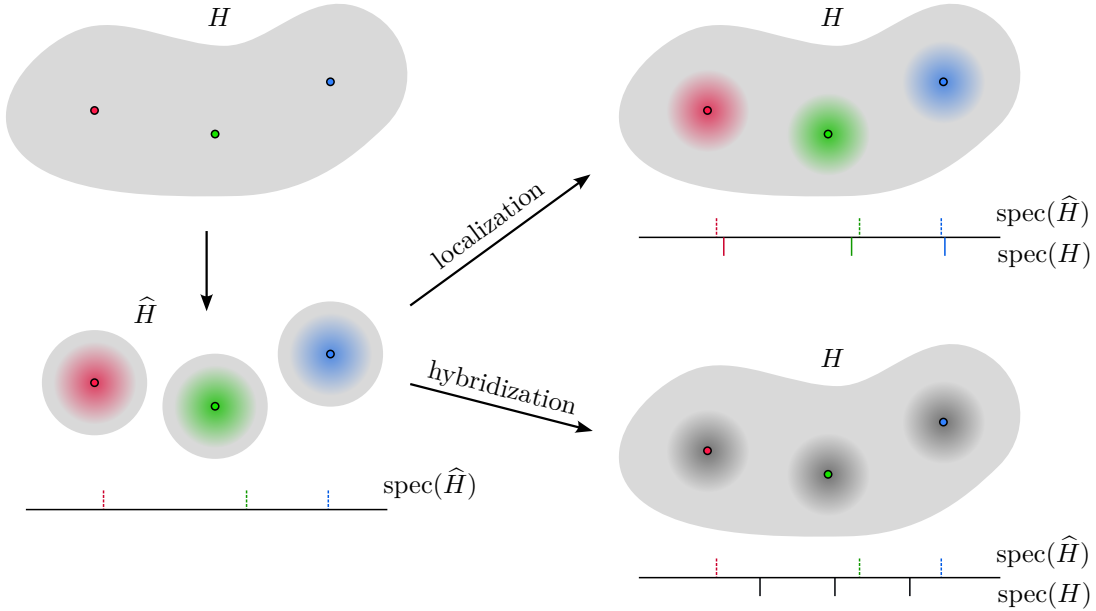}
\end{center}
\caption{A schematic illustration of the localization-hybridization dichotomy. A disordered Hamiltonian $H$ defined on a graph has three localization centres whose energies are in resonance. They are indicated in red, green, and blue. The Hamiltonian $\wh H$ is obtained from $H$ by splitting the graph into disconnected components around each localization centre. This gives rise to three eigenvalues of $\wh H$ associated with the three components. Each component associated with a centre $x$ carries a localized state $\f w(x)$ drawn as a decaying density in the colour corresponding to that of the centre. The spectrum of $\wh H$ is drawn above the real line in dotted lines that match the colour of the associated centre. Upon putting back the edges of the graph to return to $H$, we can have either localization or hybridization, depending on Mott's criterion. In each case, we draw the spectrum of $H$ below the real line. In the case of localization, the eigenvectors of $H$ remain close to $\f w(x)$ and the eigenvalues are shifted by an amount that is small compared to the eigenvalue spacing. In the case of hybridization, the eigenvectors of $H$ are delocalized over the three centres, being approximately nontrivial linear combinations of all three vectors $\f w(x)$.
\label{fig:localization}}
\end{figure}

More generally, for a disordered Hamiltonian $H$ defined on some connected graph, a simple yet instructive way to think of the rule of thumb is to suppose that the localization centres are spatially separated, and to construct the Hamiltonian $\wh H$ from $H$ by removing edges from the graph so as to disconnect the localization centres from each other. Hence, $\wh H$ is defined on a union of connected components, each of which is associated with a single localization centre. Denote by $\cal W$ the set of localization centres, which are vertices in the underlying graph. A component associated with a localization centre $x \in \cal W$ trivially gives rise to a localized state $\f w(x)$ of $\wh H$ with eigenvalue $\theta(x)$, which has the interpretation of the energy of the localization centre $x$. Upon putting back the removed edges, one can imagine two scenarios:
\begin{enumerate}[label=(\roman*)]
\item
\emph{localization}, where the eigenvectors of $H$ remain close to the eigenvectors $\f w(x)$ of $\wh H$;
\item
\emph{hybridization}, where eigenvectors $\f w(x)$ associated with several resonant localization centres $x$ with similar energies $\theta(x)$ are superimposed to form eigenvectors of $H$.
\end{enumerate}
We refer to Figure \ref{fig:localization} for an illustration of this dichotomy.

To use Mott's criterion, we note that one possible way of quantifying the tunnelling amplitude is
\begin{equation} \label{def_tau}
\tau \deq \max_{x \in \cal W} \norm{(H - \theta(x)) \f w(x)}\,.
\end{equation}
Indeed, this expression clearly generalizes the tunnelling amplitude for the two-level system given above, expressing in general the error arising from putting back the edges of $H$ missing from $\wh H$, with respect to the family $(\f w(x))_{x \in \cal W}$. It measures the extent to which the localized eigenvectors $\f w(x)$ of $\wh H$ are approximate eigenvectors of $H$ (sometimes also called quasi-modes of $H$). Under some assumptions on the localization components, $\tau$ also controls the off-diagonal part in a block-diagonal representation of $H$ in the basis $(\f w(x) \col x \in \cal W)$; see Section \ref{sec:sketch} below. Mott's criterion states that localization occurs whenever $\tau \ll \Delta$, where
\begin{equation} \label{def_Delta}
\Delta \deq \min_{x \neq y \in \cal W} \abs{\theta(x) - \theta(y)}
\end{equation}
is the eigenvalue spacing. Otherwise, we expect hybridization. The main difficulty in establishing localization is therefore to control \emph{resonances}, i.e.\ pairs of vertices $x,y$ such that $\abs{\theta(x) - \theta(y)}$ is small,
and to rule out hybridization. Although this picture is helpful in gaining intuition about localization, in many instances, such as in the proof for the Erd\H{o}s-R\'enyi graph in this paper\footnote{Or for the Anderson model e.g.\ in \cite{FroSpe}.}, it is but a rough caricature and the true picture is considerably more subtle (as we explain later in this subsection and in more detail in Section \ref{sec:sketch}).

The semilocalized phase of the Erd\H{o}s-R\'enyi graph from \cite{ADK20} is a phase where the eigenvectors are concentrated around a small number of localization centres, but where hybridization cannot be ruled out. As shown in \cite{ADK20}, the set of all possible localization centres of $\bb G$ corresponds to vertices $x \in [N]$ whose \emph{normalized degree}
\[ \alpha_x \deq \sum_{y \in [N]} H_{xy}^2 = \frac{1}{d} \sum_{y \in [N]} A_{xy}
\]  
is greater than $2$. The energy $\theta(x)$ of a localization centre $x$ is approximately equal to $\Lambda(\alpha_x)$, where we introduced the function $\Lambda \colon [2,\infty) \to [2,\infty)$ through
\begin{equation} \label{eq:defLambda}
\Lambda(\alpha) \deq \frac{\alpha}{\sqrt{\alpha-1}}\,.
\end{equation}
As shown in \cite{ADK19} (see also \cite{MR4234995}), there is a one-to-one correspondence between eigenvalues $\lambda > 2$ in the semilocalized phase and vertices $x$ of normalized degree $\alpha_x \geq 2$ given by $\lambda = \Lambda(\alpha_x) + o(1)$.
An eigenvalue $\lambda$ in the semilocalized phase has an eigenvector that is almost entirely concentrated in small balls around the set of vertices $\cal W_\lambda \deq \{x \in [N] \col \Lambda(\alpha_x) = \lambda + o(1)\}$ in resonance with the energy $\lambda$ \cite{ADK20}. The size of the set of resonant vertices $\cal W_\lambda$ is comparable to the density of states, equal to $N^{\rho_b(\lambda) + o(1)}$ for an explicit exponent $\rho_b(\lambda) < 1$ given in \eqref{lambda_max_rho} in Appendix \ref{sec:alpha_qual} below. Hence, owing to the small size of the set of resonant vertices $\cal W_\lambda$, the semilocalized phase is sharply distinguished from the delocalized phase. However,  the key issue of controlling resonances and ruling out hybridization is not addressed in \cite{ADK20} (see Figure \ref{fig:localization}).

In this paper we prove localization by ruling out hybridization among the resonant vertices. Our result holds for the largest and smallest $N^{\mu}$ eigenvalues for $\mu < \frac{1}{24}$. For small enough $\mu$, our obtained rate of exponential decay is optimal all the way up to radii of the order of the diameter of the graph. The bound $\frac{1}{24}$ is not optimal, and we expect that by a refinement of the method developed in this paper it can be improved; for the sake of keeping the argument reasonably simple, we refrain from doing so here. Heuristic arguments suggest that the optimal upper bound for $\mu$ is $\frac{1}{4}$; see Appendix \ref{sec:optimal_loc} below as well as \cite{tarzia2022}.

At this point it is helpful to review the previous works \cite{ADK20, ADK21}, which addressed the tunnelling amplitude \eqref{def_tau} in the above simple picture of localization based on disjoint neighbourhoods of localization centres. In \cite{ADK20}, the estimate $\tau \lesssim d^{-1/2}$ was established in the entire semilocalized phase, while in \cite{ADK21} it was improved to $\tau \lesssim d^{-3/2}$ at the spectral edge. In fact, here we argue that the best possible bound on $\tau$ in terms of the local approximate eigenvectors $\f w(x)$ introduced above is $\ee^{-c \frac{\log N}{\log d}}$ for some constant $c > 0$. To see this, we recall from \cite{ADK20} that the vector $\f w(x)$ is exponentially decaying around $x$ at some fixed rate $C > 0$ depending on $b$. Thus, the best possible estimate for $\tau$ arising from the exponential decay is $\ee^{-C r}$, where $r \asymp \diam(\bb G)  = \frac{\log N}{\log d} (1+ o(1))$ (see \cite{ChungLu_diameter}).

As for the eigenvalue spacing \eqref{def_Delta}, it was not addressed at all in \cite{ADK20}. In \cite{ADK21}, it was estimated at the spectral edge as $\Delta \geq d^{-1 - \epsilon}$ with high probability for any constant $\epsilon > 0$. Combined with the bound $\tau \lesssim d^{-3/2}$ at the spectral edge obtained in \cite{ADK21}, one finds that Mott's criterion is satisfied at the spectral edge. This observation was used in \cite{ADK21} to prove localization for the top $O(1)$ eigenvectors.

In the interior of the localized phase, the eigenvalue spacing at energy $\lambda$ is typically of order $N^{-\rho_b(\lambda)}$, where $\rho_b(\lambda) > 0$ is an exponent defined in \eqref{lambda_max_rho} below. This is much smaller than the best possible estimate on the tunnelling amplitude, $\ee^{-c \frac{\log N}{\log d}} = N^{-o(1)}$. Hence, Mott's criterion is never satisfied inside the localized phase, and thus the simple picture based on local approximate eigenvectors cannot be used to establish localization. In this paper we therefore introduce a new setup for proving localization.

The first key idea of our proof is to abandon the above simple picture of localization, and to replace the \emph{local} approximate eigenvectors $\f w(x)$ by \emph{global} approximate eigenvectors, denoted by $\f u(x)$, which approximate eigenvectors of $H$ much more accurately and therefore lead to a much smaller tunnelling amplitude \eqref{def_tau}. To define $\f u(x)$, we consider the graph obtained from $\bb G$ by removing all localization centres except $x$; we denote by $\lambda(x)$ the second largest eigenvalue of its adjacency matrix, and by $\f u(x)$ the associated eigenvector. The latter is localized around the vertex $x$. Crucially, the quantity \eqref{def_tau} with $\f w(x)$ replaced by $\f u(x)$ can now be estimated by a \emph{polynomial} error $N^{-\zeta}$ for some $\zeta > 0$.

The price to pay for passing from local approximate eigenvectors $\f w(x)$ to global approximate eigenvectors $\f u(x)$ is a breakdown of orthogonality. Indeed, the vectors $\f u(x)$ have a nonzero overlap, and a significant difficulty in our proof is to control these overlaps and various resulting interactions between localization centres.

To complete the verification of Mott's criterion, we need to establish a polynomial lower bound $\min_{x \neq y \in \cal W} \abs{\lambda(x) - \lambda(y)} \geq N^{-\eta}$ for some $\eta < \zeta$. Clearly, the left-hand side cannot be larger than the eigenvalue spacing $N^{-\rho_b(\lambda)}$ at the energy $\lambda$ we are considering, which yields the necessary bound $\eta > \rho_b(\lambda) > 0$. Hence, we require an anticoncentration result for the eigenvalue difference $\lambda(x) - \lambda(y)$ on a polynomial scale. Owing to the discrete law of $\bb G$ (a product of independent Bernoulli random variables), methods based on smoothness such as Wegner estimates \cite{wegner1981bounds} are not available, and obtaining strong enough anticoncentration is the most involved part of our proof. Our basic strategy is to perform a recursive resampling of neighbourhoods of increasingly large balls around $y$. At each step, we derive a \emph{concentration} bound for $\lambda(x)$ and an \emph{anticoncentration} bound for $\lambda(y)$. The key tool for the latter is a self-improving version, due to Kesten \cite{Kesten1969}, of a classical anticoncentration result of Doeblin, Lévy, Kolmogorov, and Rogozin. In order to obtain sufficiently strong anticoncentration, it is crucial to perform the recursion up to radii comparable to the diameter of $\bb G$.

We conclude this overview with a survey of related results. The eigenvalues and eigenvectors of the Erd\H{o}s-R\'enyi graph have been extensively studied in the denser regime $d \gg \log N$. Complete eigenvector delocalization for $d \gg \log N$ was established in \cite{HeKnowlesMarcozzi2018, EKYY1}. The local spectral statistics in the bulk were proved to follow the universal GOE statistics in \cite{EKYY2, UBY15} for $d \geq N^{o(1)}$. At the spectral edge, the local spectral statistics were proved to be Tracy-Widom for $d \gg N^{1/3}$ \cite{EKYY2, lee2018local}, to exhibit a transition from Tracy-Widom to Gaussian at $d \asymp N^{1/3}$ \cite{huang2020transition}, and to be Gaussian throughout the regime $N^{o(1)} \leq d \ll N^{1/3}$ \cite{huang2020transition, he2021fluctuations}. 
In fact, in the latter regime the Tracy-Widom statistics were recovered in \cite{huang2022edge, lee2021higher} after subtraction of an explicit random shift.

The random $d$-regular graph is another canonical model of sparse random graphs. Owing to the regularity constraint, it is much more homogeneous than the Erd\H{o}s-R\'enyi graph and only exhibits a delocalized phase. The eigenvectors were proved to be completely delocalized for all $d \geq 3$ in \cite{BKY15, BHY1, huang2021spectrum}, and the local spectral statistics in the bulk were shown to follow GOE statistics for $d \geq N^{o(1)}$ in \cite{BHKY15} and at the edge Tracy-Widom statistics for $N^{o(1)} \leq d \ll N^{1/3}$ in \cite{bauerschmidt2020edge, huang2023edge} or for $N^{2/3} \ll d \leq N/2$ in \cite{he2022spectral}.

Anderson transitions have been studied in a variety of models. The archetypal example is the tight-binding, or Anderson, model on $\Z^d$ \cite{abou1973selfconsistent, abou1974self, anderson1978local, aizenman2015random}. In dimensions $d \leq 2$, all eigenvectors of the Anderson model are expected to be localized, while for $d \geq 3$ a coexistence of localized and delocalized phases, separated by a mobility edge, is expected for small enough disorder. So far, only the localized phase of the Anderson model has been understood rigorously, starting from the landmark works \cite{FroSpe, AizMol}; see for instance \cite{aizenman2015random} for a recent survey.

Although a rigorous understanding of the metal-insulator transition for the Anderson tight-binding model is still elusive, some progress has been made for random band matrices. Random band matrices \cite{MFDQS,Wig,CMI1,FM1} interpolate between the Anderson model and mean-field Wigner matrices. They retain the $d$-dimensional structure of the Anderson model but have proved more amenable to rigorous analysis. They are conjectured \cite{FM1} to have a similar phase diagram as the Anderson model in dimensions $d \geq 3$. For $d = 1$ much has been understood both in the localized \cite{Sch,PSSS1,cipolloni2022dynamical, chen2022random} and the delocalized \cite{So1,EKY2,EKYY3, EK1, EK2, EK3, EK4,shcherbina2019universality,BYY1,BYY2,BEYY1,YY1,HM1} phases. For large enough $d$, recent progress in the delocalized phase has been made in \cite{xu2022bulk, yang2021delocalization, yang2022delocalization}.  A simplification of band matrices is the ultrametric ensemble \cite{FOR1}, where the Euclidean metric of $\Z^d$ is replaced with an ultrametric arising from a tree structure. For this model, a phase transition was rigorously established in \cite{VW1}.

Another modification of the $d$-dimensional Anderson model is the Anderson model on the Bethe lattice, an infinite regular tree corresponding to the case $d = \infty$. For it, the existence of a delocalized phase was shown in \cite{ASW1,FHS1,K94}. In \cite{AW1,AW2} it was shown that for unbounded random potentials the delocalized phase exists for arbitrarily weak disorder. The underlying mechanism is resonant delocalization, in which the exponentially decaying tunnelling amplitudes between localization centres are counterbalanced by an exponentially large number of possible channels through which tunnelling can occur, so that Mott's criterion is violated. As a consequence, the eigenvectors hybridize.

Heavy-tailed Wigner matrices, or Lévy matrices, whose entries have $\alpha$-stable laws for $0 < \alpha < 2$, were proposed in \cite{CB1} as a simple model that exhibits a transition in the localization of its eigenvectors; we refer to \cite{ALY1} for a summary of the predictions from \cite{CB1, TBT1}. In \cite{BG1, BG2} it was proved that eigenvectors are weakly delocalized for energies in a compact interval around the origin, and for $0 < \alpha < 2/3$ eigenvectors are weakly localized for energies far enough from the origin. In \cite{ALY1} full delocalization, as well as GOE local eigenvalue statistics, were proved in a compact interval around the origin, and in \cite{ALM1} the law of the eigenvector components was computed. Recently, by comparison to a limiting tree model, a mobility edge was established in \cite{aggarwal2022mobility} for $\alpha$ near $0$ or $1$.

\paragraph{Conventions}
Every quantity that is not explicitly called \emph{fixed} or a \emph{constant} is a sequence depending on $N$. We use the customary notations $o(\cdot)$ and $O(\cdot)$ in the limit $N \to \infty$. For nonnegative $X,Y$, if $X = O(Y)$ then we also write $X \lesssim Y$, and if $X = o(Y)$ then we also write $X \ll Y$. Moreover, we write $X \asymp Y$ to mean $X \lesssim Y$ and $Y \lesssim X$. We say that an event $\Omega$ holds with \emph{high probability} if $\P(\Omega) = 1 - o(1)$. 
Throughout this paper every eigenvector is assumed to be normalized in $\ell^2([N])$. 
Finally, we use $\kappa \in (0,1)$ to denote a small positive constant, which is used to state assumptions and definitions; a smaller $\kappa$ always results in a weaker condition.

\subsection{Results}

Let $\bb G \equiv \bb G(N,d/N)$ be the Erd{\H o}s-R\'enyi graph with vertex set $[N]$ and edge probability $d/N$ for $0 \leq d \leq N$. Let $A = (A_{xy})_{x,y \in [N]} \in \{0,1\}^{N\times N}$ be the adjacency matrix of $\bb G$. Thus, $A =A^*$, $A_{xx}=0$ for all $x \in [N]$, and  $( A_{xy} \col x < y)$ are independent $\op{Bernoulli}(d/N)$ random variables. Define the rescaled adjacency matrix
\begin{equation*}
H \deq A / \sqrt{d}\,.
\end{equation*}
We always assume that $d$ satisfies
\begin{equation} \label{eq:d_range} 
\sqrt{\log N} \, \log \log N \ll d \leq 3 \log N\,.
\end{equation}
Owing to the nonzero expectation of $H$, it is well known that the largest eigenvalue of $H$, denoted by $\lambda_1(H)$, is an outlier separated from the rest of the spectrum (see e.g.\ Proposition~\ref{pro:spectral_gap} \ref{item:spectral_gap_q} below), and we shall always discard it from our discussion. The lower bound of \eqref{eq:d_range} is made for convenience, to ensure that $\lambda_1(H)$ is separated from the bulk spectrum\footnote{In fact, to ensure the separation, the weaker lower bound 
$d \gg \sqrt{\log N / \log \log N}$ would be sufficient; see \cite{KrivelevichSudakov} and \cite[Remark~1.4]{ADK21} for detailed explanations.  We impose the slightly stronger lower bound in \eqref{eq:d_range} for convenience, as it allows us to directly import results from \cite{ADK20} that were proved under this condition.}; the upper bound of \eqref{eq:d_range} is made without loss of generality, since for $d \geq 3 \log N$ the localized phase does not exist and the entire spectrum is known to belong to the delocalized phase \cite{ADK20, ADK_delocalized}.

We denote by $B_{r}(x)$ (respectively by $S_{r}(x)$) the closed ball (respectively the sphere) of radius $r$ with respect to the graph distance of $\mathbb{G}$ around the vertex $x$.
We refer to Section \ref{sec:notations} below for a full account of notations used throughout this paper.

The localized phase is characterized by a threshold $\alpha^*$ defined, for any fixed $\kappa \in (0,1)$ and $\mu \in [0,1]$, as
\begin{equation}\label{eq:largeDegree}
\alpha^{*}\equiv\alpha^{*}(\mu) \deq \max\Big\{ \inf\big\{\alpha>0:\mathbb{P}(\alpha_{1}\geq\alpha)\leq N^{\mu-1}\big\}, 2 + \kappa \Big\}\,.
\end{equation}
We refer to Appendix \ref{sec:alpha_qual} below for the basic qualitative properties of $\alpha^*$ as well as a graph. 
We shall show exponential localization for any eigenvector with eigenvalue $\lambda$ satisfying the condition
\begin{equation} \label{lambda_condition}
\lambda \neq \lambda_1(H)\,, \qquad \abs{\lambda} \geq \Lambda(\alpha^*(\mu))+ \kappa\,,
\end{equation}
for sufficiently small $\mu>0$. In particular, the number of eigenvalues $\lambda$ satisfying \eqref{lambda_condition} is with high probability $N^{\mu + o(1)}$ as $\kappa \to 0$; see Remark \ref{rem:ev_Lambda} below.

\subsubsection{Exponential localization}
Our main result is the following theorem. We recall that by convention all eigenvectors are normalized.

\begin{theorem}[Exponential localization] 
\label{thm:main}
Suppose that \eqref{eq:d_range} holds. Fix $\kappa \in (0,1)$ and $\mu \in (0,1/24)$. Then there is a constant $c \in (0,1)$ depending on $\kappa$ such that, with high probability, for any eigenvector $\f w = (w_x)_{x \in [N]}$ of
$H$ with eigenvalue $\lambda$  satisfying \eqref{lambda_condition}, there exists a unique vertex $x\in[N]$ with 
$\alpha_x > \alpha^*(\mu)$
such that
\begin{equation} \label{eq:exponential_decay_true_eigenvector} 
w_x^2 = {\frac{\alpha_x-2}{2(\alpha_x - 1)}} + o(1)\,, \qquad \qquad 
 \norm{\f w|_{B_i(x)^c}} \lesssim \sqrt{\alpha_x} (1-c)^{i}\,, 
\end{equation} 
for all $i \in \N$ with $1 \leq i \leq \frac{1}{6} \frac{\log N}{\log d}$.  
\end{theorem}

\begin{remark}[Eigenvalue locations] \label{rem:ev_Lambda}
The eigenvalue $\lambda$ of the eigenvector $\f w$ and the associated vertex $x$ from Theorem \ref{thm:main} satisfy $\abs{\lambda} = \Lambda(\alpha_x) + o(1)$ with high probability. This follows from \eqref{eq:eigenvalue_spacing} and \eqref{eq:lambda_x_approximated_by_Lambda_alpha_x} below.

The eigenvalue locations in the localized phase were previously studied in \cite[Theorem~2.1]{ADK19} (see also \cite[Theorem~1.7]{ADK20}). In particular,
if $d > b_* \log N$ then the localized phase does not exist (see \cite[Remark~2.5]{ADK19}) and there is no eigenvalue of $H$ satisfying \eqref{lambda_condition}. 
Conversely, if $d \leq (b_* - \epsilon) \log N$ for some constant $\epsilon > 0$ then for small enough $\kappa > 0$ there is a polynomial number of eigenvalues satisfying \eqref{lambda_condition}, by \cite[Theorem~1.7]{ADK20}. By the same argument, if $d/\log N$ is small enough, then with high probability $N^{\mu + o(1)}$ eigenvalues of $H$ satisfy \eqref{lambda_condition} as $N \to \infty$ and $\kappa \to 0$.
\end{remark}

\begin{remark}[Conditions in Theorem \ref{thm:main}]
The exponential decay in Theorem~\ref{thm:main} holds up to the scale $\frac{\log N}{\log d}$ of the diameter of $\mathbb{G}$. The upper bound $1/24$ and the factor $1/6$ are not optimal (see the discussion in Appendix \ref{sec:optimal_loc}), and they can be improved with some extra effort, which we however refrain from doing here.
\end{remark}

\begin{remark} [Optimal exponential decay] \label{rem:optimal_decay}
If $\mu$ is sufficiently small then the rate $1 - c$ of exponential decay from \eqref{eq:exponential_decay_true_eigenvector} can be made explicit. 
Suppose \eqref{eq:d_range}. 
Then for each small enough constant $\eps >0$ and small enough constant $\mu >0$, depending on $\eps$, 
with high probability, 
 for any eigenvector $\f w$ of
$H$ with eigenvalue $\lambda$ satisfying \eqref{lambda_condition}, there exists a unique vertex $x\in[N]$ with $\alpha_{x} > \alpha^*(\mu)$ such that  
\begin{equation} \label{eq:improved_exponential_decay_true_eigenvector} 
\norm{ \f w|_{B_i(x)^c}} \lesssim \sqrt{\alpha_x} \bigg(\frac{1 + O(\eps)}{\sqrt{\alpha_x-1}}\bigg)^{i+1}
\end{equation} 
for each $i \in \N$ satisfying $i \ll \frac{\log N}{\log d} \frac{1}{\log\frac{10\log N}{d} }$. 
This follows from \eqref{eq:u_approx_u_x} and Proposition~\ref{prop:improved_exponential_decay} below.  

The rate of decay in \eqref{eq:improved_exponential_decay_true_eigenvector} is optimal up to the error term $O(\eps)$. Indeed, by Theorem \ref{prop:localization_profile} and the explicit form \eqref{def_ui}--\eqref{eq:DefCandidate} below, we find that
\begin{equation*}
\norm{\f w|_{B_i(x)^c}} \asymp \sqrt{\alpha_x} \pbb{\frac{1}{\sqrt{\alpha_x - 1}}}^{i+1} + o(1)\,.
\end{equation*}
for any fixed $i \in \N$. In particular, \eqref{eq:improved_exponential_decay_true_eigenvector} improves the rate $\frac{q^i}{(1-q)^2}$ with $q = (2 + o(1))\frac{\sqrt{\alpha_x - 1}}{\alpha_x}$ 
obtained in \cite[Theorem~1.7]{ADK21} at the spectral edge, corresponding to $\alpha_x = ( 1 + o(1)) \alpha^*(0)$.

\end{remark}

\subsubsection{Geometric structure of eigenvectors}
Next, we describe the precise geometric structure of the eigenvectors in the localized phase.
For any vertex $x$ with $\alpha_x > 2$ and radius $r \in \N^*$, we shall define two local vectors, $\f w_r(x)$ and $\f v_r(x)$, which depend only on $\mathbb G$ in the ball $B_r(x)$.
If $\f w$ is an eigenvector of $H$ as in Theorem~\ref{thm:main} with associated vertex $x$, then 
$\f w$ will be well approximated by $\f w_r(x)$ and $\f v_r(x)$ for suitably chosen $r \gg 1$.  

To define these local vectors, we need the following definitions.
Let $r \in \N^*$. For $\alpha > 2$ define the positive sequence $(u_i(\alpha))_{i = 0}^{r-1}$ through
\begin{equation} \label{def_ui}
u_1(\alpha) \deq \pbb{\frac{\alpha}{\alpha - 1}}^{1/2} u_0(\alpha)\,, \qquad u_i(\alpha) \deq \pbb{\frac{1}{\alpha - 1}}^{(i-1)/2} u_1(\alpha) \qquad (1 \leq i \leq r - 1)\,.
\end{equation}
We normalize the sequence by choosing $u_0(\alpha) > 0$ such that $\sum_{i = 0}^{r - 1} u_i(\alpha)^2 = 1$. 

\begin{definition}[Localization profile vectors $\f w_r(x)$ and $\f v_r(x)$] \label{def:w_x_and_v_x} 
Let $r \in \N^*$ and $x \in [N]$. 
\begin{enumerate}[label=(\roman*)] 
\item 
Denote by $\f w_r(x)$ an eigenvector of $H|_{B_r(x)}$ associated with its largest eigenvalue, chosen so that its value at $x$ is nonnegative.
Here $H|_{B_r(x)}$ denotes the matrix $H$ restricted to the vertices in $B_r(x)$ (See Section \ref{sec:notations} below.) 

\item For $\alpha_x > 2$ and $(u_i(\alpha))_{i=0}^{r-1}$ as in \eqref{def_ui}, define 
\begin{equation} \label{eq:DefCandidate}
\qquad\f{v}_{r}(x) \deq \sum_{i=0}^{r - 1}u_{i}(\alpha_x)\frac{\f 1_{S_{i}(x)}}{|S_{i}(x)|^{1/2}}\,,
\end{equation}
where $\f 1_{S_i(x)}$ denotes the indicator function of the sphere $S_i(x)$.
\end{enumerate} 
\end{definition}

Note that $\f w_r(x)$ is unique by the Perron-Frobenius theorem for irreducible matrices with nonnegative entries.

\begin{theorem}[Localization profile] \label{prop:localization_profile} 
Suppose that \eqref{eq:d_range} holds and fix $\kappa \in (0,1)$ and $\mu \in (0,1/24)$.
With high probability, for any eigenvector $\f w$ of $H$ with eigenvalue $\lambda$ satisfying \eqref{lambda_condition}, 
there exists a unique vertex $x \in [N]$ with $\alpha_x > \alpha^*(\mu)$ such that\footnote{We assume that the sign of $\f w$ is chosen suitably.}
\begin{equation} \label{w_v_approx}
\f w = \f w_r(x) + o(1) = \f v_r(x) + o(1)
\end{equation}
for each $r \in \N^*$ satisfying $\log d \ll r \leq \frac{1}{6} \frac{\log N}{\log d}$.  Here, $o(1)$ is meant with respect to the Euclidean norm on $\R^N$. 
\end{theorem} 

In particular, $\f w$ has locally the radial exponentially decaying structure of $\f v_r(x)$.

\subsubsection{Mobility edge and localization length}
Next, we combine the results of this paper with those obtained for the delocalized phase in \cite{ADK20, ADK_delocalized} to establish a mobility edge at $\pm 2$ for certain values of $d$, and analyse the structure of the eigenvectors quantitatively in the vicinity of the mobility edge.

\begin{theorem}[Mobility edge] \label{cor:mobility_edge} 
Fix $\kappa >0$ and suppose that $\frac{23}{24} + \kappa  \leq \frac{d}{b_* \log N} \leq 1 - \kappa $. Then, with high probability, for any eigenvector $\f w$ of $H$ with eigenvalue $\lambda \neq \lambda_1(H)$ we have the following dichotomy.
\begin{enumerate}[label =(\roman*)] 
\item \label{item:localized}  (Localized phase) If $\abs{\lambda} \geq 2 + \kappa$ then $\f w$ is exponentially localized as in \eqref{eq:exponential_decay_true_eigenvector} and \eqref{w_v_approx}. 
\item  \label{item:delocalized} (Delocalized phase)
If $\abs{\lambda} \leq 2 - \kappa$ then $\f w$ is completely delocalized in the sense that
\begin{equation}\label{eq:completely_delocalized_true_eigenvector} 
 \norm{\f w}_\infty^2 \leq N^{-1 + o(1)}\,.  
\end{equation} 
\end{enumerate} 
\end{theorem}

Both phases in Theorem~\ref{cor:mobility_edge} are nonempty under the assumption on $d$; see Section \ref{sec:overview} and \cite[Remark~2.5]{ADK19}. Theorem~\ref{cor:mobility_edge} establishes a dichotomy because \eqref{eq:exponential_decay_true_eigenvector} and \eqref{eq:completely_delocalized_true_eigenvector} are mutually exclusive, since $\norm{\f w}_\infty^2 \geq w_x^2 = \frac{\alpha_x - 2}{2(\alpha_x - 1)} + o(1) \gtrsim 1$ if $\alpha_x > \alpha^*(\mu) \geq 2 + \kappa$.

Next, we investigate the spatial extent of the eigenvectors near the mobility edge. To that end, we use the following notion of localization length. With each normalized vector $\f w$ we associate the length
\begin{equation} \label{def_ll}
\ell(\f w) \deq \min_{x \in [N]} \sum_{y \in [N]} \r d(x,y) \, w_y^2\,,
\end{equation}
where  $\r d(x,y)$ denotes the distance from $x$ to $y$ in the graph $\bb G$. 
Regarding $y \mapsto w_y^2$ as a probability measure on $[N]$, the quantity $\ell(\f w)$ expresses the minimal expected distance from a reference vertex $x$. The minimizing vertex $x$ has the interpretation of a localization centre for $\f w$.

Denote by $\diam(\bb G)$ the diameter\footnote{We recall that the diameter of a connected graph is the length of its longest geodesic. If the graph is disconnected, then its diameter is the maximal diameter of its connected components.} of $\bb G$. It is a classical fact \cite{ChungLu_diameter} that with high probability $\diam(\bb G) = \frac{\log N}{\log d} (1 + o(1))$ as long as $d \gg 1$.

\begin{theorem}[Localization length] \label{thm:ll}
Fix $\kappa > 0$ and suppose that $\frac{23}{24} + \kappa  \leq \frac{d}{b_* \log N} \leq 1 - \kappa$. Then, with high probability, for any eigenvector $\f w$ of $H$ with eigenvalue $\lambda \neq \lambda_1(H)$ we have
\begin{equation} \label{ll_result}
\ell(\f w) =
\begin{cases}
\frac{\abs{\lambda}}{2 \sqrt{\lambda^2 - 4}} + o(1)  & \text{if } \abs{\lambda} \geq 2 + \kappa
\\
\diam(\bb G) (1 + o(1)) & \text{if } \abs{\lambda} \leq 2 - \kappa\,.
\end{cases}
\end{equation}
\end{theorem}

\begin{remark}
By the same proof, the first estimate  of \eqref{ll_result} holds also for all eigenvectors satisfying the conditions of Theorem \ref{thm:main}. Moreover, the constant $\frac{23}{24}$ is not optimal and can be reduced with some extra effort.
\end{remark}

\begin{figure}[!ht]
\begin{center}
{\footnotesize 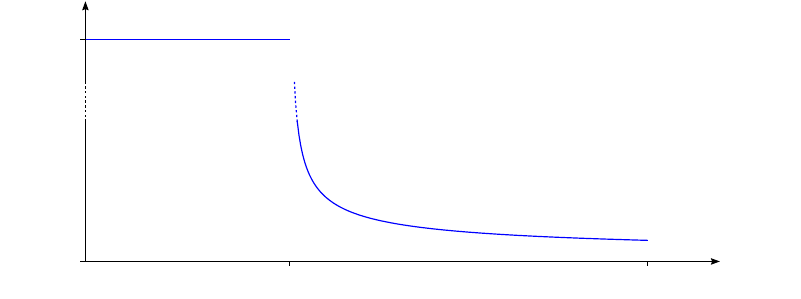}
\end{center}
\caption{An illustration of the behaviour of the localization length \eqref{def_ll} around the mobility edge, established in Theorem \ref{thm:ll}. We plot the asymptotic localization length $\ell$ of an eigenvector with eigenvalue $\lambda$ as a function of $\lambda$. Here $\frac{d}{b_* \log N}$ is a fixed number in $[\frac{23}{24} + \kappa, 1 - \kappa]$. The spectrum is asymptotically given by the interval $[-\lambda_{\max}, \lambda_{\max}]$. We only draw a portion of the spectrum near the right edge. Below the mobility edge $2$, the localization length is $\diam(\bb G) = \frac{\log N}{\log d} (1 + o(1))$. Above the mobility edge $2$, the localization length is finite and diverges as one approaches the mobility edge.
\label{fig:ll}}
\end{figure}

Theorem \ref{thm:ll} shows that the localization length diverges as one approaches the mobility edge from the localized phase, and that it equals the diameter of the graph in the delocalized phase. See Figure \ref{fig:ll} for an illustration.

\subsubsection{Eigenfunction correlator and dynamical localization}

Finally, as a consequence of Theorem~\ref{thm:main}, we control quantities commonly used to characterize Anderson localization (see e.g.\ \cite[Section~1.4]{aizenman2015random}). In particular, we establish exponential decay of the eigenfunction correlator and dynamical localization.

\begin{corollary} \label{cor:dynamical_loc}
Suppose \eqref{eq:d_range}. Then there is a constant $\mu >0$ such that 
for each fixed $\kappa \in (0,1)$, 
there exist constants $c, C > 0$ depending only on $\kappa$ such that the following holds with high probability. Let $J \subset [\Lambda(\alpha^*(\mu))+ \kappa, \sqrt{d}/2]$ be an interval with associated spectral projection $\Pi_J(H)$.
For any $x \in [N]$, any measurable function $F\colon \R \to \C$ satisfying $\norm{F}_\infty \leq 1$, and any $r \geq 0$, we have
\begin{equation} \label{eq:projection_of_evolution_with_restricted_energy} 
\norm{(\Pi_J(H) F(H) \f 1_x ) |_{B_{r}(x)^c}} \leq C \ee^{- c r}\,.
\end{equation}
In particular, denoting by $\f w_\lambda$ the normalized eigenvector of $H$ associated with $\lambda \in \spec(H)$, the \emph{eigenfunction correlator} satisfies the estimate
\begin{equation} \label{ef_correlator}
\sum_{\lambda \in \spec(H) \cap J} \abs{\scalar{\f 1_x}{\f w_\lambda} \scalar{\f w_\lambda}{\f 1_y}} \leq C \ee^{-c \, \dd(x,y)}\,, \qquad x,y \in [N]\,,
\end{equation}
and we have \emph{dynamical localization},
\begin{equation} \label{dynamical_localization}
\sup_{t \in \R} \abs{\scalar{\f 1_y}{\Pi_J(H) \, \ee^{- \ii t H} \f 1_x}} \leq C \ee^{-c \, \dd(x,y)}\,, \qquad x,y \in [N]\,.
\end{equation}
\end{corollary}

\begin{remark}
By a close inspection of the proof in Appendix \ref{sec:dynamical} (using that all error probabilities are polynomially small in $N$), we note that the estimates \eqref{ef_correlator} and \eqref{dynamical_localization} hold also in expectation, provided one multiplies both sides by the factor $\ee^{c \, \dd(x,y)}$.
\end{remark}

\paragraph{Structure of the paper}
We conclude this section with a short summary of the structure of the paper.
In Section~\ref{sec:proof_thm_main}, we collect a few basic notations, then state the three core propositions of the paper: Proposition \ref{prop:exponential_decay_u_x_Section_2}, which gives exponential decay of the approximate eigenvectors, Proposition \ref{prop:mainRigidity}, which compares the approximate eigenvalues with the true eigenvalues, and Proposition \ref{prop:MNonClustering}, which estimates the spacing between neighbouring approximate eigenvalues. After stating them, we use them to deduce Theorem~\ref{thm:main} in Section \ref{sec:exponential_localization_proof_of_theorem_main}. Then, we sketch the proofs of these three propositions in Section~\ref{sec:sketch}. Theorems \ref{prop:localization_profile}, \ref{cor:mobility_edge}, and \ref{thm:ll} are proved in the short Sections \ref{sec:localization_profile_proof_of_proposition}, \ref{sec:mobility_edge}, and \ref{sec:localization_length}, respectively. Section~\ref{sec:preliminaries} is devoted to preliminary results on the graph $\bb G$, its spectrum, and its Green function. Sections~\ref{sec:exponential_decay}, \ref{sec:rigidity} and \ref{sec:Smoothness} are devoted to the proofs of Propositions \ref{prop:exponential_decay_u_x_Section_2}, \ref{prop:mainRigidity}, and \ref{prop:MNonClustering}, respectively. In the appendices, we collect some auxiliary results and basic tools.

\section{Proof of main results} \label{sec:proof_thm_main} 

The rest of the paper is devoted to the proofs of Theorems~\ref{thm:main}, \ref{prop:localization_profile}, \ref{cor:mobility_edge}, and \ref{thm:ll}, as well as Corollary \ref{cor:dynamical_loc}. The former four are proved in this section, while Corollary \ref{cor:dynamical_loc} is proved in Appendix \ref{sec:dynamical}.

Throughout, $\kappa \in (0,1)$ denotes an arbitrary positive constant.

\subsection{Basic notations} \label{sec:notations}

We write $\N = \{0,1,2,\dots\}$. We set $[n] \defeq \{1, \ldots, n\}$
for any $n \in \N^*$ and $[0] \defeq \emptyset$.
We write $\abs{X}$ for the cardinality of a finite set $X$. For $X \subset [N]$ we write $X^c \deq [N] \setminus X$.
We use $\ind{\Omega}$ to denote the indicator function of an event $\Omega$. 

Vectors in $\R^N$ are denoted by boldface lowercase Latin letters like $\f u$, $\f v$ and $\f w$. We use the notation $\f v = (v_x)_{x \in [N]} \in \R^N$ for the entries of a vector.
We denote by $\supp \f v \deq \{x \in [N] \col v_x \neq 0\}$ the support of a vector $\f v$. We denote by $\scalar{\f v}{\f w} = \sum_{x \in [N]} v_x w_x$ the Euclidean scalar product on $\R^N$ and by $\norm{\f v} = \sqrt{\scalar{\f v}{\f v}}$ the induced Euclidean norm. For $X \subset [N]$ we set $\f v |_X \deq (v_x \ind{x \in X})_{x \in [N]}$.
For any $x \in [N]$, we define the standard basis vector $\f 1_x \defeq (\delta_{xy})_{y \in [N]} \in \R^N$, so that $w_x = \scalar{\f 1_x}{\f w}$.
To any subset $S \subset [N]$ we assign the vector $\f 1_S\in \R^N$ given by $\f 1_S \defeq \sum_{x \in S} \f 1_x$. 
In particular, $\f 1_{\{ x\}} = \f 1_x$.

We denote by $\r d(x,y)$ the distance between the vertices $x,y \in [N]$ with respect to the graph $\bb G$, i.e.\ the number of edges in the shortest path connecting $x$ and $y$. For $r \in \N$ and $x \in [N]$, we denote by $B_r(x) \deq \{y \in [N] \col \r d(x,y) \leq r\}$ the closed ball of radius $r$ around $x$, and by $S_r(x) \deq \{y \in [N] \col \r d(x,y) = r\}$ the sphere of radius $r$ around the vertex $x$.
For $X \subset [N]$ we denote by $\bb G |_X$ the subgraph on $X$ induced by $\bb G$.

For a matrix $M \in \R^{N \times N}$, $\norm{M}$ is its operator norm induced by the Euclidean norm on $\R^N$.
For an $N \times N$ Hermitian matrix $M$, we denote by $\lambda_1(M) \geq \lambda_2(M) \geq \cdots \geq \lambda_N(M)$ the ordered eigenvalues of $M$. For an $N \times N$ matrix $M\in \mathbb{R}^{N\times N}$ and a subset $X \subset [N]$, we introduce the $N \times N$ matrices $M |_X \deq (M_{xy} \ind{x,y \in X})_{x,y \in [N]}$ as well as $M^{(X)} \deq M|_{X^c}$ with entries $M^{(X)}_{xy} = M_{xy} \ind{x,y \notin X}$.

\subsection{Exponential localization -- proof of Theorem~\ref{thm:main}} 
\label{sec:exponential_localization_proof_of_theorem_main} 

In this section, after introducing some notation and stating the core propositions of the proof, we use them to prove our main result, Theorem~\ref{thm:main}. 
Recalling the definition of $\Lambda$ from \eqref{eq:defLambda}, we introduce the $\mu$-dependent sets
\begin{align}
{\cal V} &\deq \hb{x\in[N] \col \alpha_x \geq \alpha^*(\mu)}\,, \label{eq:def_cal_V} 
\\
{\cal W} &\deq \hb{x\in\cal V \col \Lambda ( \alpha_x) \geq \Lambda(\alpha^*(\mu)) + \kappa / 2}\,. \label{eq:def_cal_W} 
\end{align}
By definition, $\cal W \subset \cal V$.

The following definition introduces the fundamental approximate eigenvalues and eigenvectors underlying our proof.

\begin{definition}[$\lambda(x)$ and $\f u(x)$] \label{def:ux}
For any $x \in \cal W$, we abbreviate $\lambda(x) \deq \lambda_2(H^{(\cal V \setminus \{x\})})$.
Moreover, we denote by $\f u(x)$ a normalized eigenvector of $H^{(\cal V \setminus \{x\})}$ with eigenvalue $\lambda(x)$ satisfying $\scalar{\f 1_x}{\f u(x)} \geq 0$.
\end{definition}

As we shall see, with high probability $\lambda(x)$ is a simple eigenvalue and hence $\f u(x)$ is unique (see Corollary~\ref{cor:largest_eigenvalues_vertices_removed} below).

The proof of Theorem \ref{thm:main} consists of three main steps, which are the content of the three following propositions. The next proposition states that $\f u(x)$ has the exponential decay claimed in Theorem~\ref{thm:main} for $\f w$. It is proved in Section \ref{sec:exponential_decay} below.

\begin{proposition}[Exponential decay of $\f u(x)$] \label{prop:exponential_decay_u_x_Section_2} 
Suppose that \eqref{eq:d_range} holds. 
Then there  is a constant $c \in (0,1)$ such that, for each fixed $\mu \in [0,1/3)$,
 with high probability, for each $x \in \cal W$, 
\[ \scalar{\f 1_x}{\f u(x)} = \sqrt{\frac{\alpha_x-2}{2(\alpha_x - 1)}} + o(1)\,, \qquad \qquad \norm{\f u(x)|_{B_i(x)^c}} \lesssim \sqrt{\alpha_x} (1- c)^i \] 
for all $i \in \N$ satisfying $1 \leq i \leq \min\big\{ \frac{1}{5} - \frac{\mu}{4}, \frac{1}{3} - \mu \big\} \frac{\log N}{\log d}-2$. 
\end{proposition} 

In the proof of Theorem~\ref{thm:main}, the next two propositions will be used to conclude that 
any eigenvector of $H$ whose associated eigenvalue satisfies \eqref{lambda_condition} is close to 
$\f u(x)$ for some $x \in \cal W$. 
Given Proposition~\ref{prop:exponential_decay_u_x_Section_2}, this will directly imply 
Theorem~\ref{thm:main}.

The next proposition, Proposition~\ref{prop:mainRigidity}, states that $\lambda(x)$ and $\f u(x)$ are approximate eigenvalues and eigenvectors of $H$, 
respectively,  
with an error bounded by an inverse power of $N$. Moreover, up to such an error, 
each eigenvalue of $H$ satisfying \eqref{lambda_condition} is approximated by $\lambda(x)$ for some $x \in \cal W$. In particular, it provides an upper bound for the tunnelling amplitude, in the sense of \eqref{def_tau}, for the global approximate eigenvectors from Definition \ref{def:ux}.
Its proof is given in Section \ref{sec:rigidity} below. 

\begin{proposition}[Approximate eigenvalues] \label{prop:mainRigidity} 
Suppose \eqref{eq:d_range}. 
Fix $\mu \in [0,1/4)$ and $\zeta \in [0, 1/2 - 3\mu/2)$. 
Then, with high probability, for each $x \in \cal W$ there exists $\epsilon_x \in \R$ such that
\[
\spec( H )\cap \cal I  = \{{\lambda}(x) +\epsilon_x \col x\in\cal W\}\cap \cal I\,, 
\qquad \qquad \cal I \deq [\Lambda(\alpha^*) + 3\kappa/4,\sqrt{d}/2]
\] 
counted with multiplicity and 
\[\max_{x\in\cal W}\max\big\{ \abs{\epsilon_x}, \, \|(H-\lambda(x))\f u(x)\| \big\} \leq  N^{-\zeta}\,.\]   
\end{proposition}

The next proposition establishes a spacing of at least $N^{-\eta}$ between the approximate eigenvalues $(\lambda(x))_{x \in \cal W}$, for large enough $\eta$. 
It is proved in Section \ref{sec:Smoothness} below.

\begin{proposition}[Eigenvalue spacing] \label{prop:MNonClustering} 
Suppose \eqref{eq:d_range}. 
Fix $\mu \in (0,1/24)$ and $\eta > 8 \mu$. 
Then, with high probability, 
\[ \abs{\lambda(x) - \lambda(y)} \geq N^{-\eta} \] 
for all $x$, $y \in \cal W$ with $x \neq y$. 
\end{proposition}

We now deduce Theorem \ref{thm:main} from Propositions~\ref{prop:exponential_decay_u_x_Section_2}, 
\ref{prop:mainRigidity}, and \ref{prop:MNonClustering}.

\begin{figure}[!ht]
\begin{center}
{\footnotesize 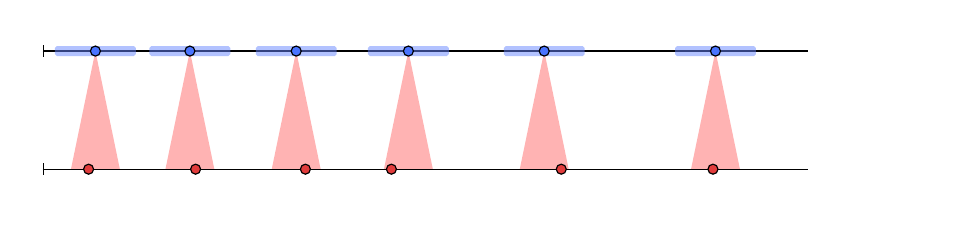}
\end{center}
\caption{An illustration of the setup for the perturbation theory in
the proof of Theorem \ref{thm:main}. We draw two instances of the
interval $[\Lambda(\alpha^*) + \kappa, \infty)$. On the top line, we draw each $\lambda(x)$ for $x
\in \cal W$ as blue dot. Each dot is surrounded by a blue buffer of width
$N^{-\eta}$. By Proposition \ref{prop:MNonClustering}, these buffers
do not intersect with high probability. On the bottom line, we draw
each eigenvalue of $H$ as red dot. Each dot on the top line gives rise to a red
region on the bottom line of width $2 N^{-\zeta}$. Since $\zeta >
\eta$, the red regions are disjoint. By
Proposition~\ref{prop:mainRigidity}, each red region contains exactly 
one eigenvalue of $H$ and each eigenvalue of $H$ is contained in a red region. Hence, the eigenvalues of $H$ are separated by
at least $N^{-\eta}/2$.
\label{fig:eigenvalues}}
\end{figure}

\begin{proof}[Proof of Theorem \ref{thm:main}]  
Let $\lambda$ be an eigenvalue of $H$ satisfying \eqref{lambda_condition}, and $\f w$ an associated, normalized eigenvector. 
Fix $\zeta \in (8\mu, 1/3)$ and $\eta \in (8\mu, \zeta)$. 
As $\mu < 1/24$, both intervals are nonempty and $\zeta < 1/2 - 3\mu/2$. 
Thus, Propositions~\ref{prop:mainRigidity} and \ref{prop:MNonClustering} are applicable with these choices of $\zeta$ and $\eta$. 
We shall show below that, on the intersection of the high-probability events of Propositions~\ref{prop:mainRigidity} and \ref{prop:MNonClustering}, there exists a unique $x \in \cal W$ such that
\begin{equation} \label{eq:u_approx_u_x} 
\f w= \f u(x)  +O(N^{\eta-\zeta})
\end{equation} 
 (under suitable choice of the sign of $\f w$). Thus, Theorem~\ref{thm:main} follows from $\eta < \zeta$ and Proposition~\ref{prop:exponential_decay_u_x_Section_2} as $\exp\big( \frac{1}{6} \frac{\log N}{\log d} \log (1-c) \big) \gg N^{-\eps}$ for any $\eps>0$. 

What remains is the proof of \eqref{eq:u_approx_u_x}.  
This is an application of perturbation theory in the form of Lemma~\ref{lem:perturbation_theory}, 
whose conditions we justify now. Note first that, combining \eqref{eq:H_minus_EH_ll_sqrt_d} below and the trivial fact $\lambda_1(\E H) = \sqrt{d} (1 + o(1))$ with rank-one eigenvalue interlacing (Lemma \ref{lem:interlacing2}), we conclude that with high probability $\lambda_1(H) = \sqrt{d} (1 + o(1))$ and $\lambda_2(H) \leq \sqrt{d}/2$. Hence, the eigenvalue $\lambda$ satisfying \eqref{lambda_condition} lies in $[\Lambda(\alpha^*(\mu))+ \kappa, d/2]$. From Propositions~\ref{prop:mainRigidity} and \ref{prop:MNonClustering} with $\eta < \zeta$, we conclude that 
\begin{equation} \label{eq:eigenvalue_spacing} 
 \dist(\lambda, \spec(H)\setminus\{\lambda\}) \geq N^{-\eta} - 2 N^{-\zeta}\,, \qquad \qquad \abs{\lambda - \lambda(x)} \leq N^{-\zeta} 
\end{equation} 
for a unique $x \in\cal W$ (see Figure~\ref{fig:eigenvalues}).
In particular, as $\eta < \zeta$, there is $\Delta \asymp N^{-\eta}$ such that $\lambda$ is the 
unique eigenvalue of $H$ in $[\lambda(x) - \Delta, \lambda(x) + \Delta]$. 
Moreover, $\norm{(H-\lambda(x))\f u(x)} \leq N^{-\zeta}$ by Proposition~\ref{prop:mainRigidity}.
Therefore, all conditions of Lemma~\ref{lem:perturbation_theory} are satisfied, and it implies \eqref{eq:u_approx_u_x}. 
This concludes the proof of Theorem~\ref{thm:main}.
\end{proof}

By choosing $\eta$ in \eqref{eq:eigenvalue_spacing} of the proof of Theorem~\ref{thm:main} sufficiently small, we conclude the following result.

\begin{corollary}[Eigenvalue spacing of $H$]  \label{cor:eigenvalue_spacing} 
Suppose \eqref{eq:d_range}. 
Fix $\mu \in (0,1/24)$ and $\eta > 8\mu$. Then, with high probability,
\begin{equation*}
\dist(\lambda, \spec(H) \setminus \{ \lambda \}) \geq N^{-\eta}\,,
\end{equation*}
for every $\lambda \in \spec(H) \cap[\Lambda(\alpha^*(\mu))+ \kappa, \infty)$.
\end{corollary}

\begin{remark}[Eigenvalue spacing in critical regime] 
In the critical regime, i.e.\ when $d \asymp \log N$, the lower bound on $\eta$ in the conditions of 
Proposition~\ref{prop:MNonClustering} and Corollary~\ref{cor:eigenvalue_spacing} can be weakened to $\eta > 4\mu$. 
Details can be found in Remark~\ref{rem:critical61} below. 
\end{remark} 

\begin{remark}[Eigenvector mass on vertices in $\cal V \setminus \{x \}$] 
With high probability the following holds. 
Let $\f w$ be an eigenvector of $H$ associated with the vertex $x$ as in Theorem~\ref{thm:main}. Then, 
from $\f u(x)|_{\cal V \setminus \{x \}} = 0$ and \eqref{eq:u_approx_u_x}, we conclude 
\[  
\norm{\f w|_{\cal V\setminus \{x \}}} \lesssim N^{-\eps} 
\] 
for any small enough $\eps >0$. 
\end{remark}

\subsection{Sketch of the proof} \label{sec:sketch} 
In this subsection we sketch the proof of Theorem~\ref{thm:main}. We use the definitions and notations from Sections \ref{sec:notations} and \ref{sec:exponential_localization_proof_of_theorem_main}.

The basic strategy is to find an orthogonal matrix $U$, a diagonal $n\times n$
matrix $\Theta = \diag(\theta_1, \dots, \theta_n)$, a symmetric $(N-n) \times (N-n)$ matrix $X$, and a symmetric $N\times N$ matrix $E$
such that the following holds. 
In the basis of the columns $\f u_1, \dots, \f u_N$ of $U$, the matrix $H$ has the form  
\begin{equation}\label{eq:blockSketch}
U^* H U = \begin{pmatrix} \Theta & 0 \\ 0 & X \end{pmatrix}  + E\,, 
\end{equation} 
where the matrices $E$ and $X$ satisfy 
\begin{subequations} \label{eq:norm_E_spec_X_away} 
\begin{align} 
& \|E\|\ll \min \big\{ \abs{\theta_i - \theta_j} \col i \neq j \big\}, \label{eq:norm_E_smaller_level_spacing} \\ 
& \dist(\spec(X) , I) \geq \kappa\,; \label{eq:spec_X_away_from_I}   
\end{align}
\end{subequations} 
here $I$ denotes the interval containing the eigenvalues of $H$ that we are interested in (cf.\ \eqref{lambda_condition}). We call the first $n$ columns $\f u_1, \dots, \f u_n$ of $H$ \emph{profile vectors}.

If $\norm{E} = o(1)$ then, for each $i \in [n]$, the vector $\f u_i$ is an approximate eigenvector of $H$ 
with approximate eigenvalue $\theta_i$. Unlike approximate eigenvalues, in general approximate eigenvectors have nothing to do with the actual eigenvectors. For $\f u_i$ to be close to an eigenvector of $H$, we require the stronger estimates \eqref{eq:norm_E_spec_X_away}, which can be regarded as a version of Mott's criterion in terms of the profile vectors encoded by $\f u_1, \dots, \f u_n$. Localization then follows provided the $\f u_i$ are shown to be localized.

In \cite{ADK19}, it was showed that there is a one-to-one correspondence between eigenvalues of $H$ in the semilocalized phase $[2 + o(1), \infty) \setminus \{ \lambda_1(H)\}$ 
and vertices $x$ of $\mathbb{G}$ with normalized degree $\alpha_x \geq 2 + o(1)$. 
Subsequently, in \cite{ADK20,ADK21}, the eigenvectors of $H$ in the semilocalized phase were investigated using the decomposition \eqref{eq:blockSketch}. There, the profile vectors $\f u_i$ were supported in balls $B_r(x)$ around vertices $x$ 
of sufficiently large $\alpha_x$, where $r \gg 1$. We refer to such vectors as \emph{local profile vectors}: they are spatially localized (in the graph distance) and their supports are disjoint. Examples of such local profile vectors are $\f w_r(x)$ and $\f v_r(x)$ for $x \in \cal W$, defined in Definition~\ref{def:w_x_and_v_x} (others were defined in \cite{ADK20,ADK21}).

The local profile vectors are exponentially decaying with a rate $c > 0$ depending on $b$ and the energy. The best possible error estimate for $\norm{(H - \theta_i) \f u_i}$ under the condition that $\f u_i$ is supported in $B_r(x)$ is obtained by choosing $\f u_i = \f w_r(x)$, the top eigenvector of $H|_{B_r(x)}$; in that case the error is purely a boundary term of order $\ee^{-cr}$. If the supports of the profile vectors are separated by more than 1 (in the graph distance) then it is easy to see that $\norm{E} \leq \max_{i} \norm{(H- \theta_i) \f u_i}$ (see also our formulation of Mott's criterion \eqref{def_tau}). Hence, the best possible error estimate for $\norm{E}$ is $\ee^{-cr}$. Since the diameter of $\mathbb{G}$ is $\frac{\log N}{\log d}(1 + o(1))$ 
with high probability, the best bound resulting from this approach is $\norm{E} \lesssim N^{-c/\log d} =N^{-o(1)}$ for some constant $c>0$. However, inside the semilocalized phase this bound is always much larger than the typical eigenvalue spacing $N^{-\eta}$ for some $\eta > 0$. Recalling the condition \eqref{eq:norm_E_smaller_level_spacing}, we conclude that any approach to prove localization in the semilocalized phase that uses local profile vectors is doomed to fail.

The reason why any approach based on local profile vectors fails is that local profile vectors (such as $\f w_r(x)$) are supported on balls containing a comparatively small set of vertices, and the mass of the true eigenvectors outside of such balls is not small enough to be fully negligible.
This leads us to introduce the \emph{global profile vectors} $\f u(x)$, $x \in \cal W$, from Definition \ref{def:ux}, with associated approximate eigenvalues $\lambda(x)$. They are defined as the second eigenvector-eigenvalue pair of the matrix $H^{(\cal V\setminus \{x\})}$. Thus, $\Theta = \diag( (\lambda(x))_{x \in \cal W})$ and the first $n = \abs{\cal W}$ columns of $U$ are given by the orthonormalization of the family $(\f u(x))_{x \in\cal W}$. The global profile vector $\f u(x)$ and the best possible local profile vector $\f w_r(x)$ are each defined as eigenvectors of the graph after removal of a set of vertices, $\abs{\cal V} - 1 \sim N^{\mu} \ll N$ vertices for the former and $\abs{B_r(x)^c} \sim N$ vertices for the latter. This suggests that $\f u(x)$ is a better approximation of a true eigenvector of $H$. The price to pay is that its definition is less explicit, and, crucially, the family $(\f u(x))_{x \in \cal W}$ is not orthogonal owing to the global profile vectors having nonzero overlaps. As explained below, the need to control the overlaps of the global profile vectors presents a serious complication.

The proof of Theorem \ref{thm:main} consists of three main steps:
\begin{enumerate}[label=(\roman*)] 
\item \label{item:sketch_0}
exponential decay for $\f u(x)$ around $x$,
\item \label{item:sketch_i}
$\|E\|\leq N^{-\zeta}$ and $\spec(X)$ is separated from $I$,
\item  \label{item:sketch_ii} 
$\min_{x\neq y\in \cal W} |\lambda(x)-\lambda(y)| \geq N^{-\eta}$, 
\end{enumerate} 
for some constants $\zeta>\eta>0$ (see also Section \ref{sec:exponents} below). These items corresponds to Propositions \ref{prop:exponential_decay_u_x_Section_2}, \ref{prop:mainRigidity}, and  \ref{prop:MNonClustering}, respectively. We outline their proofs in Sections \ref{sec:sketch_decay_u_x}, \ref{sec:sketch_approximate_eigenvalues}, and \ref{sec:sketch_eigenvalue_spacing}, respectively.

\subsubsection{Exponential decay of $\f u(x)$} \label{sec:sketch_decay_u_x}
First we explain the need to introduce the two different vertex sets $\cal W \subset \cal V$. In the definition of $\lambda(x)$ and $\f u(x)$, all vertices in $\cal V\setminus \{x \}$ are removed, while the profile vectors $\f u(x)$ are only considered for $x$ in the smaller set $\cal W$. The difference between $\cal W$ and $\cal V$ is used precisely to obtain a spectral gap for $H^{(\cal V \setminus \{x \})}$ around $\lambda(x)$, $x \in \cal W$.

To show exponential decay of $\f u(x)$, we use its definition, a simple computation, and a truncated 
Neumann series expansion to obtain
\begin{equation} \label{eq:sketch_exp_decay_u_x_resolvent_expansion} 
\begin{aligned} 
\f u(x)|_{\cal V^c} & = c_x \bigg(1-\frac{H^{(\cal V)}}{\lambda(x)}\bigg)^{-1}\f 1_{S_1(x)}\\ 
 & =c_x\sum_{k=0}^{n-1}\bigg(\frac{H^{(\cal V)}}{\lambda(x)}\bigg)^{k}\f 1_{S_1(x)}+c_x\bigg(1-\frac{H^{(\cal V)}}{\lambda(x)}\bigg)^{-1}\bigg(\frac{H^{(\cal V)}}{\lambda(x)}\bigg)^{n}\f 1_{S_1(x)}\,, 
\end{aligned} 
 \end{equation}  
where $c_x \deq \frac{\scalar{\f 1_x}{\f u(x)}}{\lambda(x) \sqrt{d}}$. 
Each term of the sum is supported in $B_n(x)$ and, thus, vanishes when restricting to $B_n(x)^c$.
Hence, as $\abs{\scalar{\f 1_x}{\f u(x)}} \lesssim 1$, $\norm{\f 1_{S_1(x)}}= \sqrt{d \alpha_x}$, and $\norm{(\lambda(x) - H^{(\cal V)})^{-1}} \lesssim 1$ by the spectral gap of $H^{(\cal V\setminus \{ x\})}$ around $\lambda(x) = \lambda_2(H^{(\cal V\setminus \{x\})})$ mentioned above, we obtain $\|\f u(x)|_{B_n(x)^c}\|\lesssim \sqrt{\alpha_x} \, q^n $ with $q=\lambda(x)^{-1}\|H^{(\cal V)}|_{B_{n+1}(x)}\|<1$. 
This is the desired exponential decay. 
We remark that in \eqref{eq:sketch_exp_decay_u_x_resolvent_expansion} 
we establish exponential decay of the Green function of $H^{(\cal V)}$ evaluated at $\lambda(x)$, 
using that $\lambda(x)$ is away from the spectrum of $H^{(\cal V)}$. 
This is an instance of a Combes-Thomas estimate, and we translate it to an exponential decay for the eigenvector $\f u(x)$. 
Furthermore, we show that $\lambda(x)$ is isolated in the spectrum of $H^{(\cal V\setminus\{x\})}$ and, thus,
 perturbation theory implies that $\f u(x) = \f w_r(x) + o(1)$. 

The rate obtained from the above argument is far from optimal, but an extension of this argument does yield the optimal rate of decay for $\norm{\f u(x)|_{B_i(x)^c}}$ for small enough $\mu$. To that end, we choose $n = C i$ for a large constant $C>0$ in \eqref{eq:sketch_exp_decay_u_x_resolvent_expansion}, which makes the remainder term in \eqref{eq:sketch_exp_decay_u_x_resolvent_expansion} with $n = C i$ 
 subleading. 
It remains to estimate the terms $(H^{(\cal V)})^k{\f 1_{S_1(x)}}$ for $k = i$, \ldots, $ C i$, since they vanish for $k < i$ when restricted to $B_i(x)^c$. For $k \geq i$, we relate $(H^{(\cal V)})^k{\f 1_{S_1(x)}}$ to the number of a family of walks on the graph $\mathbb{G}$. 
We obtain optimal bounds on this number by a path counting argument, exploiting the tree structure of $\mathbb{G}|_{B_r(x)}$, a precise bound on the degrees in $B_r(x)\setminus \{x \}$ and the concentration of the sphere sizes $\abs{S_i(x)}$ for $i \leq r$.

\subsubsection{Approximate eigenvalues} \label{sec:sketch_approximate_eigenvalues} 
We now sketch how \ref{item:sketch_i} is proved. 
In contrast to the case of local profile vectors discussed above, the proof of $\norm{E} \lesssim \max_{x \in \cal W} \norm{(H- \lambda(x)) \f u (x)}$ requires also a control of the nonzero overlaps $\scalar{\f u(x)}{H \f u(y)}$ for $x \neq y$. By exponential decay of $\f u(x)$, it is easy to see that these overlaps are $N^{-o(1)}$, but, as explained above, a polynomial bound $N^{-c}$ is required to prove localization.
The construction of $\f u(x)$ and $\lambda(x)$ and a simple computation reveal that 
\begin{equation} \label{eq:sketch_computation_H_minus_lambda_x_u_x} 
(H-\lambda(x))\f u(x)=\sum_{y\in{\cal V} \setminus \{x \}}\epsilon_{y}(x) \f 1_{y}\,, \qquad \qquad  
\epsilon_{y}(x) \deq \frac{1}{\sqrt{d}}\sum_{t\in S_1(y)} \langle \f 1_t, \f u(x) \rangle\,.
\end{equation} 
Then the main idea to estimate $\eps_y(x)$ is the following elementary bound. 
Let $\cal T$ be a finite set. For any $(u_t)_{t\in \cal T} \in \C^{\cal T}$ and any random $T \subset \cal T$ we have
\begin{equation}
\E \qBB{\sum_{t \in T} \abs{u_t}^2}= \sum_{t \in \cal T} \E [\ind{t \in T}] \abs{u_t}^2 \leq \max_{t \in \cal T} \P(t \in T) \sum_{t \in \cal T} \abs{u_t}^2\,.\label{eq:magic_lemma} 
\end{equation}
Heuristically, by the independence of the edges in the Erd{\H o}s-R\'enyi graph, the edges between $\cal V \setminus \{x\}$ and $(\cal V \setminus \{x\})^c$ are sampled independently of the
subgraphs $\mathbb{G}|_{({\cal V}\setminus\{x\})^c}$ and $\mathbb{G}|_{\cal V\setminus\{x\}}$ and are, therefore, independent of $\f u(x)$. 
Hence, \eqref{eq:sketch_computation_H_minus_lambda_x_u_x}, \eqref{eq:magic_lemma} and $\abs{S_1(y)} \lesssim \log N$ yield
\[
\mathbb{E}\big[|\epsilon_y(x)|^2\big|A^{(y)}\big] \leq \frac{1}{d} \,\mathbb{E}\bigg[|S_1(y)|\sum_{t\in S_1(y)}\langle \f 1_t, \f u(x) \rangle^2\bigg|A^{(y)}\bigg] \lesssim \frac{(\log N)^2}{Nd}\|\f u(x)\|^2=N^{-1+o(1)}\,. 
\]
From $|{\cal V}|=N^{\mu+o(1)}$ and Chebyshev's inequality, we therefore conclude $\|(H-\lambda(x))\f u(x)\|^2\leq N^{2\mu-1+o(1)}$ with high probability.

Since $(\f u(x))_{x \in \cal W}$ is not an orthogonal family, we choose the first columns of $U$ in \eqref{eq:blockSketch} as the Gram-Schmidt orthonormalization $( \f u^\perp(x))_{x \in \cal W}$ 
of $(\f u(x))_{x \in \cal W}$ (with respect to a fixed order on $\cal W$), i.e.\ 
\begin{equation} \label{eq:sketch_def_u_perp}  
\f u ^\perp (x) \deq \frac{\f u(x) - \Pi_{<x} \f u(x)}{\norm{\f u(x) - \Pi_{<x} \f u(x)}}\,, 
\end{equation} 
where $\Pi_{<x}$ denotes the orthogonal projection onto $\Span\{ \f u(y) \colon y \in \cal W, \, y < x\}$. 
It remains to show that, for any $x \in \cal W$, $\norm{(H-\lambda(x))\f u^\perp(x)}$ is also bounded by an inverse power of $N$. 
The denominator in \eqref{eq:sketch_def_u_perp} is $\gtrsim 1$ since $\scalar{\f 1_x}{\f u(x)} \gtrsim 1$ and $\scalar{\f 1_y}{\f u(x)} = 0$ for all $y \neq x$. 
Moreover, $H-\lambda(x)$ applied to the numerator of \eqref{eq:sketch_def_u_perp} is bounded by a negative power of $N$, since 
$\norm{(H-\lambda(x))\f u(x)}\leq N^{\mu - 1/2 + o(1)}$ as shown above, and $\norm{\ol{\Pi}_{< x} H {\Pi}_{<x}} \leq N^{\mu - 1/2 + o(1)}$ where $\ol{\Pi}_{<x} \deq 1 - \Pi_{<x}$. The latter bound is proved using the above estimate on $\norm{(H -\lambda(x))\f u(x)}$.  Given this construction and the bounds explained above, we extend $(\f u^\perp(x))_{x \in \cal W}$ to an orthonormal basis and choose these vectors 
as columns of $U$. In particular, the first $n = \abs{\cal W}$ columns of $U$ are given by $(\f u^\perp(x))_{x \in\cal W}$. 

What remains is to show that $\spec(X)$ is separated from $I$. To that end, we decompose its domain of definition into the span of $\f w_1$, the eigenvector of $H$ associated with its largest eigenvalue, 
and its orthogonal complement. This largest eigenvalue and the overlaps between $\f w_1$ and $\f u^\perp(x)$ for all $x \in \cal W$ can be controlled relatively 
precisely, and we omit $\f w_1$ from the remaining explanations below.
It suffices to show that $\lambda_1(X) \leq \Lambda(\alpha^*) + o(1)$, by the definitions of $\cal V$ and $\cal W$.
This upper bound on $\lambda_1(X)$ is equivalent to
\begin{equation} \label{X_estimate_sketch}
\lambda_1(\ol{\Pi}_{\cal W} H \ol{\Pi}_{\cal W}) \leq \Lambda(\alpha^*) + o(1)\,,
\end{equation}
where $\ol{\Pi}_{\cal W}\deq 1 - \Pi_{\cal W}$ 
and $\Pi_{\cal W}$ is the orthogonal projection onto $\Span\{ \f u^\perp(x) \colon x \in \cal W\}$. An estimate of the form \eqref{X_estimate_sketch} was first derived in \cite{ADK20}, except that there the projection was defined in terms local profile vectors. Since we are using the global profile vectors $\f u^\perp(x)$, in our case this estimate is considerably more involved. The rough strategy is to make a link between \eqref{X_estimate_sketch} and a corresponding estimate for local profile vectors, which was already established in \cite{ADK20} using bounds on the non-backtracking matrix of $H$, an Ihara-Bass type identity, and a local delocalization result for approximate eigenvectors.

To that end, let $Q$ be the orthogonal projection onto the complement of $\bigcup_{y \in \cal V \setminus \cal W} B_{2r_\star - 1}(y)$, 
where $r_\star \asymp \sqrt{\log N}$ is chosen as in \cite[eq.~(1.8)]{ADK20}. In particular, the local profile vectors $\f v(x)$ from \cite{ADK20} 
satisfy $\supp \f v(x)\subset B_{r_\star}(x)$ for every $x \in \cal V$. 
We denote by $\Pi_Q$ the projection onto $\Span\{ Q \f u^\perp(x) \colon x \in \cal W\}$. 
If we obtain a small enough upper bound on $\norm{\Pi_{\cal W} - \Pi_{Q}}$ then 
it suffices to show $\lambda_1(\ol{\Pi}_Q H \ol{\Pi}_Q) \leq \Lambda(\alpha^*) + o(1)$, where $\ol{\Pi}_{Q} \deq 1 -  \Pi_{Q}$. 
To get a sufficient estimate for $\norm{\Pi_{\cal W} - \Pi_{Q}}$, we need that $\f u^\perp(x)|_{B_r(y)}$ is polynomially small in $N$ for $x \neq y$ 
as a summation over $x,y \in \cal W$ is required. This is achieved through an argument motivated by \eqref{eq:magic_lemma}. 
Let $\Pi_{\f v}$ be the orthogonal projection onto $\Span\{ \f v(x) \colon x \in \cal V \setminus \cal W\}$. 
By definition of $Q$ and $\supp \f v(x) \subset B_{r_\star}(x)$, $\Pi_Q$ and $\Pi_{\f v}$ commute. 
Thus, 
$\lambda_1(\ol{\Pi}_Q H \ol{\Pi}_Q) \leq \max\{\lambda_1(\Pi_{\f v} H \Pi_{\f v}), \lambda_1(\ol{\Pi}_Q \ol{\Pi}_{\f v} H\ol{\Pi}_{\f v}\ol{\Pi}_Q ) \} + o(1)$, 
where $\ol{\Pi}_{\f v} \deq 1 -  \Pi_{\f v}$. 
By \cite{ADK20}, $\lambda_1(\Pi_{\f v} H \Pi_{\f v}) \leq \Lambda(\alpha^*) + o(1)$. 
For eigenvectors of $H$ which are orthogonal to the local profile vectors $\f v(x)$ and whose associated eigenvalues are large enough, 
we obtain a weak delocalization estimate by following an argument in \cite{ADK20}. 
This weak delocalization estimate shows that 
$\lambda_1(\ol{\Pi}_Q \ol{\Pi}_{\f v} H\ol{\Pi}_{\f v}\ol{\Pi}_Q ) \leq \Lambda(\alpha^*) + o(1)$ if 
 $\lambda_1(\ol{\Pi}_{\cal W} H^{(\cal V\setminus \cal W)} \ol{\Pi}_{\cal W}) \leq \Lambda(\alpha^*) + o(1)$. 
The last bound is finally obtained by a careful analysis of the spectrum of $H^{(\cal V\setminus \cal W)}$, which is based 
on viewing it as a perturbation of ${\Pi}_{\cal W} H^{(\cal V\setminus \cal W)} {\Pi}_{\cal W} + \ol{\Pi}_{\cal W} H^{(\cal V\setminus \cal W)} \ol{\Pi}_{\cal W}$
and analysing $H^{(\cal V\setminus \cal W)}$ on $\ran\Pi_{\cal W}$ in detail.

\subsubsection{Eigenvalue spacing} \label{sec:sketch_eigenvalue_spacing} 
We now sketch how to prove \ref{item:sketch_ii}. To that end, we fix $a \neq b \in \cal W$. 
To prove that $\lambda(a)$ and $\lambda(b)$ are not too close to each other, we choose an appropriate radius $r$ 
on the scale $\frac{\log N}{\log d}$ of the diameter of $\mathbb{G}$. 
Then we fix the two subgraphs $\mathbb{G}|_{B_{r}(b)}$ and $\mathbb{G}|_{B_{r}(b)^{c}}$ and show that resampling the edges between $S_r(b)$ 
and $B_r(b)^c$ results in a substantial change of $\lambda(b)$ while $\lambda(a)$ remains almost unchanged: we establish simultaneous \emph{anticoncentration} for $\lambda(b)$  and \emph{concentration} for $\lambda(b)$, which yields anticoncentration for their difference.
The edges between $S_r(b)$ and $B_r(b)^c$ form an independent family of Bernoulli random variables 
by definition of the Erd{\H o}s-R\'enyi graph. 

On a more formal level, we work conditionally on $\cal F \deq \sigma (B_r(b), A|_{B_r(b)},A|_{B_r(b)^c})$ and prove the following two statements in order to obtain a lower bound on $|\lambda(a)-\lambda(b)|$. 
\begin{enumerate}[label=(\alph*)]
\item \label{item:lambda_a} $\lambda(a)$ fluctuates little under resampling of the edges between $S_{r}(b)$ and $B_{r}(b)^{c}$, i.e.\ the \emph{concentration} estimate 
\[\mathbb{P}(|\lambda(a)-z|\ll N^{-\eta} | \cal F)\geq 1-N^{-\eta/2+o(1)}\]
holds if $z$ is the second largest eigenvalue of $H^{(\cal V \cup B_r(b)\setminus \{a\})}$, which is $\cal F$-measurable.
\item \label{item:lambda_b} $\lambda(b)$ fluctuates a lot under resampling of the edges between $S_{r}(b)$ and $B_{r}(b)^{c}$, i.e.\ the \emph{anticoncentration} estimate 
\[\mathbb{P}(|\lambda(b)-z|\geq N^{-\eta}| \cal F)\geq 1- N^{-\eta/2+o(1)}\]  
holds for any $\cal F$-measurable spectral parameter $z$ in $I$ (see the definition of $I$ after \eqref{eq:spec_X_away_from_I}). 
\end{enumerate}
We justify \ref{item:lambda_a} by replacing $\f u(a)$ with an $\cal F$-measurable version. 
This allows for a use of \eqref{eq:magic_lemma} in a similar fashion as in the first part of Section~\ref{sec:sketch_approximate_eigenvalues} 
and reveals that, conditionally on $\cal F$, $\lambda(a)$ is concentrated around the second largest eigenvalue of $H^{(\cal V \cup B_r(b)\setminus \{ a\})}$ 
since $\abs{B_r(b)}$ is not too large due to our choice of $r$.

The proof of \ref{item:lambda_b} is much more elaborate. We start by noting that $\lambda(b)$ is characterized by the equation
\begin{equation} \label{eq:sketch_Schur_lambda_b} 
\lambda(b) + \frac{1}{d} \sum_{x,y \in S_1(b)} (H^{(\cal V)} - \lambda(b))^{-1}_{xy} = 0\,, 
\end{equation} 
as follows from Schur's complement formula. The main strategy is to derive a recursive family of equations for the Green function, starting from \eqref{eq:sketch_Schur_lambda_b} and extending to increasingly large spheres around $b$, to which Kesten's self-improving anticoncentration result can be applied.
To quantify anticoncentration, we use \emph{L\'evy's concentration function} 
\begin{equation} \label{eq:def_concentration_function_Q} 
 Q(X,L) \deq \sup_{t\in \mathbb{R}}\mathbb{P}(X\in[t-L,t+L])\,, 
\end{equation} 
where $X$ is a random variable and $L >0$ is deterministic.
\begin{proposition}[Theorem~2 of \cite{Kesten1969}] \label{prop:ConcentrationSum}  
There exists a universal constant $K$ such that for any independent random variables $X_1,\ldots, X_n$ satisfying $Q(X_i,L) \leq 1/2$ we have 
\begin{equation} \label{kesten_estimate}
Q\bigg(\sum_{i \in [n]} X_i,L\bigg) \leq \frac{K}{\sqrt{n}} \max_{i \in [n]} Q(X_i,L)\,. 
\end{equation}
\end{proposition}
This result is an improvement due to Kesten \cite{Kesten1969} of a classical anticoncentration result of Doeblin, Lévy, Kolmogorov, and Rogozin. Kesten's insight was that such an estimate can be made self-improving, as manifested by the factor $\max_{i \in [n]} Q(X_i,L)$ on the right-hand side.
This factor is crucial for our argument, as it allows us to successively improve the upper bound on $Q$.

We now explain more precisely how the expression of $\lambda(b)$ in terms of a large number of Green function entries is obtained. We shall tacitly use that $\mathbb{G}|_{B_{r}(b)}$ is a tree, which can be easily shown to be true with high probability.
Applying Schur's complement formula at $x\in S_i(b)$, using standard resolvent identities, and arguing similarly as in \eqref{eq:magic_lemma} to control errors yields
 \begin{equation}\label{eq:treeResolvent}
\frac{1}{G_{xx}(i-1,z)}= - z-\frac{1}{d}\sum_{y\in S_{1}^+(x)}G_{yy}(i,z)+o(1)\,, 
\end{equation} 
where $G_{xx}(i,z) \deq (H^{(B_{i}(b) \cup \cal V )}-z)_{xx}^{-1}$ and $S_1^+(x) = S_1(x) \cap S_{i+1}(b)$ is the set of children of $x$ in the tree $\mathbb{G}|_{B_r(x)}$ rooted at $b$. The error $o(1)$ is polynomially small in $N$; it comprises error terms arising from removing vertices from $H$ and neglecting all off-diagonal Green function entries.

Setting the error term in \eqref{eq:treeResolvent} to zero, we obtain a recursive equation for the idealized Green function entries $(g_x(z))_{x \in B_r(b)\setminus \{b\}}$, given by
\begin{equation} \label{eq:sketch_def_g} 
\frac{1}{g_{x}(z)} \deq \begin{cases} - z - \frac{1}{d} \sum_{y \in S_1^+(x)} G_{yy}(r,z) & \text{ if } x \in S_r(b) \\ 
- z-\frac{1}{d}\sum_{y\in S_1^+(x)}g_y(z)& \text{ if } x \in B_{r-1}(b) \setminus \{ b \}\,, \end{cases}  
\end{equation} 
which is an approximate version of the recursion \eqref{eq:treeResolvent} for the actual Green function entries. The recursion begins at the boundary of the ball $B_r(b)$ and propagates inwards.
We note that, for any $1 \leq i \leq r$, conditioned on $\cal F$, the family $(g_x(z))_{x \in S_i(b)}$ is independent if $\mathbb{G}|_{B_r(b)}$ 
is a tree and $z$ is $\cal F$-measurable. 
From \eqref{eq:treeResolvent}, \eqref{eq:sketch_def_g} and $r \asymp \frac{\log N}{\log d}$, it is not hard to conclude by induction 
that, for large enough $\eta >0$, with high probability, $g_x(z) = G_{xx}(0,z) + o(N^{-\eta}) = ( H^{(\cal V)}- z)^{-1}_{xx} + o(N^{-\eta})$ for all $x \in S_1(b)$.
Hence, if we can prove that $Q(g_x(z) , N^{-\eta}) \leq N^{-c}$ for all $x \in S_1(b)$ and some constant $c>0$, 
then 
a union bound over $a \neq b \in \cal W$, $\abs{\cal W} \asymp N^{\mu}$, 
the smallness of the off-diagonal entries of $(H^{(\cal V)} - \lambda(b))^{-1}$ as argued after \eqref{eq:treeResolvent}, 
and 
\eqref{eq:sketch_Schur_lambda_b} 
imply anticoncentration for $\lambda(b) - z$. This is \ref{item:lambda_b}, 
which together with \ref{item:lambda_a} implies that  
 $\min_{a\neq b \in \cal W} \abs{\lambda(a) - \lambda(b)} \geq N^{-\eta}$ with high probability, i.e.\ 
\ref{item:sketch_ii}.

Therefore, to complete the sketch of \ref{item:sketch_ii}, what remains is to prove $Q(g_x(z) , N^{-\eta}) \leq N^{-c}$ for all $x \in S_1(b)$, whose 
proof we sketch now. 
Throughout the entire argument we condition on $\cal F$ and use that $\mathbb{G}|_{B_r(x)}$ is a tree. To begin the recursion, we first show that $Q(g_x, d^{-1}) \leq 1/2$ for any $x \in S_r(b)$. This follows from the first case of \eqref{eq:sketch_def_g}, using a weak lower bound on the entries $G_{yy}(r,z)$ and anticoncentration from the fact that the size of $S_1^+(x)$ is a binomial random variable conditioned on $\cal F$.

Next, let $x \in S_i(b)$ for $1 \leq i \leq r - 1$.
By using the second case in \eqref{eq:sketch_def_g} and Proposition~\ref{prop:ConcentrationSum}, 
we iteratively refine the resolution, i.e.\ decrease the second argument of $Q$, and decrease the upper bound on $Q$, which are $d^{-1}$ and $1/2$, respectively, at the starting point. 
Indeed, conditioning on $A |_{B_i(b)}$, the second case in \eqref{eq:sketch_def_g} and rescaling the second argument of $Q$ yield
\begin{align} 
Q\big(g_x(z), (T^2d)^{-r + i} \big) &\leq Q\bigg( \sum_{y \in S_1^+(x)} g_y(z), (T^2d)^{-r + i +1}\bigg)
\notag \\ \label{eq:sketch_iteration_Q_without_robust}  
&\leq \frac{K}{\sqrt{\abs{S_1^+(x)}}} 
\max_{y \in S_1^+(x)} Q \big( g_y(z), (T^2d)^{-r + i + 1} \big)\,, 
\end{align} 
where we applied Proposition~\ref{prop:ConcentrationSum} using the independence of $(g_y(z))_{y \in S_1^+(x)}$ 
in the second step. 
Here, we also used that $Q(f(X), L) \leq Q (X, T^{-2}L)$ if $f(t) \deq \frac{1}{t}$ and $X \in [T^{-1},T]$ 
since the derivative of $f$ is bounded from below by $T^{-2}$ on this interval, 
and that, with high probability, $g$ lies in $[T^{-1}, T]$ for $T \asymp \sqrt{\frac{\log N}{d}}$. 

The estimate \eqref{eq:sketch_iteration_Q_without_robust} yields the desired self-improvement provided that $\abs{S_1^+(x)}$ is large enough.
However, $\abs{S_1^+(x)}$ is not large enough for all vertices $x$ in $B_r(b)$ (and in fact consistently applying \eqref{eq:sketch_iteration_Q_without_robust} at all vertices yields an anticoncentration bound at the root $b$ that is far from optimal and too weak to conclude localization). Sometimes, a better bound than \eqref{eq:sketch_iteration_Q_without_robust} can be obtained by replacing Kesten's estimate \eqref{kesten_estimate} with the trivial estimate
\begin{equation} \label{concentration_trivial}
Q\bigg(\sum_{i \in [n]} X_i,L\bigg) \leq \min_{i \in [n]} Q(X_i,L)\,,
\end{equation}
which follows immediately from the independence of the random variables $X_i$. Although this estimate lacks the factor $K / \sqrt{n}$ from \eqref{kesten_estimate}, it replaces the maximum with a minimum. Thus, an important ingredient in our recursive self-improving anticoncentration argument is an algorithm that determines which of \eqref{kesten_estimate} or \eqref{concentration_trivial} is to be used at any given vertex $x \in B_r(b)$. It relies on the notion of \emph{robust vertices}.

Recursively, a vertex $x \in B_{r}(b)$ is called \emph{robust} if $x \in S_r(b)$ or $S_1^+(x)$ contains at least $d/2$ robust vertices. 
We denote the set of robust vertices by $\cal R$. 
An important auxiliary result is that the root $b$ is robust with high probability, which in particular implies that 
$S_i(b) \cap \cal R$ is large for any $i \leq r$. 
Therefore, we can restrict to $x \in S_i(b) \cap \cal R$ and proceed similarly as in \eqref{eq:sketch_iteration_Q_without_robust} to obtain
\[ 
Q(g_x(z), (T^2d)^{-r +i  }) \leq Q\bigg( \sum_{y \in S_1^+(x)\cap \cal R} g_y(z), (T^2d)^{-r + i + 1}\bigg) \leq \frac{K\sqrt{2}}{\sqrt{d}} 
\max_{y \in S_1^+(x)\cap \cal R} Q \big( g_y(z), (T^2d)^{-r + i + 1} \big)\,, 
\] 
where in the first step we used \eqref{concentration_trivial}.
Thus, we obtain $ Q(g_x(z), (T^2d)^{-r}) \leq (K \sqrt{2d^{-1}})^{r-1} $ for all $x \in S_1(b)$ 
and, therefore, by choosing $r$ such that $(T^2d)^{r+1} = N^{\eta}$, we arrive at $Q(g_x(z), N^{-\eta}) \leq N^{-\eta/2 + o(1)}$ for all $x \in S_1(b)$. 
This is the desired anticoncentration bound, \ref{item:lambda_b}, which, as explained above, implies \ref{item:sketch_ii}.

\subsubsection{The three main exponents of $N$} \label{sec:exponents}
Throughout this paper we use the three exponents $\mu$, $\zeta$, $\eta>0$ to control three central quantities of the argument. We summarize their roles here for easy reference.
\begin{itemize}
\item $N^\mu$ is the typical size of the vertex sets $\cal V$ and $\cal W$ (cf.\ Proposition \ref{prop:graphProperties} \ref{item:size_cal_V}), 
as the parameter $\mu$ is introduced to control $\alpha^*$ (see \eqref{eq:largeDegree}). 
Consequently, $N^\mu$ is also the typical number of eigenvalues of $H$ satisfying \eqref{lambda_condition}. 
Therefore, the factor $N^\mu$ emerges from union bounds when a property is required for all $x \in \cal V$ or $x \in \cal W$. 
\item $N^{-\zeta}$ is the upper bound on the eigenvalue approximation that we establish in Proposition~\ref{prop:mainRigidity}, i.e.\ on the distance between $\lambda(x)$ and the eigenvalues of $H$ or, more precisely, the upper bound on $\norm{(H-\lambda(x))\f u(x)}$. Proposition~\ref{prop:mainRigidity} requires the condition $\zeta < 1/2-3\mu/2$. 
\item $N^{-\eta}$ is the lower bound on the eigenvalue spacing, or correspondingly $\min_{x \neq y \in \cal W} \abs{\lambda(x)-\lambda(y)}$, that we prove in Proposition~\ref{prop:MNonClustering}. Since the typical eigenvalue spacing is $N^{-\mu}$ in the interval we consider, we clearly need $\mu  < \eta$. 
In fact, our proof of Proposition~\ref{prop:MNonClustering} requires the stronger condition $8\mu < \eta$ for technical reasons 
as well as $\mu <1/24$. 
\end{itemize}
To apply perturbation theory in the proof of Theorem~\ref{thm:main}, we need that the error in the eigenvalue approximation $N^{-\zeta}$ be smaller
than the eigenvalue spacing $N^{-\eta}$. This means that $\eta < \zeta$.

\subsection{Localization profile -- proof of Theorem~\ref{prop:localization_profile}} 
\label{sec:localization_profile_proof_of_proposition}

Theorem~\ref{prop:localization_profile} is an immediate consequence of the following result.

\begin{proposition}[Local approximation for $\f u(x)$] 
\label{prop:IsolatedEigenvalue}
Let $\mu \in [0,1/3)$. 
Then, with high probability, the following holds for all $x\in \cal W$. 
If $r \in \N$ satisfies 
\begin{equation} \label{eq:condition_r2} 
\log d  \ll r \leq \min\bigg\{ \frac{1}{6} \frac{\log N}{\log d}, \min\bigg\{ \frac{1}{5} - \frac{\mu}{4}, \frac{1}{3} - \mu \bigg\} \frac{\log N}{\log d} - 2 \bigg\}\,,  
\end{equation} 
then 
\begin{enumerate}[label=(\roman*)] 
\item \label{item:u_x_approx_v_r_x}  $\f{u}(x) = \f{v}_r(x)+o(1)$,  
\item \label{item:u_x_approx_w_r_x}  $\f u(x) = \f w_r(x) + o(1)$. 
\end{enumerate} 
\end{proposition}

Part \ref{item:u_x_approx_v_r_x} is proved in Section~\ref{sec:properties_spectrum} below. 
Part \ref{item:u_x_approx_w_r_x} follows from Corollary~\ref{cor:comparison_u_x_w_x} below, since $\Lambda(\alpha_x) \geq 2$.

\begin{proof}[Proof of Theorem~\ref{prop:localization_profile}] 
From \eqref{eq:u_approx_u_x} in the proof of Theorem~\ref{thm:main}, we know that $\f w = \f u(x) + O(N^{-\eps})$ 
for a unique $x \in\cal W$ and some small enough $\eps>0$. 
Therefore, Theorem~\ref{prop:localization_profile} follows from Proposition~\ref{prop:IsolatedEigenvalue}. 
\end{proof}  

\subsection{Mobility edge -- proof of Theorem \ref{cor:mobility_edge}} \label{sec:mobility_edge}
Part \ref{item:localized} follows from Theorem~\ref{thm:main}, since 
$\inf\{ \alpha >0 \col \P(\alpha_1 \geq \alpha) \leq N^{-23/24} \} \leq 2$ if $d \geq \frac{23}{24} b_* \log N$, by Lemma~\ref{lem:Bennett} below. 
Part \ref{item:delocalized} was proved in \cite[Theorem~1.1 (ii)]{ADK_delocalized}.

\subsection{Localization length -- proof of Theorem \ref{thm:ll}} \label{sec:localization_length} 

For any eigenvector $\f w$ we introduce the function
\begin{equation*}
q(u) \deq \sum_{y} \r d(u,y) \scalar{\f 1_y}{\f w}^2\,,
\end{equation*}
so that $\ell(\f w) = \min_u q(u)$.

For the localized phase, suppose that $\lambda$ is an eigenvalue satisfying $\abs{\lambda} \geq 2 + \kappa$ with associated eigenvector $\f w$. Denote by $x$ the unique vertex associated with $\f w$ from Theorems \ref{thm:main} and \ref{prop:localization_profile}. Then the following estimates hold on the intersection of the high-probability events of these two theorems. By Theorem \ref{thm:main}, there exists a constant $R \equiv R_\kappa$ such that
\begin{equation} \label{ll_ux_far}
q(u) \geq \r d(u,x) \scalar{\f 1_x}{\f w}^2 \geq R \frac{\alpha_x - 2}{2 (\alpha_x - 1)} + o(1) \geq \frac{\alpha_x}{\alpha_x - 2} \qquad \text{for all } u \in B_R(x)^c\,,
\end{equation}
where the third inequality holds for large enough constant $R$.
Moreover, for any $\epsilon > 0$ there exists a constant $R' \in \N$ such that for $u \in B_R(x)$ we have
\begin{align}
q(u) &= \sum_{y \in B_{R'}(x)} \r d(u,y) \scalar{\f 1_y}{\f w}^2  + \sum_{y \in B_{R'}(x)^c} \r d(u,y) \scalar{\f 1_y}{\f w}^2
\notag \\
&= \sum_{y \in B_{R'}(x)} \r d(u,y) \scalar{\f 1_y}{\f w}^2 + O(\epsilon)
\notag \\  \label{qu_computation}
& = \sum_{y \in B_{R'}(x)} \r d(u,y) \scalar{\f 1_y}{\f v_r(x)}^2 + O(\epsilon) + o(1)\,,
\end{align}
where the second step follows from the estimate $\r d(u,y) \leq R + \r d(x,y)$ and the exponential decay from Theorem \ref{thm:main}, and the third step from Theorem \ref{prop:localization_profile} with some $\log d \ll r \ll \frac{\log N}{\log d}$.

To analyse the sum, we abbreviate $k \deq \r d(x,u)$ and introduce the set  $T_i(u,x) \deq S_i(x) \cap S_{i+k}(u)$ for $0 \leq i \leq r$,  which is the set of vertices in $S_i(x)$ whose geodesic to $x$ does not pass through $u$.
By Proposition \ref{prop:graphProperties} \ref{item:B_r_tree} below, the graph $\bb G \vert_{B_r(x)}$ is a tree, which implies
\begin{equation} \label{dux}
\r d(u,y) =
\begin{cases}
k + i & \text{if } y \in T_i(u,x)
\\
\abs{k - i} & \text{if } y \in S_i(x) \setminus T_i(u,x)\,.
\end{cases}
\end{equation}
Next, we estimate $\abs{S_i(x) \setminus T_i(u,x)}$. For $1 \leq i \leq k$, the set $S_i(x) \setminus T_i(u,x)$ consists of the unique vertex on the geodesic from $x$ to $u$ at distance $i$ from $x$. For $i > k$, we have 
$S_i(x) \setminus T_i(u,x) = S_i(x) \cap S_{i - k}(u)$. 
Here we used the tree structure of $\bb G \vert_{B_r(x)}$. Hence, we conclude that for $i > k$ we have $\abs{S_i(x) \setminus T_i(u,x)} \leq \abs{S_{i - k}(u)} \lesssim (\log N) d^{i - k - 1}$, by Proposition~\ref{prop:graphProperties_every_vertex}.

Next, by \cite[Lemma 5.4]{ADK20}, with high probability we have $\abs{S_i(x)} = \abs{S_1(x)} d^{i - 1} (1 + o(1))$ for $1 \leq i \leq r$. 
Since $\abs{S_1(x)} \geq \alpha^* d \geq 2d$, we conclude that $\abs{S_i(x)} \geq d^i$ for $i \leq r$. Putting all of these estimates together, we conclude for all $1 \leq i \leq r$ that
\begin{equation*}
\frac{\abs{S_i(x) \setminus T_i(u,x)}}{\abs{S_i(x)}} \lesssim
\begin{cases}
d^{-i} & \text{if } 1 \leq i \leq k
\\
(\log N) d^{-k-1} & \text{if } k < i \leq r\,.
\end{cases}
\end{equation*}
Since $T_i(u,x) = S_i(x)$ for $k = 0$ or $i = 0$, using the condition \eqref{eq:d_range} we conclude that
\begin{equation} \label{TS_estimate}
\frac{\abs{S_i(x) \setminus T_i(u,x)}}{\abs{S_i(x)}} = o(1)
\end{equation}
for all $u \in B_R(x)$ and $0 \leq i \leq r$.

Next, using \eqref{dux} and recalling the definition \eqref{eq:DefCandidate}, we write the above sum as
\begin{align*}
\sum_{y \in B_{R'}(x)} \r d(u,y) \scalar{\f 1_y}{\f v_r(x)}^2 &=
\sum_{i = 0}^{R'} \pBB{\sum_{y \in T_i(u,x)} (k+i) \scalar{\f 1_y}{\f v_r(x)}^2 + \sum_{y \in S_i(x) \setminus T_i(u,x)} \abs{k-i} \scalar{\f 1_y}{\f v_r(x)}^2}
\\
&= \sum_{i = 0}^{R'}\pbb{ (k+i) \, u_i(\alpha_x)^2 \frac{\abs{T_i(u,x)}}{\abs{S_i(x)}} + \abs{k - i} \, u_i(\alpha_x)^2 \frac{\abs{S_i(x)  \setminus T_i(u,x)}}{\abs{S_i(x)}}}
\\
&=\sum_{i = 0}^{R'} (k + i) u_i(\alpha_x)^2 + o(1)\,,
\end{align*}
where in the last step we used \eqref{TS_estimate}, the fact that $R'$ is constant, and $\sum_{i = 0}^{r-1} u_i(\alpha_x)^2 = 1$.
The latter sum is clearly minimized for $k = 0$.
Recalling \eqref{qu_computation}, we therefore conclude that for any $u \in B_{R}(x)$ we have
\begin{equation*}
q(u) \geq q(x) + O(\epsilon) + o(1)\,, \qquad q(x) = \sum_{i = 1}^\infty i u_i(\alpha_x)^2 + O(\epsilon) + o(1)
\end{equation*}
for large enough $R'$ depending on $\epsilon$.
Since $\epsilon > 0$ was arbitrary, and recalling \eqref{ll_ux_far}, by taking the minimum over $u \in [N]$, it therefore suffices to show that
\begin{equation} \label{ll_geom_series}
\sum_{i = 1}^\infty i u_i(\alpha_x)^2 = \frac{\alpha_x}{2(\alpha_x - 2)} + o(1) = \frac{\abs{\lambda}}{2 \sqrt{\lambda^2 - 4}} + o(1)\,.
\end{equation}
The first equality of \eqref{ll_geom_series} is an elementary computation using the definition \eqref{def_ui}, recalling the normalization $\sum_{i = 0}^{r - 1} u_i(\alpha)^2 = 1$ and that $r \gg 1$. The second equality of \eqref{ll_geom_series} follows from the estimate $\abs{\lambda} = \Lambda(\alpha_x) + o(1)$ by Remark \ref{rem:ev_Lambda}, which can be inverted to obtain $\frac{1}{\alpha_x} = \frac{1}{2}\pb{1 - \frac{\sqrt{\lambda^2 - 4}}{\abs{\lambda}}} + o(1)$.

For the delocalized phase, suppose that $\lambda$ is an eigenvalue satisfying $\abs{\lambda} \leq 2 - \kappa$ with associated eigenvector $\f w$. We use  \cite[Theorem~1.1 (ii)]{ADK_delocalized} to deduce that with probability $1 - O(N^{-10})$ we have $\norm{\f w}_\infty^2 \leq N^{-1 + o(1)}$. Hence, with probability $1 - O(N^{-9})$ we have, for any $x \in [N]$ and $r \geq 0$,
\begin{multline*}
q(x) \geq \sum_{y \in B_r(x)} \r d(x,y) \, \scalar{\f 1_y}{\f w}^2 + r \sum_{y \in B_r(x)^c} \scalar{\f 1_y}{\f w}^2
\\
\geq r \pbb{1 - \sum_{y \in B_r(x)} \scalar{\f 1_y}{\f w}^2} \geq r \pbb{1 - N^{-1 + o(1)} \abs{B_r(x)}}\,.
\end{multline*}
Next, we deduce from \cite[Lemma 1]{ChungLu_diameter} that for any constant $\epsilon > 0$ there is a constant $\delta > 0$ such that, with high probability, if $r \leq (1 - \epsilon) \frac{\log N}{ \log d}$ then $\abs{B_r(x)} \leq N^{1 - \delta}$ for all $x \in [N]$. Choosing $r \deq \floor{(1 - \epsilon) \frac{\log N}{ \log d}}$, we conclude that, for any $\epsilon > 0$, with high probability, for all eigenvectors $\f w$ with eigenvalue $\lambda$ satisfying $\abs{\lambda} \leq 2 - \kappa$, we have $\ell(\f w) \geq (1 - \epsilon + o(1)) \frac{\log N}{\log d}$.
Since $\epsilon > 0$ was an arbitrary constant, we conclude the stronger lower bound $\ell(\f w) \geq (1 - o(1)) \frac{\log N}{\log d} = \diam(\bb G) (1 + o(1))$, where the last step follows from \cite{ChungLu_diameter}. The complementary upper bound $\ell(\f w) \leq  \diam(\bb G)$ follows by definition of $\diam(\bb G)$, since $\r d(x,y) \leq \diam(\bb G)$ for all $x,y \in [N]$ (here we used that under our assumption on $d$, the graph $\bb G$ is with high probability connected). This concludes the proof.

\section{Preliminaries}  \label{sec:preliminaries}

The rest of this paper is devoted to the proofs of Propositions~\ref{prop:exponential_decay_u_x_Section_2}, 
\ref{prop:mainRigidity} and \ref{prop:MNonClustering}. We begin with a short section that collects some basic properties of the graph $\bb G$ and its spectrum.

\subsection{Properties of the graph} 
In this subsection, we collect some basic local properties of the Erd{\H o}s-R\'enyi graph $\bb G$ 
around vertices in $\cal V$.

\begin{proposition} 
\label{prop:graphProperties_every_vertex} 
Suppose that $\sqrt{\log N} \ll d \leq 3 \log N$.
With high probability, the following holds. 
\begin{enumerate}[label=(\roman*)] 
\item \label{item:upper_bound_alpha_x} 
$\max_{x\in [N]}|S_1(x)|\leq  10 \log N$. 
\item \label{item:upper_bound_size_B_r} 
$|B_i(x)|\lesssim \max\{\abs{S_1(x)},d\} d^{i-1}$ for all $x \in [N]$ and all $i \in \N$ with $i \leq \frac{1}{3}\frac{\log N}{\log d}$. 
\end{enumerate}
\end{proposition} 

Item \ref{item:upper_bound_alpha_x} is a simple application of Bennett's inequality and 
\ref{item:upper_bound_size_B_r} follows from \cite{ADK19,ADK20}; a detailed proof is given in Section~\ref{sec:proof_prop_graphProperties} below. In particular, by the assumption $d \gg \sqrt{\log N}$, we get from Proposition \ref{prop:graphProperties_every_vertex} \ref{item:upper_bound_alpha_x} that 
\begin{equation} \label{Lambda_alpha_bound}
\Lambda(\alpha_x) \lesssim \sqrt{\alpha_x} \lesssim  \sqrt{\frac{\log N}{d}} \ll \sqrt{d}
\end{equation}
for all $x \in [N]$ on the high-probability event from Proposition~\ref{prop:graphProperties_every_vertex}.

\begin{proposition}
\label{prop:graphProperties}
Let $\mu  \in [0,1/3)$ be a constant. 
Suppose that $\sqrt{\log N} \ll d \leq 3 \log N$.
With high probability, for any $r \in \N$ satisfying  
\begin{equation} \label{eq:condition_r} 
1\leq r \leq  \min\bigg\{\frac{1}{5}- \frac{\mu}{4},\frac{1}{3}-\mu \bigg\}\frac{\log N}{\log d}\,, 
\end{equation} 
the following holds. 
\begin{enumerate}[label=(\roman*)] 
\item \label{item:size_cal_V} $|\cal V|\leq N^{\mu+o(1)}$. 
\item \label{item:B_r_tree} $\mathbb{G}|_{B_{r}(x)}$ is a tree for all $x \in \cal V$. 
\label{enu:1-1}
\item \label{item:B_r_disjoint} $B_{r}(x)\cap B_{r}(y)=\emptyset$ for all $x, y\in{\cal V}$ satisfying $x \neq y$. 
\label{enu:2-1}
\item \label{item:upper_degree_in_balls_around_vertices_in_V} 
Let $\nu \in [0,1]$. If $1- \mu - \nu - r \frac{\log d}{\log N} \gtrsim 1$ then $\abs{S_1(y)} \leq \alpha^*(\nu) d $ for all $y \in \bigcup_{x \in  \cal V} (B_r(x)\setminus \{x \})$.  
\end{enumerate}
\end{proposition}

These statements are all consequences of \cite{ADK19,ADK21}; we explain the details of the proof of Proposition~\ref{prop:graphProperties} in Section~\ref{sec:proof_prop_graphProperties} below.

\subsection{Properties of the spectrum}\label{sec:properties_spectrum} 
In this subsection we collect basic spectral properties of $H$ and some of its submatrices.
 
\begin{definition}\label{def:w_1_eigenvector_H} 
Let $\f w_1$ be a normalized eigenvector of $H$ with nonnegative entries associated with the largest eigenvalue $\lambda_1(H)$ of $H$. 
\end{definition} 
Note that, with high probability, $\f w_1$ is unique and coincides with the Perron-Frobenius eigenvector of the 
giant component of $\mathbb{G}$.

\begin{proposition} \label{pro:spectral_gap} 
Suppose \eqref{eq:d_range}. With high probability the following holds.
\begin{enumerate}[label=(\roman*)]
\item \label{item:L2}
For any $\mu \in [0,1]$ we have $\max\{ \lambda_2(H^{(\cal V)}), -\lambda_N(H^{(\cal V)})\} \leq \Lambda(\alpha^*) + o(1)$.
\item \label{item:L1}
Fix $\mu \in [0,1/3)$.
If $X \subset \bigcup_{x\in \cal V}B_r(x)$ with $r \in \N$  as in \eqref{eq:condition_r},
then $\lambda_1(H^{(X)})$ and the corresponding eigenvector $\f w$ of $H^{(X)}$ satisfy  
\[ 
\lambda_1(H^{(X)}) = \sqrt{d}(1 + o(1))\,, \qquad \qquad  
\normbb{\f w - \frac{\f 1_{X^c}}{\abs{X^c}^{1/2}}} = o(1)\,.
\]  
\item \label{item:Candidate}
Fix $\mu \in [0,1)$. 
If $x \in \cal V$, $r \in \N$ satisfies $\log d \ll r \leq \frac{1}{6} \frac{\log N}{\log d}$
and $X \subset [N]$ satisfies $X \cap B_r(x) = \emptyset$, then 
\[ \norm{(H^{(X)}  - \Lambda(\alpha_x)) \f v_r (x) }  = o(1)\,,
\] 
where $\f v_r(x)$ was defined in \eqref{eq:DefCandidate}.
\item \label{item:spectral_gap_q} 
Fix $\mu \in [0,4/5)$. For any $r \in \N$ satisfying $r \ll \frac{d}{\log \log N}$, 
there is a normalized vector $\f q$ with $\supp \f q \subset \big( \bigcup_{x \in \cal V} B_{r+1}(x) \big)^c$ such that 
\[ \norm{(H - \sqrt{d})\f q } \lesssim d^{-1/2}\,, \qquad \norm{\f w_1 - \f q } \lesssim d^{-1}\,, \qquad \norm{\f q - N^{-1/2}\f 1_{[N]}} \lesssim d^{-1/2}\,. \]  
\end{enumerate} 
\end{proposition}

These results follow essentially from \cite{ADK19,ADK20,ADK21}; the detailed proof is 
presented in Section~\ref{sec:proof_pro_spectral_gap} below.

\begin{definition} \label{def:Omega}
We denote by $\Omega$ the intersection of the high-probability events of Propositions \ref{prop:graphProperties_every_vertex}, \ref{prop:graphProperties} and \ref{pro:spectral_gap}. 
\end{definition}

In particular, on $\Omega$ the estimate \eqref{Lambda_alpha_bound} holds.

\begin{corollary} \label{cor:largest_eigenvalues_vertices_removed} 
Fix $\mu \in [0,1/3)$. On $\Omega$, the following holds for all $x \in \cal W$ and all $X \subset [N]$.
If 
$X \cap B_r(x) = \emptyset$ and $\cal V \setminus \{x \} \subset X \subset \bigcup_{y\in \cal V}B_r(y)$ for some $r \in \N$ satisfying 
$\log d \ll r \leq \min \big\{ \frac{1}{6}, \frac{1}{5} - \frac{\mu}{4}, \frac{1}{3} - \mu \big\} \frac{\log N}{\log d}$
then 
\[ \lambda_1(H^{(X)}) = \sqrt{d} ( 1+ o(1))\,, \qquad \lambda_2(H^{(X)}) = \Lambda(\alpha_x) + o(1)\,, \qquad\lambda_3(H^{(X)}) \leq \Lambda(\alpha^*) + o(1)\,. \] 
\end{corollary} 

Corollary~\ref{cor:largest_eigenvalues_vertices_removed} with $X = \cal V \setminus \{x \}$ directly implies that, on $\Omega$, 
\begin{equation} \label{eq:lambda_x_approximated_by_Lambda_alpha_x} 
\lambda(x) = \lambda_2(H^{(\cal V\setminus \{ x\})}) = \Lambda(\alpha_x) + o(1)\,. 
\end{equation} 

\begin{proof}[Proof of Corollary \ref{cor:largest_eigenvalues_vertices_removed}] 
The statement about $\lambda_1(H^{(X)}) $ is identical to Proposition~\ref{pro:spectral_gap} \ref{item:L1}. By eigenvalue interlacing (Lemma \ref{lem:interlacing1}), we have
\[
\lambda_{3}(H^{(X)})\leq\lambda_{3}(H^{({\cal V}\setminus\{x\})})\leq\lambda_{2}(H^{({\cal V})})\leq\Lambda(\alpha^{*})+o(1)\,,
\] where for the last inequality we used Proposition~\ref{pro:spectral_gap} \ref{item:L2}.
Finally, Proposition \ref{pro:spectral_gap} \ref{item:Candidate} implies that there exists an eigenvalue of $H^{(X)}$ at distance $o(1)$ from $\Lambda(\alpha_x)$. Because of the estimates on $\lambda_{1}(H^{(X)})$ and $\lambda_{3}(H^{(X)})$ just proven, and because $\Lambda(\alpha_x) \geq \Lambda(\alpha^*) + \kappa/2$ for $x \in \cal W$ (recall \eqref{eq:def_cal_W}) as well as $\Lambda(\alpha_x) \ll \sqrt{d}$ by \eqref{Lambda_alpha_bound}, this eigenvalue has to be $\lambda_{2}(H^{(X)})$.
\end{proof}

\begin{proof}[Proof of Proposition~\ref{prop:IsolatedEigenvalue} \ref{item:u_x_approx_v_r_x}]
We show that the conclusion of Proposition \ref{prop:IsolatedEigenvalue} 
\ref{item:u_x_approx_v_r_x} holds on $\Omega$. 
The proof uses a spectral gap of $H^{(\cal V \setminus \{x \})}$ around $\lambda (x) = \lambda_2(H^{(\cal V \setminus \{x \})})$, that $\f v_r(x)$ is an approximate eigenvector 
of $H^{(\cal V \setminus \{ x\})}$ by Proposition \ref{pro:spectral_gap} \ref{item:Candidate}, and perturbation theory.
Indeed, from Corollary~\ref{cor:largest_eigenvalues_vertices_removed} with $X = \cal V \setminus \{ x\}$, recalling the definition of $\cal W$ (see \eqref{eq:def_cal_W}), we obtain that $\lambda_2(H^{(\cal V \setminus \{x \})})$ is separated from the other eigenvalues of $H^{(\cal V \setminus \{x \})}$ by a positive constant.
Owing to \eqref{eq:lambda_x_approximated_by_Lambda_alpha_x}, 
Proposition~\ref{pro:spectral_gap} \ref{item:Candidate} with $X = \cal V \setminus \{ x\}$ 
and Lemma~\ref{lem:perturbation_theory} below imply  $\norm{\f u(x) -\f v_r(x)} =o(1)$, i.e.\ 
Proposition~\ref{prop:IsolatedEigenvalue} \ref{item:u_x_approx_v_r_x}.
\end{proof}

We conclude the following corollary from Proposition~\ref{prop:IsolatedEigenvalue} \ref{item:u_x_approx_v_r_x} and its proof. 

\begin{corollary} \label{cor:lower_bound_u_x_in_x} 
Fix $\mu \in [0,1/3)$.
On $\Omega$ we have $\scalar{\f 1_x}{\f u(x)} = \sqrt{\frac{\alpha_x-2}{2(\alpha_x - 1)}} + o(1) \gtrsim 1$ for all $x \in \cal W$. 
\end{corollary} 

\begin{proof}
We first note that for any $r \to \infty$ as $N\to \infty$, we have $u_0(\alpha_x) = \sqrt{\frac{\alpha_x - 2}{2(\alpha_x - 1)}} + o(1)$ by \eqref{def_ui}.
By Proposition~\ref{prop:IsolatedEigenvalue} \ref{item:u_x_approx_v_r_x} and its proof, $\scalar{\f 1_x}{\f u(x)}  = \scalar{\f 1_x}{\f v_r(x)} + o(1) = u_0(\alpha_x) + o(1) \gtrsim 1$ on $\Omega$, where in second step we used the definition \eqref{eq:DefCandidate}, and in the last step we used that $\alpha_x \geq 2 + \kappa$, so that the sequence $(u_i(\alpha_x))_{i}$ from \eqref{def_ui} is exponentially decaying in $i$, uniformly in $r$.
\end{proof}

\subsection{Properties of the Green function} 
In this subsection, we fix $\mu \in [0,1/3)$.
Define
\begin{equation*}
\cal J = [\Lambda(\alpha^*) + \kappa / 4, \sqrt{d} /2]\,.
\end{equation*}
We shall use that whenever $z \in \cal J$, all Green functions appearing in our proof are bounded, which is the content of the following result.

\begin{lemma} \label{lem:G_bounded}
Suppose \eqref{eq:d_range} and \eqref{eq:condition_r}.
On $\Omega$, for $z,z' \in \cal J$ and $X \subset [N]$ satisfying $\cal V \subset X \subset \bigcup_{x \in \cal V} B_r(x)$, we have
\begin{align}
\norm{(H^{(X)} - z)^{-1}} & \leq {8} / \kappa \,, \label{eq:bound_resolvent} \\  
\norm{(H^{(X)} - z)^{-1}-(H^{(X)} - z’)^{-1}}& \leq (8 / \kappa)^2  \, |z-z’|\,.  \label{eq:G_error_z_zi} 
\end{align}
\end{lemma}

\begin{proof}[Proof of Lemma \ref{lem:G_bounded}]
Eigenvalue interlacing (Lemma \ref{lem:interlacing1}) and 
 Proposition \ref{pro:spectral_gap} \ref{item:L2} and \ref{item:L1} imply 
\[ \lambda_2(H^{(X)})\leq  \lambda_2(H^{(\cal V)}) \leq \Lambda(\alpha^*)+ \kappa/8 \qquad \text{and} \qquad \lambda_1(H^{(X)})= \sqrt{d}(1+o(1))\geq \sqrt{d}/2+ \kappa/8\,.  
\]
Therefore, $\text{dist}(z,\text{Spec}(H^{(X)}))\geq \kappa/8$ for any $z\in \cal J$, which proves \eqref{eq:bound_resolvent}.
The Lipschitz bound \eqref{eq:G_error_z_zi} follows from \eqref{eq:bound_resolvent} and the resolvent identity. 
\end{proof}

\begin{lemma} \label{lem:lowboundGreen}
Suppose \eqref{eq:d_range} and \eqref{eq:condition_r}.
On $\Omega$, for any $\cal V \subset X\subset \bigcup_{x\in \cal V}B_r(x)$, $z\in \cal J$, and  $y\notin X$, we have   \[-( H^{(X)}-z)^{-1}_{yy}\geq (3z)^{-1}\,.\]
\end{lemma}
\begin{proof}[Proof of Lemma \ref{lem:lowboundGreen}]
Denoting by $\lambda_1\geq \lambda_2\geq \dots $ and $\f w_1,\f w_2,\dots$ the eigenvalues and eigenvectors of $H^{(X)}$, respectively, we have 
\begin{equation} \label{eq:lower_bound_resolvent_proof} 
-(H^{(X)}-z)^{-1}_{yy} = \frac{\scalar{\f 1_y}{\f w_1}^2}{z- \lambda_1}+ \sum_{2\leq i\leq N}\frac{\scalar{\f 1_y}{\f w_i}^2}{z- \lambda_i}\geq \frac{1}{z- \lambda_N}(1-\scalar{\f 1_y}{\f w_1}^2)-\left|\frac{\scalar{\f 1_y}{\f w_1}^2}{ \lambda_1-z}\right|\,, 
\end{equation}
where in the second step we used $z - \lambda_N \geq z- \lambda_i > 0$ for all $i \geq 2$, 
which follows from eigenvalue interlacing (Lemma \ref{lem:interlacing1}), Proposition~\ref{pro:spectral_gap} \ref{item:L2}, and the condition $z \in \cal J$. 

To estimate the right-hand side of \eqref{eq:lower_bound_resolvent_proof}, we use Proposition~\ref{prop:graphProperties} \ref{item:size_cal_V} and Proposition~\ref{prop:graphProperties_every_vertex} \ref{item:upper_bound_size_B_r} as well as \eqref{eq:condition_r} to estimate $\abs{X} \leq \sum_{x \in \cal V} \abs{B_r(x)} \leq N^{1/3 + o(1)}$ on $\Omega$. From Proposition~\ref{pro:spectral_gap} \ref{item:L1} we therefore deduce that
\begin{equation} \label{eq:lower_bound_resolvent_proof2} 
\scalar{\f 1_y}{\f w_1}^2 = \abs{X^c}^{-1/2} + o(1) = o(1)\,.
\end{equation}
We conclude that the first term on the right-hand side of \eqref{eq:lower_bound_resolvent_proof} is bounded from below by $((2+o(1))z)^{-1}$, 
as $z - \lambda_N \leq 2z$ from Proposition~\ref{pro:spectral_gap} \ref{item:L2}. 
The second term on the right-hand side of \eqref{eq:lower_bound_resolvent_proof} 
is estimated using \eqref{eq:lower_bound_resolvent_proof2} as well as $\abs{\lambda_1 - z} \gtrsim \sqrt{d} \geq z$ 
by Proposition~\ref{pro:spectral_gap} \ref{item:L1} and $z \in \cal J$. 
\end{proof}

\section{Exponential decay of $\f u(x)$ and proof of Proposition~\ref{prop:exponential_decay_u_x_Section_2}} 
\label{sec:exponential_decay} 

In this section we establish the exponential decay of $\f u(x)$ around the vertex $x$. 
In particular, Proposition~\ref{prop:exponential_decay_u_x_Section_2} is a direct consequence of Proposition~\ref{prop:exponential_decay_u_x} below. 
Moreover, we prove in Corollary~\ref{cor:comparison_u_x_w_x} below that $\f u(x)$ is well approximated by $\f w_r(x)$, the eigenvector of $H|_{B_r(x)}$ 
corresponding to its largest eigenvalue. This implies 
Proposition~\ref{prop:IsolatedEigenvalue} \ref{item:u_x_approx_w_r_x}.

Throughout this section we use the high-probability event $\Omega$ from Definition \ref{def:Omega}.

\subsection{Simple exponential decay of $\f u(x)$}
In this subsection we establish exponential decay at some positive but not optimal rate.

\begin{proposition}
\label{prop:exponential_decay_u_x} 
Suppose that \eqref{eq:d_range} holds. Then there  is a constant $c \in (0,1)$ such that,
for each fixed $\mu \in [0,1/3)$, on $\Omega$, for each $x \in \cal W$ there exists $q_x > 0$ such that
\begin{equation} \label{eq:exponential_decay_u_x_bounded_normalized_degree} 
\norm{\f u(x)|_{B_i(x)^c}} \lesssim \sqrt{\alpha_x} \, q_x^{i}\,, \qquad \qquad q_x = \frac{\Lambda(\alpha^*(1/2)) + o(1)}{\lambda(x)}  \leq 1 - c\,,  
\end{equation}
for all $i \in \N$ satisfying $1 \leq i \leq \min\big\{ \frac{1}{5} - \frac{\mu}{4}, \frac{1}{3} - \mu \big\} \frac{\log N}{\log d}-2$. 
\end{proposition} 

\begin{proof} 
We note that $\supp \f u(x) \subset ([N] \setminus \cal V)\cup \{x \}$ and decompose $\f u(x) = \scalar{\f 1_x}{\f u(x)} \f 1_x + Q \f u(x)$, 
where $Q \deq \sum_{y \in [N] \setminus \cal V } \scalar{\f 1_y}{\,\cdot\,} \f 1_y$ is the orthogonal projection on the coordinates in $[N] \setminus \cal V$. 
We apply the projection $Q$ to the eigenvalue-eigenvector relation 
\[\lambda(x) \f u (x)=H^{(\cal V \setminus \{x \})} \f u(x) = \scalar{\f 1_x}{\f u(x)} \f 1_{S_1(x)}/\sqrt{d} + H^{(\cal V\setminus \{x \})} Q \f u(x)\,, \]  
solve for $Q \f u(x)$, and obtain 
\begin{equation} \label{eq:Q_u_x_resolvent_1_S_1} 
Q \f u(x) = \frac{\scalar{\f 1_x}{\f u(x)}}{\sqrt{d}} \big(\lambda(x) -  H^{(\cal V)} \big)^{-1} \f 1_{S_1(x)}\,, 
\end{equation} 
where we used that $H^{(\cal V \setminus \{x \})}  \f 1_x = \f 1_{S_1(x)}/\sqrt{d}$ and 
$Q H^{(\cal V \setminus \{x \})} Q = H^{(\cal V)}$. 
We also used that $\lambda(x)  - H^{(\cal V)}$ is invertible, which can be seen as follows. 
From Proposition~\ref{pro:spectral_gap} \ref{item:L1} with $X = \cal V$ and $r=0$, we conclude $\lambda_1(H^{(\cal V)}) = \sqrt{d}(1 + o(1))$ 
and, hence, 
\eqref{eq:lambda_x_approximated_by_Lambda_alpha_x} and  
Proposition~\ref{pro:spectral_gap} \ref{item:L2} yield
that $\lambda(x)$ is not an eigenvalue of $H^{(\cal V)}$ and 
\begin{equation} \label{eq:dist_lambda_x_spec_H_cal_V} 
\dist(\lambda(x), \spec H^{(\cal V)}) \gtrsim 1\,. 
\end{equation} 

Let
\begin{equation*}
Q_i \deq \sum_{y \in [N] \setminus (\cal V \cup  B_i(x))} \scalar{\f 1_y}{\,\cdot\,} \f 1_y
\end{equation*}
the projection onto 
the coordinates in $[N] \setminus (\cal V \cup B_i(x))$. As $\supp \f u(x) \subset ([N] \setminus \cal V)\cup \{x \}$ 
and $Q_i Q = Q_i$, 
we conclude from \eqref{eq:Q_u_x_resolvent_1_S_1} that 
\begin{equation} \label{eq:u_x_restricted_B_i_complement} 
\f u(x)|_{B_i(x)^c} = Q_i \f u(x) = Q_i Q \f u(x)  
= \frac{\scalar{\f 1_x}{\f u(x)}}{\lambda(x) \sqrt{d}} Q_i \bigg(1  -  \frac{H^{(\cal V)}}{\lambda(x)}  \bigg)^{-1} \f 1_{S_1(x)}\,. 
\end{equation}
For any $n \in \N$ we have
\begin{equation} \label{eq:resolvent_identity_telescopic} 
 \bigg(1  -  \frac{H^{(\cal V)}}{\lambda(x)}  \bigg)^{-1} = \sum_{k=0}^n \bigg(\frac{H^{(\cal V)}}{\lambda(x)}  \bigg)^{k} + \bigg(1  -  \frac{H^{(\cal V)}}{\lambda(x)}  \bigg)^{-1} \bigg( \frac{H^{(\cal V)}}{\lambda(x)}  \bigg)^{n+1}\,. 
\end{equation} 
Since $H^{(\cal V)}$ is a local operator, we conclude that $Q_i (H^{(\cal V)})^k \f 1_{S_1(x)} = 0$ 
if $k + 1 \leq i$. 
Hence, fixing $i \geq 1$ and applying \eqref{eq:resolvent_identity_telescopic} with $n = i - 1$ 
to \eqref{eq:u_x_restricted_B_i_complement}, we get 
\begin{equation} \label{eq:representation_u_x_outside_B_i} 
 \f u(x)|_{B_i(x)^c} = \frac{\scalar{\f 1_x}{\f u(x)}}{\sqrt{d}} Q_i \big(\lambda(x)  -  H^{(\cal V)}\big)^{-1} \bigg( \frac{H^{(\cal V)}}{\lambda(x)}  \bigg)^{i} \f 1_{S_1(x)} = O \bigg( \frac{\norm{(H^{(\cal V)})^{i} \f 1_{S_1(x)}}}{\sqrt{d}\lambda(x)^{i}} \bigg)\,. 
\end{equation} 
Here, in the last step, we employed 
$\abs{\scalar{\f 1_x}{\f u(x)}} \leq 1$, 
$\norm{Q_i} \leq 1$, $\norm{(\lambda(x)  - H^{(\cal V)})^{-1}} \lesssim 1$ by \eqref{eq:dist_lambda_x_spec_H_cal_V}. 

Let $r \in \N$ be the largest integer satisfying \eqref{eq:condition_r}.  We now claim that, on $\Omega$, 
\begin{equation}\label{eq:norm_H_cal_V_support_in_B_r} 
 \norm{H^{(\cal V)}\f v} \leq (\Lambda(\alpha^*(1/2)) + o(1))\norm{\f v} 
\end{equation} 
for any $\f v$ such that $\supp \f v \subset B_r(x)$. Before proving \eqref{eq:norm_H_cal_V_support_in_B_r}, we use it conclude the proof of \eqref{eq:exponential_decay_u_x_bounded_normalized_degree}.

For $i \leq r - 2$ we have $\supp (H^{(\cal V)})^{i} \f 1_{S_1(x)} \subset B_r(x)$, so that iterative applications of \eqref{eq:norm_H_cal_V_support_in_B_r}
yield 
\begin{equation} \label{eq:norm_H_power_i} 
\norm{(H^{(\cal V)})^{i} \f 1_{S_1(x)}} \lesssim (\Lambda(\alpha^*(1/2)) + o(1))^{i} \sqrt{\alpha_x} \sqrt{d}\,.
\end{equation} 
Hence, the first bound in \eqref{eq:exponential_decay_u_x_bounded_normalized_degree} follows from \eqref{eq:representation_u_x_outside_B_i}. 
We now show the second bound in \eqref{eq:exponential_decay_u_x_bounded_normalized_degree}. 
If $\frac{\log N}{d} \geq T$ for a sufficiently large constant $T$ then $q_x \leq 1 - c$ for some constant $c > 0$ 
by using Corollary~\ref{cor:Lambda_alpha_1_2_Lambda_alpha_mu} with $\nu = 1/2$, $\lambda(x) \geq \Lambda(\alpha^*(\mu))$, $\mu \leq 1/3$ and possibly increasing its $T$. 
If $\frac{\log N}{d} < T \lesssim 1$ then $\alpha_x \lesssim 1$ by Proposition~\ref{prop:graphProperties_every_vertex} \ref{item:upper_bound_alpha_x} 
and the estimate $q_x \leq 1- c$ for some constant $c > 0$ follows from 
$\Lambda(\alpha_x) \leq 2 \sqrt{\alpha_x} \lesssim 1$ and $\lambda(x) = \Lambda(\alpha_x) + o(1) \geq 
\Lambda(\alpha^*) + \kappa/4 \geq \Lambda(\alpha^*(1/2)) + \kappa/4$, 
which is a consequence of \eqref{eq:lambda_x_approximated_by_Lambda_alpha_x} and the definition of $\cal W$.  
This completes the proof of the second bound in~\eqref{eq:exponential_decay_u_x_bounded_normalized_degree}.

What remains is the proof of \eqref{eq:norm_H_cal_V_support_in_B_r}. 
For any $\f v \in \R^{[N]}$,  the Cauchy-Schwarz inequality implies 
$\norm{(\E H)^{(\cal V)}\f v} \leq \sqrt{\frac{d \abs{\supp \f v}}{N}}\norm{\f v}$. By Proposition~\ref{prop:graphProperties_every_vertex} \ref{item:upper_bound_size_B_r} and \eqref{eq:condition_r}, if $\supp \f v \subset B_{r}(x)$ 
then $\abs{\supp \f v} \leq N^{1/5+o(1)}$. Therefore, 
\begin{multline*}  
\norm{H^{(\cal V)} \f v} \leq \norm{ (H - \E H)^{(\cal V)} \f v} + o(\norm{\f v})  \\ 
\leq \norm{ (H - \E H)^{(\{ y \,:\, \alpha_y \geq \alpha^*(1/2) \})} \f v} + o(\norm{\f v})
\leq (\Lambda(\alpha^*(1/2)) + o(1))\norm{\f v}\,, 
\end{multline*}  
where in the second step, we used 
$\cal V^c \cap B_r(x) \subset \{ y \col \alpha_y \geq \alpha^*(1/2) \} ^c$ by Proposition~\ref{prop:graphProperties} \ref{item:upper_degree_in_balls_around_vertices_in_V} 
and, in the third step, 
$\norm{(H - \E H)^{(\{y \,:\, \alpha_y \geq \alpha^*(1/2)\})}} \leq \Lambda(\alpha^*(1/2)) + o(1)$ by 
Lemma~\ref{lem:bound_A_minus_EA}\footnote{Note that Lemma~\ref{lem:bound_A_minus_EA} holds on $\Omega$, since Lemma~\ref{lem:bound_A_minus_EA} is used to prove Proposition~\ref{pro:spectral_gap} \ref{item:L2}.}. 
This concludes the proof of \eqref{eq:norm_H_cal_V_support_in_B_r}.
\end{proof}

\subsection{Approximating $\f u(x)$ by $\f w_r(x)$} 

From Proposition~\ref{prop:exponential_decay_u_x} and its proof, we deduce the following result, which compares $\f u(x)$ and $\f w_r(x)$ from Definition~\ref{def:w_x_and_v_x}.

\begin{corollary} \label{cor:comparison_u_x_w_x} 
Suppose that \eqref{eq:d_range} holds and fix $\mu \in [0,1/3)$. 
Then, on $\Omega$, for all $x \in \cal W$ and $r \in \N$ satisfying \eqref{eq:condition_r2}, 
\[ 
 \f u(x) = \f w_r(x) + O\big(\alpha_x q_x^{r} \Lambda(\alpha_x)^{-1} \big)
= \f w_r(x) + o\big(\Lambda(\alpha_x)^{-1} \big) \,,  
\] 
where $q_x$ is the same as in Proposition~\ref{prop:exponential_decay_u_x}. 
\end{corollary} 

\begin{proof} 
We shall apply Lemma~\ref{lem:perturbation_theory} with $M = H|_{B_r(x)}$, $\wh{\lambda} = \lambda(x)$ 
and $\f v = \f u(x)$. First, we check its conditions by studying the spectral gap 
of $H|_{B_r(x)}$ around its largest eigenvalue. 
As $r \leq \frac{1}{6} \frac{\log N}{\log d}$, we conclude from 
 Proposition~\ref{pro:spectral_gap} \ref{item:Candidate} with $X = B_r(x)^c$  
 that $H|_{B_r(x)} = H^{(B_r(x)^c)}$ has an eigenvalue $\Lambda(\alpha_x) + o(1)$. 
 The bound \eqref{eq:norm_H_cal_V_support_in_B_r} implies that $\norm{H|_{B_r(x) \setminus \{ x\}}} \leq \norm{H^{(\cal V)}}
\leq \Lambda(\alpha^*(1/2)) + o(1)$, where in the first step we used Proposition \ref{prop:graphProperties} \ref{enu:2-1}. By eigenvalue interlacing (Lemma \ref{lem:interlacing1}), we therefore deduce that
\[ \lambda_1(H|_{B_r(x)}) = \Lambda(\alpha_x) + o(1)\,, \qquad \qquad \lambda_2(H|_{B_r(x)}) \leq \lambda_1(H|_{B_r(x) \setminus \{ x\}}) 
\leq \Lambda(\alpha^*(1/2)) + o(1)\,. \] 
Thus, the definition of $\cal W$ in \eqref{eq:def_cal_W} implies
\begin{equation*}
\lambda_1(H|_{B_r(x)}) - \lambda_2(H|_{B_r(x)}) \gtrsim \Lambda(\alpha_x) \big( 1 - \frac{\Lambda(\alpha^*(1/2))}{\Lambda(\alpha^*(\mu)) +\kappa/2}\big) \gtrsim \Lambda(\alpha_x)\,,
\end{equation*}
where the last inequality follows from Corollary~\ref{cor:Lambda_alpha_1_2_Lambda_alpha_mu} if $\frac{\log N}{d} \geq T$ for some large enough constant $T$, and from $\Lambda(\alpha^*(1/2)) + \kappa/2 
\leq \Lambda(\alpha^*(\mu)) + \kappa/2 \leq \Lambda(\alpha_x) \lesssim 1$ (see \eqref{Lambda_alpha_bound}) otherwise. 
Hence, owing to \eqref{eq:lambda_x_approximated_by_Lambda_alpha_x}, 
 there is $\Delta \gtrsim \Lambda(\alpha_x)$ such that $H|_{B_r(x)}$ has precisely one eigenvalue in $[\lambda(x) - \Delta, 
\lambda(x) + \Delta]$.  

Let $P_r$ be the orthogonal projection onto the coordinates in $B_r(x)$. 
The eigenvalue-eigenvector relation $(H^{(\cal V \setminus \{x \})} - \lambda(x)) \f u(x) = 0$
and 
$P_r H^{(\cal V\setminus \{x \})} P_r = H|_{B_r(x)}$
imply 
\begin{equation*}
(H|_{B_r(x)} - \lambda(x)) \f u (x) = - P_r H^{(\cal V\setminus \{x \})} (1 - P_r)\f u (x) 
- \lambda(x) (1 - P_r) \f u(x)\,.
\end{equation*}
Therefore, since $(1 - P_r)\f u(x) = \f u(x)|_{B_r(x)^c}$, we get
\begin{equation} \label{eq:approximate_eigenvector_wh_u_r_x} 
\norm{(H|_{B_r(x)} - \lambda(x)) \f u (x)} \leq \lambda(x) \norm{\f u(x)|_{B_r(x)^c}} 
+ \norm{P_r H^{(\cal V\setminus \{x \})} (\f u(x)|_{B_r(x)^c})} \lesssim \alpha_x q_x^{r}\,. 
\end{equation} 
In the last step, we used $\lambda(x) \lesssim \alpha_x^{1/2}$ and Proposition~\ref{prop:exponential_decay_u_x} to estimate the first term. 
For the second term, we used that $P_r H^{(\cal V\setminus \{x \})} (\f u(x)|_{B_r(x)^c}) 
= P_r H^{(\cal V\setminus \{x \})} ( \f u(x)|_{S_{r+1}(x)})$ by the locality of $H^{(\cal V\setminus \{x \})}$ as well as the identity
$H^{(\cal V\setminus \{x \})} (\f u(x)|_{S_{r+1}(x)}) = H^{(\cal V)} (\f u(x)|_{S_{r+1}(x)})$, which yield
\[
\norm{H^{(\cal V\setminus \{x \})} (\f u(x)|_{S_{r+1}(x)})} \lesssim \Lambda(\alpha^*(1/2)) \norm{\f u(x)|_{S_{r+1}(x)}} \lesssim \alpha_x^{1/2} \norm{\f u(x)|_{B_r(x)^c}}  \lesssim \alpha_x q_x^{r} 
\] 
due to \eqref{eq:norm_H_cal_V_support_in_B_r} as $r + 1 \leq \min\big\{\frac{1}{5} - \frac{\mu}{4} , 
 \frac{1}{3} - \mu \big\} \frac{\log N}{\log d} $,   
$\Lambda(\alpha^*(1/2))^2 \leq \Lambda(\alpha^*)^2 \lesssim \alpha_x$, and Proposition~\ref{prop:exponential_decay_u_x}.  

Since $\alpha_x \lesssim \frac{\log N}{d}$ by Proposition~\ref{prop:graphProperties_every_vertex} \ref{item:upper_bound_alpha_x},  
 $q_x \leq 1 - c $  for some constant $c   > 0$ by Proposition~\ref{prop:exponential_decay_u_x}, 
$r \gg \log d \gtrsim \log \big(\frac{10\log N}{d} \big)$ by \eqref{eq:d_range} 
and $\alpha_x \geq 2$, we obtain 
\begin{equation} \label{eq:approximate_eigenvectors_control_error_term}  
\alpha_x q_x^r \ll 1 \lesssim \alpha_x^{1/2}\,.
\end{equation} 
Therefore, owing to \eqref{eq:approximate_eigenvector_wh_u_r_x} and 
$\alpha_x^{1/2} \asymp \Lambda(\alpha_x) \lesssim \Delta$, 
we can now apply Lemma~\ref{lem:perturbation_theory} with $M = H|_{B_r(x)}$, $\wh{\lambda}   = \lambda(x)$ 
and $\f v = \f u(x)$  
to conclude $\f u(x) = \f w_r(x) + O ( \alpha_x q_x^{r} \Lambda(\alpha_x)^{-1})$ from \eqref{eq:approximate_eigenvector_wh_u_r_x} and $\Delta \gtrsim \Lambda(\alpha_x)$ as well as $\scalar{\f 1_x}{\f u(x)} \geq 0$ and 
$\scalar{\f 1_x}{\f w_r(x)} \geq 0$.   
This proves the first equality in Corollary~\ref{cor:comparison_u_x_w_x}.

Since $\alpha_x q_x^r \ll 1$ by \eqref{eq:approximate_eigenvectors_control_error_term}, 
the last equality in Corollary~\ref{cor:comparison_u_x_w_x} 
follows immediately.
\end{proof}

\subsection{Optimal exponential decay of $\f u(x)$} 

In this subsection we establish an explicit rate of exponential decay of $\f u(x)$. It holds only for a smaller set of eigenvectors near the spectral edge, requiring $\mu$ to be small enough. As pointed out in Remark \ref{rem:optimal_decay}, up to the error $O(\epsilon)$, this rate is optimal. 

\begin{proposition}\label{prop:improved_exponential_decay} 
Suppose \eqref{eq:d_range}. Then the following holds.
\begin{enumerate}[label=(\roman*)]
\item \label{item:improved_exp_decay_sub} (Subcritical regime) There are constants $T \geq 1$ and $c>0$ such that 
if $\frac{\log N}{d} \geq T$ then, on $\Omega$, for any small enough constant  $\epsilon > 0$  and 
each $\mu \in [0,\eps]$, 
\[ \norm{\f u(x)|_{B_i(x)^c}} \lesssim \sqrt{\alpha_x} \bigg( \frac{1 + O(\eps)}{\sqrt{\alpha_x - 1}} \bigg)^{i+1} \] 
for all $x \in \cal W$ and $i \leq c \frac{\log N}{(\log d) \log\frac{10\log N}{d} }$.
\item \label{item:improved_exp_decay_crit} 
(Critical regime) 
There is a constant $c > 0$ such that, for any constants $T  \geq 1$ and $\epsilon > 0$ with $\eps \leq c \kappa$, 
if $\frac{\log N}{d} \leq T$ and $\mu \in [0,c T^{-1} \eps^2]$ then, with high probability, 
\[ \norm{ \f u(x)|_{B_i(x)^c}} \lesssim \sqrt{\alpha_x} \bigg(\frac{1 + O(\eps)}{\sqrt{\alpha_x-1}}\bigg)^{i+1}\]  
for all $x \in \cal W$ and $i \leq \frac{c \eps^2}{T} \frac{\log N}{\log d}$. 
\end{enumerate} 
\end{proposition}

\begin{proof} 
Let $x \in \cal W$. 
For any $i$, $n \in \N$, we conclude from \eqref{eq:u_x_restricted_B_i_complement} 
and \eqref{eq:resolvent_identity_telescopic} (see also \eqref{eq:representation_u_x_outside_B_i}) that 
\[ \norm{\f u(x)|_{B_i(x)^c}} \lesssim \frac{1}{\lambda(x)\sqrt{d}} \normbb{Q_i \sum_{k=0}^n \bigg( \frac{H^{(\cal V)}}{\lambda(x)} \bigg)^k 
\f 1_{S_1(x)}} + \frac{\norm{(H^{(\cal V)})^{n+1} \f 1_{S_1(x)}}}{\sqrt{d} \lambda(x)^{n+1}} \] 
since 
$\abs{\scalar{\f 1_x}{\f u(x)}} \lesssim 1$, $\norm{Q_i} \leq 1$ and $\norm{(\lambda(x) - H^{(\cal V)})^{-1}} 
\lesssim 1$ by \eqref{eq:dist_lambda_x_spec_H_cal_V}.  
We denote the right-hand side of \eqref{eq:condition_r} by $r$. In the following, we always assume that 
$n + 1\leq r$ and tacitly use the graph properties listed in Proposition~\ref{prop:graphProperties}. 
As in the proof of Proposition~\ref{prop:exponential_decay_u_x} (see \eqref{eq:norm_H_power_i} and use $q_x \leq 1 - \eps$ by \eqref{eq:exponential_decay_u_x_bounded_normalized_degree}),
we find some constant $\epsilon > 0$ such that if $n + 1 \leq r $ then
\begin{equation}\label{eq:proof_improved_exp_decay_error_term}
\frac{\norm{(H^{(\cal V)})^{n+1} \f 1_{S_1(x)}}}{\sqrt{d} \lambda(x)^{n+1}} \lesssim \sqrt{\alpha_x} (1- \eps)^{n+1}\,. 
\end{equation}
We shall use the following result, whose proof is given at the end of this subsection.

\begin{claim} \label{lem:exponential_decay_main_term} 
Suppose that all vertices in $B_{r - 1}(x)\setminus \{x \}$ have degree at most $\tau d$, where $2 \leq 2 \sqrt{\tau} < \lambda(x)$. Then for all $n \leq r - 1$ and $i \in \N$ we have on $\Omega$
\begin{equation}\label{eq:bound_Q_i_sum_k_to_n_A_power_k_1_S_1} 
\frac{1}{\lambda(x)\sqrt{d}} \normbb{ Q_i \sum_{k = 0}^n \bigg( \frac{H^{(\cal V)}}{\lambda(x)} \bigg)^k \f 1_{S_1(x)} } 
\lesssim \sqrt{\alpha_x} \bigg( \frac{2}{\lambda(x)+ \sqrt{\lambda(x)^2 - 4\tau}} \bigg)^{i + 1}\,. 
\end{equation} 
\end{claim}

We now explain how we choose $\tau$ and $n$ in the subcritical and the critical regimes in order to deduce 
Proposition~\ref{prop:improved_exponential_decay} from Claim~\ref{lem:exponential_decay_main_term}. 
In the subcritical regime, i.e.\ for the proof of \ref{item:improved_exp_decay_sub}, we 
choose $\tau \defeq \alpha^*(1-3\eps)$ and $n \leq \eps \frac{\log N}{\log d}$ for some small enough constant $\eps > 0$.
By arguing similarly as in the proof of Corollary~\ref{cor:Lambda_alpha_1_2_Lambda_alpha_mu} below, we 
find a constant $T \geq 1$ such that 
\begin{equation} \label{eq:proof_improved_exp_decay_ratio_tau_lambda_sub} 
\frac{4 \tau}{\lambda(x)^2} \leq \frac{4 \alpha^*(1-\eps) (1 + o(1))}{\Lambda(\alpha^*(\mu))^2} = O(\eps) < 1
\end{equation}
by \eqref{eq:lambda_x_approximated_by_Lambda_alpha_x} and \eqref{eq:def_cal_W} if $\frac{\log N}{d} \geq T$. 
Here, the last inequality holds if $\eps$ is sufficiently small. Hence, the assumption on $\tau$ in Claim \ref{lem:exponential_decay_main_term} holds.
Moreover, with our choice of $\tau$, the degrees in $B_n(x)\setminus\{x\}$ are bounded by $\tau d$ 
due to Proposition~\ref{prop:graphProperties} \ref{item:upper_degree_in_balls_around_vertices_in_V} 
and our assumption on $n$. 
Hence, the assumptions of Claim \ref{lem:exponential_decay_main_term} hold.
As $\lambda(x)\lesssim \sqrt{\frac{\log N}{d}}$ by \eqref{eq:lambda_x_approximated_by_Lambda_alpha_x}  
and Proposition~\ref{prop:graphProperties_every_vertex} \ref{item:upper_bound_alpha_x}, 
 the right-hand side of \eqref{eq:proof_improved_exp_decay_error_term} 
is bounded by the right-hand side of \eqref{eq:bound_Q_i_sum_k_to_n_A_power_k_1_S_1} 
if $n = i C \log \big(\frac{10\log N}{d} \big)$ for some sufficiently large constant $C$. Since we need that $n +1 \leq r$, this yields the upper bound on $i$ in \ref{item:improved_exp_decay_sub}. 
From \eqref{eq:proof_improved_exp_decay_ratio_tau_lambda_sub} and \eqref{eq:lambda_x_approximated_by_Lambda_alpha_x}, we deduce that 
$\lambda(x) + \sqrt{\lambda(x)^2 - 4\tau} = \lambda(x) \big( 1 + \sqrt{1 - \frac{4\tau}{\lambda(x)^2}}\big)
= \Lambda(\alpha_x) ( 1 +o(1))(2 + O(\eps)) \geq 2\sqrt{\alpha_x - 1} (1 + O(\eps))^{-1}$, 
which completes the proof of \ref{item:improved_exp_decay_sub}.

For the proof of \ref{item:improved_exp_decay_crit}, we note that $\mu \leq c \eps^2 T^{-1} < 1/3$ for a sufficiently small constant $c>0$ since $T \geq 1$, $\eps \leq c \kappa$ and $\kappa \leq 1$. 
We choose $\tau = 1 + \eps$ for some constant $\eps >0$. 
As $\lambda(x)^2 - 4 \tau = ( \lambda(x)^2 - 4)(1 - \frac{4\eps}{\lambda(x)^2 - 4}) > 0$ 
due to \eqref{eq:lambda_x_approximated_by_Lambda_alpha_x} and \eqref{eq:def_cal_W} 
if $\eps \leq c \kappa$ with a small enough constant $c > 0$. This establishes the first condition for \eqref{eq:bound_Q_i_sum_k_to_n_A_power_k_1_S_1}. 
We set $R \defeq \frac{c\eps^2}{T} \frac{\log N}{\log d}$ and conclude from
Bennett's inequality (Lemma \ref{lem:Bennett} below) and Proposition~\ref{prop:graphProperties} 
that
\[ 
\P\Big( \alpha_y \leq \tau \text{ for all } x \in \cal W \text{ and for all } y \in B_{R}(x)\setminus \{x \} \Big) 
\leq \exp \bigg( \bigg( \mu + R \frac{\log d}{\log N} - \frac{d}{\log N} h(\eps)\bigg)\log N \bigg)\,,
\]  
where $h(\epsilon) \deq (1 + \epsilon) \log (1 + \epsilon) - \epsilon$. 
Since $h(\eps) \geq 3c \eps^2$ for some small enough constant $c >0$, the upper bound on $\mu$ imposed in the 
statement, the definition of $R$ and $\frac{d}{\log N} \geq \frac{1}{T}$ imply that the factor in front of $\log N$ is negative. 
Therefore, with high probability, we can apply \eqref{eq:bound_Q_i_sum_k_to_n_A_power_k_1_S_1} simultaneously 
for all $x \in \cal W$. 
We choose $n = C ( i + 1)$  
and deduce  
that if $C > 0$ is a large enough constant then 
the error in \eqref{eq:proof_improved_exp_decay_error_term} is dominated by the right-hand side of \eqref{eq:bound_Q_i_sum_k_to_n_A_power_k_1_S_1} 
as $\lambda(x) \lesssim 1$ for $\frac{\log N}{d} \leq T$. 
We recall that $r$ denotes the right-hand side of \eqref{eq:condition_r} and note that the condition of \eqref{eq:proof_improved_exp_decay_error_term}, 
$n + 1 = C (i + 1) +1 \leq r$,  
can be satisfied by possibly decreasing the constant $c > 0$. 
Finally, we obtain $\lambda(x) + \sqrt{\lambda(x)^2 - 4\tau} = \lambda(x) + \sqrt{\lambda(x)^2 - 4\tau} 
+ O( \frac{\eps}{\sqrt{\lambda(x)^2 - 4}} ) = 2 \sqrt{\alpha_x - 1} ( 1 + O(\eps))^{-1}$ 
similarly as argued above and in the proof of \ref{item:improved_exp_decay_sub} using $\eps \leq c \kappa$ 
and $\kappa \leq 1$. 
This proves \ref{item:improved_exp_decay_crit} and, thus, Proposition~\ref{prop:improved_exponential_decay}.
\end{proof} 

\begin{proof}[Proof of Claim~\ref{lem:exponential_decay_main_term}]
For $y \in S_1(x)$ we denote by $B_n^+(y)$ the ball of radius $n$ around $y$ in the graph $\bb G |_{[N] \setminus \{x\}}$, and we write $S_n^+(y) \deq B_n^+(y) \setminus B_{n-1}^+(y)$.
Since $\mathbb{G}|_{B_{n+1}(x) \setminus \cal V}$ is a forest on $\Omega$ and for each $y \in S_1(x)$, $\mathbb{G}|_{B_{n}^+(y)}$ is a tree, we obtain
\[ (A^{(\cal V)})^k \f 1_{S_1(x)} = \sum_{y \in S_1(x)} \sum_{j=1}^{k+1} \sum_{z \in S_{j-1}^+(y)} \f 1_z N_k(y,z)\,, \] 
where $N_k(y,z)$ denotes the number of walks on $\mathbb{G}|_{B_{n}^+(y)}$ of length $k$ between $y$ and 
$z$.

Thus, for $i \geq 0$, we conclude 
\begin{align*}  
Q_i \sum_{k=0}^n \frac{1}{\lambda(x)^k d^{k/2}} (A^{(\cal V)})^k  \f 1_{S_1(x)}
&=   Q_i \sum_{y \in S_1(x)} \sum_{j=1}^{n+1} \sum_{z \in S_{j-1}^+(y)} \sum_{k= j-1}^n\f 1_z \frac{N_k(y,z)}{\lambda(x)^k d^{k/2}} \\ 
&=  \sum_{y \in S_1(x)} \sum_{j=i + 1}^{n+1} \sum_{z \in S_{j-1}^+(y)}\f 1_z \bigg(  \sum_{k= j-1}^n 
\frac{N_k(y,z)}{\lambda(x)^k d^{k/2}}  \bigg)\,,
\end{align*} 
where the second step follows from the definition of $Q_i$. Since the sets $\{B_n^+(y) \col y \in S_1(x)\}$, are disjoint, we conclude that
\begin{equation} \label{eq:norm_squared_power_H} 
\normbb{Q_i \sum_{k = 0}^n \bigg( \frac{H^{(\cal V)}}{\lambda(x)} \bigg)^k \f 1_{S_1(x)}}^2 
= \sum_{y \in S_1(x)} \sum_{j=i + 1}^{n+1} \sum_{z \in S_{j-1}^+(y)} \bigg( \sum_{k=j-1}^n 
\frac{N_k(y,z)}{\lambda(x)^k d^{k/2}} \bigg)^2\,.  
\end{equation}

Next, let $z \in S_{j-1}^+(y)$ for some $y \in S_1(x)$. 
We note that  $N_k(y,z) = 0$ if $k - (j-1)$ is odd due to the bipartite structure of a tree. 
If this difference is even, then we now show that 
\begin{equation} \label{eq:number_paths} 
N_k(y,z) \leq (M^k)_{1j} (\tau d)^{(k - (j-1))/2}\,, 
\end{equation}  
where $M$ is the adjacency matrix of $\N^*$ (regarded as a graph where consecutive numbers are adjacent), see~\eqref{eq:def_M_adjacency_matrix_N} below for a precise definition, 
under the assumption that the degree of each vertex in $\mathbb{G}|_{B_{n}^+(y)}$ is bounded by $\tau d$. 

For the proof of \eqref{eq:number_paths}, we introduce the set of walks on $\N^*$
\begin{equation*}
W_k (j) \deq \{ \text{walks } \gamma \text{ on } \N^* \text{ of length } k \text{ such that } \gamma(0) = 1 
\text{ and } \gamma(k) = j\}\,,
\end{equation*}
and for each $\gamma \in W_k(j)$ we introduce the set of walks on $\mathbb{G}|_{B_{n}^+(y)}$ that project down to $\gamma$,
\begin{equation*}
 W_k(y,z;\gamma) \deq  \{ \text{walks }\Gamma \text{ on } \mathbb{G}|_{B_{n}^+(y)} \text{ of length } k \text{ such that } \Gamma(0) = y, \Gamma(k) = z \text{ and } \r d(x, \Gamma(\,\cdot\,)) 
= \gamma \} \,.
\end{equation*}
By definition, any walk $\Gamma \in  W_k(y,z;\gamma)$ projects down to a walk $\gamma \in W_k(j)$, which implies
\begin{equation} \label{Nk_W}
N_k(y,z) \leq \sum_{\gamma \in W_k(j)} \abs{W_k(y,z;\gamma)}\,. 
\end{equation}
In order to prove \eqref{eq:number_paths}, we fix $\gamma \in W_k(j)$ and estimate $\abs{W_k(y,z;\gamma)}$. 
For $i \in [j]$, let $T_i \defeq \max\{ t \in \{0, \dots, k\} \col \gamma(t) = i \}$, i.e.\ $T_i$ is the last time 
when $\gamma$ hits $i$. Clearly, $T_1 < T_2 < \ldots < T_{j-1} < T_j = k$.  See Figure \ref{fig:paths} for an illustration of the walks $\gamma$ and $\Gamma$. At each time $T_i$ with $i \in [j-1]$, the walk $\gamma$ takes a step to the right, and by definition of the times $T_i$, any walk $\Gamma \in W_k(y,z;\gamma)$ takes a step outwards on the geodesic from $y$ to $z$. This means that $j - 1$ of the $k$ steps of $\Gamma$ are fixed. Of the remaining $k - (j - 1)$ steps, half correspond to steps to the left of $\gamma$, which again correspond to a uniquely determined step of $\Gamma$ along the unique path back towards $y$. Hence, only $(k - (j - 1))/2$ of the $k$ steps of $\Gamma$ are free to choose. Since the degrees in $\mathbb{G}|_{B_{n}^+(y)}$ are bounded by $\tau d$, we obtain  
$\abs{W_k(y,z; \gamma)} \leq (\tau d)^{(k-(j-1))/2}$. Plugging this estimate into \eqref{Nk_W} implies \eqref{eq:number_paths}, since $\abs{W_k(j)} = (M^k)_{1j}$ by definition of $M$.

\begin{figure}[!ht]
\begin{center}
{\footnotesize 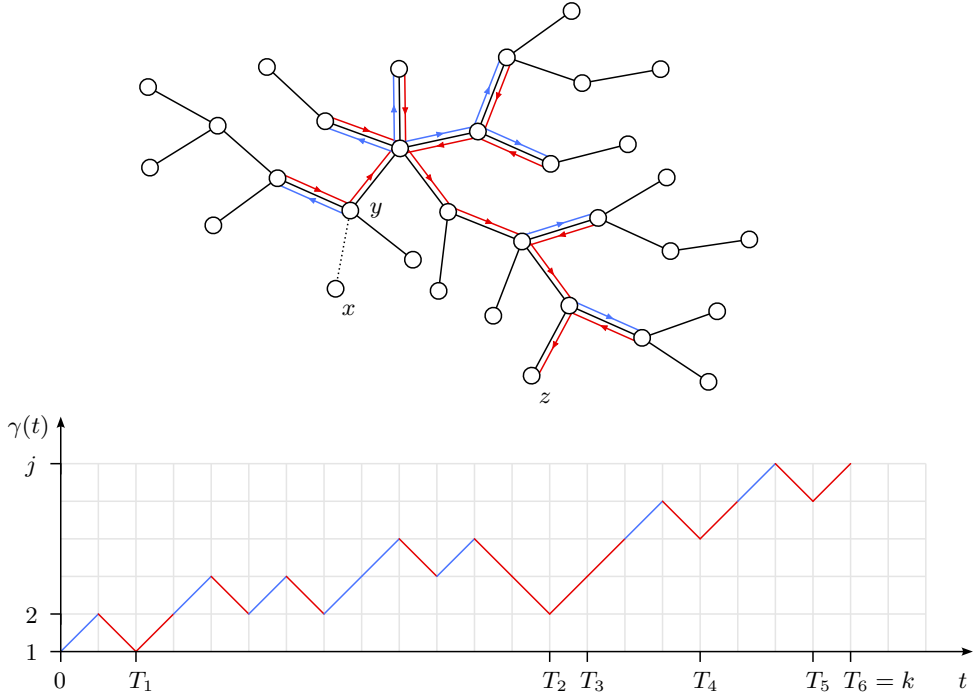}
\end{center}
\caption{An illustration of two walks $\gamma \in W_k(j)$ (bottom) and $\Gamma \in W_k(y,z;\gamma)$ (top). Here $k = 21$ and $j = 6$. By definition of $W_k(y,z;\gamma)$, for each $t \in \{0, \dots, k\}$ we have $\gamma(t) = \r d(x, \Gamma(t))$. The time $T_i$ is the last time when $\gamma$ hits $i$. We draw an edge $\{\Gamma(t), \Gamma(t+1)\}$ in red if the choice of the vertex $\Gamma(t+1)$ is uniquely determined by $\gamma$ and in blue otherwise. In the latter case, there are at most $\tau d$ possible choices for $\Gamma(t+1)$. Red edges arise in two ways: (i) a step to the right in $\gamma$ following a time $T_i$ ($j - 1$ in total, in $\Gamma$ corresponding to a step towards $z$ along the geodesic from $y$ to $z$); (ii) a step to the left in $\gamma$ ($(k - (j - 1))/2$ in total, in $\Gamma$ corresponding to a step towards $y$).
\label{fig:paths}}
\end{figure}

Next, applying \eqref{eq:number_paths} to \eqref{eq:norm_squared_power_H}, we obtain 
\begin{align*}  
\normbb{Q_i \sum_{k = 0}^n \bigg( \frac{H^{(\cal V)}}{\lambda(x)} \bigg)^k \f 1_{S_1(x)}}^2
 &\leq \sum_{j = i + 1}^{n+1} \abs{S_j(x)} \bigg( \sum_{k = j -1}^n \frac{(M^k)_{1j} (\tau d)^{(k - (j-1))/2}}{\lambda(x)^k d^{k/2}} \bigg)^2 \\ 
 &\leq \sum_{j = i + 1}^{n+1} \abs{S_j(x)} \frac{1}{(\tau d)^{j - 1}} 
\bigg( \sum_{k = j - 1}^n \frac{(M^k)_{1j} \tau^{k/2}}{\lambda(x)^k} \bigg) ^2 \\ 
&\lesssim \sum_{j = i + 1}^{n+1} \abs{S_1(x)} \frac{1}{\tau^{j-1}}  
\bigg( \sum_{k = 0}^\infty \frac{(M^k)_{1j} \tau^{k/2}}{\lambda(x)^k} \bigg) ^2 \\ 
&= \sum_{j = i + 1}^{n+1} \abs{S_1(x)} \frac{1}{\tau^{j-1}}  
\bigg(\bigg( 1 - \frac{\sqrt{\tau}}{\lambda(x)} M\bigg)^{-1}_{1j} \bigg)^2\\ 
&=  \sum_{j = i+1}^{n+1} \abs{S_1(x)} \lambda(x)^2 \bigg(\frac{2}{ \lambda(x) +  \sqrt{\lambda(x)^2 - 4\tau}} \bigg)^{2j} 
\\ 
&\lesssim \abs{S_1(x)} \lambda(x)^2  \bigg( \frac{2}{\lambda(x)+ \sqrt{\lambda(x)^2 - 4\tau}} \bigg)^{2 (i + 1)}\,. 
\end{align*}  
Here, in the third step, we used that $\abs{S_j(x)} \lesssim \abs{S_1(x)} d^{j-1}$ by Proposition~\ref{prop:graphProperties_every_vertex} \ref{item:upper_bound_size_B_r}
and that $(M^k)_{1j} \geq 0$ for all $k \in \N$. 
The fourth step follows from the condition $2\sqrt{\tau}/\lambda(x) < 1$, the invertibility and the Neumann series representation of $M$ in Lemma~\ref{lem:resolvent_adjacency_matrix_positive_integers} below with $t = \lambda(x) / \sqrt{\tau}$. 
The fifth step is a consequence of the representation of $( 1 - \frac{\sqrt{\tau}}{\lambda(x)}M)^{-1}_{1j}$ 
in Lemma~\ref{lem:resolvent_adjacency_matrix_positive_integers} below. 
In the last step, we used that $\lambda(x) = \Lambda(\alpha_x) + o(1) \geq 2 + \kappa/4$ to sum up the geometric series and conclude that the series is $\lesssim 1$. 
This completes the proof of \eqref{eq:bound_Q_i_sum_k_to_n_A_power_k_1_S_1} and, thus, the one of 
 Claim~\ref{lem:exponential_decay_main_term}. 
\end{proof}

\section{Approximate eigenvalues -- proof of Proposition~\ref{prop:mainRigidity}} \label{sec:rigidity}

In this section we prove Proposition~\ref{prop:mainRigidity}, by showing that in the interval $\cal I$ there is a one-to-one correspondence between eigenvalues of $H$ and the points $\lambda(x)$ for $x \in \cal W$, up to a polynomially small error term.

\subsection{Proof of Proposition~\ref{prop:mainRigidity}} 
\label{sec:proof_Prop_StrongRigidityLambda}

We recall $\f w_1$ from Definition~\ref{def:w_1_eigenvector_H}.

\begin{definition}  \label{def:Pi}
Let $\Pi$ be the orthogonal projection onto $\Span\big(\{\f w_1\} \cup \{\f u(x) \col x\in{ \cal W} \}\big)$ and $\ol{\Pi} \deq 1 - \Pi$. 
\end{definition} 

Note that the set $\{\f w_1\} \cup \{\f u(x) \col x\in{\cal W} \}$ in the definition of $\Pi$ is not orthogonal. Throughout the proof, we regard $H$ as a block matrix associated with the orthogonal sum decomposition $\ran\Pi \oplus (\ran\Pi)^\perp$. 

\begin{proposition}
\label{prop:StrongRigidity} 
Suppose \eqref{eq:d_range}. 
Fix $\mu \in [0,1/3)$ and $\zeta \in [0,1/2 - \mu)$.
With high probability, the following holds. 
\begin{enumerate}[label=(\roman*)] 
\item \label{item:StrongRigidity_approx_eigenvector} 
$\norm{(H- \lambda(x))\f u(x)} \leq N^{-\zeta}$ for all $x \in \cal W$. 
\item \label{item:StrongRigidity_upper_block} 
If $\zeta < 1/2 - 3 \mu/2$ then 
\[\spec(\Pi H \Pi) \setminus \{ 0 \} =\{\lambda_1(H)\}\cup\{{\lambda}(x) +\epsilon_x \col x\in \cal W\}\] 
counted with multiplicity, where $|\epsilon_x|\leq N^{-\zeta}$ for all $x\in \cal W$ and $\lambda_1(H)=\sqrt{d}(1+o(1))$. 
\item \label{item:StrongRigidity_offdiag_blocks} 
If $\zeta < 1/2 - 3 \mu/2$ then 
$\|\overline{\Pi} H \Pi\|\leq N^{-\zeta}$.  
\item \label{item:StrongRigidity_lower_block} 
If $\mu < 1/4$ then $\lambda_1(\overline{\Pi} H \overline{\Pi})\leq \Lambda(\alpha^*) + \kappa/2 +o(1)$.  
\end{enumerate}
\end{proposition}

\begin{proof}[Proof of Proposition~\ref{prop:mainRigidity}] 
Owing to the block decomposition $H = \Pi H \Pi + \ol{\Pi}H \ol{\Pi} + \ol{\Pi} H \Pi + \Pi H \ol{\Pi}$, 
Proposition~\ref{prop:StrongRigidity} \ref{item:StrongRigidity_offdiag_blocks} yields
\[
\spec(H) \setminus \{0\}  = \big\{ \lambda+\epsilon_\lambda \col \lambda \in (\spec(\Pi H \Pi) \cup \spec(\overline{\Pi} H \overline{\Pi}) ) \big\} \setminus \{0\}
\]
counted with multiplicities, where $|\epsilon_\lambda|\leq 2 \|\overline{\Pi} H \Pi\|\leq N^{-\zeta}$ for all $\lambda$. 
Therefore, Proposition~\ref{prop:mainRigidity} follows from the definition of $\cal I$ as well as Proposition~\ref{prop:StrongRigidity} \ref{item:StrongRigidity_upper_block}, \ref{item:StrongRigidity_lower_block}, and \ref{item:StrongRigidity_approx_eigenvector}. 
\end{proof}

The rest of this section is devoted to the proof of Proposition \ref{prop:StrongRigidity}. We assume the condition \eqref{eq:d_range} throughout. We recall the definition of the high probability event  $\Omega$  from Definition~\ref{def:Omega}. For any event $A$ and random variable $X$ we write
\begin{equation} \label{def_POmega}
\P_\Omega(A) \deq \P(\Omega \cap A)\,, \qquad \E_\Omega[X] \deq \E[X \ind{\Omega}]\,.
\end{equation}

\subsection{Proof of Proposition \ref{prop:StrongRigidity} \ref{item:StrongRigidity_approx_eigenvector}}

The proof of Proposition \ref{prop:StrongRigidity} \ref{item:StrongRigidity_approx_eigenvector} relies on the following result, whose proof is given at the end of this subsection.

\begin{proposition} \label{prop:EstimateExy}
Let $\mu \in [0,1/3)$.
For any $x\in{\cal W}$, we have the decomposition  
\begin{equation} \label{eq:represenation_error_H_u_x} 
(H-\lambda(x))\f u(x)=\sum_{y\in{\cal V} \setminus \{x \}}\epsilon_{y}(x) \f 1_{y}\,.
\end{equation}
Moreover, for any $x,y\in[N]$, we have the estimate 
\[
\mathbb{E}_{\Omega}[\ind{x \in \cal W} \ind{y\in{\cal V \setminus \{x \} }}\,\epsilon_{y}(x)^{2}]\leq d^{-1}(10 \log N)^{2}N^{2\mu - 3}\,.
\]
\end{proposition}

From Proposition~\ref{prop:EstimateExy}, for any $x\in [N]$, we obtain 
\begin{equation} \label{eq:expectation_norm_H_minus_lambda_u} 
\mathbb{E}_{\Omega} \qb{\ind{x\in \cal W}\|(H-\lambda(x))\f u(x)\|^2} \leq d^{-1}(10 \log N)^{2} N^{2\mu-2}\,.
\end{equation}

\begin{proof}[Proof of Proposition~\ref{prop:StrongRigidity} \ref{item:StrongRigidity_approx_eigenvector}] 
From \eqref{eq:expectation_norm_H_minus_lambda_u}, a union bound, and Chebyshev's inequality, 
 we conclude that 
\[
\mathbb{P} \pb{\exists x\in \cal W,\, \|(H-\lambda(x))\f u(x)\| > N^{-\zeta}} \leq   \mathbb{P}(\Omega^c) + N^{2\zeta+2\mu-1+o(1)}
\]
for any $\zeta >0$. 
This proves Proposition~\ref{prop:StrongRigidity} \ref{item:StrongRigidity_approx_eigenvector} since $\zeta < 1/2 -\mu$ by assumption, and $\P(\Omega^c) = o(1)$ by Propositions~\ref{prop:graphProperties_every_vertex}, 
\ref{prop:graphProperties} and \ref{pro:spectral_gap}.  
\end{proof} 

We shall need modifications of the sets $\cal V$ and $\cal W$ defined in \eqref{eq:def_cal_V} and \eqref{eq:def_cal_W}, respectively. For $X \subset [N]$ we define
\begin{align} \label{eq:def_mathcal_V_X} 
\cal V^{(X)} &\deq \hbb{ y \in [N]\setminus X  \col \frac{\abs{S_1(x) \setminus X}}{d} \geq \alpha^*(\mu)} \,,
\\ \label{eq:def_mathcal_W_X}
\cal W^{(X)} &\deq \hbb{ y \in \cal V^{(X)}  \col \Lambda\pbb{\frac{\abs{S_1(x) \setminus X}}{d}} \geq  \Lambda(\alpha^*(\mu)) + \kappa / 2} \,.
\end{align}
The point of these definitions is that $\cal V^{(X)}$ and $\cal W^{(X)}$ depend only on the edges in $\bb G \vert_{X^c}$.
The following remark states that on the event $\Omega$ the effect of the upper index in these definitions amounts simply to excluding vertices.

\begin{remark} \label{rem:cal_V_without_y_cal_W_without_y} 
On ${\Omega}\cap \{ y \in \mathcal V \}$, owing to Proposition~\ref{prop:graphProperties} \ref{item:B_r_disjoint}, 
we have $\cal V^{(y)} = \cal V \setminus \{ y\}$ and $\cal W^{(y)} = \cal W \setminus \{ y\}$. 
\end{remark}

\begin{proof}[Proof of Proposition \ref{prop:EstimateExy}]
Note that $H - H^{(\cal V \setminus \{x \})}$ is $d^{-1/2}$ times the adjacency matrix of the subgraph of $\bb G$ containing the edges incident to $\cal V \setminus \{x\}$.
Hence, because $(H^{({\cal V}\setminus\{x\})}-\lambda(x))\f u(x)=0$ and $\supp \f u(x) \subset (\cal V \setminus \{x \})^c$, we find
\begin{equation} \label{eps_splitting}
(H-\lambda(x))\f u(x) = (H - H^{(\cal V \setminus \{x \})}) \f u (x)  = \sum_{y\in{\cal V} \setminus \{x \}}\epsilon_{y}(x)\f 1_{y}\,,
\qquad \epsilon_{y}(x)\deq \frac{1}{\sqrt{d}}\sum_{t\in S_{1}(y)}\scalar{\f 1_t}{\f u(x)}\,.
\end{equation}
See Figure \ref{fig:epsilon} for an illustration.
\begin{figure}[!ht]
\begin{center}
{\footnotesize 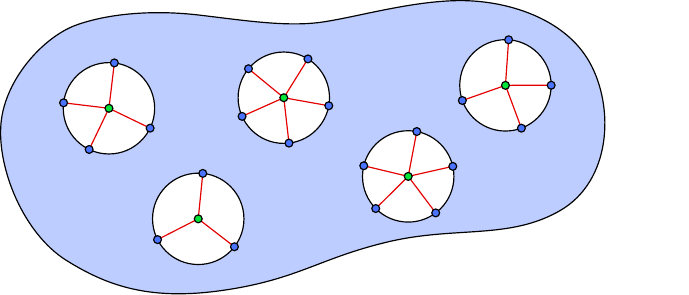}
\end{center}
\caption{An illustration of the identity \eqref{eps_splitting}. The set $\cal V \setminus \{x\}$ is drawn in green and its complement $(\cal V \setminus \{x\})^c$ in blue. The blue neighbours of the green vertices are drawn explicitly, and the remaining blue vertices are represented by the shaded blue region. The edges of incident to $\cal V \setminus \{x\}$ are drawn in red. The adjacency matrix of these red edges is $A - A^{(\cal V \setminus \{x\})}$. The vector $\f u(x)$ is supported on the blue vertices, and hence the vector $(A - A^{(\cal V \setminus \{x\})}) \f u(x)$ is supported on the green vertices. Its value at a green vertex $y$ equals the sum of the entries of $\f u(x)$ at the blue vertices adjacent to $y$.
\label{fig:epsilon}}
\end{figure}
The Cauchy-Schwarz inequality implies 
\begin{equation}
\label{eq:EstimateExy}
\epsilon_{y}(x)^{2}\leq \frac{|S_1(y)|}{d}\sum_{t\in S_{1}(y)}\scalar{\f 1_t}{\f u(x)}^{2}
\leq \frac{10 \log N}{d} \sum_{t \in S_1(y)} \scalar{\f 1_t}{\f u(x)}^2 
\end{equation}
on ${\Omega}$ due to Proposition~\ref{prop:graphProperties_every_vertex} \ref{item:upper_bound_alpha_x}. 
From now on we fix $x,y\in[N]$. 
Moreover, let 
 $\tilde{\f u }(x)$ be the eigenvector of $H^{(\cal V^{(y)} \cup \{ y \} \setminus \{ x \})}$ associated 
with its second largest eigenvalue $\lambda_2(H^{(\cal V^{(y)} \cup \{ y \} \setminus \{ x \})})$ and satisfying $\scalar{\f 1_x}{\tilde{\f u}(x)} \geq 0$.
On ${\Omega}\cap \{ y \in \mathcal V \}$, we have ${\cal V}={\cal V}^{(y)} \cup \{ y\}$ by Remark~\ref{rem:cal_V_without_y_cal_W_without_y} and, thus, $\f u(x)=\tilde{\f u}(x)$.  As $x \neq y$, on $\Omega \cap \{ y \in \cal V\}$, the events $\{ x \in \cal W \}$ and $\{ x \in \cal W^{(y)} \}$ coincide by Remark~\ref{rem:cal_V_without_y_cal_W_without_y}. 
Therefore, since $\cal W^{(y)}$ and $\wt{\f u}(x)$ are $\sigma(A^{(y)})$-measurable, we obtain 
\begin{align}
&\mspace{-20mu}\mathbb{E}_{{\Omega}}\bigg[\ind{x\in{\cal W}}\ind{y\in \cal V \setminus \{x \}}\sum_{t\in S_{1}(y)}\scalar{\f 1_t}{\f u(x)}^{2}\bigg] 
\nonumber \\
& \leq \mathbb{E}\bigg[\ind{x\in{{\cal W}^{(y)}}} \sum_{t \in [N]} \scalar{\f 1_t}{\tilde{\f u}(x)}^{2} \, \mathbb{E}\big[\ind{
\alpha^*(\mu) d \leq \abs{S_1(y)} \leq 10 \log N}\ind{t\in S_{1}(y)}\big| A^{(y)} \big]\bigg] \nonumber \\
& \leq \P ( x \in \cal W^{(y)}) \max_{t \in [N]} \P \big( t \in S_1(y), \, \alpha^*(\mu) d \leq \abs{S_1(y)} \leq 10 \log N \big) \nonumber \\ 
 & \leq \mathbb{P}(x\in{\cal W}^{(y)})\sum_{ \alpha^*(\mu)d\leq k\leq 10 \log N}\mathbb{P}(\abs{S_1(y)}=k)\frac{k}{N} \nonumber \\
  & \leq \mathbb{P}(x\in{\cal V})\mathbb{P}(y\in \cal V)\frac{10 \log N}{N}\,. \label{eq:bound_individual_term_eps_y_x} 
\end{align}
Here, in the first step, in addition, we spelled out the condition $y \in \cal V$ as $\abs{S_1(y)} \geq \alpha^*(\mu)d$, used that $\abs{S_1(y)} \leq 10 \log N$ on $\Omega$ 
by Proposition~\ref{prop:graphProperties_every_vertex} \ref{item:upper_bound_alpha_x} 
and then dropped the indicator function $\ind{\Omega}$. 
The second step follows from the independence of the event $\{ \alpha^*(\mu) d \leq \abs{S_1(y)} \leq 10 \log N, t\in S_{1}(y) \}$ from $\sigma(A^{(y)})$ and 
 $\sum_{t\in [N]}\scalar{\f 1_t}{\tilde{\f u}(x)}^{2}=\|\tilde{\f u}(x)\|^{2}=1$. 
(For these two steps, see also \eqref{eq:magic_lemma}.)  
For the third step, we conditioned on $\abs{S_1(y)}$ and used that if $\abs{S_1(x)} = k$ then $t$ lies in a uniformly 
distributed subset of $[N]\setminus \{y \}$ with $k$ elements. For the last step, we used that $\cal W^{(y)} \subset \cal W \subset \cal V$.

Finally, applying \eqref{eq:bound_individual_term_eps_y_x} 
 to \eqref{eq:EstimateExy} and using the estimate $\P(x \in \cal V) \leq N^{\mu - 1}$ (by the definitions \eqref{eq:def_cal_V} and \eqref{eq:largeDegree}) concludes the proof of Proposition~\ref{prop:EstimateExy}.  
\end{proof}

\subsection{Proof of Proposition \ref{prop:StrongRigidity} \ref{item:StrongRigidity_upper_block}, \ref{item:StrongRigidity_offdiag_blocks}} 

In this section, we conclude Proposition~\ref{prop:StrongRigidity} \ref{item:StrongRigidity_upper_block} 
and \ref{item:StrongRigidity_offdiag_blocks} from the following result, which is also proved in this section.

\begin{definition} \label{def:u_perp}
Order the elements of $\cal W$ in some arbitrary fashion, and denote by $(\f u^{\perp}(x))_{x\in \cal W}$ the Gram-Schmidt orthonormalization of $(\f u(x))_{x\in \cal W}$.
\end{definition}

\begin{proposition} \label{prop:RigidityUPerp}
Let $\mu \in [0,1/3)$. Then the following holds with high probability. 
For any $x\in \cal W$, we have 
\begin{equation} \label{eq:H-lambda}
\|(H-\lambda(x))\f u^{\perp}(x)\|  \lesssim N^{\mu-1/2 + o(1)}\,.
\end{equation}
More generally, denoting $\cal D = \sum_{x \in \cal W} \lambda(x) \f u^\perp (x) (\f u^\perp(x))^*$ we have 
\begin{equation} \label{eq:H_approx_diagonal_u_basis} 
\|(H-\cal D)\f u\| \lesssim  N^{3\mu/2-1/2 + o(1)}\|\f u\|
\end{equation} 
for all $\f u\in \Span\{ \f u(x) \col x\in \cal W\}$.
\end{proposition}

\begin{proof}[Proof of Proposition \ref{prop:StrongRigidity} \ref{item:StrongRigidity_upper_block}]
By Definition \ref{def:Pi}, $\f w_1$ is an eigenvector of $\Pi H \Pi$ with eigenvalue $\lambda_1(H)$. 
Let $\zeta  < 1/2 - 3\mu / 2$. 
From \eqref{eq:H_approx_diagonal_u_basis}, we conclude that, for each $x \in \cal W$, there is $\eps_x \in [-N^{-\zeta}, N^{-\zeta}]$  
such that $\{ \lambda(x) + \eps_x \col x \in \cal W\} \subset \spec(\Pi H \Pi)\setminus \{ 0\}$ counted with multiplicity. 
By Proposition~\ref{pro:spectral_gap} \ref{item:L1}, $\lambda_1(H)  = \sqrt{d}(1 + o(1))$. 
Hence, by \eqref{Lambda_alpha_bound} and \eqref{eq:lambda_x_approximated_by_Lambda_alpha_x},
\begin{equation} \label{eq:lambda_1_H_compared_lambda_x} 
\lambda_1(H) \gg \Lambda(\alpha_x)  + o(1) = \lambda(x)
\end{equation} 
for any $x \in\cal W$.
Therefore, we have found $1 + \abs{\cal W}$ non-zero eigenvalues of $\Pi H \Pi$ (counted with multiplicity). 
Since the dimension of $\ran \Pi$ is at most $1 + \abs{\cal W}$, this completes the proof of Proposition~\ref{prop:StrongRigidity} \ref{item:StrongRigidity_upper_block}. 
\end{proof}

\begin{proof}[Proof of Proposition \ref{prop:StrongRigidity} \ref{item:StrongRigidity_offdiag_blocks}]
In order to estimate $\norm{\ol{\Pi}H \Pi}$, let $\f v \in \ran \Pi$. 
We decompose $\f v = \alpha \f w_1 + \f u$ from some $\alpha \in \R$ and $\f u \in \Span \{ \f u(x) \col x \in \cal W\}$. Let $\zeta < 1/2- 3\mu/2$. From \eqref{eq:H_approx_diagonal_u_basis}, we obtain that 
\begin{equation} \label{eq:norm_bound_ol_Pi_H_v} 
\ol{\Pi} H \f v = \ol\Pi \alpha \lambda_1(H) \f w_1 + \ol \Pi \cal D \f u + o (N^{-\zeta} \norm{\f u }) 
= o(N^{-\zeta} \norm{\f u})\,,
\end{equation}
where the last steps follows from $\f w_1 \in \ran \Pi$ and $\ran \cal D \subset \ran \Pi$ due to the definitions of $\cal D$ and $\Pi$. 
It remains to show that $\norm{\f u} \lesssim \norm{\f v}$. 
From \eqref{eq:H_approx_diagonal_u_basis}, we also conclude 
\begin{equation} \label{eq:overlap_u_1_span_u_x}  
\lambda_1(H) \scalar{\f w_1}{\f u} = \scalar{\f w_1}{H \f u } = \scalar{\f w_1}{\cal D \f u } + o(N^{-\zeta} \norm{\f u})\,.
\end{equation} 
Since 
$\norm{\cal D} \leq \max_{x \in \cal W} \lambda(x) \ll \lambda_1(H)$ by \eqref{eq:lambda_1_H_compared_lambda_x}, we obtain from \eqref{eq:overlap_u_1_span_u_x} that $\abs{\scalar{\f w_1}{\f u}} \ll \norm{\f u}$ 
and, thus, 
$2 \abs{\alpha} \abs{\scalar{\f w_1}{\f u}} \leq \alpha^2 + \abs{\scalar{\f w_1}{\f u}}^2 = \alpha^2 +  o(\norm{\f u}^2)$. 
Therefore, 
\begin{equation} \label{eq:norm_v_span_w_1_and_u_x}  
\norm{\f v}^2 = \alpha^2 + \norm{\f u}^2 + 2 \alpha \scalar{\f w_1}{\f u} \geq \alpha^2 + \norm{\f u}^2 
- 2 \abs{\alpha}\abs{\scalar{\f w_1}{\f u}} \geq \alpha^2/2 + (1- o(1)) \norm{\f u}^2\,. 
\end{equation} 
Hence, $\norm{\f u} \leq (1 - o(1))\norm{\f v}$ which completes the proof of Proposition~\ref{prop:StrongRigidity} \ref{item:StrongRigidity_offdiag_blocks} due to \eqref{eq:norm_bound_ol_Pi_H_v}. 
\end{proof}

The main ingredient in the proof of Proposition~\ref{prop:RigidityUPerp} is the following result. It uses the following orthogonal projection.

\begin{definition}[$\Pi_{\cal X}$] \label{def:Pi_X}
For any $\cal X \subset \cal W$, we denote by $\Pi_{\cal X}$ the orthogonal projection onto $\Span\{ \f u(x) \col x\in\cal X \}$ and we define $\overline{\Pi}_{\cal X}\deq I-\Pi_{\cal X}$.
\end{definition}

\begin{proposition} \label{prop:RigidityOffDiagonal}
If $\mu \in [0,1/3)$ then
\[
\mathbb{E}_{\Omega}\bigg[\max_{\cal X\subset \cal W}\|\overline{\Pi}_{\cal X}H\Pi_{\cal X}\|^{2}\bigg]\lesssim d^{-1}(\log N)^{2}N^{2\mu-1}\,. 
\]
\end{proposition}

Before proving Proposition~\ref{prop:RigidityOffDiagonal} we deduce Proposition~\ref{prop:RigidityUPerp} from it. 

\begin{proof}[Proof of Proposition \ref{prop:RigidityUPerp}]
Recall from Definition \ref{def:u_perp} that the elements of $\cal W$ are ordered in an arbitrary fashion. 
For $x \in \cal W$, let $\Pi_{<x}$ be the orthogonal projection onto $\Span\{ \f u(y) \col y \in \cal W, \, y < x \}$ and $\ol{\Pi}_{<x} \deq 1 - \Pi_{<x}$. 
Then for any $x \in \cal W$ we have 
\begin{equation} \label{eq:representation_Gram_Schmidt} 
\f u^\perp (x) = \frac{\f u(x) -\Pi_{<x}\f u(x) }{\|\f u(x) -\Pi_{<x}\f u(x) \|}= \frac{\overline{\Pi}_{<x}\f u(x) }{\|\overline{\Pi}_{<x}\f u(x) \|}\,. 
\end{equation} 
In order to estimate $(H-\lambda(x))\ol{\Pi}_{<x}\f u(x)$, we conclude 
from the definitions of $\Pi_{<x}$ and $\ol{\Pi}_{<x}$ that 
\[
H\overline{\Pi}_{<x} = \overline{\Pi}_{<x}H\overline{\Pi}_{<x} + \Pi_{<x}H\overline{\Pi}_{<x}  = \overline{\Pi}_{<x}H-\overline{\Pi}_{<x}H\Pi_{<x} +\Pi_{<x}H\overline{\Pi}_{<x}\,,
  \] 
which then yields 
\[(H-\lambda(x))\overline{\Pi}_{<x}\f u(x)  = \overline{\Pi}_{<x} (H-\lambda(x)) \f u(x)  +\Pi_{<x}H\overline{\Pi}_{<x}\f u(x)  -\overline{\Pi}_{<x}H\Pi_{<x}\f u(x)\,. 
 \]
Hence,
\[\|(H-\lambda(x))\overline{\Pi}_{<x}\f u(x) )\| \leq  \|(H-\lambda(x)) \f u(x) \| +2 \|\overline{\Pi}_{<x}H\Pi_{<x} \|\,. \]
From \eqref{eq:expectation_norm_H_minus_lambda_u} combined with a union bound over $x$, Proposition \ref{prop:RigidityOffDiagonal}, and Chebyshev's inequality, we deduce that with high probability
\begin{equation*}
\max_{x \in \cal W} \|(H-\lambda(x))\overline{\Pi}_{<x}\f u(x) )\| \leq N^{\mu-1/2+o(1)}\,.
\end{equation*}
Moreover, since $\scalar{\f 1_x}{\f u(y)} = 0$ for $x$, $y \in \cal W$ satisfying $y < x$, we find that with high probability, for all $x \in \cal W$,
\begin{equation*}
\|\overline{\Pi}_{<x}\f u(x) \|^2 = \|\f u(x) -\Pi_{<x}\f u(x) \|^2 \geq \scalar{\f 1_x}{\f u(x) -\Pi_{<x}\f u(x)}^2 = \scalar{\f 1_x}{\f u(x)}^2\gtrsim 1\,,
\end{equation*}
where the last step follows from Corollary~\ref{cor:lower_bound_u_x_in_x} and the fact that $\P(\Omega) = 1 - o(1)$. Plugging these two estimates into \eqref{eq:representation_Gram_Schmidt} yields \eqref{eq:H-lambda}.

For the proof of \eqref{eq:H_approx_diagonal_u_basis}, we write 
 $\f u = \sum_{x\in \cal W} a_x \f u^\perp(x) $ and apply \eqref{eq:H-lambda} and Cauchy-Schwarz,
\[
\|(H-\cal D)\f u \| = \normbb{\sum_{x\in \cal W} a_x (H-\lambda(x))\f u^\perp(x) } \leq N^{\mu-1/2+o(1)}|\cal W|^{1/2}\bigg(\sum_{x\in \cal W} a_x^2\bigg)^{1/2}\,. 
\]
Hence, \eqref{eq:H_approx_diagonal_u_basis} 
follows from $\cal W \subset \cal V$ and Proposition~\ref{prop:graphProperties} \ref{item:size_cal_V}. 
\end{proof}

\begin{proof}[Proof of Proposition \ref{prop:RigidityOffDiagonal}]
Let $\f v\in\ran \Pi_{\cal X}$. We write $\f v=\sum_{x\in{\cal X}}a_{x}\f u(x)$ with
$a_{x}\in\mathbb{R}$, use \eqref{eq:represenation_error_H_u_x} and $\ol{\Pi}_{\cal X} \f u (x) = 0$ 
for each $x \in \cal X$ to obtain 
\begin{align*}
\|\overline{\Pi}_{\cal X}H\f v\|^{2} & =\normbb{\sum_{x\in{\cal X}}a_{x}\overline{\Pi}_{\cal X}\bigg(\lambda(x) \f u(x)+\sum_{y\in{\cal V \setminus\{ x\}}}\epsilon_{y}(x)\f 1_{y}\bigg)}^{2}\\
 & =\normbb{\overline{\Pi}_{\cal X}\sum_{x\in \cal X,y\in{\cal V \setminus \{x \}}}\epsilon_{y}(x)a_{x}\f 1_{y}}^{2}\\
 & \leq\sum_{y\in{\cal V}}\bigg(\sum_{x\in{\cal X \setminus \{y \}}}\epsilon_{y}(x)a_{x}\bigg)^{2}\\
 & \leq\bigg(\sum_{x\in \cal X}a_{x}^{2} \bigg)\bigg(\sum_{x\in{\cal W},\, y \in \cal V \setminus\{x\}}\epsilon_{y}(x)^{2}\bigg)\,.
\end{align*}
Moreover, since $\scalar{\f 1_y}{\f u(x)} = 0$ for any $y \in \cal V \setminus \{ x\}$, on $\Omega$ we have
\[
\|\f v\|^{2}\geq\sum_{x\in{\cal X}}|\scalar{\f 1_x}{\f v}|^{2}=\sum_{x\in{\cal X}}a_{x}^{2}\abs{\scalar{\f 1_x}{\f u(x)}}^{2}\geq c\sum_{x\in{\cal X}}a_{x}^{2}
\]
 for some constant $c>0$ by Corollary~\ref{cor:lower_bound_u_x_in_x}. 
Therefore,  on $\Omega$ we have
\begin{align*}
\|\overline{\Pi}_{\cal X}H\f v\|^{2} & \leq\frac{1}{c}\bigg(\sum_{x \in \cal W,\, y\in{\cal V \setminus \{x\}}}\epsilon_{y}(x)^{2}\bigg)\|\f v\|^{2}
\end{align*}
and, in particular, $\|\overline{\Pi}_{\cal{X}}H\Pi_{\cal X}\|^{2}\leq\frac{1}{c}\left(\sum_{x \in \cal W,\, y\in{\cal V \setminus \{x\}}}\epsilon_{y}(x)^{2}\right)$ for any $\cal X \subset \cal W$. 
By Proposition~\ref{prop:EstimateExy}, we therefore conclude  
\[
\mathbb{E}_{\Omega}\left[\sum_{x \in \cal W,\, y\in{\cal V\setminus\{x\}}}\epsilon_{y}(x)^{2}\right]=\sum_{x,y\in [N]}  \mathbb{E}_{\Omega}
\big[ \ind{x \in \cal W} \ind{y\in{\cal V \setminus \{x \}}} \, \epsilon_{y}(x)^{2}\big]\lesssim d^{-1}(\log N)^{2}N^{2\mu-1}\,,
\]
as claimed.
\end{proof}

\subsection{Proof of Proposition~\ref{prop:StrongRigidity} \ref{item:StrongRigidity_lower_block}} 
\label{sec:proof_Prop_StrongRigidity} 

For the proof of Proposition~\ref{prop:StrongRigidity} \ref{item:StrongRigidity_lower_block}, we shall need several notions from the works \cite{ADK19, ADK20}.
As in \cite[eq.~(1.8)]{ADK20}, we set   
\begin{equation}\label{eq:def_r_star} 
 r_\star = \floor{c \sqrt{\log N}} 
\end{equation} 
for the constant $c>0$ from \cite{ADK20}. 
Following \cite[eq.~(1.9)]{ADK20}, we define 
\begin{equation} \label{eq:def_xi_xi_u} 
\xi \deq \frac{\sqrt{\log N}}{d} \log d\,, \qquad \qquad  
\xi_u \deq \frac{\sqrt{\log N}}{d} \frac{1}{u} 
\end{equation}
for $u >0$. 
For any $\tau \in [1 + \xi^{1/2}, 2]$, we denote by $\mathbb{G}_\tau$ the \emph{pruned graph} introduced in 
\cite[Proposition~3.1]{ADK20}. 
We denote the balls and spheres in $\mathbb{G}_\tau$ around a vertex $x \in [N]$ by $B_i^\tau(x)$ and $S_i^\tau(x)$, respectively. 
The pruned graph $\mathbb{G}_\tau$ is a subgraph of $\mathbb{G}$, which possesses a number of useful properties listed 
in \cite[Proposition~3.1]{ADK20}. In particular, the balls $B_{2r_\star}^\tau(x)$ and $B_{2r_\star}^\tau(y)$ 
in $\mathbb{G}_\tau$ are disjoint if $x$, $y \in [N]$ satisfy $x \neq y$ and $\min\{\alpha_x, \alpha_y \} \geq \tau$.

Recalling the definition of $u_i(\alpha)$ from \eqref{def_ui}, for any $x \in [N]$ with 
$\alpha_x \geq 2 + \xi^{1/4}$ and $\sigma = \pm$, as in \cite[eq.~(3.5)]{ADK20}, we define 
\begin{equation} \label{eq:def_v_tau} 
\f v_\sigma^\tau (x) \deq \sum_{i=0}^{r_\star} \sigma^i u_i(\alpha_x) \frac{\f 1_{S_i^\tau(x)}}{\norm{\f 1_{S_i^\tau(x)}}}\,,
 \end{equation}  
where, for the last coefficient $u_{r_\star}(\alpha_x)$ we make the special choice $u_{r_\star}(\alpha_x) \deq u_{r_\star - 1}(\alpha_x) / \sqrt{\alpha_x}$, and $u_0(\alpha_x)> 0$ is chosen such that $\f v_\sigma^\tau(x)$ is normalized, i.e.\ $\sum_{i=0}^{r_\star} u_i^2(\alpha_x) = 1$.

\begin{remark} \label{rem:v_tau_orthonormal}
The family $(\f v_\sigma^\tau(x) \col x \in [N], \alpha_x \geq 2 + \xi^{1/4}, \, \sigma = \pm)$ is orthonormal. See \cite[Remark~3.3]{ADK20}
\end{remark}
As in \cite[Definition~3.6]{ADK20}, we denote the adjacency matrix of $\mathbb{G}_\tau$ by $A^\tau$ and define the 
matrix 
\begin{equation} \label{eq:def_H_tau} 
 H^\tau \deq ( A^\tau - \chi^\tau ( \E A) \chi^\tau)/\sqrt{d}\,, 
\end{equation} 
where $\chi^\tau$ is the orthogonal projection onto $\Span\{ \f 1_y \col y \notin \bigcup_{x \col \alpha_x \geq \tau} B_{2r_\star}^\tau(x) \}$. 
Moreover, we recall \cite[Definition~3.10]{ADK20}.

\begin{definition}[$\Pi^\tau$, $\wh H^\tau$]
Define the orthogonal projections (see Remark \ref{rem:v_tau_orthonormal})
\[ \Pi^\tau \deq \sum_{x \,:\, \alpha_x \geq 2 + \xi^{1/4}} \sum_{\sigma = \pm} \f v_\sigma^\tau(x) \f v_\sigma^\tau(x)^*\,, \qquad \qquad \ol{\Pi}^\tau \deq 1 - \Pi^\tau \] 
and the associated block matrix (recall \eqref{eq:def_H_tau})
\begin{equation} \label{eq:def_wh_H_tau} 
\wh{H}^\tau \deq \sum_{x \,:\, \alpha_x \geq 2 + \xi^{1/4}} \sum_{\sigma = \pm} \sigma \Lambda(\alpha_x) \f v_\sigma^\tau(x) \f v_\sigma^\tau(x)^* + \ol{\Pi}^\tau H^\tau \ol{\Pi}^\tau\,. 
\end{equation}
\end{definition}

We note that, for any $\tau \in [1+ \xi^{1/2},2]$, by \eqref{eq:d_range} we have   
\begin{equation} \label{eq:xi_xi_tau_minus_1_are_o_1} 
\xi = o(1)\,, \qquad \qquad \xi_{\tau-1} = o(1)\,. 
\end{equation} 
For any $\tau \in [1 + \xi^{1/2}, 2]$, the definition of $\wh{H}^\tau$ in \eqref{eq:def_wh_H_tau} and \cite[Proposition~3.12]{ADK20} yield that, 
with high probability, 
\begin{equation} \label{eq:bound_norm_wh_H_tau}  
\norm{\wh{H}^\tau} \leq \max \bigg\{ \max_{x \col \alpha_x \geq 2 + \xi^{1/4}} \Lambda(\alpha_x), 
\, 2 \tau + C \xi \bigg\} \lesssim \sqrt{\frac{\log N}{d}} \ll \sqrt{d}\,, 
\end{equation}  
where we used \eqref{Lambda_alpha_bound}, 
\eqref{eq:xi_xi_tau_minus_1_are_o_1} and \eqref{eq:d_range} in the last two steps. 

Owing to \eqref{eq:xi_xi_tau_minus_1_are_o_1}, \cite[Lemmas~3.8, 3.11]{ADK20}\footnote{We stress that the definition of $H$ in \cite{ADK20} differs from that in the current paper by $\E H$; see \cite[Definition~3.6]{ADK20}.}, 
with high probability, we have 
\begin{equation} \label{eq:norm_H_minus_EH_minus_wh_H_tau} 
\norm{H - \E H - \wh{H}^\tau} = o(1)\,. 
\end{equation} 
From \eqref{eq:norm_H_minus_EH_minus_wh_H_tau} and \eqref{eq:bound_norm_wh_H_tau}, 
we conclude that, with high probability, 
\begin{equation} \label{eq:H_minus_EH_ll_sqrt_d}   
 \norm{H - \E H} \ll \sqrt{d}\,. 
\end{equation}

After these preparations, we can start the proof of Proposition~\ref{prop:StrongRigidity} \ref{item:StrongRigidity_lower_block}. We begin with the following definition.

\begin{definition}[$Q_r$] \label{def:Q_r} 
For $r \in \N$, denote by $Q_r$ the orthogonal projection defined by restriction to the set $(\bigcup_{y\in{\cal V}\setminus{\cal W}}B_{r}(y))^{c}$. 
\end{definition}

\begin{definition}[$Q$, $\Pi_Q$] \label{def:Q}
Let $r_\star$ be as in \eqref{eq:def_r_star}. 
Set $Q \deq Q_{2r_\star - 1}$ and 
define
$\Pi_{Q}$ as the orthogonal projection onto $\Span(\{Q\f  u^{\perp}(x) \col x\in{\cal W}\} \cup \{ Q \f w_1 \} )$, and write $\ol \Pi_Q \deq 1 - \Pi_Q$.
\end{definition}
Then
\begin{align}
\lambda_{1}(\overline{\Pi}H\overline{\Pi}) & \leq\lambda_{1}(\overline{\Pi}(H-\mathbb{E}H)\overline{\Pi})+\|\overline{\Pi}(\mathbb{E}H)\overline{\Pi}\| \nonumber \\
& \leq\lambda_{1}(\overline{\Pi}_{Q}\widehat{H}^{\tau}\overline{\Pi}_{Q})+ 2 \norm{H -\E H} \|\Pi-\Pi_{Q}\|+\|H-\mathbb{E}H-\widehat{H}^{\tau}\|+\norm{\ol{\Pi} (\E H) 
\ol{\Pi}} \nonumber \\  
& \leq\lambda_{1}(\overline{\Pi}_{Q}\widehat{H}^{\tau}\overline{\Pi}_{Q})+ 2 \norm{H -\E H} \|\Pi-\Pi_{Q}\|+o(1)\,, 
\label{eq:lambda_1_ol_Pi_H_ol_Pi_first_bound} 
\end{align}
whose last step follows from \eqref{eq:norm_H_minus_EH_minus_wh_H_tau} and $\norm{\ol{\Pi} (\mathbb E H) \ol{\Pi}} = o(1)$. The latter bound is a consequence of    
\begin{equation} \label{eq:bound_ol_Pi_EH_ol_Pi} 
\|\overline{\Pi}(\mathbb{E}H)\overline{\Pi}\|= \sqrt{d} \|\overline{\Pi}  (\f  e\f e^{*} - 1/N ) \overline{\Pi}\|=\sqrt{d}\|\overline{\Pi}(\f e-\f w_1)(\f e^{*}-\f w_1^{*})\overline{\Pi}\| + o(d^{-1/2}) \lesssim d^{-1/2}\,, 
\end{equation} 
where we introduced $\f e \deq N^{-1/2}\f 1_{[N]}$, used $\Pi \f w_1 = \f w_1$ in the second step, and in the last step we used Proposition~\ref{pro:spectral_gap} \ref{item:spectral_gap_q} to estimate $\norm{\f e - \f w_1} \leq \norm{\f w_1 - \f q} + \norm{\f q - \f e} \lesssim d^{-1/2}$.

Fix $\mu \in [0,1/4)$. We then claim that with high probability
\begin{equation} \label{eq:projections_Pi_Pi_Q_approximation} 
\|\Pi-\Pi_Q\| \lesssim d^{-1}
\end{equation} 
and
\begin{equation} \label{eq:lambda_1_ol_Pi_Q_hat_H_ol_Pi_Q} 
\lambda_1(\ol{\Pi}_Q \wh{H}^\tau \ol{\Pi}_Q) \leq \Lambda(\alpha^*) + \kappa/2\,. 
\end{equation}
Using \eqref{eq:projections_Pi_Pi_Q_approximation} and \eqref{eq:lambda_1_ol_Pi_Q_hat_H_ol_Pi_Q}, Proposition~\ref{prop:StrongRigidity} \ref{item:StrongRigidity_lower_block} follows immediately from  \eqref{eq:lambda_1_ol_Pi_H_ol_Pi_first_bound} and \eqref{eq:H_minus_EH_ll_sqrt_d}.  What remains to prove Proposition~\ref{prop:StrongRigidity} \ref{item:StrongRigidity_lower_block}, therefore, is the proof of \eqref{eq:projections_Pi_Pi_Q_approximation} and \eqref{eq:lambda_1_ol_Pi_Q_hat_H_ol_Pi_Q}.

\subsection{Proof of \eqref{eq:projections_Pi_Pi_Q_approximation}}

Clearly,
\begin{equation} \label{eq:proof_pi_minus_pi_Q_small_1}  
\norm{\Pi - \Pi_Q} = \norm{\Pi\ol{\Pi}_Q - \ol{\Pi} \Pi_Q} \leq \norm{\Pi\ol{\Pi}_Q} + \norm{\ol{\Pi} \Pi_Q} 
= \norm{\ol{\Pi}_Q\Pi} + \norm{\ol{\Pi} \Pi_Q}\,,
\end{equation}
where in the last step we used that $(\Pi\ol{\Pi}_Q)^* = \ol{\Pi}_Q \Pi$.
In order to estimate the terms on the right-hand side, we continue with 
\begin{equation} \label{eq:proof_pi_minus_pi_Q_small_2}  
\norm{\ol{\Pi}_Q \Pi}^2 = \sup_{\substack{\f v \in \ran \Pi \\ \norm{\f v} = 1}} \norm{\ol{\Pi}_Q \f v}^2 
\leq \sup_{\substack{\f v \in \ran \Pi \\ \norm{\f v} = 1}} \inf_{\substack{\f u \in \ran \Pi_Q \\ \norm{\f u} = 1}} \norm{\f v - \f u}^2\,. 
\end{equation}
We now apply the next lemma, whose proof is given in Section \ref{sec:pf_auxiliary} below.

\begin{lemma}\label{lem:UXBallY} 
Fix $\mu \in [0,1/3)$. 
With high probability, 
 for all $x\in \cal W$, $y \in \cal V \setminus \cal W$ and $r \in\N$ satisfying $r \ll \frac{\log N}{\log d}$, 
we have 
\[
\|\f u^\perp(x)|_{B_r(y)}\| \lesssim N^{-1/2 + \mu + o(1)}\,.  
\]
\end{lemma}

Owing to Lemma~\ref{lem:UXBallY} as well as parts \ref{item:B_r_disjoint} and \ref{item:size_cal_V} of Proposition~\ref{prop:graphProperties}, we obtain 
\begin{equation} \label{eq:proof_pi_minus_pi_Q_small_3}  
\norm{\f u^\perp(x) - Q \f u^\perp(x)} = \Bigg(\sum_{y \in \cal V\setminus \cal W} \norm{\f u^\perp(x)|_{B_{2r_\star - 1}(y)}}^2\Bigg)^{1/2} \leq N^{-1/2 + 3 \mu/2 + o(1)}\,. 
\end{equation}
Let $\f q$ be as in Proposition~\ref{pro:spectral_gap} \ref{item:spectral_gap_q} with $r = 2 r_\star - 2 \ll \frac{d}{\log \log N}$ by \eqref{eq:d_range}. 
Since $\supp \f q \subset \big( \bigcup_{x \in \cal V} B_{r + 1} (x) \big)^c =  \big( \bigcup_{x \in \cal V} B_{2r_\star-1} (x) \big)^c$,  
we conclude from Proposition~\ref{pro:spectral_gap} \ref{item:spectral_gap_q} and the definition of $Q$ that 
\begin{equation} \label{eq:proof_pi_minus_pi_Q_small_4}  
\|\f w_1 -Q \f w_1 \| = \norm{(1-Q) ( \f w_1 - \f q)} \leq \|\f w_1-\f q\|\lesssim d^{-1}\,. 
\end{equation}

For $\gamma_1, \gamma_x \in \R$ for $x \in \cal W$ write
\begin{equation} \label{def_u_v}
\f v \deq \gamma_1 \f w_1 + \sum_{x \in \cal W} \gamma_x \f u^\perp(x) \in \ran \Pi\,,
\qquad
\f u \deq \gamma_1 Q \f w_1 + \sum_{x \in \cal W} \gamma_x Q \f u^\perp(x) \in \ran \Pi_Q\,.
\end{equation}
Then
\[ 
\norm{\f v - \f u} \lesssim d^{-1} \abs{\gamma_1} + N^{-1/2 + 3\mu/2 + o(1)} \sum_{x \in \cal W} \abs{\gamma_x} 
\lesssim d^{-1} \abs{\gamma_1} + N^{-1/2 + 2\mu +o(1)} \bigg(\sum_{x\in \cal W} \abs{\gamma_x}^2 \bigg)^{1/2} 
\lesssim d^{-1} \norm{\f v}\,. 
\] 
Here, we used in \eqref{eq:proof_pi_minus_pi_Q_small_4} and \eqref{eq:proof_pi_minus_pi_Q_small_3} in the 
first step, Proposition~\ref{prop:graphProperties} \ref{item:size_cal_V} in the second step
and, in the fourth step, $\mu < 1/4$ as well as 
$\norm{\f v}^2 \asymp \abs{\gamma_1}^2 + \sum_{x \in \cal W} \abs{\gamma_x}^2$ 
(the inequality $\gtrsim$ follows from \eqref{eq:norm_v_span_w_1_and_u_x} and the orthogonality of $(\f u^\perp(x))_{x \in \cal W}$; the inequality $\lesssim$ is trivial). 
Hence, if $\norm{\f v} = 1$ then $\norm{\f u} = 1 + O(d^{-1})$ and, thus, $\norm{\f v - \frac{\f u}{\norm{\f u}}} \lesssim d^{-1}$. 
Therefore, $\norm{\ol{\Pi}_Q \Pi} \lesssim d^{-1}$ by \eqref{eq:proof_pi_minus_pi_Q_small_2}.

Finally, similarly to \eqref{eq:proof_pi_minus_pi_Q_small_2}, we have
\begin{equation*}
\norm{\ol{\Pi} \Pi_Q}^2 \leq \sup_{\substack{\f u \in \ran \Pi_Q \\ \norm{\f u} = 1}} \inf_{\substack{\f v \in \ran \Pi \\ \norm{\f v} = 1}} \norm{\f v - \f u}^2\,,
\end{equation*}
and the same argument as above, with the representation \eqref{def_u_v}, implies that the right-hand side is $O(d^{-1})$.
By \eqref{eq:proof_pi_minus_pi_Q_small_1}, we therefore conclude \eqref{eq:projections_Pi_Pi_Q_approximation}.

\subsection{Proof of \eqref{eq:lambda_1_ol_Pi_Q_hat_H_ol_Pi_Q}}
We begin by introducing another orthogonal projection.

\begin{definition}[$\Pi_{\f v}$]
Let $\Pi_{\f v}$ be the orthogonal projection onto $\Span\{\f v^\tau_+(z) \col z\in \cal V\setminus \cal W\}$.
\end{definition}
For any $z\in \cal V \setminus \cal W$, we have $\supp\f v^\tau_\pm(z)\subset B_{r_\star}(z)$  
and, thus, by definition of $Q$, $\supp \f v^\tau_\pm(z) \cap \supp Q\f u^\perp (x) = \emptyset$ for any $x \in \cal W$ and 
$\supp \f v^\tau_\pm(z) \cap \supp Q \f w_1 = \emptyset$. 
Therefore, $\f v^\tau_\pm(z)$ is orthogonal to $\ran \Pi_Q$, i.e.\ $\Pi_Q \f v_\pm^\tau(x) = 0$. 
That implies $\Pi_{\f v} \Pi_Q = 0 = \Pi_Q \Pi_{\f v}$ and, in particular, $\Pi_{\f v}$ 
and $\Pi_{Q}$ commute. 
Since $\Pi_{\f v}$ and $\wh{H}^\tau$ commute and $\Pi_{\f v}$ and $\Pi_Q$ commute, we obtain 
\begin{equation} \label{eq:lambda_1_ol_Pi_Q_hat_H_ol_Pi_Q_second} 
\begin{aligned}
\lambda_{1}(\overline{\Pi}_{Q}\wh{H}^{\tau}\overline{\Pi}_{Q}) & =\lambda_{1}(\overline{\Pi}_{Q}\Pi_{\f v}\wh{H}^{\tau}\Pi_{\f v}\overline{\Pi}_{Q}+\overline{\Pi}_{Q}\overline{\Pi}_{\f v}\wh{H}^{\tau}\overline{\Pi}_{\f v}\overline{\Pi}_{Q})\\
 & =\max\{\lambda_{1}(\overline{\Pi}_{Q}\Pi_{\f v}\wh{H}^{\tau}\Pi_{\f v}\overline{\Pi}_{Q}),\, \lambda_{1}(\overline{\Pi}_{Q}\overline{\Pi}_{\f v}\wh{H}^{\tau}\overline{\Pi}_{\f v}\overline{\Pi}_{Q})\}\,. 
\end{aligned}
\end{equation}
By the definition of $\wh{H}^\tau$ in \eqref{eq:def_wh_H_tau}, we have 
\begin{equation} \label{eq:lambda_1_Pi_V_hat_H_Pi_V} 
\lambda_1(\Pi_{\f v}\widehat{H}^\tau \Pi_{\f v})= \max_{x\in{\cal V}\setminus{\cal W}}\Lambda(\alpha_{x}) \leq\Lambda(\alpha^{*})+\kappa/2\,. 
\end{equation}
What remains, therefore, is to estimate $\lambda_{1}(\overline{\Pi}_{Q}\overline{\Pi}_{\f v}\wh{H}^{\tau}\overline{\Pi}_{\f v}\overline{\Pi}_{Q})$.\

Suppose that $\f w$ is a normalized eigenvector of $\overline{\Pi}_{Q}\overline{\Pi}_{\f v}\widehat{H}^\tau\overline{\Pi}_{\f v}\overline{\Pi}_{Q}$ such that the associated eigenvalue $\lambda$ satisfies $\lambda\geq \Lambda(\alpha^*)+\kappa/2$. 
We now check that the next lemma, whose proof is given in Section \ref{sec:pf_auxiliary} below, is applicable to $\f w$ and $\lambda$ for any $x \in \cal V \setminus \cal W$.

\begin{lemma}\label{lem:weak_delocalization}  
Let $r_\star \in \N$ be as in \eqref{eq:def_r_star}. In particular, $r_\star \asymp \sqrt{\log N}$. 
Suppose $\tau \in [1 + \xi^{1/2},2]$. 

There is a constant $C>0$ such that the following holds with high probability\footnote{We note that the statement actually holds with \emph{very high probability}, meaning that for each $\nu>0$, there is a constant $C \equiv C_\nu >0$ such that \eqref{eq:weak_delocalization_conclusion} holds with probability at least $1- N^{-\nu}$ for all sufficiently large $N$.}. 
Let $x \in [N]$, $\lambda > 2 \tau + C \xi$ and $\f w$ satisfy 
\begin{equation} \label{eq:condition_weak_delocalization} 
(\widehat{H}^\tau \f w)|_{B_{2r_\star-1}^\tau(x)}=\lambda \f w|_{B_{2r_\star-1}^\tau(x)}\,. 
\end{equation}
If $\alpha_x \geq 2 + \xi^{1/4}$ and $\f v_-^\tau(x) \perp \f w \perp \f v_+^{\tau}(x)$ or   
$2 + \xi^{1/4} > \alpha_x \geq \tau$ then
\begin{equation} \label{eq:weak_delocalization_conclusion} 
 \frac{\abs{\scalar{\f 1_x}{\f w}}}{\norm{\f w |_{B_{2r_\star}^\tau(x)}}} \lesssim  \frac{\lambda^2}{(\lambda - 2 \tau - C \xi)^2} \bigg( \frac{2 \tau + C \xi}{\lambda} \bigg)^{r_\star}\,. 
\end{equation} 
An analogous result holds if $\lambda < - 2\tau - C \xi$. 
\end{lemma}

Choosing $\tau = 1 + \xi^{1/2}$ and taking $x \in \cal V \setminus \cal W$, we now verify the conditions of Lemma~\ref{lem:weak_delocalization} for $\f w$ and $\lambda$. 
Owing to \eqref{eq:xi_xi_tau_minus_1_are_o_1}, there is a constant $c\equiv c_\kappa>0$ such that 
$\frac{2 \tau + C \xi}{\lambda} \leq 1 - c$ as $\lambda \geq \Lambda(\alpha^*) + \kappa/2 \geq 2 + 2c$ due to the 
definition of $\alpha^*(\mu)$ in \eqref{eq:largeDegree}. 
In particular, $\lambda > 2 \tau + C \xi$. 
From $\overline{\Pi}_{Q}\overline{\Pi}_{\f v}\widehat{H}^\tau\overline{\Pi}_{\f v}\overline{\Pi}_{Q}\f w = \lambda \f w$ we conclude that $\ol\Pi_Q \f w = \f w$ and $\ol\Pi_{\f v} \f w = \f w $.  
Thus, as $\Pi_{\f v}$ and $\wh{H}^\tau$ commute, we get $\ol{\Pi}_Q \wh{H}^\tau \f w = \lambda \f w$. 
Restricting both sides in the last identity to $B_{2r_\star - 1}^\tau(x)$ yields $(\wh{H}^\tau \f w)|_{B_{2r_\star-1}^\tau(x)} = \lambda \f w |_{B^\tau_{2r_\star - 1}(x)}$ as $Q$ is the restriction to $\big(\bigcup_{y \in \cal V \setminus \cal W} B_{2r_\star-1}(y) \big)^c$ and $B_{2r_\star - 1}^\tau(x) \subset B_{2 r_\star -1}(x)$ by \cite[Proposition~3.1 (iv)]{ADK20}. 
This proves \eqref{eq:condition_weak_delocalization}. 
Note that $\alpha_x \geq 2 + \xi^{1/4}$. 
\ From $\ol{\Pi}_{\f v} \f w = \f w$, we conclude $\f w \perp \f v^\tau_+(x)$. 
Since $(\f v_\sigma^\tau(y) \col \alpha_y \geq 2 + \xi^{1/4}, \, \sigma = \pm)$ is an orthonormal family by Remark \ref{rem:v_tau_orthonormal}, 
the definition of $\Pi_{\f v}$ implies $\Pi_{\f v} \f v_-^\tau(y) = 0$ for all $y \in \cal V$. 
Therefore, as moreover $\Pi_Q \f v_-^\tau(x) = 0$, we have  $\ol{\Pi}_Q \ol{\Pi}_{\f v} \wh{H}^\tau \ol{\Pi}_{\f v} \ol{\Pi}_Q \f v_-^\tau(x) = -\Lambda(\alpha_x) \f v_-^\tau(x)$. 
Hence, $\f w \perp \f v_-^\tau(x)$. 
Therefore, we have verified all assumptions of Lemma~\ref{lem:weak_delocalization} with $\tau =1 + \xi^{1/2}$ 
for $\f w$, $\lambda$, and any $x \in \cal V \setminus \cal W$.

Since $\frac{2 \tau + C \xi}{\lambda} \leq 1 - c$ as shown above, for the right-hand side of \eqref{eq:weak_delocalization_conclusion}, we get 
\[ \frac{\lambda^2}{(\lambda - 2 \tau - C\xi)^2} \bigg( \frac{2\tau + C \xi}{\lambda} \bigg)^{r_\star} 
= \bigg( 1 - \frac{2 \tau + C \xi}{\lambda} \bigg)^{-2} \bigg( \frac{2 \tau + C \xi}{\lambda} \bigg)^{r_\star} 
\leq c^{-2} (1 - c)^{r_\star} \ll d^{-1/2} 
\] 
as $r_\star \asymp \sqrt{\log N}$. 
Therefore, Lemma~\ref{lem:weak_delocalization}, the disjointness of the balls $(B_{2r_\star}^\tau(x))_{x \in \cal V \setminus \cal W}$ (see the paragraph after \eqref{eq:def_xi_xi_u} as well as \cite[Proposition~3.1 (i)]{ADK20}) and $\norm{\f w}^2 = 1$ imply 
\begin{equation} \label{eq:norm_w_restricted_V_setminus_W}
\norm{\f w|_{\cal V\setminus \cal W}}^2 = \sum_{x \in \cal V \setminus \cal W} \scalar{\f 1_x}{\f w}^2 \ll d^{-1} \sum_{x \in \cal V \setminus \cal W} 
\norm{\f w|_{B_{2r_\star}^\tau(x)}}^2 \leq d^{-1}\,. 
\end{equation} 
 
We recall from Definition~\ref{def:Q_r} that $Q_0$ denotes the orthogonal projection defined by restriction to the set $(\cal V\setminus \cal W)^c$. 
Since $\f w = \f w |_{\cal V\setminus \cal W} + Q_0 \f w$, we obtain 
from \eqref{eq:norm_w_restricted_V_setminus_W}, $Q_0  \wh{H}^\tau Q_0  = (\wh{H}^\tau)^{(\cal V \setminus \cal W)}$, and $\norm{\wh{H}^\tau} \ll \sqrt{d}$ by \eqref{eq:bound_norm_wh_H_tau} that  
\begin{align}
 \lambda = \scalar{\f w}{\widehat{H}^\tau \f w} & =\scalar{\f w}{(\widehat{H}^\tau)^{({\cal V}\setminus{\cal W})} \f w}+ o(1) 
\nonumber \\
 & =\scalar{\f w}{\overline{\Pi}_{Q}(\widehat{H}^\tau)^{({\cal V}\setminus{\cal W})}\overline{\Pi}_{Q} \f w}+ o(1)  
\nonumber \\
 & \leq\lambda_{1}(\overline{\Pi}_Q (\widehat{H}^\tau)^{({\cal V}\setminus{\cal W})}\overline{\Pi}_Q)+o(1) 
\nonumber \\
 & \leq\lambda_{1}(\overline{\Pi} (\widehat{H}^\tau+(\mathbb{E}H))^{({\cal V}\setminus{\cal W})}\overline{\Pi})+o(1) 
\nonumber \\
 & \leq\lambda_{1}(\overline{\Pi}H^{({\cal V}\setminus{\cal W})}\overline{\Pi})+o(1)\,. 
\label{eq:lambda_1_ol_Pi_Q_hat_H_ol_Pi_Q_third} 
\end{align}
Here, the third step is a consequence of $\ol\Pi_Q \f w = \f w$. 
In the fifth step, we used \eqref{eq:bound_ol_Pi_EH_ol_Pi} and \eqref{eq:projections_Pi_Pi_Q_approximation} and the 
sixth step follows from \eqref{eq:norm_H_minus_EH_minus_wh_H_tau} and $\norm{M^{(\cal V\setminus\cal W)}} \leq \norm{M}$ for any matrix $M$.

We now apply the next result, whose proof is given in Section \ref{sec:pf_auxiliary} below.
\begin{lemma} \label{prop:PiBoundVminusW}
Fix $\mu \in [0,1/3)$.  With high probability, 
$\lambda_{1}(\overline{\Pi}H^{({\cal V}\setminus{\cal W})}\overline{\Pi})\leq\Lambda(\alpha^{*})+o(1).$
\end{lemma}

From Lemma~\ref{prop:PiBoundVminusW} and \eqref{eq:lambda_1_ol_Pi_Q_hat_H_ol_Pi_Q_third}, 
we deduce that $\lambda \leq \Lambda(\alpha^*) + o(1)$, in contradiction with the assumption $\lambda\geq \Lambda(\alpha^*)+\kappa/2$. We conclude that $\lambda_1(\overline{\Pi}_{Q}\overline{\Pi}_{{\f v}}\widehat{H}^\tau\overline{\Pi}_{{\f v}}\overline{\Pi}_{Q})\leq
\Lambda(\alpha^*) + \kappa/2$. 
Owing to \eqref{eq:lambda_1_ol_Pi_Q_hat_H_ol_Pi_Q_second} and \eqref{eq:lambda_1_Pi_V_hat_H_Pi_V}, 
this proves \eqref{eq:lambda_1_ol_Pi_Q_hat_H_ol_Pi_Q}.

\subsection{Proofs of auxiliary results} \label{sec:pf_auxiliary}

In this final subsection we prove Lemmas~\ref{lem:UXBallY}, \ref{lem:weak_delocalization}, and \ref{prop:PiBoundVminusW}. 

\begin{proof}[Proof of Lemma \ref{lem:UXBallY}]

For fixed $y \in [N]$, on the event $\Omega\cap \{ y \in \cal V\setminus \cal W\}$, if $x \in \cal W^{(y)}$, let $\f u^{(y)}(x)$ be the eigenvector of $H^{((\cal V^{(y)}\cup \{ y \} ) \setminus \{ x\})}$ with eigenvalue $\lambda_2(H^{((\cal V^{(y)}\cup \{ y \} ) \setminus \{ x\})})$ and satisfying $\scalar{\f 1_x}{\f u^{(y)}(x)} > 0$. See Corollaries \ref{cor:largest_eigenvalues_vertices_removed} and \ref{cor:lower_bound_u_x_in_x} for the existence and uniqueness of $\f u^{(y)}(x)$.

In analogy to Definition \ref{def:u_perp}, let $((\f u^{(y)})^\perp(x))_{x \in \cal W^{(y)}}$ be the Gram-Schmidt orthonormalization of $(\f u^{(y)}(x))_{x \in \cal W^{(y)}}$. 
On $\Omega \cap \{y \in \cal V \setminus \cal W\}$, we have $\cal V^{(y)} \cup \{y \} = \cal V$ and $\cal W^{(y)} = \cal W$ by Remark~\ref{rem:cal_V_without_y_cal_W_without_y}. 
Therefore, $\f u^{(y)}(x) = \f u(x)$ and, thus, $(\f u^{(y)})^\perp(x) = \f u^\perp(x)$ for all $x \in \cal W^{(y)} = \cal W$ on $\Omega \cap \{y \in \cal V \setminus \cal W\}$. 
Hence, for fixed $x$, $y \in [N]$ with $x \neq y$, we estimate
\begin{align*}
&\E_\Omega \bigg[ \ind{x \in \cal W,\, y \in \cal V \setminus \cal W} \norm{\f u^\perp(x)|_{B_r(y)}}^2 \bigg] \\
& \hspace*{.5cm}  = \E_\Omega\bigg[ \ind{x \in \cal W^{(y)}} \sum_{a \in [N]\setminus \{y\}} \scalar{\f 1_a}{(\f u^{(y)})^\perp(x)}^2 \ind{a \in B_r(y)}\ind{y \in \cal V \setminus \cal W}  \bigg] \\ 
& \hspace*{.5cm}  \leq  \E\bigg[ \ind{x \in \cal W^{(y)}} \ind{\Omega_r^{(y)}} \sum_{a \in [N]\setminus \{ y\}} \scalar{\f 1_a}{(\f u^{(y)})^\perp(x)}^2 \sum_{b \in B_{r-1}^{(y)}(a)} \E\Big[ \ind{b \in S_1(y)}\ind{\alpha^*(\mu)d \leq \abs{S_1(y)} \leq 10 \log N} \condB A^{(y)}\Big]  \bigg] \\ 
&  \hspace*{.5cm} \leq \E \bigg[ \ind{x \in \cal W^{(y)}} \ind{\Omega_r^{(y)}} \sum_{a \in [N]\setminus \{y \}} \scalar{\f 1_a}{( \f u^{(y)})^\perp(x)}^2 \abs{B_{r-1}^{(y)}(a)} \bigg]
\\
&  \hspace*{2cm}
\times \max_{b \in [N] \setminus \{y\}} \P ( b \in S_1(y), \alpha^*(\mu) d \leq \abs{S_1(y)} \leq 10 \log N)  \\ 
& \hspace*{.5cm} \lesssim (\log N)d^{r-1} \P (x \in \cal W^{(y)}) \max_{b \in [N]\setminus \{ y\}} 
\P (b \in S_1(y), \, \alpha^*(\mu) d \leq \abs{S_1(y)} \leq 10 \log N) \\ 
& \hspace*{.5cm} \lesssim \frac{(\log N)^2 d^{r-1}}{N} \P (x \in \cal V) \P(y \in \cal V) \\ 
& \hspace*{.5cm} \leq N^{-3 + 2\mu + o(1)}\,. 
\end{align*} 
Here, in the first step, we also used that $\scalar{\f 1_y}{(\f u^{(y)})^\perp(x)} = 0$ for all $y \in \cal V \setminus \cal W$ and $x \in \cal W^{(y)}$. 
In the second step, we conditioned on $A^{(y)}$, 
employed the notations 
\begin{align*} 
B_i^{(y)}(a) & \deq \text{ball of radius $i$ around $a$ in the graph $\bb G |_{[N] \setminus \{y\}}$}
\\ 
\Omega^{(y)}_r & \deq \{ \abs{B_r^{(y)}(z)} \lesssim (\log N) d^{r-1} \text{ for all } z \in [N] \setminus \{y \} \} 
\end{align*} 
and used that $\cal W^{(y)}$, $(\f u^{(y)})^\perp(x)$, and $B_{r-1}^{(y)}(a)$ are $A^{(y)}$-measurable, 
that $\Omega \subset \Omega_r^{(y)}$ by Proposition~\ref{prop:graphProperties_every_vertex}, as 
$r \ll \frac{\log N}{\log d}$, and 
that $ a \in B_r(y)$ is equivalent to $b \in S_1(y)$ for some $b \in B_{r-1}^{(y)}(a)$. 
The third step follows from the independence of $b \in S_1(y)$ and $\abs{S_1(y)}$ from $A^{(y)}$. 
The normalization of $(\f u^{(y)})^{\perp}(x)$ and the definition of $\Omega_r^{(y)}$ imply the fourth step. 
In the fifth step, we argued as in the last steps of \eqref{eq:bound_individual_term_eps_y_x} and, finally, 
we used $r \ll \frac{\log N}{\log d}$ as well as the definitions of $\cal V$ in \eqref{eq:def_cal_V} and of $\alpha^*$ in \eqref{eq:largeDegree}. 

Therefore, a union bound over $x$, $y \in [N]$ with $x \neq y$ and Chebyshev's inequality complete 
the proof of Lemma~\ref{lem:UXBallY}. 
\end{proof}

\begin{proof}[Proof of Lemma~\ref{lem:weak_delocalization}] 
In order to prove Lemma~\ref{lem:weak_delocalization}, we follow \cite[Proof of Proposition~3.14 (i)]{ADK20}\footnote{Note that ${\cal V}_\tau = \{ x \in [N] \col \alpha_x \geq \tau\}$ and $\cal V = {\cal V}_{2 + \xi^{1/4}}$ in \cite{ADK20}, see \cite[eq.~(3.3)]{ADK20}, which differs from the definition of $\cal V$ in the present paper, 
see~\eqref{eq:def_cal_V}.},
whose assumptions are all satisfied apart from the eigenvalue-eigenvector relation $\wh{H} \f w = \lambda \f w$. 

We now explain the necessary, minor, modifications. 
We start with the case $\alpha_x \geq 2 + \xi^{1/4}$ and $\f w\perp \f v_\pm^\tau(x)$ which corresponds to the case 
$x \in \cal V$ in \cite[Proposition 3.14]{ADK20}. 
The eigenvalue-eigenvector relation is used in \cite[Proof of Proposition~3.14 (i)]{ADK20} only in \cite[eq.~(3.55)]{ADK20}. 
Using \eqref{eq:condition_weak_delocalization} instead of the eigenvalue-eigenvector relation and the notation of \cite[Proof of Proposition~3.14 (i)]{ADK20}, 
we now verify the first two steps in \cite[eq.~(3.55)]{ADK20}. 
We write $P_{2r_\star - 1}$ for the orthogonal defined by restriction to the set $B_{2r_\star - 1}^\tau(x)$. 
For any $i < r_\star$, since $\supp \f g_i \subset B_{2r_\star - 1}^\tau(x)$ by \cite[eq.~(3.52)]{ADK20}, 
 we have $P_{2r_\star-1} \f g_i = \f g_i$.  
Therefore, for any $i < r_\star$, 
we obtain  
\begin{align*} 
\lambda u_i = \scalar{\f g_i}{\lambda P_{2r_\star-1} \f w} 
 = \scalar{( \wh{H}^\tau - \Lambda(\alpha_x) \f v_+^\tau(x) \f v_+^\tau(x) + \Lambda(\alpha_x) \f v_-^\tau(x) \f v_-^\tau(x)) \f g_i}{\f w}  
= \scalar{\wh{H}^{\tau,x} \f g_i}{\f w}\,. 
\end{align*} 
Here, we used \eqref{eq:condition_weak_delocalization} and $\f w \perp \f v_\pm^\tau(x)$ in the second step 
and $\supp \f g_i \subset B_{2r_\star-1}^\tau(x)$
as well as the definition of $\wh{H}^{\tau,x}$ from \cite[eq.~(3.48)]{ADK20} in the third step. 

The remaining steps in \cite[eq.~(3.55)]{ADK20} and the remainder of \cite[Proof of Proposition~3.14 (i)]{ADK20} 
including the case $2 + \xi^{1/4} > \alpha_x \geq \tau$, which corresponds to the case $x \in \cal V_\tau \setminus \cal V$ 
in \cite{ADK20}, 
are obtained in the same way as in \cite{ADK20}. 
\end{proof}

\begin{proof}[Proof of Lemma \ref{prop:PiBoundVminusW}]
We write $H^{(\cal V \setminus \cal W)}$ in the block decomposition  
\[
H^{({\cal V}\setminus{\cal W})}=\Pi H^{({\cal V}\setminus{\cal W})}\Pi+\overline{\Pi}H^{({\cal V}\setminus{\cal W})}\overline{\Pi}+\overline{\Pi}H^{({\cal V}\setminus{\cal W})}\Pi+\Pi H^{({\cal V}\setminus{\cal W})}\overline{\Pi}\,.
\]
The nonzero eigenvalues of the block diagonal arise as the eigenvalues of the individual diagonal blocks, i.e.\ 
\[
\spec(\Pi H^{({\cal V}\setminus{\cal W})}\Pi+\overline{\Pi}H^{({\cal V}\setminus{\cal W})}\overline{\Pi}) \setminus \{ 0\} 
=\big( \spec(\Pi H^{({\cal V}\setminus{\cal W})}\Pi)\cup\spec(\overline{\Pi}H^{({\cal V}\setminus{\cal W})}\overline{\Pi}) \big) \setminus \{ 0\} \,,
\]
counted with multiplicity. 
Therefore, for any $1\leq i$, $j\leq N$ there are at least
$i+j$ eigenvalues of the block diagonal larger than $\min\{\lambda_{i}(\overline{\Pi}H^{({\cal V}\setminus{\cal W})}\overline{\Pi}),\, \lambda_{j}(\Pi H^{({\cal V}\setminus{\cal W})}\Pi)\}$ provided this number is positive. 
Hence, we conclude 
\begin{equation} \label{eq:proof_PiBoundVminusW_aux1} 
\min\{\lambda_{1}(\overline{\Pi}H^{({\cal V}\setminus{\cal W})}\overline{\Pi}),\,\lambda_{1+|{\cal W}|}(\Pi H^{({\cal V}\setminus{\cal W})}\Pi)\}
\leq\lambda_{2+|{\cal W}|}(H^{({\cal V}\setminus{\cal W})})+2\|\ol{\Pi} H^{({\cal V}\setminus{\cal W})} {\Pi} \|\,.
\end{equation}
Moreover, using eigenvalue interlacing (Lemma \ref{lem:interlacing1}) and Proposition \ref{pro:spectral_gap} \ref{item:L2}, we obtain
\begin{equation} \label{eq:proof_PiBoundVminusW_aux2} 
\lambda_{2+|{\cal W}|}(H^{({\cal V}\setminus{\cal W})})\leq \lambda_{2}(H^{({\cal V})}) \leq \Lambda(\alpha^{*})+o(1)\,.
\end{equation}
Lemma~\ref{prop:PiBoundVminusW} follows from \eqref{eq:proof_PiBoundVminusW_aux1} and \eqref{eq:proof_PiBoundVminusW_aux2}
provided that we show that there is a constant $c >0$ such that 
\begin{equation} \label{eq:proof_PiBoundVminusW_aux3} 
\lambda_{1+|{\cal W}|}(\Pi H^{({\cal V}\setminus{\cal W})}\Pi) \geq \Lambda(\alpha^*) + c\,, 
\end{equation}
and 
\begin{equation} \label{eq:proof_PiBoundVminusW_aux4} 
\norm{\ol{\Pi} H^{(\cal V\setminus \cal W)} \Pi}=o(1)\,.
\end{equation}
 
For the proof of \eqref{eq:proof_PiBoundVminusW_aux3}, we recall the projections $\Pi_{\cal W}$ from Definition \ref{def:Pi_X} and $Q_0$ from Definition \ref{def:Q_r}.
By the definition of $\f u(x)$ for $x \in \cal W$, we have $\Pi_{\cal W} = Q_0 \Pi_{\cal W}$. 
Hence, \eqref{eq:H_approx_diagonal_u_basis} implies $(H - \cal D)Q_{0} \Pi_{\cal W} = o(1)$. 
Applying $Q_0$ to the last relation yields 
\begin{equation} \label{eq:H_cal_V_setminus_cal_W_Pi_cal_W_equal_cal_D} 
H^{(\cal V\setminus \cal W)} \Pi_{\cal W} = \cal D \Pi_{\cal W} + o(1)
\end{equation} 
 as $Q_0 H Q_0 = H^{(\cal V \setminus \cal W)}$ 
and $Q_0 \cal DQ_{0} = \cal D$. 
Since $\ran \Pi_{\cal W} \subset \ran \Pi$ and $\Pi_{\cal W} \cal D \Pi_{\cal W} = \cal D$, the definition of $\cal D$ and \eqref{eq:H_cal_V_setminus_cal_W_Pi_cal_W_equal_cal_D} imply that 
$\Pi H^{(\cal V\setminus \cal W)} \Pi$ has at least $\abs{\cal W}$ many eigenvalues in 
$[\min_{x\in\cal W} \lambda(x) - o(1), \max_{x \in \cal W} \lambda(x) + o(1)]$. 
Note that $\min_{x \in \cal W} \lambda(x) \geq \Lambda(\alpha^*) + c $ for some constant $c>0$ by
 \eqref{eq:lambda_x_approximated_by_Lambda_alpha_x} and \eqref{eq:def_cal_W}. 

Furthermore, let $\f q$ be as in Proposition~\ref{pro:spectral_gap} \ref{item:spectral_gap_q} with $r = 1$. 
From $\norm{H^{(\cal V\setminus \cal W)} } \leq \norm{H} \lesssim \sqrt{d}$ by \eqref{eq:H_minus_EH_ll_sqrt_d} 
and $\norm{\E H} \lesssim \sqrt{d}$, 
Proposition~\ref{pro:spectral_gap} \ref{item:spectral_gap_q} and 
$H^{(\cal V \setminus \cal W)} \f q = H \f q$, we deduce 
\begin{equation} \label{eq:H_cal_V_setminus_cal_W_w_1} 
H^{(\cal V \setminus \cal W)} \f w_1 = H^{(\cal V\setminus \cal W)} \f q + o(1) 
= H \f q + o(1) = H\f w_1 + o(1)  = \lambda_1(H) \f w_1 + o(1)\,. 
\end{equation} 
Hence, $\Pi H^{(\cal V \setminus \cal W)}\Pi$ has an eigenvalue at $\lambda_1(H) + o(1) = \sqrt{d}(1 + o(1)) 
\gtrsim \sqrt{d} \gg \max_{x\in \cal W} \lambda(x) + o(1)$ 
by Proposition~\ref{pro:spectral_gap} \ref{item:L1}, \eqref{Lambda_alpha_bound} and \eqref{eq:lambda_x_approximated_by_Lambda_alpha_x}. 
Therefore, $\Pi H^{(\cal V\setminus \cal W)} \Pi$ has $1 + \abs{\cal W}$ many eigenvalues 
larger or equal to $\Lambda(\alpha^*) + c$ for some constant $c>0$. 
This proves \eqref{eq:proof_PiBoundVminusW_aux3}. 

Finally, \eqref{eq:proof_PiBoundVminusW_aux4} follows from  
\eqref{eq:H_cal_V_setminus_cal_W_Pi_cal_W_equal_cal_D}, $\ran \cal D \subset \ran \Pi_{\cal W} \subset \ran \Pi$ and 
\eqref{eq:H_cal_V_setminus_cal_W_w_1}. 
This completes the proof of Lemma~\ref{prop:PiBoundVminusW}. 
\end{proof}

\section{Eigenvalue spacing -- proof of Proposition \ref{prop:MNonClustering}}
\label{sec:Smoothness}

We recall the definition of the high-probability event $\Omega$ from 
Definition~\ref{def:Omega}.
In this section we use the notation from \eqref{def_POmega}, as well as the conditional versions
\begin{equation*}
\P_\Omega(A \,|\, \cal F) \deq \P(\Omega \cap A \,|\, \cal F)\,,\qquad \E_\Omega[X \,|\, \cal F] \deq \E[X \ind{\Omega} \,|\, \cal F]\,.
\end{equation*}
Throughout this section, we assume that 
$d$ satisfies \eqref{eq:d_range}, $\mu \in (0,1/3)$ and that  $\eta$ satisfies
\begin{equation} \label{eta_assumption}
0 < \frac{\eta}{2} < \min \bigg\{ \frac{1}{6}, \frac{1}{5} - \frac{\mu}{4}, \frac{1}{3} - \mu \bigg\}\,. 
\end{equation}
Proposition~\ref{prop:MNonClustering} follows directly from the following result.

\begin{proposition}\label{prop:NonConcentrationFixA}  
For any $a\neq b\in [N]$, we have 
\begin{equation} \label{P_conc_61}
\mathbb{P}_\Omega(a, b \in \cal W ,\, |\lambda(a) -\lambda(b)|\leq N^{-\eta} )\leq N^{-2 + 2\mu -\eta/4+o(1)}\,. 
\end{equation}
\end{proposition}

\begin{remark} \label{rem:critical61}
If we restrict ourselves to the critical regime $d \asymp \log N$, then Proposition \ref{prop:NonConcentrationFixA} can be improved by replacing the factor $N^{-\eta}$ inside the probability in \eqref{P_conc_61} by $N^{-\eta/2}$. See Remark \ref{rem:critical_eta} below for more details.
\end{remark}

\begin{proof}[Proof of Proposition \ref{prop:MNonClustering}]
The condition \eqref{eta_assumption} holds by assumptions on $\mu$ and $\eta$. 
A union bound and Proposition~\ref{prop:NonConcentrationFixA} yield 
\begin{align*}
\mathbb{P}(\exists x \neq y\in \mathcal W \col |\lambda(x) -\lambda(y)|\leq N^{-\eta} ) & \leq \mathbb{P}(\Omega^c)+\binom{N}{2} \sup_{a \neq b\in[N]}\P_\Omega \big( a,b\in \mathcal W, \abs{\lambda(a) - \lambda(b)} \leq N^{-\eta}  \big) \\ 
& \leq o(1) + N^{2\mu-\eta/4+o(1)}\,, 
\end{align*}
where we used that $\P(\Omega^c) = o(1)$ by definition of $\Omega$ (recall Definition \ref{def:Omega}). 
As $\eta > 8 \mu$, we conclude that the right-hand side is $o(1)$.
\end{proof}

\subsection{Key tools of the proof of Proposition \ref{prop:NonConcentrationFixA}}
The rest of this section is devoted to the proof of Proposition \ref{prop:NonConcentrationFixA}.
Throughout, we fix deterministic vertices $a \neq b \in [N]$ and suppose that $\eta$ satisfies \eqref{eta_assumption}. We use the following definitions. Let
\begin{equation} \label{eq:choice_r} 
r \deq \bigg\lfloor \frac{\eta}{2} \, \frac{\log N}{\log d}\bigg\rfloor - 1\,. 
\end{equation}
In particular, $d^{r + 1} \leq N^{\eta/2} < d^{r+2}$.
Note that $r$ from \eqref{eq:choice_r} satisfies the condition of Proposition~\ref{prop:graphProperties_every_vertex} \ref{item:upper_bound_size_B_r} and \eqref{eq:condition_r}, i.e.\ the condition of Proposition~\ref{prop:graphProperties}. 

\begin{figure}[!ht]
\begin{center}
{\footnotesize 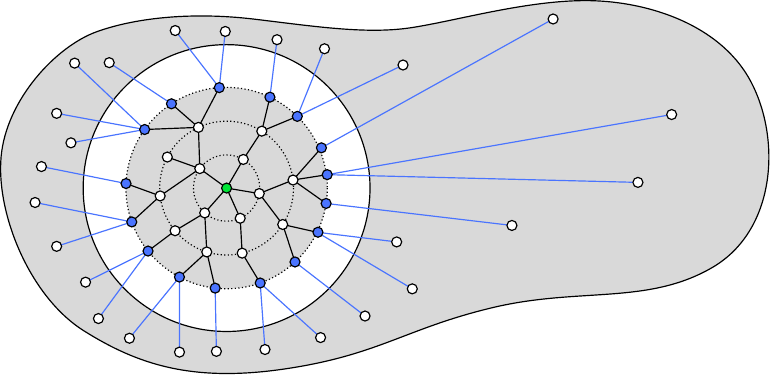}
\end{center}
\caption{An illustration of the $\sigma$-algebra $\cal F_i$. Here $i = 3$, and the vertex $b$ is drawn in green. Conditioning on $\cal F_i$ means that the graph is fixed in the ball $B_i(b)$ and its complement, drawn in grey. The only randomness is the choice of the edges from $S_i(b)$ to $B_i(b)^c$, drawn in blue. By Lemma \ref{lem:FconditionLaw}, these edges are chosen independently with probability $d/N$.
\label{fig:calF}}
\end{figure}

\begin{definition} \label{def:calFi}
We define the $\sigma$-algebra
\begin{equation*}
\cal F_i \deq \sigma(B_{i}(b), A |_{B_{i}(b)}, A |_{B_{i}(b)^c})
\end{equation*}
for $0 \leq i\leq r$, and we abbreviate $\cal F \equiv \cal F_r$. 
\end{definition}

More explicitly, we define $\cal F_i$ inductively through the filtration $\cal G_1 \subset \cal G_2 \subset \cdots \subset \cal G_i$, where $\cal G_1 = \sigma(A |_{\{b\}})$ and $\cal G_{k+1} = \sigma(\cal G_k, A |_{S_k(b)})$, since by construction $S_k(b)$ is $\cal G_k$-measurable. Then we set $\cal F_i = \sigma(\cal G_i, A |_{B_{i}(b)^c})$, using that $B_i(b)$, and hence $B_i(b)^c$, is $\cal G_i$-measurable. See Figure \ref{fig:calF} for an illustration of $\cal F_i$. The following lemma is an immediate consequence of the definition of $\cal F_i$ and the independence of the family $(A_{xy})$.

\begin{lemma} \label{lem:FconditionLaw}
Conditionally on $\cal F_i$, the random variables $(A_{xy} \col x \in S_i(b), y \in B_i(b)^c)$ are independent Bernoulli random variables with mean $d/N$.
\end{lemma}

\begin{definition} \label{def:H(i)}
For $0 \leq i \leq r$ define $H(i) \deq H^{(B_i(b) \cup \cal V^{(B_i(b))})}$ and $G(i,z) \deq (H(i) - z)^{-1}$.
\end{definition}

By definition, $H(i)$ is $\cal F_i$-measurable.
We shall need the following regularization of the function $t \mapsto t^{-1}$.

\begin{definition} \label{def_iota}
Define $\iota \colon \R \to \R$ through  
\begin{equation} \label{eq:def_iota} 
\iota (t) \deq \begin{cases}t^{-1} & \text{if }t\in [T^{-1},T] \\ -t+ T+T^{-1}  & \text{otherwise}\,,
\end{cases}
\end{equation}
with $T \deq 10 \max\{\sqrt{d^{-1}\log N} , \kappa^{-1} \}$.
\end{definition}

\begin{remark} \label{rem:properties_of_iota}  
The function $\iota$ is an involution on $\R$ with Lipschitz constant $T^2$. 
\end{remark} 

\begin{definition}
For $x \in S_i(b)$ we define $S_1^+(x) \deq S_1(x) \cap S_{i+1}(b)$.
\end{definition}

\begin{definition} \label{def:gx}
For each $z \in \R$, we define the family $(g_x(z))_{x\in B_r(b)\setminus \{b\}}$ recursively through
\begin{equation}\label{eq:idealDef}
\begin{cases}  g_x(z) = - \iota \Big( z + \frac{1}{d} \sum_{y \in S_1^+(x) \setminus \cal V^{(B_r(b))}} G_{yy}(r,z)\Big) & \text{ if }x \in S_r(b) \\
g_x(z) = -\iota\Big(z + \frac{1}{d} \sum_{y \in S_1^+(x)} g_y(z) \Big)  &  \text{ if }x \in B_{r - 1}(b)\setminus \{b\}\,. 
\end{cases}
\end{equation} 
\end{definition}

The following remark contains a crucial independence property of the family $(g_x(z))_x$. 

\begin{remark} \label{rem:g_independence}
Conditioned on $\cal F$, 
if $\mathbb{G}|_{B_r(b)}$ is a tree and $z$ is $\cal F$-measurable, then for  any $1 \leq i \leq r$ the family $(g_x(z))_{x \in S_i(b)}$ is independent.
To see this, we first note that $(g_x(z))_{x \in S_r(b)}$ is independent conditioned on $\cal F$ because of Lemma \ref{lem:FconditionLaw}. 
For $i \leq r - 1$, the statement follows inductively by using the tree structure of $\mathbb{G}|_{B_r(b)}$.
\end{remark}

\begin{remark}
The function $\iota \colon \R \to \R$ is a regularized version of the function $t \mapsto t^{-1}$ as a function $\R_+ \to \R_+$. The regularization acts on both small values of $t$ (for $t < T^{-1}$) and large values of $t$ (for $t > T$). The former regularization is needed to ensure the Lipschitz continuity of the function $\iota$, which is used in the proof of Proposition \ref{prop:idealRigidity} to ensure the stability of $g_x(z)$ under change of the argument $z$. The latter regularization is needed to ensure the Lipschitz continuity of the function $\iota^{-1}$, which is used in the proof of Proposition \ref{prop:InsideRecursion} below to ensure that anticoncentration of a random variable is preserved, up to a factor $T^2$, by applying $\iota$. In the proof of Proposition~\ref{prop:idealRigidity} below, we show and use that with high probability the argument of $\iota$ is always contained in the interval $[T^{-1}, T]$ where $\iota$ coincides with $t \mapsto t^{-1}$. Moreover, we choose the lower and upper bounds, $T^{-1}$ and $T$, to be each other's inverses for convenience, since in that case $\iota$ is an involution. Actually, an inspection of our proof shows that the lower bound $T^{-1}$ could be replaced with the larger value $\kappa/10$. Finally, we note that in the critical regime $d \asymp \log N$, the parameter $T$ is of order one. This observation can be used to improve Proposition \ref{prop:NonConcentrationFixA} somewhat in that regime; see Remarks~\ref{rem:critical61} and~\ref{rem:critical_eta}.
\end{remark}

We now state the three key propositions that underlie the proof of Proposition \ref{prop:NonConcentrationFixA} -- Propositions \ref{prop:idealConcentration}, \ref{prop:comparison_lambda_x_lambda_x_complement_B_i}, and \ref{prop:idealRigidity}. Their proofs are postponed to Sections \ref{sec:proof_idealconcentration}, \ref{sec:proof_prop_comparison_lambda_x_lambda_x_complement_B_i}, and \ref{sec:proofidealrigidity}, respectively.

Using the independence from Remark \ref{rem:g_independence}, we obtain the following anticoncentration estimate for the family $(g_x(z))_{x\in B_r(b)\setminus\{b\}}$.

\begin{proposition}[Anticoncentration of $g_x$] \label{prop:idealConcentration} 
Let $z \in \cal J$ be $\cal F$-measurable and $(g_x(z))_x$ be defined as in Definition \ref{def:gx}. Then for any $a,b \in [N]$ we have
\begin{equation} \label{anticoncent_estimate}
\mathbb{P}_\Omega\bigg(a, b \in \cal W, \, \absbb{\frac{1}{d}\sum_{x\in S_1(b)} g_{x}(z) +z} \leq N^{-\eta} \bigg) \leq N^{-2 + 2\mu-\eta/4+o(1)}\,.
\end{equation}
\end{proposition}

For each $0 \leq i \leq r$ we shall need the following $\cal F_i$-measurable approximation of $\lambda(a)$.

\begin{definition}
For $i \in [r]$ we abbbreviate $\lambda(a,i) \deq \lambda_2(H^{((\cal V^{(B_{i}(b))} \cup B_{i}(b)) \setminus \{a\})})$.
\end{definition}

By definition, $\lambda(a,i)$ is $\cal F_i$-measurable. The next results states that $\lambda(a,i)$ is with high probability close to $\lambda(a)$.

\begin{proposition}[Comparison of $\lambda(a)$ and $\lambda(a,i)$] \label{prop:comparison_lambda_x_lambda_x_complement_B_i} 
For any small enough $\epsilon>0$, we have  
\[
\mathbb{P}_\Omega \pb{a , b \in \cal W, \, \exists  i \in [r] ,\, \abs{  \lambda(a)-\lambda(a,i) }\geq \epsilon}\leq \epsilon^{-2} N^{-3 +  2\mu+\eta/2+o(1)}\,.  
\]
\end{proposition}

The next result states that, when choosing the spectral parameter $z = \lambda(a,r)$, the Green function entries of $H^{(\cal V)}$ on $S_1(b)$ are well approximated by the family $(g_x)$ from Definition \ref{def:gx}.

\begin{proposition}[Approximation of Green function by $g_x$] \label{prop:idealRigidity} 
Let $(g_x(z))_x$ be defined as in Definition~\ref{def:gx} with $z \deq \lambda(a,r)$. 
For any constant $c > 0$ and any $\epsilon \leq N^{-c}$, we have 
\[\mathbb{P}_\Omega \pB{a,b \in \cal W,\, \exists x,y\in S_1(b), \, |(H^{(\cal V)}-z)^{-1}_{xy} -g_x(z)\ind{x=y}|\geq \epsilon}\leq \epsilon^{-2} N^{-3 + 2\mu+ 2 \eta + o(1)}\,.\] 
\end{proposition}

\subsection{Proof of Proposition \ref{prop:NonConcentrationFixA}}

In this subsection we prove Proposition \ref{prop:NonConcentrationFixA}. We begin by introducing the following events that we use throughout this section. Recall the definitions of the sets $\cal V^{(X)}$ and $\cal W^{(X)}$ from \eqref{eq:def_mathcal_V_X} and \eqref{eq:def_mathcal_W_X}.

\begin{definition}
We define
\begin{equation*}
\Xi \deq \{ a, b \in \cal W \}\,,
\end{equation*}
and, $1 \leq i \leq r$,
$$ \Xi_i \deq \{a\in \cal W^{(B_i(b))},b\in \cal W\}\,. $$
\end{definition}

\begin{remark} \label{rem:b_in_cal_W_measurable_Xi_i}
We record the following straightforward properties of $\Xi$ and $\Xi_i$.
\begin{enumerate}[label=(\roman*)]
\item \label{itm:W1}
$\Xi_i$ is $\cal F_i$-measurable.
\item \label{itm:W2}
$\Xi_i \subset \Xi$ (by definition of $\cal W^{(B_i(b))}$).
\item \label{itm:W3}
On $\Omega$, for any $b \in \cal W$ and $1 \leq i \leq r$ we have $\cal V^{(B_i(b))} = \cal V \setminus \{b\}$ and $\cal W^{(B_i(b))} = \cal W\setminus \{b\}$ (see Proposition~\ref{prop:graphProperties} \ref{item:B_r_disjoint}). In particular,  $\Omega \cap \Xi_i = \Omega \cap \Xi$ for all $1 \leq i \leq r$.
\end{enumerate}
\end{remark}

Next, we state a basic result that, together with Lemma~\ref{lem:G_bounded}, is used throughout this section to establish the boundedness of the Green function for certain spectral parameters.

\begin{lemma} \label{lem:lambda_i_in_D} 
On $\Omega \cap \Xi$ we have $\lambda(a,i) = \Lambda(\alpha_a) + o(1)$ and $\lambda(a,i) \in \cal J$ for all $1 \leq i \leq r$.
\end{lemma}

\begin{proof}
By Remark \ref{rem:b_in_cal_W_measurable_Xi_i} \ref{itm:W3}, on $\Omega$ the assumptions of Corollary~\ref{cor:largest_eigenvalues_vertices_removed} are satisfied for $x = a$ and 
$X =( \cal V^{(B_{i}(b))} \cup B_{i}(b)) \setminus \{a\}$. 
Therefore, $\lambda(a,i) = \Lambda(\alpha_a)+o(1)$, from which we conclude that $\lambda(a,i) \geq \Lambda(\alpha^*) + \kappa/4$ by the definition \eqref{eq:def_cal_W} of $\cal W$ as well as $\lambda(a,i) \ll \sqrt{d}$ by \eqref{Lambda_alpha_bound}.
\end{proof}

\begin{proof}[Proof of Proposition \ref{prop:NonConcentrationFixA}]
By spectral decomposition of $H^{(\cal V \setminus \{b\})}$, we have $\im (H^{(\cal V \setminus \{b\})} - \lambda(b) - \ii t)_{bb}^{-1} \geq \scalar{\f 1_b}{\f u(b)}^2 / t$ for any $t > 0$.
By Corollary \ref{cor:lower_bound_u_x_in_x}, on $\Omega$ we have $\scalar{\f 1_b}{\f u(b)} \neq 0$, which implies
\begin{equation*}
\lim_{t \downarrow 0} \frac{1}{(H^{(\cal V \setminus \{b\})} - \lambda(b) - \ii t)^{-1}_{bb}} = 0\,,
\end{equation*}
and hence Schur's complement formula yields
\[
\lambda(b) + \frac{1}{d} \sum_{x,y \in S_1(b)} ( H^{(\cal V)}-\lambda (b))_{xy}^{-1} =0\,. \] 
Therefore, with $z=\lambda (a,r)$, we obtain from the definition of the family $(g_x(\lambda(a,r)))_x$ in Definition~\ref{def:gx} that 
\begin{align*}
 \absbb{\lambda (a,r) + \frac{1}{d}\sum_{x\in S_1(b)} g_{x}(\lambda(a,r))} 
  & \leq |\lambda(a,r) - \lambda(b)|    
+ \frac{1}{d} \sum_{x,y \in S_1(b)} \absb{(H^{(\cal V)}-\lambda(a,r))^{-1}_{xy} - (H^{(\cal V)}-\lambda(b))^{-1}_{xy}} 
\\ &\quad  + \frac{1}{d} \sum_{x,y \in S_1(b)} |g_x(\lambda(a,r))\ind{x=y}- (H^{(\cal V)}-\lambda(a,r))^{-1}_{xy}| 
\\  &\leq  C\kappa^{-2}(\log N)^2 (|\lambda (a)-\lambda(b)|+ |\lambda (a,r)-\lambda(a)|) \\ 
 & \quad  + \frac{1}{d} \sum_{x,y \in S_1(b)} |g_x(\lambda(a,r))\ind{x=y} -  ( H^{(\cal V)}-\lambda(a,r))^{-1}_{xy}|  
\end{align*}
on the event $\Omega$, where we used \eqref{eq:G_error_z_zi} and Proposition \ref{prop:graphProperties} \ref{item:upper_bound_alpha_x} for the second inequality. Here $C$ is some positive constant.
Thus, for any $\gamma\geq 0$, we obtain 
\begin{align*}
&\mspace{-20mu}\mathbb{P}_{\Omega \cap \Xi}(C\kappa^{-2}(\log N)^2|\lambda(b)-\lambda(a)|\leq \gamma) 
\\
  &\leq  \mathbb{P}_{\Omega \cap \Xi}\bigg(\absbb{\frac{1}{d}\sum_{x\in S_1(b)} g_{x}(\lambda(a,r)) +\lambda(a,r)}\leq 3\gamma\bigg) 
\\
&\quad  + \mathbb{P}_{\Omega\cap \Xi}\big(C\kappa^{-2}(\log N)^2|\lambda(a)-\lambda(a,r)|\geq \gamma\big) 
\\ 
&\quad
+ \mathbb{P}_{\Omega\cap \Xi}\bigg(\frac{1}{d} \sum_{x,y \in S_1(b)} \absB{( H^{(\cal V)}-\lambda(a,r))^{-1}_{xy}-g_x(\lambda(a,r))\ind{x=y}}\geq \gamma\bigg)\,.   
\end{align*} 
Now Proposition \ref{prop:NonConcentrationFixA} follows with the choice $\gamma \deq N^{-\eta}/3$, applying Proposition \ref{prop:idealConcentration} to the first line, Proposition \ref{prop:comparison_lambda_x_lambda_x_complement_B_i} to the second line, and Proposition \ref{prop:idealRigidity} combined with Proposition \ref{prop:graphProperties} \ref{item:upper_bound_alpha_x} on $\Omega$ to the third line. Here we used Lemma \ref{lem:lambda_i_in_D} to ensure that $\lambda(a,i) \in \cal J$.
\end{proof}

\subsection{Proof of Proposition~\ref{prop:comparison_lambda_x_lambda_x_complement_B_i}} 
\label{sec:proof_prop_comparison_lambda_x_lambda_x_complement_B_i} 

In this subsection we prove Proposition~\ref{prop:comparison_lambda_x_lambda_x_complement_B_i}.

\begin{proof}[Proof of Proposition~\ref{prop:comparison_lambda_x_lambda_x_complement_B_i}]
We follow the proof of Proposition \ref{prop:EstimateExy}. 
Let $\f u(a,i)$ be a normalized eigenvector of $H^{((\cal V^{(B_i(b))}\cup B_i(b))\setminus \{ a \})}$ associated with the 
eigenvalue $\lambda(a,i)$. 
Then, as $\supp \f u(a,i) \subset (\cal V^{(B_i(b))} \cup B_i(b))^c \cup \{a \}$, 
on the event $\Omega \cap \{ a, b \in \cal W\}$, using Remark \ref{rem:b_in_cal_W_measurable_Xi_i} \ref{itm:W3} we obtain 
\begin{equation} \label{Hlambda_split}
(H^{(\cal V\setminus{\{a\}})}-\lambda(a,i))\f u(a,i)= (H^{(\cal V \setminus \{a \})} - H^{((\cal V^{(B_i(b))}\cup B_i(b))\setminus \{ a \})})\f u(a,i)=\sum_{y\in S_{i}({b})}\epsilon_{ay}(i) \f 1_{y}\,, 
\end{equation}
where $\epsilon_{ay}(i) \deq \frac{1}{\sqrt{d}} \sum_{v\in S_1^+(y)}\scalar{\f 1_v}{\f u(a,i)}$.
From Lemma \ref{lem:FconditionLaw} we conclude, using Cauchy-Schwarz,
\begin{align}
\mathbb{E}_\Omega\Bigg[\sum_{y\in S_{i}({b})}(\epsilon_{ay}(i))^{2} \,\bigg|\,{\cal F}_{i}\Bigg]
&\lesssim \sum_{y \in S_i(b)} \frac{\log N}{d} \sum_{v \in B_i(b)^c} \E \qb{ \ind{v \in S_1^+(y) } \scalar{\f 1_v}{\f u(a,i)}^2 \, \big|\, \cal F_i}
\notag \\
&=
\sum_{y \in S_i(b)} \frac{\log N}{d} \sum_{v \in B_i(b)^c} \P(v \in S_1^+(y) \,|\, \cal F_i) \,  \scalar{\f 1_v}{\f u(a,i)}^2
\notag \\ \label{eq:eps_sq_estimate}
&\leq \frac{\log N}{N} \abs{S_i(b)}\,,
\end{align}
where we used that on $\Omega$ we have $\abs{S_1^+(y)} \lesssim \log N$, that $\f u(a,i)$ is $\cal F_i$-measurable, and that by Lemma \ref{lem:FconditionLaw} we have $\P(v \in S_1^+(y) \,|\, \cal F_i) = \frac{d}{N}$ for any $y \in S_i(b)$ and $v \notin B_i(b)$.

By the definition \eqref{eq:def_cal_W} of $\cal W$, we have $\Lambda(\alpha_a) \geq \Lambda(\alpha^*) + \kappa/2$, and on $\Omega$ we have $\Lambda(\alpha_a) \ll \sqrt{d}$ (recall \eqref{Lambda_alpha_bound}).
Hence, by Corollary~\ref{cor:largest_eigenvalues_vertices_removed}, for any small enough $\epsilon > 0$, if there are a scalar $\wh{\lambda}$ and a normalized vector $\wh{\f u}$ such that, 
for some $a \in \cal W$,  
$\norm{(H^{(\cal V\setminus \{ a\})} - \wh{\lambda})\wh{\f u}} \leq \epsilon$ and  
 $\wh{\lambda} = \Lambda(\alpha_a) + o(1)$, then 
$\abs{\lambda(a) - \wh{\lambda}} \leq \epsilon$ (as $\lambda(a) = \lambda_2(H^{(\cal V\setminus \{a \})})$).

We apply this observation to the choices $\wh{\lambda} = \lambda(a,i)$ and $\wh{\f u} = \f u(a,i)$, for which Lemma~\ref{lem:lambda_i_in_D} yields $\lambda(a,i) = \Lambda(\alpha_a) + o(1)$. 
Hence, for any small enough $\epsilon > 0$, we have $\abs{\lambda(a) - \lambda(a,i)} \leq \epsilon$ provided that 
$\norm{(H^{(\cal V\setminus \{a\})}- \lambda(a,i)) \f u(a,i)} \leq \epsilon$. 
We can then estimate
\begin{align*}
 \mathbb{P}_\Omega(a,b \in \cal W, \,| \lambda(a)-\lambda(a,i)|> \epsilon) 
& \leq \mathbb{P}_\Omega\left(\Xi_i \cap \hB{\normb{( H^{(\cal V\setminus{\{a\}})}-\lambda(a,i) )\f u(a,i)}^2> \epsilon^2}\right)
\\ & \leq \epsilon^{-2} \,\mathbb{E}\Bigg[\ind{\Xi_i} \ind{\abs{S_i(b)} \leq N^{\eta/2 + o(1)}} \mathbb{E}_\Omega\Bigg[\sum_{y\in S_{i}({b})}(\epsilon_{ay}(i))^{2}\,\bigg|\,{\cal F}_{i}\Bigg] \Bigg]
\\
&\leq \epsilon^{-2} \, \P(\Xi_i) \, \frac{\log N}{N} N^{\eta/2 + o(1)}
\\
&\leq \epsilon^{-2} \, N^{2\mu+\eta/2-3+o(1)}\,, 
\end{align*}
where in the first step we used Remark \ref{rem:b_in_cal_W_measurable_Xi_i} \ref{itm:W3}, in the second step \eqref{Hlambda_split}, Remark \ref{rem:b_in_cal_W_measurable_Xi_i} \ref{itm:W1}, and the estimate $\abs{S_i(b)}  \leq N^{\eta/2 + o(1)}$ on $\Omega$ by Proposition~\ref{prop:graphProperties} \ref{item:upper_bound_size_B_r} and \eqref{eq:choice_r}, in the third step \eqref{eq:eps_sq_estimate}, and the fourth step Remark \ref{rem:b_in_cal_W_measurable_Xi_i} \ref{itm:W2} and Lemma \ref{lem:Wab_estimate} below.
Now Proposition~\ref{prop:comparison_lambda_x_lambda_x_complement_B_i} follows from a union bound over $i \in [r]$ with $r \lesssim \log N$.
\end{proof}

The following result is used throughout the rest of this section.

\begin{lemma} \label{lem:Wab_estimate}
For any $a \neq b \in [N]$ we have $\P(a,b \in \cal W) \leq N^{- 2 + 2 \mu}$.
\end{lemma}

\begin{proof} 
Since $0\leq \Lambda'(\alpha)=\frac{\alpha-2}{2(\alpha-1)^{3/2}} \leq \frac{1}{2}$ for all $\alpha \geq 2$, we have 
\begin{align*}
P(a,b \in \cal W)  &\leq \P\Big( \min\{|S_1(a)|,|S_1(b)| \} \geq (\alpha^*+\kappa) d \Big)
\\
&= \E \qB{\P\Big( \min\{|S_1(a)|,|S_1(b)| \} \geq (\alpha^*+\kappa) d \,\Big|\, A_{ab}\Big)}
\\ 
&\leq \P\big(|S_1(a)| \geq (\alpha^*+\kappa) d -1\big)^2
\\
&\leq \P \big( \abs{S_1(a)} \geq \alpha^* d \big)^2
\\
&\leq N^{-2+2\mu}\,, 
\end{align*}
where we used in the third step that conditionally on $A_{ab}$, $S_1(a)$ and $S_1(b)$ are independent and in the last step the definition of $\alpha^*$.
\end{proof}

\subsection{Proof of Proposition~\ref{prop:idealRigidity}}  \label{sec:proofidealrigidity}

This subsection is devoted to the proof of Proposition \ref{sec:proofidealrigidity}. We begin with the following result, which contains two estimates. The first one is an approximate version of Schur's complement formula, where $G(i - 1,z)$ is related to $G(i,z)$ at the cost of an error term; this amounts to removing not just the vertex $x \in S_i(b)$ at which the Green function is evaluated, but the entire ball $B_i(b)$. The second estimate provides an upper bound on the off-diagonal entries of the Green function.

\begin{lemma} 
\label{lem:ResolvantBalls}  
{
\begin{enumerate}[label=(\roman*)]
\item \label{itm:resballs1}
For $1 \leq i \leq r$ let $z_i \in \cal J$ be $\cal F_i$-measurable and $x \in S_i(b)$. Then
\begin{equation}\label{eq:schur_balls}
-\frac{1}{G_{xx}(i-1,z_i)} = z_i + \frac{1}{d} \sum_{y \in S_1^+(x)} G_{yy}(i,z_i)  +\cal E_i(x)
\end{equation}
where the error term $\cal E_i(x)$ satisfies, for any $\epsilon > 0$,
\begin{equation} \label{eq:bound_error_term_schur} 
\mathbb{P}_\Omega\left(\left|\cal E_i(x) \right|> \epsilon \,|\, \cal F_i \right)
\ind{b \in \cal W} 
\lesssim   d^2 \kappa^{-4} (\log N)^2 \abs{S_i(b)}^2 N^{-1} \epsilon^{-2} \,.  
\end{equation}
\item \label{itm:resballs2}
Let $z \in \cal J$ be $\cal F_1$-measurable. For any $x \neq y \in [N]$ and $\epsilon > 0$, we have 
\begin{equation} \label{eq:bound_off_diagonal} 
\P_\Omega \big( \abs{ (H^{(\cal V)} - z)^{-1}_{xy}} \geq \epsilon \,|\, \cal F_1 \big) \ind{b \in \cal W} \ind{x,y \in S_1(b)} \lesssim \kappa^{-4} (\log N) N^{-1} \epsilon^{-2}\,. 
\end{equation} 
\end{enumerate}
}
\end{lemma}

\begin{proof}[Proof of Proposition \ref{prop:idealRigidity}]
We choose $z_i \deq \lambda(a,i)$, which is $\cal F_i$-measurable, and set $z = z_r$. 
For $\epsilon > 0$ we introduce the event $\Theta \defeq \Theta_1 \cap \Theta_2 \cap \Theta_3$ given by  
\begin{align*}
\Theta_1 & \defeq \{\forall i\leq r, \forall x\in S_i (b) ,  |\cal E_i(x)|\leq \epsilon\}\,, \\ 
\Theta_2 & \defeq \{\forall i\leq r, |z-z_i|\leq \epsilon\}\,, \\ 
\Theta_3 & \defeq \{ \forall x\neq y \in S_1(b), \abs{(H^{(\cal V)}-z_1)^{-1}_{xy}} \leq \epsilon \}  \,,
\end{align*}
with $\cal E_i(x)$ defined as in Lemma \ref{lem:ResolvantBalls}. We estimate the probability of $\Theta_1^c$ using $\Xi_i \in \cal F_i$, Lemmas \ref{lem:ResolvantBalls} \ref{itm:resballs1} and \ref{lem:Wab_estimate}, Proposition~\ref{prop:graphProperties} \ref{item:upper_bound_size_B_r}, we well as Remark \ref{rem:b_in_cal_W_measurable_Xi_i} as
\begin{equation*}
\mathbb{P}_{\Omega\cap \Xi_i}(\exists x\in S_i (b) ,  |\cal E_i(x)| > \epsilon) \leq \mathbb{P}_{\Omega}(\Xi) N^{-1 + 3\eta/2 + o(1)} \epsilon^{-2} \leq N^{-3 + 2 \mu + 3\eta/2 + o(1)} \epsilon^{-2}\,.
\end{equation*}
Similarly, we find using Lemma \ref{lem:ResolvantBalls} \ref{itm:resballs2} that
\begin{equation*}
\P_{\Omega \cap \Xi_1}(\Theta_3^c) \leq \P_\Omega(\Xi) N^{-1 + \eta + o(1)} \epsilon^{-2} \leq N^{-3 + 2 \mu + \eta/2 + o(1)} \epsilon^{-2}\,.
\end{equation*}
Hence, using Remark~\ref{rem:b_in_cal_W_measurable_Xi_i} and Proposition~\ref{prop:comparison_lambda_x_lambda_x_complement_B_i}, we have
\begin{align}
\mathbb{P}_\Omega(\Theta^c\cap \Xi ) &\leq \sum_{i=1}^r \bigg( \mathbb{P}_{\Omega\cap \Xi_i}(\exists x\in S_i (b) ,  |\cal E_i(x)| > \epsilon)+\mathbb{P}_{\Omega\cap \Xi} (|z-z_i| > \epsilon) \bigg) + 
\P_{\Omega \cap \Xi_1} \big ( \Theta_3^c\big ) \nonumber\\ 
& \leq N^{-3 + 2 \mu +2\eta+o(1)} \epsilon^{-2}\,.
\label{eq:P_Omega_Theta_complement_cap_Xi} 
\end{align} 

We shall show below that there is a constant $C >0$ such that on the event $\Omega \cap \Theta \cap \Xi$ we have
\begin{equation} \label{eq:g_x_approx_G_xx_induction}
|g_{x}(z)-G_{xx}(i  - 1 ,z)|\leq C \epsilon \sum_{j=0 }^{r-i+1} (C\kappa^{-4}d^{-1})^{j}|S_{j}^{+}(x)|\leq N^{o(1)} \epsilon
\end{equation} 
for all $1 \leq  i \leq r$ and all $x \in S_i(b)$.
The second inequality in \eqref{eq:g_x_approx_G_xx_induction} follows from $r\ll \log N$ and $d^{-j}|S_j^+(x)|\leq d^{-j}|B_j(x)|\lesssim \log N$ for all $x\in [N]$ by  Proposition~\ref{prop:graphProperties} \ref{item:upper_bound_size_B_r}.

Before proving \eqref{eq:g_x_approx_G_xx_induction}, we conclude Proposition~\ref{prop:idealRigidity} 
from \eqref{eq:g_x_approx_G_xx_induction} (after renaming $\epsilon \mapsto N^{-o(1)}\epsilon$) with $i = 1$, the definition of $\Theta_3$ 
and \eqref{eq:P_Omega_Theta_complement_cap_Xi}, 
where we used that $G(0,z) = (H^{(\cal V)}-z)^{-1}$ on $\Omega \cap \{ b \in \cal W\}$ as $\cal V^{(B_i(b))} = \cal V \setminus \{ b\}$.  

What remains, therefore, is the proof of \eqref{eq:g_x_approx_G_xx_induction}. We prove it by inductively decreasing $i$ starting from $i = r + 1$. 
By convention, for any $x\in S_{r+1}(b)$ we denote $g_x(z)\defeq G_{xx}(r,z)$, so that \eqref{eq:g_x_approx_G_xx_induction} trivially holds for $i = r+1$. 
Therefore, we can assume throughout the following argument that $g_x(z)$ is defined by the second case in Definition \ref{def:gx} for all $x \in B_r(b) \setminus \{ b\}$ (note that on the event $\Omega \cap \Xi$ we have $S_1^+(x) \setminus \cal V^{(B_r(b))} = S_1^+(x)$ for any $x \in S_r(b)$, by Proposition \ref{prop:graphProperties} \ref{item:B_r_disjoint}).

To verify the induction step, we assume that \eqref{eq:g_x_approx_G_xx_induction} holds on $S_j(b)$ for all $i+1 \leq j \leq r + 1$ and consider $x\in S_i(b)$. We first show that on $\Omega \cap \Xi \cap \Theta$
\begin{equation} \label{eq:g_x_defined_by_Schur} 
 g_x(z) = -\bigg(z+\frac{1}{d}\sum_{y\in S_{1}^{+}(x)}g_{y}(z)\bigg)^{-1}\,. 
\end{equation} 
To that end, we conclude from the induction hypothesis, \eqref{eq:G_error_z_zi}, the definition of $\Theta_2$ 
and Proposition~\ref{prop:graphProperties} \ref{item:upper_bound_alpha_x} 
that 
 \begin{align}
 \absbb{\frac{1}{d}\sum_{y\in S_{1}^{+}(x)}G_{yy}(i,z_{i})-\frac{1}{d}\sum_{y\in S_{1}^{+}(x)}g_{y}(z)}
 & \leq\frac{1}{d}\sum_{y\in S_{1}^{+}(x)}\left(|G_{yy}(i,z)-g_{y}(z)|+|G_{yy}(i,z_{i})-G_{yy}(i,z)|\right) \nonumber \\
 & \leq d^{-1}\abs{S_1(x)} \bigg(\max_{y \in S_1^+(x)} |G_{yy}(i,z)-g_{y}(z)|+(8/\kappa)^{2}|z-z_{i}|\bigg) \nonumber \\
 & \leq N^{o(1)} \epsilon\,. \label{eq:upper_bound_difference_sum_G_sum_g} 
\end{align}
By Lemma~\ref{lem:lambda_i_in_D}, on $\Omega \cap \Xi$ we have $z_i \in \cal J$.
We recall from Lemma~\ref{lem:lowboundGreen} that all Green function entries $G_{xx}(i-1,z_i)$ and $G_{yy}(i,z_{i})$ are negative for $z_i \in \cal J$. We use the upper bound \eqref{eq:bound_resolvent} for $-G_{xx}(i-1,z_i)$ as well as \eqref{eq:schur_balls} to obtain, on $\Omega \cap \Xi \cap \Theta_1 \cap \Theta_2$,
\begin{equation} \label{eq:upper_lower_bound_Schur_complement_expansion} 
\frac{10}{9T}\leq  \frac{\kappa}{9}\leq -\frac{1}{G_{xx}(i-1,z_i)}-\cal E_i(x)= z_{i}+\frac{1}{d}\sum_{y\in S_{1}^{+}(x)}G_{yy}(i,z_{i})\leq\Lambda(\alpha_{a})+ o(1) \leq \sqrt{2 \alpha_a} + o(1) \leq \frac{T}{2}\,. 
\end{equation}
Here, we also used Corollary~\ref{cor:largest_eigenvalues_vertices_removed}, which is applicable as 
$\cal V^{(B_i(b))} = \cal V \setminus \{b \}$ on $\Omega\cap \{b\in \cal W\}$ by Proposition~\ref{prop:graphProperties} \ref{item:B_r_disjoint}, 
 $\alpha_a \leq 10 d^{-1}\log N$ on $\Omega$ by Proposition~\ref{prop:graphProperties} \ref{item:upper_bound_alpha_x}, and the definitions of $\Lambda$ and $T$. 
Then, by the assumption $\epsilon \leq N^{-c}$, we find that \eqref{eq:upper_lower_bound_Schur_complement_expansion} and \eqref{eq:upper_bound_difference_sum_G_sum_g} 
imply $T^{-1}\leq z+\frac{1}{d}\sum_{y\in S_{1}^{+}(x)}g_{y}(z)\leq T$, which yields 
\eqref{eq:g_x_defined_by_Schur} by the definitions of $\iota$ and $g_x(z)$ in \eqref{eq:def_iota} 
and Definition~\ref{def:gx}, respectively. This concludes the proof of \eqref{eq:g_x_defined_by_Schur}.

Finally, we find on $\Omega \cap \Xi \cap \Theta$
\begin{align*}
& |g_{x}(z)-G_{xx}(i-1,z)|
\\ & \leq g_{x}(z)G_{xx}(i-1,z_i)\left|g_{x}(z)^{-1}-G_{xx}(i-1,z_i)^{-1}\right|+(8/\kappa)^{2}|z-z_i|\\
 & \leq (8/\kappa)^2 \bigg(\absbb{\frac{1}{d}\sum_{y\in S_{1}^{+}(x)}G_{yy}(i,z_{i})-\frac{1}{d}\sum_{y\in S_{1}^{+}(x)}g_{y}(z)}+|z-z_{i}|+\abs{{\cal E}_{i}(x)}\bigg) \\
 & \leq (8/\kappa)^2\bigg(\frac{1}{d}\sum_{y\in S_{1}^{+}(x)} \sum_{j=0 }^{r-i} (C\kappa^{-4}d^{-1})^{j}|S_{j}^{+}(y)|C \epsilon  +((8/\kappa)^2 d^{-1}|S_{1}^{+}(x)|+1)|z-z_{i}|+\abs{{\cal E}_{i}(x)} \bigg)\\
 & \leq C^{-1}\bigg(64 \sum_{j=1}^{r-i+1} (C\kappa^{-4}d^{-1})^{j}|S_{j}^{+}(x)|+64^2\kappa^{-4}(d^{-1}|S_{1}^{+}(x)|+2)\bigg)C \epsilon\,.
\end{align*}
Here we used \eqref{eq:G_error_z_zi}
in the first inequality, 
 \eqref{eq:schur_balls} and \eqref{eq:g_x_defined_by_Schur} in the second, \eqref{eq:bound_resolvent} and  the iteration hypothesis  in the third and that $|S_{j+1}^+(x)| =\sum_{y\in S_{1}^{+}(x)} |S_{j}^{+}(y)|$ and the definitions of $\Theta_1$ and $\Theta_2$ in the last one. We conclude 
\eqref{eq:g_x_approx_G_xx_induction} for large enough $C$.
\end{proof}

What remains is the proof of Lemma \ref{lem:ResolvantBalls}. 

\begin{proof}[Proof of Lemma \ref{lem:ResolvantBalls}]
We begin with \ref{itm:resballs1}.
Throughout the proof, we fix $i$ and we condition on $\cal F_i$. 
We always assume that $b \in \cal W$, which is an $\cal F_i$-measurable event since $i \geq 1$. 

We start by observing that the event
\begin{equation*}
\Gamma_i \deq \{ \cal V^{(B_{i-1}(b))} = \cal V^{(B_i(b))} \} \cap \{A |_{S_i(b)} = 0\}
\end{equation*}
satisfies $\Omega \subset \Gamma_i$, by  Proposition~\ref{prop:graphProperties}  \ref{item:B_r_tree}, \ref{item:B_r_disjoint}. See Figure \ref{fig:Gamma} for an illustration of $\Gamma_i$.

\begin{figure}[!ht]
\begin{center}
{\footnotesize 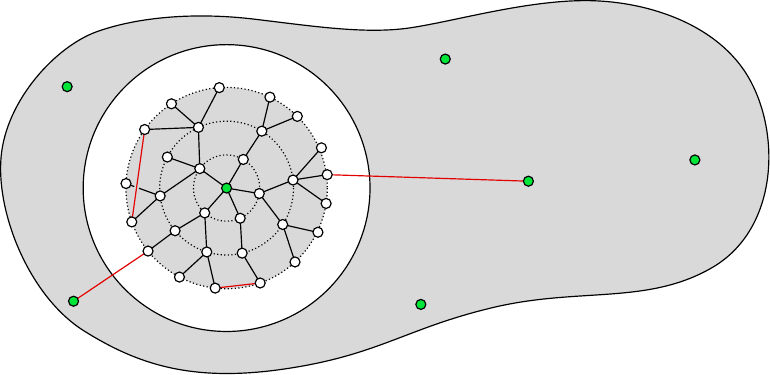}
\end{center}
\caption{An illustration of the event $\Gamma_i$ for $i = 3$. The vertices of $\cal V$ are drawn in green. By definition of $\Gamma_i$, the red edges are forbidden.
\label{fig:Gamma}}
\end{figure}

For the proof we abbreviate
\begin{equation} \label{def_calTi}
\cal T_i \deq  B_{i-1}(b) \cup  \cal V^{(B_{i}(b))}\,,
\end{equation}
so that on the event $\Gamma_i$ we have $H(i-1) = H^{(\cal T_i)}$.
From Schur's complement formula  we get on the event $\Gamma_i$
\begin{equation} \label{Gi_schur}
-\frac{1}{G_{xx}(i-1,z_i)} = z_i + \frac{1}{d} \sum_{u,v \in S_1^+(x)} (H^{(\cal T_i \cup \{x\}
)}-z_i)^{-1}_{uv} \,,  
\end{equation}
where we used that $x$ has no neighbours in $S_i(b)$ by definition of $\Gamma_i$.

We now decompose the error term $\cal E_i(x)$ from \eqref{eq:schur_balls} into several summands estimated separately.
To that end, let $\{ x, y_1, \ldots, y_{\abs{S_i(b)}-1} \} = S_i(b)$ be an enumeration of $S_i(b)$. 
We set $y_0 \deq x$ and 
\begin{align*} 
 \cal E_i^{(0)} & \deq  \frac{1}{d}  \sum_{u,v \in S_1^+(x), \, u \neq v} \big ( H^{(\cal T_i \cup \{ x\})} - z_i\big)^{-1}_{uv}
, \\ 
 \cal E_i^{(j)} & \deq \frac{1}{d}   \sum_{u \in S_1^+(x)} \Big( \big(H^{(\cal T_i \cup \{y_0,  \ldots, y_j\})} - z_i\big)^{-1}_{uu} 
- \big( H^{(\cal T_i \cup \{ y_0, \ldots, y_{j-1}\})} - z_i\big)^{-1}_{uu} \Big)\,.  
 \end{align*} 
for $j = 1, \ldots, \abs{S_i(b)} - 1$. 
Then, from \eqref{Gi_schur} we conclude that \eqref{eq:schur_balls} holds with
\[ 
\cal E_i(x) = \cal E_i^{(0)} - \sum_{j=1}^{\abs{S_i(b)}-1} \cal E_i^{(j)}\,.  
\] 
Chebyshev's and the Cauchy-Schwarz inequalities yield
\begin{equation} \label{eq:bound_prob_cal_E_i} 
\P_\Omega( \abs{\cal E_i(x)} > \epsilon \,|\, \cal F_i) 
\leq \epsilon^{-2} \, \abs{S_i(b)}\bigg( \E_\Omega[\abs{\cal E_i^{(0)}}^2 \,|\,\cal F_i] +  \sum_{j=1}^{\abs{S_i(b)}-1} 
\E_\Omega [ \abs{\cal E_i^{(j)}}^2 \,|\, \cal F_i] \bigg)\,. 
\end{equation}  

Therefore, it remains to estimate $\E_{\Omega} [ \abs{\cal E_i^{(j)}}^2 \,|\, \cal F_i]$ for all $j = 0, \ldots, \abs{S_i(b)}-1$. 
To that end, we introduce a $\sigma$-algebra refining $\cal F_i$ from Definition \ref{def:calFi}.
For any subset $X \subset S_i(b)$ we define the $\sigma$-algebra\footnote{Somewhat more carefully (since the set $X$ is random), the precise definition of $\cal F_i(X)$ is as follows. As the underlying probability space is finite, any $\sigma$-algebra, in particular $\cal F_i$, is atomic. Conditioning on $\cal F_i$ means that we restrict ourselves to a single atom $\cal A$ of $\cal F_i$. On this atom $S_i(b)$ and $X$ are deterministic. Then $\cal F_i(X)$ is by definition the smallest $\sigma$-algebra on $\cal A$ such that $S_1^+(y)$ is measurable for all $y \in S_i(b) \setminus X$.
}
\begin{equation*}
\cal F_i(X) \deq \sigma\pb{\cal F_i, (S_1^+(y))_{y \in S_i(b) \setminus X}}\,,
\end{equation*}
using that $S_i(b)$ is $\cal F_i$-measurable.
See Figure \ref{fig:calFX} for an illustration of $\cal F_i(X)$.
Moreover, for $X \subset S_i(b)$ we define the event 
\begin{equation*}
\Delta_i(X) \deq \hB{\normb{(H^{(\cal T_i \cup X)}-z_i)^{-1}}\leq  8 \kappa^{-1}  }\,.
\end{equation*}
We note that for any $X \subset S_i(b)$, the event $\Delta_i(X)$ lies in $\cal F_i(X)$, since $H^{(\cal T_i \cup X)}$ is $\cal F_i(X)$-measurable (see the definition \eqref{def_calTi} and Figure \ref{fig:calFX}).
Furthermore, by Lemma~\ref{lem:G_bounded} for $z_i \in \cal J$ we have
\begin{equation} \label{Omega_Delta}
\Omega \subset \bigcap_{X \subset S_i(b)} \Delta_i(X)\,.
\end{equation}

\begin{figure}[!ht]
\begin{center}
{\footnotesize 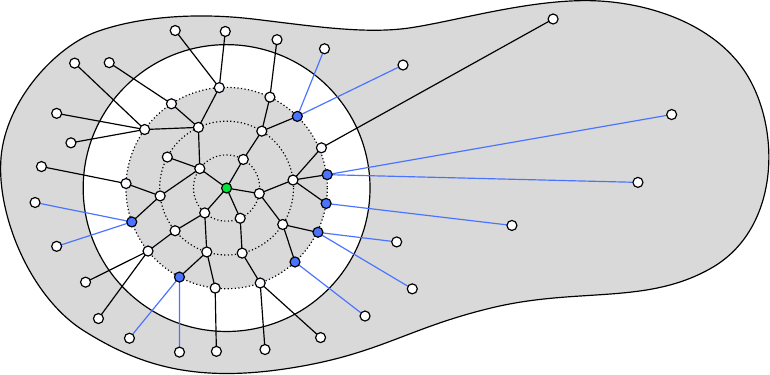}
\end{center}
\caption{An illustration of the $\sigma$-algebra $\cal F_i(X)$. Here $i = 3$, the vertex $b$ is drawn in green, and the set $X \subset S_i(b)$ is drawn in blue. Conditioning on $\cal F_i(X)$ means that we fix all edges within $B_i(b)$, within $B_i(b)^c$, and connecting $S_i(b) \setminus X$ with $B_i(b)^c$. The only randomness is the choice of the edges from $X$ to $B_i(b)^c$, drawn in blue. After removal of the vertices in $\cal T_i \cup X$ (see \eqref{def_calTi}), only black edges and edges within $B_i(b)^c$ remain, which shows that $H^{(\cal T_i \cup X)}$ is $\cal F_i(X)$-measurable. Note that for $X = S_i(b)$ we have $\cal F_i(X) = \cal F_i$, and we recover the illustration from Figure \ref{fig:calF}.
\label{fig:calFX}}
\end{figure}

For the estimate of $\cal E_i^{(0)}$, we introduce the sets 
$\cal Q \deq \{ (u,v) \in (B_i(b)^c)^2 \col u \neq v \}$ and $Q \deq \cal Q \cap (S_1^+ (x))^2$ and 
the family $(Z_q)_{q \in \cal Q}$ defined by 
$Z_{(u,v)} \deq (H^{(\cal T_i \cup \{ x \})} - z_i)^{-1}_{uv}$.  
We first note that, for any $ q = (u,v) \in \cal Q$, we have 
\begin{equation} \label{eq:bound_E_i_x_aux1} 
 \P ( q \in Q \,|\, \cal F_i(\{x \})) = \P (u \in S_1^+(x), v \in S_1^+(x) \,|\, \cal F_i ) = \frac{d^2}{N^2} \,.
\end{equation} 
In the last step we used that, by Lemma~\ref{lem:FconditionLaw}, for any $X \subset S_i(b)$, conditionally on $\cal F_i(X)$, the random variables $(A_{xy} \col x \in X, y \in B_i(b)^c)$ are independent Bernoulli-$\frac{d}{N}$ random variables. Moreover, we get 
\begin{equation} \label{eq:bound_E_i_x_aux2} 
 \sum_{q \in \cal Q} \abs{Z_q}^2 \ind{\Delta_i(\{x \})} 
\leq \op{Tr} ( (H^{(\cal T_i \cup \{ x \})} - z_i)^{-2} ) \ind{\Delta_i(\{ x\})}   \lesssim N \kappa^{-2}\,. 
\end{equation} 
Since $\cal E_i^{(0)} = \frac1d \sum_{q \in Q} Z_q$, the inclusion $\Omega \subset \Gamma_i \cap \Delta_i(\{x\})$ and the Cauchy-Schwarz inequality imply  
\begin{equation} \label{eq:bound_cal_E_i_0} 
\begin{aligned}  
\E_\Omega \big[ \abs{\cal E_i^{(0)}}^2 \big| \cal F_i(\{ x\})  \big] 
& \leq \frac{1}{d^2}   \E_\Omega \bigg[ \abs{Q} \sum_{q \in Q} \abs{Z_q}^2 \bigg| \cal F_i(\{ x\}) \bigg] \\ 
& \lesssim \frac{(\log N)^2}{ d^2}\, \E \bigg[ \sum_{q \in Q} \abs{Z_q}^2  \ind{\Delta_i(\{x \})} \bigg| \cal F_i(\{ x\}) \bigg] \\ 
& \lesssim \frac{(\log N)^2}{ d^2}  \max_{q \in \cal Q} \P ( q \in Q \,|\, \cal F_i(\{x\})) \sum_{q \in \cal Q} \abs{Z_q}^2 \ind{\Delta_i(\{x \})}  \\ 
& \lesssim \frac{ (\log N)^2}{ \kappa^2 N}\,,
\end{aligned} 
\end{equation} 
where we used that $\abs{S_1(x)} \lesssim \log N$ on $\Omega$, by  Proposition~\ref{prop:graphProperties} \ref{item:upper_bound_alpha_x}, as well 
as the $\cal F_i(\{x \})$-measurability of $\abs{Z_q}^2 \ind{\Delta_i(\{x \})}$.
The last step follows from \eqref{eq:bound_E_i_x_aux1} and \eqref{eq:bound_E_i_x_aux2}. 

To bound $\cal E_i^{(j)}$ for a fixed $j  \in \{ 1, \ldots, \abs{S_i(b)}-1\}$, we conclude 
on the event $\Gamma_i$ from the resolvent identity that 
\[ 
\cal E_i^{(j)} = \frac{1}{ d^{3/2}} \sum_{u \in S_1^+(x)} \big( H^{(\cal T_i\cup \{ y_0, \ldots, y_{j-1}\})} -z_i\big)^{-1}_{uy_j} 
\sum_{v \in S_1^+(y_j)} \big( H^{(\cal T_i \cup  \{ y_0, \ldots, y_{j}\})} - z_i\big)^{-1}_{vu}\,. 
\]
Therefore, by applying the Cauchy-Schwarz inequality twice and  using \eqref{Omega_Delta} with $X = \{y_0, \dots, y_{j}\}$, we obtain
\begin{align} 
 \abs{\cal E_i^{(j)}}^2 \ind{\Omega}  
& \leq 
\ind{\Omega} \frac{\abs{S_1^+(x)}}{d^3} \sum_{u \in S_1^+(x)} \abs{(H^{(\cal T_i \cup \{ y_0, \ldots, y_{j-1}\})} - z_i)^{-1}_{uy_j}}^2 
\absbb{\sum_{v \in S_1^+(y_j)} (H^{(\cal T_i \cup \{y_0,\ldots ,y_j\})} - z_i)^{-1}_{vu}}^2 
\nonumber \\ 
& \lesssim \ind{\Omega}  \frac{\abs{S_1^+(x)} \abs{S_1^+(y_j)}}{d^3\kappa^2}  \sum_{u \in S_1^+(x)} \abs{(H^{(\cal T_i \cup \{ y_0, \ldots, 
y_{j-1}\})} - z_i)^{-1}_{uy_j}}^2 \nonumber \\ 
& \lesssim \ind{\Omega} \frac{(\log N)^2}{d^3\kappa^2} \sum_{u \in S_1^+(x)} \abs{(H^{(\cal T_i \cup \{ y_0, \ldots, 
y_{j-1}\})} - z_i)^{-1}_{uy_j}}^2\,,   \label{eq:bound_cal_E_i_j_aux1} 
\end{align} 
where in the last step we used $\abs{S_1(x)} + \abs{S_1(y_j)} \lesssim \log N$ from 
 Proposition~\ref{prop:graphProperties} \ref{item:upper_bound_alpha_x}  on $\Omega$. 

We now set $Y_u \deq (H^{(\cal T_i \cup \{ y_0, \ldots, y_{j-1}\})} - z_i)^{-1}_{uy_j}$ for $u \in B_i(b)^c$, 
 $U \deq S_1^+(x)$ and $\cal U \deq B_i(b)^c$. 
We apply $\E[ \, \cdot \,|\, \cal F_i(\{y_0, \ldots, y_{j-1}\})]$ to \eqref{eq:bound_cal_E_i_j_aux1} 
and, similarly as in \eqref{eq:bound_cal_E_i_0}, obtain
\begin{align} 
\E_{\Omega}[\abs{\cal E_i^{(j)}}^2\,|\, \cal F_i(\{y_0, \ldots, y_{j-1}\})] 
& \lesssim \frac{(\log N)^2}{ d^3 \kappa^2} \max_{u \in \cal U} \P ( u\in U \,|\, \cal F_i(\{ y_0, \ldots, y_{j-1}\})) 
\sum_{u \in \cal U} \abs{Y_{u}}^2 \ind{\Delta_i(\{ y_0, \ldots, y_{j-1}\})} \nonumber \\ 
& \lesssim \frac{(\log N)^2}{ d^2 \kappa^4  N}, \label{eq:bound_cal_E_i_j} 
\end{align}  
where we used that $(Y_u)_{u \in \cal U}$ is $\cal F_i(\{y_0, \dots, y_{j-1}\})$-measurable, and that $\sum_{u \in \cal U} \abs{Y_{u}}^2\ind{\Delta_i(\{ y_0, \ldots, y_{j-1}\})}  =  \norm{(Y_u)_{u \in B_i(b)^c}}^2 \ind{\Delta_i(\{ y_0, \ldots, y_{j-1}\})}  \lesssim \kappa^{-2}$. We also used the remark following \eqref{eq:bound_E_i_x_aux1}.

Finally, using the estimates \eqref{eq:bound_cal_E_i_0} and \eqref{eq:bound_cal_E_i_j} 
in \eqref{eq:bound_prob_cal_E_i} together with the 
tower property of the conditional expectation complete the proof of \eqref{eq:schur_balls} 
and \eqref{eq:bound_error_term_schur}.  This concludes the proof of \ref{itm:resballs1}.

Next, we prove \ref{itm:resballs2}.
For the proof of \eqref{eq:bound_off_diagonal}, we fix $x,y \in S_1(b)$ and conclude from the resolvent identity 
 \cite[eq.~(3.5)]{BenyachKnowles2017} that 
\begin{align*}
(H^{(\cal V)}- z)^{-1}_{xy} &= - (H^{(\cal V)}-z)^{-1}_{xx} \sum_{u \notin \cal V \cup \{x \}} 
H_{xu} (H^{(\cal V \cup \{x\})} - z)^{-1}_{uy}  
\\
&= - \frac{1}{\sqrt{d}} (H^{(\cal V)}-z)^{-1}_{xx} \sum_{u \notin \cal V \cup \{x \}} \ind{u \in S_1^+(x)} 
(H^{(\cal V\cup \{x\})} - z)^{-1}_{uy}\,. 
\end{align*}
Therefore, the Cauchy-Schwarz inequality and Lemma~\ref{lem:G_bounded} imply 
\[ \abs{(H^{(\cal V)}- z)^{-1}_{xy}}^2 \ind{\Omega} \lesssim \ind{\Omega} 
 \frac{\abs{S_1^+(x)}}{\kappa^2 d}  
\sum_{u \notin \cal T_1 \cup \{x \}} \ind{u \in S_1^+(x)}\abs{(H^{(\cal T_1 \cup \{x\})} - z)^{-1}_{uy}}^2\,, 
\] 
where we also used that 
$\cal V = \{b \} \cup \cal V^{(B_1(b))}=\cal T_1$ on $\Omega$ by Proposition~\ref{prop:graphProperties} \ref{item:B_r_disjoint}. 
Hence, since $\abs{S_1(x)} \lesssim \log N$ on $\Omega$ by Proposition~\ref{prop:graphProperties} 
\ref{item:upper_bound_alpha_x}, $\Omega \subset \Delta_1(\{x\})$ and $\Delta_1(\{x\}) \in \cal F_1(\{x\})$, we obtain 
\begin{multline*}
\E\big[ \abs{(H^{(\cal V)}- z)^{-1}_{xy}}^2 \ind{\Omega} \,|\, \cal F_1(\{x \}) \big] 
\\
\lesssim \frac{\log N}{\kappa^2 d} \max_{u \notin \cal T_1 \cup \{ x\}} \P (u \in S_1^+(x) \,|\, \cal F_1(\{x \})) \sum_{u \notin \cal T_1 \cup \{ x\}} \abs{(H^{(\cal T_1 \cup \{x\})} - z)^{-1}_{uy}}^2 \ind{\Delta_1(\{x\})}
\lesssim \frac{\log N}{\kappa^4 N}\,. 
\end{multline*}
Here, in the last step, we used $\P (u \in S_1^+(x) \,|\, \cal F_1(\{x \})) \leq d/N$ and $\sum_{u \notin \cal T_1 \cup \{ x\}} \abs{(H^{(\cal T_1 \cup \{x\})} - z)^{-1}_{uy}}^2 \ind{\Delta_1(\{x\})} \lesssim \kappa^{-2}$. 
Thus, Chebyshev's inequality and the tower property of the conditional expectation complete the proof of \eqref{eq:bound_off_diagonal} and, 
therefore, the one of Lemma~\ref{lem:ResolvantBalls}.    
\end{proof}

\subsection{Proof of Proposition~\ref{prop:idealConcentration}} \label{sec:proof_idealconcentration}

This subsection is devoted to the proof of Proposition \ref{prop:idealConcentration}. We start by introducing the notion of a robust vertex.

\begin{figure}[!ht]
\begin{center}
{\footnotesize 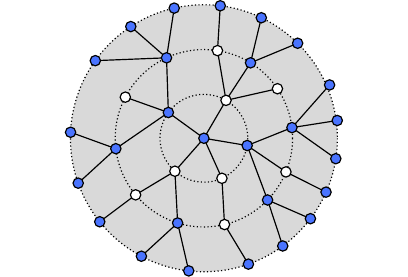}
\end{center}
\caption{An illustration of the set $\cal R \equiv \cal R_r(\bb G, b)$ of robust vertices. Here $r = 3$ and $d = 4$. We draw the ball $B_r(b)$ around the root $b$. The robust vertices are drawn in blue and the non-robust vertices in white. In this example, the root $b$ is robust. 
\label{fig:robust}}
\end{figure}

\begin{definition} \label{def:NiceVertex}
Let $b \in [N]$ and $r \in \N^*$.
We call a vertex $y \in B_{r}(b)$ \emph{robust} if 
\begin{enumerate}
\item[(a)] $y\in S_{r}(b)$ or
\item[(b)] $y\in B_{r-1}(b)$ and at least $d/2$ vertices in $S_1^{+}(y)$ are robust.
\end{enumerate}
We denote by $\cal R \equiv \cal R_r(\bb G, b) \subset B_r(b)$ the set of robust vertices.
\end{definition}

Note that $\cal R$ is an $\cal F$-measurable random set. See Figure \ref{fig:robust} for an illustration of Definition \ref{def:NiceVertex}. The following result states that with high probability the root $b$ is robust, conditioned on $S_1(a)$ and $S_1(b)$.

\begin{proposition}[The root is robust] \label{prop:goodCenter}
Suppose that $\sqrt{\log N} \ll d \lesssim \log N$ and that $r$ satisfies \eqref{eq:condition_r}.  Then $\P_\Omega(b \notin \cal R \,|\, S_1(a), S_1(b)) \lesssim N^{-1/2}$ whenever $a,b \in \cal W$.
\end{proposition}

The proof of Proposition \ref{prop:goodCenter} is given at the end of this subsection. From now on, we choose $r$ as in \eqref{eq:choice_r}.

\begin{definition} \label{def:Upsilon} 
Let $z \in \cal J$ be an $\cal F$-measurable real random variable. 
We introduce the event $\Upsilon$ on which the following conditions hold.
\begin{enumerate}[label=(\Alph*)] 
\item \label{item:UpsilonA}
$\bb G |_{B_r(b)}$ is a tree. 
\item $b \in \cal R$.
\item \label{item:Upsilon_condition_G}  $-G_{yy}(r,z) \geq (3z)^{-1}$ for all $y \in (B_r(b) \cup \cal V^{(B_r(b))})^c$.  
\item \label{item:UpsilonD}
$\abs{B_r(b) \cup \cal V^{(B_r(b))}} \leq N^{1/2}$.
\end{enumerate}
\end{definition} 

\begin{lemma} \label{lem:Upsilon_estimate}
We have $\Upsilon \in \cal F$ and $\P_\Omega(\Upsilon^c \,|\, S_1(a), S_1(b)) = O(N^{-1/2})$ whenever $a,b \in \cal W$.
\end{lemma}
\begin{proof}
That $\Upsilon \in \cal F$ follows from the Definitions \ref{def:H(i)} and \ref{def:NiceVertex}. The estimate follows from Proposition \ref{prop:goodCenter} and the facts that on $\Omega$ the conditions \ref{item:UpsilonA}, \ref{item:Upsilon_condition_G}, and \ref{item:UpsilonD} hold surely. That \ref{item:Upsilon_condition_G} holds surely on $\Omega$ follows from Lemma \ref{lem:lowboundGreen} and the observation that $\cal V = \cal V^{(B_r(b))}$ on $\Omega$ (see Proposition~\ref{prop:graphProperties} \ref{item:B_r_disjoint}). That \ref{item:UpsilonD} holds surely on $\Omega$ follows from the statements \ref{item:size_cal_V}, \ref{item:upper_bound_alpha_x}, \ref{item:upper_bound_size_B_r} of Proposition \ref{prop:graphProperties} (recall the choice \eqref{eq:choice_r}).
\end{proof}

We recall the definition of L\'evy's concentration function $Q(X,L)$ from \eqref{eq:def_concentration_function_Q}.

\begin{remark} \label{rem:Q}
The concentration function has the following obvious properties.
\begin{enumerate}[label=(\roman*)]
\item \label{itm:Q1}
For any $u > 0$ we have $Q(u X, uL) = Q(X,L)$. More generally, if $f$ is a continuous bijection on $\R$ such that $f^{-1}$ is $K$-Lipschitz, then $Q(f(X), L) \leq Q(X,KL)$.
\item \label{itm:Q2}
If $X$ and $Y$ are independent then $Q(X + Y, L) \leq \min \{Q(X,L), Q(Y,L)\}$.
\end{enumerate}
Property \ref{itm:Q2} is in general not sharp, and in some situations it can be improved considerably; see Proposition \ref{prop:ConcentrationSum}. Nevertheless, in some situations 
 \ref{itm:Q2} gives a better bound than Proposition \ref{prop:ConcentrationSum}; this is due to the minimum in \ref{itm:Q2} as opposed to the maximum in Proposition \ref{prop:ConcentrationSum}. An important theme in the proof of Proposition \ref{prop:idealConcentration} is a judicious mix of \ref{itm:Q2} and Proposition \ref{prop:ConcentrationSum}. How to do this mix is encoded by the set of robust vertices from Definition \ref{def:NiceVertex}.
\end{remark}

We denote by $Q^{\cal F}$ Lévy's concentration function with respect to the probability measure $\mathbb{P}(\, \cdot\, |\,\cal F)$. The main tool behind the proof of Proposition \ref{prop:idealConcentration} is the following anticoncentration estimate for $g_x(z)$.

\begin{proposition}\label{prop:InsideRecursion}
Let $z \in \cal J$ be a $\cal F$-measurable real random variable.
There exists a constant $\chi>0$ such that, on the event $\Upsilon$, for any $1 \leq i \leq r$ and $x\in S_{i}(b)\cap \cal R$  we have 
\begin{align*}
Q^{\cal F}\bigg(g_x(z),\frac{1}{8 z (T^2 d)^{r-i+1}}\bigg) \leq \frac{1}{2 (\chi d)^{(r-i)/2}} \,. 
\end{align*}
\end{proposition}

Before proving Proposition \ref{prop:InsideRecursion}, we use it to conclude the proof of Proposition \ref{prop:idealConcentration}.

\begin{proof}[Proof of Proposition~\ref{prop:idealConcentration}] 
We estimate
\begin{align} 
&\mspace{-20mu} \mathbb{P}_\Omega\bigg(a,b\in \cal W, \absbb{\frac{1}{d}\sum_{x\in S_1(b)} g_{x}(z) +z} 
 \leq  N^{-\eta} \bigg)
\notag  \\
& \leq  \mathbb{P}_\Omega\bigg(\Xi_r \cap \Upsilon \cap \hbb{ \absbb{\frac{1}{d}\sum_{x\in S_1(b)} g_{x}(z) +z} 
 \leq  N^{-\eta}} \bigg) + \P_\Omega(\Xi \cap \Upsilon^c)
\notag  \\ \label{P58_proof}
 & = \E \qbb{ \ind{\Xi_r} \, \ind{\Upsilon} \, \P \pbb{\absbb{\frac{1}{d}\sum_{x\in S_1(b)} g_{x}(z) +z} 
 \leq  N^{-\eta} \, \bigg| \, \cal F}}  + \P_\Omega(\Xi \cap \Upsilon^c)\,,
\end{align}
where we used Remark \ref{rem:b_in_cal_W_measurable_Xi_i} \ref{itm:W1} and \ref{itm:W3} as well as $\Upsilon \in \cal F$ by Lemma \ref{lem:Upsilon_estimate}. The second term on the right-hand side of \eqref{P58_proof} is estimated as
\begin{equation*}
\P_\Omega(\Xi \cap \Upsilon^c) = \E [\ind{a \in \cal W} \, \ind{b \in \cal W}\, \P_\Omega(\Upsilon^c \,|\, S_1(a), S_1(b))] \lesssim N^{-1/2} \, \P(a, b \in \cal W)\,,
\end{equation*}
by Lemma \ref{lem:Upsilon_estimate}.

To estimate the first term on the right-hand side of \eqref{P58_proof}, we use Proposition \ref{prop:InsideRecursion} with $i = 1$ and the estimate
\begin{equation} \label{Tr_estimate}
8z d (T^2 d)^{r} \leq  d^{2r+3/2}\leq N^{\eta}
\end{equation}
where we used the definitions of $r$ and $T$ from \eqref{eq:choice_r} and Definition \ref{def_iota}, as well as $d \gg \sqrt{\log N}$. On $\Upsilon$ we have $b \in \cal R$ and hence, by Definition \ref{def:NiceVertex}, there exists $x_* \in S_1(b) \cap \cal R$. Using that $(g_x(z))_{x \in S_1(b)}$ is an independent family, on the event $\Upsilon$ conditionally on $\cal F$, this yields, on the event $\Upsilon$,
\begin{align*}
\P \pbb{\absbb{\frac{1}{d}\sum_{x\in S_1(b)} g_{x}(z) +z} 
 \leq  N^{-\eta} \, \bigg| \, \cal F}
 &\leq
Q^\cal F \pbb{\frac{1}{d}\sum_{x\in S_1(b)} g_{x}(z) \,,\,
  \frac{1}{8 z d (T^2 d)^{r}}}
 \\
  &\leq
Q^\cal F \pbb{\frac{1}{d} \,g_{x_*}(z) \,,\,
  \frac{1}{8 z d (T^2 d)^{r}}}
  \\
  &=
  Q^\cal F \pbb{ g_{x_*}(z) \,,\,
  \frac{1}{8 z (T^2 d)^{r}}}\,,
\end{align*}
where in the second step we used Remark \ref{rem:Q} \ref{itm:Q2} and in the last step Remark \ref{rem:Q} \ref{itm:Q1}. From Proposition \ref{prop:InsideRecursion} we therefore conclude that the first term on the right-hand side of \eqref{P58_proof} is bounded by
\begin{equation*}
\frac{1}{2 (\chi d)^{(r-i)/2}} \,\P(\Xi_r) \leq N^{-\eta / 4 + o(1)} \, \P(a,b \in \cal W)\,,
\end{equation*}
where we used the definition \ref{eq:choice_r} and Remark \ref{rem:b_in_cal_W_measurable_Xi_i} \ref{itm:W2}.
The claim now follows from Lemma \ref{lem:Wab_estimate}.
\end{proof} 

\begin{remark} \label{rem:critical_eta}
If we restrict ourselves to the critical regime $d \asymp \log N$, then the factor $N^{-\eta}$ inside the probability in \eqref{anticoncent_estimate} can be improved to $N^{-\eta/2}$. To see this, we note that in this regime the parameter $T$ from Definition \ref{def_iota} satisfies $T = 10 \kappa^{-1}$ since, in the critical regime, the estimate \eqref{eq:upper_lower_bound_Schur_complement_expansion} with small enough $\kappa$ remains valid for this smaller choice of $T$. Thus, the estimate \eqref{Tr_estimate} in the proof of Proposition~\ref{prop:idealConcentration} can be replaced with $8z d (T^2 d)^{r} \leq  (C d)^r$, which is bounded by $N^{\eta/2 + o(1)}$.
\end{remark}

The key tool behind the proof of Proposition~\ref{prop:InsideRecursion} is Proposition~\ref{prop:ConcentrationSum} due to Kesten. 

\begin{proof}[Proof of Proposition \ref{prop:InsideRecursion}] 
Throughout the proof, the argument $z$ of $g_x(z)$ for any $x \in B_r(b)$ will always be the random variable 
$z$ from Definition~\ref{def:Upsilon}. Therefore, we omit this argument from our notation and write $g_x \equiv g_x(z)$. 

For $i \in [r]$ we define
\begin{equation*}
L_i \deq \frac{1}{8zT^2d} \pbb{\frac{1}{T^2 d}}^{r - i}\,, \qquad P_i \deq \frac{1}{2} \pbb{\frac{\sqrt{2} K}{d^{1/2}}}^{r - i}\,.
\end{equation*}
where $K$ is the universal constant from Proposition \ref{prop:ConcentrationSum}. We prove Proposition~\ref{prop:InsideRecursion} by showing that, for all $i \in [r]$ and $x \in S_i(b)\cap \cal R$ we have
\begin{equation}\label{eq:proof_InsideRecursion_induction} 
Q^{\cal F}\big(g_x,L_{i}\big) \leq P_{i}\,.
\end{equation} 
We show \eqref{eq:proof_InsideRecursion_induction} by induction on $i = r, r - 1, \dots, 1$.

We start the induction at $i = r$. Abbreviate $\cal X \deq (B_r(b) \cup \cal V^{(B_r(b))})^c$, which is an $\cal F$-measurable set. For $x\in S_r(b)$, conditioned on $\cal F$, $(\ind{y\in S_1^+(x)})_{y \in \cal X}$ are i.i.d.\ Bernoulli random variables. Hence, conditioned on $\cal F$ we have 
\begin{equation} \label{G_sum_eqdist}
\sum_{y\in S_1^+(x) \setminus \cal V^{(B_r(b))}} G_{yy}(r,z) \eqdist \sum_{k=0}^{|\cal X|}\ind{|S_1^+(x) \cap \cal X|=k} \sum_{i=1}^{k}G_{\sigma(i)\sigma(i)}(r,z)\,,
\end{equation}
where $\sigma$ is a uniform random enumeration of $\cal X$ (i.e.\ a bijection $[\abs{\cal X}] \to \cal X$) that is independent of $|S_1^+(x) \cap \cal X|$. Because of the condition \ref{item:Upsilon_condition_G} in Definition~\ref{def:Upsilon}, for any $k\neq l$, $|\sum_{i=1}^{l}G_{\sigma(i)\sigma(i)}(r,z)-\sum_{i=1}^{k}G_{\sigma(i)\sigma(i)}(r,z)|\geq (3z)^{-1} $. 
Therefore, for any $t\in \mathbb{R}$ we get on $\Upsilon$
\begin{equation*}
\P \pBB{\absBB{\sum_{k=0}^{|\cal X|}\ind{|S_1^+(x) \cap \cal X|=k} \sum_{i=1}^{k}G_{\sigma(i)\sigma(i)}(r,z) - t } \leq \frac{1}{8z} \, \bigg|\, \sigma, \cal F}
\leq \max_{0 \leq k\leq |\cal X|} \mathbb{P}(|S_1^+(x) \cap \cal X|=k \,|\, \cal F)\leq \frac12\,,
\end{equation*}
where in the last step we used that $|S_1^+(x) \cap \cal X| \eqdist \op{Binom}(\abs{\cal X}, d/N)$ conditioned on $\cal F$, that $\abs{\cal X} \geq N - N^{1/2}$ by the condition \ref{item:UpsilonD} in Definition \ref{def:Upsilon}, and that $d \gg 1$. From \eqref{G_sum_eqdist} we therefore conclude that
\begin{equation*}
Q^{\cal F} \pBB{\sum_{y\in S_1^+(x) \setminus \cal V^{(B_r(b))}} G_{yy}(r,z) , \frac{1}{8z}} \leq \frac12\,.
\end{equation*}
Hence, Remarks \ref{rem:properties_of_iota} and \ref{rem:Q} \ref{itm:Q1} imply
\begin{equation*}
Q^{\cal F}\pbb{g_x,\frac{1}{8z T^2 d}} \leq \frac12\,,
\end{equation*}
which is \eqref{eq:proof_InsideRecursion_induction} for $i = r$.

For the induction step, we let $i < r$, choose $x \in S_i(b) \cap \cal R$, and assume that 
\[Q^{\cal F}(g_y,L_{i+1})\leq P_{i+1}\] 
for all $y\in S_{i+1}(b)\cap \cal R$. Note that $S_1^+(x)$ and $\cal R$ are $\cal F$-measurable, and that the family $(g_y)_{y \in S_1+(x)}$ is independent on $\Upsilon$ conditioned on $\cal F$, by Remark \ref{rem:g_independence} and Definition \ref{def:Upsilon} \ref{item:UpsilonA}. Hence, we can apply Proposition \ref{prop:ConcentrationSum} to the concentration function $Q^{\cal F}$ to obtain
\begin{equation} \label{QF_Li1}
Q^{\cal F}\pBB{\sum_{y\in S_1^+(x)\cap \cal R} g_y,L_{i+1}}\leq \frac{K}{\sqrt{\abs{ S_1^+(x)\cap \cal R}}} P_{i+1} \leq \frac{K\sqrt{2} P_{i+1}}{d^{1/2}}\,,
\end{equation}
where the last inequality follows from $\abs{ S_1^+(x)\cap \cal R} \geq d/2$, by Definition \ref{def:NiceVertex}.
Moreover, the conditional independence of the sums $\sum_{y \in S_1^+(x) \cap \cal R} g_y$ and $\sum_{y\in S_1^+(x) \setminus \cal R} g_y$ combined with Remark \ref{rem:Q} \ref{itm:Q1} and \ref{itm:Q2} yields 
\begin{align*}
Q^{\cal F}\pBB{\frac{1}{d}\sum_{y\in S_1^+(x)} g_y,\frac{L_{i+1}}{d}}
&= Q^{\cal F}\pBB{\sum_{y\in S_1^+(x)} g_y, L_{i+1}}
\\
&= Q^{\cal F}\pBB{\sum_{y\in S_1^+(x) \cap\cal R} g_y+\sum_{y\in S_1^+(x) \setminus \cal R} g_y,L_{i+1}}
\\
&\leq Q^{\cal F}\pBB{\sum_{y\in S_1^+(x) \cap \cal R} g_y,L_{i+1}}\,. 
\end{align*}
Hence, by Remark \ref{rem:properties_of_iota}, Remark \ref{rem:Q} \ref{itm:Q1}, Definition~\ref{def:gx}, and \eqref{QF_Li1}, we obtain
\begin{align*}
Q^{\cal F}\pBB{g_x,\frac{L_{i+1}}{T^2 d}} \leq \frac{K\sqrt{2} P_{i+1}}{d^{1/2}}\,,
\end{align*}
 which is \eqref{eq:proof_InsideRecursion_induction}.
This completes the proof of \eqref{eq:proof_InsideRecursion_induction} and, hence, the one of Proposition~\ref{prop:InsideRecursion}. 
\end{proof}

\begin{proof}[Proof of Proposition \ref{prop:goodCenter}]
The proof proceeds in two steps: first by establishing the claim for the root of a Galton-Watson branching process with Poisson offspring distribution with mean $d$, and then concluding by a comparison argument.

Denote by $\cal P_s$ a Poisson random variable with expectation $s$. Let $W$ denote the Galton-Watson branching process with Poisson offspring distribution $\cal P_d$, which we regard as a random rooted ordered tree\footnote{A rooted ordered tree (also called plane tree) is a rooted tree in which an ordering is specified among the children of each vertex.} whose root we call $o$. We use the graph-theoretic notations (such as $S_i(x)$) from Section \ref{sec:notations} also on rooted ordered trees. Moreover, we extend Definition \ref{def:NiceVertex} to a rooted ordered tree $T$ in the obvious fashion, and when needed we use the notation $\cal R \equiv \cal R_r(T,o)$ to indicate the radius $r$, the tree $T$, and the root $o$ explicitly.

We define the parameter $\delta \deq \mathbb{P}(\cal P_{3d/4}\leq \frac{d}{2})$.
By Bennett's inequality (see Lemma \ref{lem:Bennett} below), we find that $\delta \leq \ee^{-cd}$ for some universal constant $c > 0$. We shall show by induction on $i$ that
\begin{equation} \label{GW_R_induction}
\P(o \notin \cal R_i(W,o)) \leq \delta
\end{equation}
for all $i \geq 0$.
For $i = 0$ we have $\P(o \notin \cal R_i(W,o))  = 0$ 
since $o \in \cal R_0(W,o)$ by (the analogue of) Definition \ref{def:NiceVertex}, and \eqref{GW_R_induction} is trivial.

To advance the induction, we suppose that \eqref{GW_R_induction} holds for some $i \geq 0$. By Definition \ref{def:NiceVertex},
\begin{equation*}
\P(o \notin \cal R_{i+1}(W,o)) = \P \pBB{\sum_{x \in S_1(o)} \ind{x \in  \cal R_{i+1}(W,o)} < \frac{d}{2}}\,.
\end{equation*}
By definition of the branching process $W$, conditioned on $S_1(o)$, the random variables $(\ind{x \in  \cal R_{i+1}(W,o)})_{x \in S_1(o)}$ are independent Bernoulli random variables with expectation
\begin{equation*}
\P(x \in  \cal R_{i+1}(W,o) | S_1(o)) = \P(o \in \cal R_{i}(W,o)) \eqd 1 - \zeta_i\,,
\end{equation*}
where $x \in S_1(o)$.
We conclude that $\sum_{x \in S_1(o)} \ind{x \in  \cal R_{i+1}(W,o)} \eqdist \cal P_{d (1 - \zeta_i)}$. Using the induction assumption $\zeta_i \leq \delta$ from \eqref{GW_R_induction} and the bound $\delta \leq \ee^{-c d} < 1/4$ for large enough $d$, we therefore conclude that
\begin{equation*}
\P(o \notin \cal R_{i+1}(W)) = \P \pbb{\cal P_{d (1 - \P(o \notin \cal R_{i}(W,o)))} < \frac{d}{2}} \leq \P \pbb{\cal P_{3d/4} \leq \frac{d}{2}} = \delta\,.
\end{equation*}
This concludes the proof of \eqref{GW_R_induction} for all $i \geq 0$.

Hence, denoting by $\cal B_{n,p}$ a random variable with law $\op{Binom}(n,p)$, we conclude that if $\abs{S_1(o)} \geq d$ then
\begin{multline} \label{robust_GW_estimate}
\P(o \notin \cal R_r(W,o) \,|\, S_1(o)) = \P \pBB{\sum_{x \in S_1(o)} \ind{x \in  \cal R_{r - 1}(W,o)} < \frac{d}{2} \, \bigg|\, S_1(o)} 
= \P\pbb{\cal B_{\abs{S_1(o)}, 1 - \zeta_{r-1}} < \frac{d}{2} \, \bigg|\, S_1(o)}
\\
\leq \P\pbb{\cal B_{d, 1 - \delta} < \frac{d}{2}}
=
\P\pbb{\cal B_{d, \delta} \geq \frac{d}{2}} \leq \ee^{-c d \delta \frac{1}{2 \delta} \log \frac{1}{2 \delta}}\leq \ee^{-c d^2} \leq N^{-1}
\end{multline}
for some universal constant $c > 0$,
where in the third step we used that $\abs{S_1(o)} \geq d$ and $\zeta_{r - 1} \leq \delta$ by \eqref{GW_R_induction}, in the fifth step Bennett's inequality (see Lemma \ref{lem:Bennett} below), in the sixth step that $\delta \leq \ee^{-cd}$, and in the last step the assumption $d \gg \sqrt{\log N}$. This concludes the estimate for the Galton-Watson process $W$.

Next, we analyse $\P_\Omega(b \notin \cal R_r(\bb G, b) \,|\, S_1(a), S_1(b))$. We note first that we can assume that $\abs{B_1(a)} \leq 1 + 10 \log N$ and that $B_1(a)$ and $B_1(b)$ are disjoint, for otherwise the above probability vanishes by definition of $\Omega$ and  Proposition~\ref{prop:graphProperties}.

We observe that a rooted ordered tree can be regarded as an equivalence class of (labelled) rooted trees up to a relabelling of the vertices that preserves the ordering of the children of each vertex. We denote by $[\bb T,x]$ the equivalence class of the labelled rooted tree $(\bb T, x)$, where $x$ is the root. By convention, if $\bb T$ is not a tree (i.e.\ if it contains a cycle) then its equivalence class is the empty tree. We denote by $\fra T_r$ the set of rooted ordered trees of depth $r$.  Moreover, we denote by $\fra T_r^* \subset \fra T_r$ the subset of rooted ordered trees with at most $N^{1/5}$ vertices and whose root is a robust vertex with $\abs{S_1(b)}$ children. Abbreviating $\Delta \deq \{B_r(b) \subset [N] \setminus B_1(a)\}$, we can write
\begin{align} 
\P_\Omega(b \notin \cal R_r(\bb G, b) \, |\, S_1(a), S_1(b))
&= \P_\Omega([\bb G|_{B_r(b)},b] \notin \fra T_r^* \,|\, S_1(a), S_1(b))
\notag\\ \label{P_robust_Omega}
&\leq
\P(\{[\bb G|_{B_r(b)},b] \notin \fra T_r^* \}\cap \Delta\,|\, S_1(a), S_1(b))\,,
\end{align}
where we used that, since $a,b \in \cal W \subset \cal V$, Proposition~\ref{prop:graphProperties} implies that on the event $\Omega$ the graph $\bb G|_{B_r(b)}$ is a tree with at most $N^{1/5}$ vertices, and that $\Omega \subset \Delta$.

For the following, let $T \in \fra T_r^*$ and denote by $o$ its root. For $1 \leq i \leq r$ we introduce the event
\begin{equation*}
\Theta_i \deq \{[\bb G|_{B_{i}(b)},b] =  T |_{B_i(o)}\} \cap \{B_i(b) \subset [N] \setminus B_1(a)\}\,.
\end{equation*}
In particular, $\P(\Theta_1 \,|\, S_1(a), S_1(b)) = 1$ because of the assumed disjointness of $S_1(a)$ and $S_1(b)$.
We now estimate
\begin{equation} \label{Theta_r_proba}
\P(\{[\bb G|_{B_r(b)},b] = T\}\cap \Delta \,|\, S_1(a), S_1(b))
= \P(\Theta_r \, | \, S_1(a), S_1(b))
\end{equation}
recursively, for $1 \leq i \leq r - 1$, using the expression
\begin{multline*}
\frac{\P(\Theta_{i+1}| \, S_1(a), S_1(b))}{\P(\Theta_{i}| \, S_1(a), S_1(b))}
\\
=  \frac{(N' - \abs{B_i(o)})!}{(N' - \abs{B_{i+1}(o)})! \prod_{x \in S_i(o)} \abs{S_1^+(x)}!} \prod_{x \in S_i(o)} \pbb{\frac{d}{N}}^{\abs{S_1^+(x)}} \pbb{1 - \frac{d}{N}}^{N - \abs{B_i(o)} - \abs{S_1^+(x)}}\,,
\end{multline*}
where $N' \deq N - \abs{B_1(a)}$,  and the graph-theoretic quantities on the left-hand side are in terms of $\bb G$ and on the right-hand side in terms of the deterministic rooted ordered tree $T$.
Here, the multinomial factor in front arises from a choice of the $\abs{S_i(o)}$ disjoint subsets representing the children of the vertices in $S_i(o)$ from $N' - \abs{B_i(o)}$ available vertices, and the remaining product follows by independence of the edges in $\bb G$. We deduce that
\begin{align*}
\frac{\P(\Theta_{i+1}| \, S_1(a), S_1(b))}{\P(\Theta_{i}| \, S_1(a), S_1(b))} &= (1 + O(N^{-3/4})) \frac{(N' - \abs{B_i(o)})!}{(N' - \abs{B_{i+1}(o)})! \, N^{\abs{S_{i+1}(o)}}}  \prod_{x \in S_i(o)} \frac{d^{\abs{S_1^+(x)}}}{\abs{S_1^+(x)}!} \ee^{-d}
\\
&= (1 + O(N^{-3/5}))  \prod_{x \in S_i(o)} \frac{d^{\abs{S_1^+(x)}}}{\abs{S_1^+(x)}!} \ee^{-d}\,,
\end{align*}
where we used that $N - N' \leq 1 + 10 \log N$ and $\abs{B_r(o)} \leq N^{1/5}$. By induction on $i$ and comparison with the Galton-Watson tree $W$, using that
\[ \frac{\P ( W|_{B_{i+1}(o)} = T|_{B_{i+1}(o)} \,|\, S_1(o))} {\P ( W|_{B_{i}(o)} = T|_{B_{i}(o)} \,|\, S_1(o))} 
= \prod_{x \in S_i(o)} \frac{d^{\abs{S_1^+(x)}}}{\abs{S_1^+(x)}!} \ee^{-d}\,, \] 
as well as $\P(\Theta_1 | S_1(a), S_1(b)) = 1$, we therefore conclude from \eqref{Theta_r_proba} that if $\abs{S_1(b)} = \abs{S_1(o)}$ then
\begin{equation*}
\P(\{[\bb G|_{B_r(b)},b] = T\} \cap \Delta  \,|\, S_1(a), S_1(b)) = (1 + O(N^{-1/2})) \, \P(W |_{B_r(o)} = T \,|\, S_1(o)) 
\end{equation*}
for all $T \in \fra T_r^*$. Thus,
\begin{align*}
\P(\{[\bb G|_{B_r(b)},b] \notin \fra T_r^*\} \cap \Delta \,|\, S_1(a), S_1(b)) &\leq 1 - \sum_{T \in \fra T_r^*} \P(\{[\bb G|_{B_r(b)},b] = T\} \cap \Delta \,|\, S_1(a), S_1(b))
\\
&= 1 - \sum_{T \in \fra T_r^*} \P(W |_{B_r(o)} = T \,|\, S_1(o)) + O(N^{-1/2}) 
\\
&= \P(o \notin \cal R_r(W, o) \,|\, S_1(o)) + O(N^{-1/2})\,.
\end{align*}
The claim now follows from \eqref{robust_GW_estimate} and \eqref{P_robust_Omega}, noting that if $b \in \cal W$ then $\abs{S_1(b)} = \abs{S_1(o)} \geq d$.
\end{proof}

\appendix

\section{Quantitative behaviour of $\alpha^*$} 
\label{app:alpha_star}

We recall the definition of $\alpha^*(\mu)$ from \eqref{eq:largeDegree} and remark that, besides $\mu$, 
it depends on $N$ and $d$. 
Our analysis of $\alpha^*(\mu)$ is based on quantitatively approximating the distribution of $\alpha_1$ by
a Poisson distribution. 
Owing to the Poisson approximation of binomial random variables (see e.g.\ \cite[Lemma~A.6]{ADK20}) 
 and the Stirling approximation of factorials, 
we have 
\begin{equation} \label{eq:binomial_Poisson_Stirling_approximation}  
\P(d \alpha_1 = k ) = \frac{d^k}{k!} \ee^{-d} \bigg( 1 + O \bigg( \frac{k^2}{N^2} + \frac{d^2}{N} \bigg) \bigg) = \exp(-f_d(k/d) + O(k^{-1})) 
\end{equation}
for $k \leq \sqrt{N}$, 
where $f_d(\alpha) \deq d ( \alpha \log \alpha - \alpha + 1) + \frac{1}{2} \log (2 \pi \alpha d )$  for $\alpha >0$.

\begin{lemma} \label{lem:expansion_alpha_star} 
If $\mu \in [0,1-\eps]$ for some constant $\eps \in (0,1)$ then there is a constant $T \geq 2$ such that if $d \geq 1$ and $t \defeq \frac{\log N}{d} \geq T$ then 
\begin{equation} \label{eq:alpha_star_expansion} 
 \alpha^*(\mu) = \frac{(1-\mu)t}{\log t} \bigg( 1+ O\bigg( \frac{\log \log t}{\log t} \bigg) \bigg)\,. 
\end{equation} 
\end{lemma} 

Lemma \ref{lem:expansion_alpha_star} is proved at the end of this appendix. We use it to derive two simple consequences.

\begin{lemma} \label{lem:lower_bound_alpha_x} 
Let $1 \leq d \leq 3 \log N$ and $\mu \in [0,1- \eps]$ for some constant $\eps \in (0,1)$.  
If $x \in \cal V$ then 
\[ \alpha_x \gtrsim \frac{\log N}{d \log \big( \frac{10\log N}{d} \big)}\,. \] 
\end{lemma} 

\begin{proof}
Let $T$ be as in Lemma~\ref{lem:expansion_alpha_star}. If $\frac{\log N}{d}\geq T$ then Lemma~\ref{lem:lower_bound_alpha_x} follows directly from Lemma~\ref{lem:expansion_alpha_star} and the definition of $\cal V$ 
in \eqref{eq:def_cal_V}. 
If $\frac{\log N}{d} < T \lesssim 1$ then $\alpha_x \geq \alpha^*(\mu) \geq 2 + \kappa$ directly implies 
the lemma as $\frac{\log N}{d \log \big( \frac{10\log N}{d} \big)} \asymp 1$.  
\end{proof}

\begin{corollary} \label{cor:Lambda_alpha_1_2_Lambda_alpha_mu} 
Let $\mu$, $\nu \in [0,1-\eps]$ for some constant $\eps >0$. 
Then there is a constant $T  \equiv T(\eps) \geq 2 $ such that if $t\deq \frac{\log N}{d} \geq T$ then 
\[ 
\frac{\Lambda(\alpha^*(\nu))}{\Lambda(\alpha^*(\mu))} 
= \sqrt{\frac{1-\nu}{1-\mu}} \bigg(1 + O\bigg( \frac{\log \log t}{\log t} \bigg) \bigg)\,. 
\] 
\end{corollary} 

\begin{proof} 
Choose $T$ as in Lemma~\ref{lem:expansion_alpha_star}, which then yields that 
$\frac{\alpha^*(\nu)}{\alpha^*(\mu)} = \frac{1-\nu}{1-\mu} \big( 1 + O\big( \frac{\log \log t}{\log t}\big) \big)$. 
Hence, using $\min\{\alpha^*(\nu),\alpha^*(\mu)\} \gtrsim \frac{t}{\log t}$, we get 
\[ 
\frac{\Lambda(\alpha^*(\nu))}{\Lambda(\alpha^*(\mu))} = \sqrt{\frac{\alpha^*(\nu)}{\alpha^*(\mu)}} \bigg( 1 + O\bigg ( \frac{1}{\alpha^*(\nu)} + \frac{1}{\alpha^*(\mu)}\bigg)\bigg) = \sqrt{\frac{1-\nu}{1-\mu}} \bigg(1 + O\bigg( \frac{\log \log t}{\log t} \bigg) \bigg)\,. \qedhere
\] 
\end{proof}

\begin{proof}[Proof of Lemma \ref{lem:expansion_alpha_star}]
Fix $\mu \in [0, 1 - \eps ]$ for some constant $\eps \in (0,1)$.
Throughout the proof we abbreviate $\alpha^* \equiv \alpha^*(\mu)$. 
A simple application of Bennett's inequality (see Lemma \ref{lem:Bennett} below) shows that $\alpha^* \lesssim (\log N)/d$. 

If $Z$ is a random variable with law $\mathrm{Binom}(n,d/n)$ then there is a constant $C>0$ such that
\[ 
\P(Z \geq  k )  = \P( Z = k ) \bigg( 1 + O \bigg( \frac{d}{k} \bigg) \bigg) 
\] 
for any $k \in \N$ with  $2d \leq k \leq \sqrt{n}/C$ (see e.g.\ \cite[eq.~(3.4)]{BBK2}). 
Therefore, for any $\alpha$ satisfying $2 d \leq \alpha d \lesssim  \log N$,  we obtain 
\begin{multline}\label{eq:probability_alpha_1_geq_alpha_Poisson_approx} 
 \P( \alpha_1 \geq \alpha) = \P ( \alpha_1 = \alpha) ( 1 + O(\alpha^{-1})) \\ 
 = \exp ( - f_d(\alpha)
 + O( (\alpha d)^{-1})) ( 1+O(\alpha^{-1})) = \exp ( -f_d(\alpha) + O( \alpha^{-1}))\,, 
\end{multline} 
where we used \eqref{eq:binomial_Poisson_Stirling_approximation} in the second step and 
incorporated the error term into the argument of the exponential in the last step. 
As $\alpha^* \in \N/d$, the definition of $\alpha^*$ in \eqref{eq:largeDegree} implies  
\begin{equation} \label{eq:alpha_star_N_mu} 
 \P ( \alpha_1 \geq \alpha^*) \leq N^{\mu - 1} \leq \P ( \alpha_1 \geq \alpha^* - d^{-1})\,.   
\end{equation} 
Since $f_d(\alpha^* - d^{-1}) = f_d(\alpha^*) + O( (\log \alpha^*)^{-1})$ we conclude from 
\eqref{eq:probability_alpha_1_geq_alpha_Poisson_approx} and \eqref{eq:alpha_star_N_mu} 
that
$f_d(\alpha^*) = (1 - \mu) \log N + O \big ( (\log \alpha^*)^{-1} \big)$ 
which yields 
\begin{equation} \label{eq:alpha_star_first_expansion} 
\alpha^* = \frac{(1-\mu)t }{- 1 + \log \alpha^*} + O \big( (\log \alpha^*)^{-1} \big) 
\end{equation} 
with $t \defeq \frac{\log N}{d}$.  
We replace $\alpha^*$ in the first term on the right-hand side of \eqref{eq:alpha_star_first_expansion} 
by the entire right-hand side of \eqref{eq:alpha_star_first_expansion} 
and choose $T \geq 2$ a sufficiently large constant to arrive at  
\begin{equation} \label{eq:alpha_star_second_expansion} 
\alpha^* = \frac{(1-\mu)t}{\log t + \log( 1- \mu) - \log ( \log \alpha^* - 1) + O( ( \alpha^* \log \alpha^*)^{-1})  -1 } 
 + O( (\log \alpha^*)^{-1})\,.   
\end{equation} 
By possibly increasing the constant $T$, we see that $\alpha^* \gtrsim \frac{t}{\log t}$. 
Using this lower bound to replace $\alpha^*$ in the error terms of \eqref{eq:alpha_star_second_expansion} 
by $t$ yields \eqref{eq:alpha_star_expansion} and, thus, 
completes the proof of Lemma~\ref{lem:expansion_alpha_star}. 
\end{proof}

\section{Qualitative behaviour of degree sequence and of $\alpha^*$} \label{sec:alpha_qual}
In this appendix we describe the qualitative behaviour of the normalized degree sequence $(\alpha_x)_{x \in [N]}$ and apply it to $\alpha^*(\mu)$ from \eqref{eq:largeDegree}. For definiteness, we focus on the critical regime, where $d = b \log N$ for some constant $b$, and consider the limit $N \to \infty$ with $\kappa = o(1)$. For detailed proofs, we refer to \cite[Appendix A.4]{ADK20}.

\begin{figure}[!ht]
\begin{center}
{\footnotesize 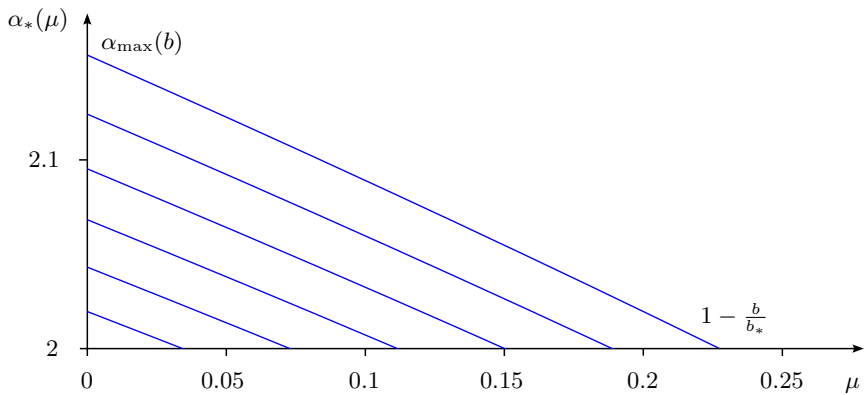}
\end{center}
\caption{An illustration of the function $\mu \mapsto \alpha_*(\mu)$, the asymptotic version of $\alpha^*(\mu)$. We plot it for the values $b = 2.0, 2.1, \dots, 2.5$, corresponding to the graphs from top to bottom. Each graph crosses the horizontal axis at $\mu = 1 - \frac{b}{b_*}$ and the vertical axis at $\alpha = \alpha_{\max}(b)$.
\label{fig:alpha}}
\end{figure}

For $b \geq 0$ and $\alpha \geq 2$ define
\begin{equation} \label{def_theta}
\theta_b(\alpha) \deq 1 - b (\alpha \log \alpha - \alpha + 1)\,.
\end{equation}
For any $b \leq b_*$, it is easy to see that for any $\mu \in \qb{0,1 - \frac{b}{b_*}}$ the equation $\mu = \theta_b(\alpha)$ has a unique solution $\alpha \geq 2$, which we denote by $\alpha_*(\mu)$. By Poisson approximation from \eqref{eq:binomial_Poisson_Stirling_approximation}, we deduce that, with high probability, for any $\alpha \geq 0$ we have
\begin{equation} \label{degree_counting}
\abs{\h{x \in [N] \col \alpha_x \geq \alpha}} = N^{\theta_b(\alpha) + o(1)} + o(1)\,.
\end{equation}
Recalling the definition \eqref{eq:largeDegree}, we hence conclude that $\alpha^*(\mu) = \alpha_*(\mu) + o(1)$. In other words, $\alpha_*(\mu)$ is the asymptotic value of $\alpha^*(\mu)$. We refer to Figure \ref{fig:alpha} for an illustration of the function $\alpha_*(\mu)$.

Next, we define
\begin{equation} \label{def_lambda_max}
\alpha_{\max}(b) \deq \alpha_*(0) \,, \qquad \lambda_{\max}(b) \deq \Lambda(\alpha_{\max}(b))\,,
\end{equation}
which, by \eqref{degree_counting} and Remark \ref{rem:ev_Lambda}, 
have the interpretation of the asymptotic largest normalized degree and largest nontrivial eigenvalue of $H$, respectively. Moreover, for $b < b_*$ we define
\begin{equation} \label{lambda_max_rho}
\rho_b(\lambda) \deq
\begin{cases}
\theta_b(\Lambda^{-1}(\lambda))_+ & \text{if } \abs{\lambda} \geq 2
\\
1 & \text{if } \abs{\lambda} < 2\,,
\end{cases}
\end{equation}
where $\Lambda^{-1}(\lambda) = \frac{\lambda^2}{2}(1 + \sqrt{1 - 4/\lambda^2})$ for $\abs{\lambda} \geq 2$. Then with high probability the density of states around energy $\lambda \in \R$ equals $N^{\rho_b(\lambda) + o(1)}$. For $\lambda > 2$, this follows from Remark \ref{rem:ev_Lambda} and \eqref{degree_counting}. The function $\rho_b$ is illustrated in Figure \ref{fig:rho} below.

\begin{figure}[!ht]
\begin{center}
{\small 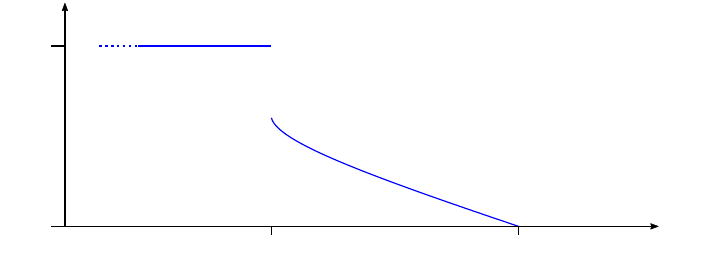}
\end{center}
\caption{The behaviour of the exponent $\rho_b$ of the density of states from \eqref{lambda_max_rho}, as a function of the energy $\lambda$. Here $d = b \log N$ with $b = 1$ and $\lambda_{\max}(b) \approx 2.0737$. We only plot a neighbourhood of the threshold energy $2$. The jump at $2$ of $\rho_b$ is from $\rho_b(2^-) = 1$ to $\rho_b(2^+) = 1 - b / b_* = 2 - 2 \log 2$.
\label{fig:rho}}
\end{figure}

\section{Largest eigenvalues of $H$ and some submatrices -- proof of Proposition~\ref{pro:spectral_gap}} \label{sec:proof_pro_spectral_gap}

In this section we prove Proposition~\ref{pro:spectral_gap}. We start with the following auxiliary result.

\begin{lemma} \label{lem:bound_A_minus_EA} 
Let $d$ satisfy \eqref{eq:d_range}.
With high probability, for all $\nu \in [0,1]$, we have  
\begin{align} \label{eq:proof_spectral_gap_aux1} 
  \norm{(H - \E H)^{(\{x \,:\, \alpha_x \geq \alpha^*(\nu)\})}} 
 \leq \Lambda(\alpha^*(\nu)) + o(1)\,.  
\end{align}
\end{lemma} 

Before giving the proof of Lemma~\ref{lem:bound_A_minus_EA}, we use it to establish Proposition~\ref{pro:spectral_gap}. 
Lemma~\ref{lem:bound_A_minus_EA} is the core ingredient in the proof of Proposition~\ref{pro:spectral_gap} \ref{item:L2}.
The main ideas for the proofs of Proposition~\ref{pro:spectral_gap} \ref{item:L2} and Lemma~\ref{lem:bound_A_minus_EA} 
stem from \cite{ADK19}, where a version of Proposition~\ref{pro:spectral_gap} \ref{item:L2} was shown in 
\cite[eq.~(8.11)]{ADK19}. However, \cite[eq.~(8.11)]{ADK19} does not cover the full range of $d$ required in Proposition~\ref{pro:spectral_gap} 
\ref{item:L2} and does not have precise enough error bounds. Therefore, we apply the results from \cite{ADK20}, 
which refined the ideas from \cite{ADK19}.   

\begin{proof}[Proof of Proposition~\ref{pro:spectral_gap}] 
For the proof of \ref{item:L2}, 
we choose $\nu = \mu$ in \eqref{eq:proof_spectral_gap_aux1}, observe that $\{ x \col \alpha_x \geq \alpha^*(\mu)\} = \cal V$ and use eigenvalue interlacing (Lemma \ref{lem:interlacing2}) to obtain \ref{item:L2}  since $(\E H)^{(\cal V)}$ is a positive semidefinite rank-one matrix.

For the proof of \ref{item:L1}, we regard $H^{(X)}$ as perturbation of $(\E H)^{(X)}$. 
We note that $\norm{H^{(X)} - (\E H)^{(X)}} \leq \norm{H - \E H}\ll \sqrt{d}$ by \eqref{eq:H_minus_EH_ll_sqrt_d}, that $ (\E H)^{(X)}$ is a rank-one matrix with nonzero eigenvalue $\sqrt{d}\big(1 - \frac{\abs{X}}{N}\big)$ and associated eigenvector $\frac{\f 1_{X^c}}{\abs{X^c}^{1/2}}$, and that $\abs{X} \leq \sum_{x \in \cal V} \abs{B_r(x)} \leq N^{1/3 + o(1)}$ by Proposition~\ref{prop:graphProperties} 
\ref{item:size_cal_V} and Proposition~\ref{prop:graphProperties_every_vertex} \ref{item:upper_bound_size_B_r}, recalling the condition \eqref{eq:condition_r}. We hence conclude \ref{item:L1} using Lemma~\ref{lem:perturbation_theory}.

For the proof of \ref{item:Candidate}, we first note that 
$H^{(X)} \f v_r(x) = H \f v_r(x)$ since $X \cap B_r(x) = \emptyset$ and $\supp \f v \subset B_{r - 1}(x)$. To estimate $(H - \E H - \Lambda(\alpha_x))\f v_r(x)$, we now apply \cite[Proposition~5.1]{ADK19}. 
From the assumption $\log d \ll r \leq \frac{\log N}{6 \log d}$, $\abs{S_1(x)} \leq 10 \log N$ by Proposition~\ref{prop:graphProperties_every_vertex} \ref{item:upper_bound_alpha_x} and the lower bound in \eqref{eq:d_range}, 
we conclude that the condition on $r$ in \cite[Proposition~5.1]{ADK19} is satisfied\footnote{See \cite[eq.~(5.1)]{ADK19} for the definition of $r_x$.}. 
The condition $\log d \ll r \leq \frac{\log N}{6 \log d}$, $2 + \kappa \leq  \alpha^* \leq \alpha_x$ for $x \in \cal V$ and Proposition~\ref{prop:graphProperties_every_vertex} 
\ref{item:upper_bound_alpha_x} imply that the event from \cite[eq.~(5.4)]{ADK19} occurs with high probability. 
Therefore, using the lower bound on $d$ from \eqref{eq:d_range} and Lemma~\ref{lem:lower_bound_alpha_x}, 
we obtain from \cite[Proposition~5.1]{ADK19} that 
\[ 
\norm{(H - \E H- \Lambda(\alpha_x))\f v_r(x)} \lesssim (\log \log N)^{-1}\bigg( \log \bigg( \frac{10 \log N}{d} \bigg) \bigg)^{1/2} 
  = o(1) 
\] 
with high probability. 
Using that $\supp \f v_r(x) \subset B_{r - 1}(x)$ and $\E H_{xy} = \frac{\sqrt{d}}{N} \ind{x \neq y}$, we deduce that $\norm{(\E H)\f v_r (x) } \leq \sqrt{\frac{d}{N}\abs{B_r(x)}} = o(1)$, where in the last step we used the upper bound on $r$ and Proposition~\ref{prop:graphProperties_every_vertex} \ref{item:upper_bound_size_B_r}. We hence conclude that  $\norm{(H^{(X)} -\Lambda(\alpha_x))\f v_r(x)} = o(1)$
by recalling $H^{(X)} \f v_r(x) = H \f v_r(x)$.

We now turn to the proof of \ref{item:spectral_gap_q}. 
For the given $r \in \N$, we choose $\f q$ as in \cite[Proposition~6.1]{ADK21}. 
Then \cite[Proposition~6.1]{ADK21} yields $\norm{\f q - \f e} \lesssim d^{-1/2}$ and $\norm{(H- \sqrt{d})\f q}\lesssim d^{-1/2}$. 
In order to estimate $\norm{\f w_1 - \f q}$, we use Lemma~\ref{lem:perturbation_theory} with $M = H$, $\wh{\lambda} = \sqrt{d}$ and $\f v = \f q$ 
whose assumptions we check next. 
From \eqref{eq:H_minus_EH_ll_sqrt_d} and eigenvalue interlacing (Lemma \ref{lem:interlacing2}), we conclude that $\lambda_2(H) \ll \sqrt{d}$.  
Therefore, there is $\Delta \asymp \sqrt{d}$ such that $H$ has a unique eigenvalue in $[\sqrt{d} - \Delta, \sqrt{d} + \Delta]$. 
Since $\norm{(H-\sqrt{d})\f q} \lesssim d^{-1/2}$ as explained above, the conditions of Lemma~\ref{lem:perturbation_theory} with $\Delta \asymp \sqrt{d}$ 
and $\eps  \asymp d^{-1/2}$ are satisfied 
and we obtain $\norm{\f w_1 - \f q} \lesssim d^{-1}$.

Finally, from \cite[Proposition~6.1]{ADK21}, we know the inclusion $\supp \f q \subset \big( \bigcup_{x \in \cal U} B_{r+1}(x) \big)^c$, where  
\begin{equation} \label{eq:def_cal_U} 
\cal U \deq \begin{cases}  \big\{ x \in [N] \col \alpha_x \geq 2 + \xi^{1/4} \big\} & \text{ if } d > (\log N)^{3/4} \\  
 \big\{x \in [N] \col \alpha_x \geq \mathfrak a/5 \big\} & \text{ if }d \leq (\log N)^{3/4} 
\end{cases} 
\end{equation}
with $\xi$ from \eqref{eq:def_xi_xi_u} and $\mathfrak a$ the solution of $h(\mathfrak a - 1) = \frac{\log N}{d}$ (compare with \eqref{eq:def_h} below). 

We now show that $\cal V \subset \cal U$. If $d > (\log N)^{3/4}$, the inclusion $\cal V \subset \cal U$ is 
obvious as $\xi = o(1)$ by \eqref{eq:xi_xi_tau_minus_1_are_o_1} and $\alpha_x \geq 2 + \kappa$ by \eqref{eq:def_cal_V} and \eqref{eq:largeDegree}. 
For $d \leq (\log N)^{3/4}$, we note that a short analysis of the definition of $\mathfrak a$ reveals that 
 $\mathfrak a = \frac{t}{\log t}( 1 + o(1))$ with $t \deq \frac{\log N}{d}$ if $t \geq C$ for some sufficiently large constant $C$. 
By possibly increasing $C$, we conclude from Lemma~\ref{lem:expansion_alpha_star} that $\alpha^*(\mu) \geq \mathfrak a /5$ 
if $\mu \in [0,1/5)$ and $t \geq C$. 
This implies $\cal V \setminus \cal W$ if $d \leq (\log N)^{3/4}$ 
and completes the proof of Proposition~\ref{pro:spectral_gap} \ref{item:spectral_gap_q}.  
\end{proof} 

\begin{remark} 
We note that in Proposition~\ref{pro:spectral_gap} \ref{item:spectral_gap_q},  
actually $\supp \f q \subset \big( \bigcup_{x \in \cal U} B_{r+1}(x) \big)^c$ holds 
as shown in the above proof. 
\end{remark}

\begin{proof}[Proof of Lemma~\ref{lem:bound_A_minus_EA}]  

We first recall $\xi$ and $\xi_u$ from \eqref{eq:def_xi_xi_u}. 
Throughout the proof, we choose $\tau \deq 1 + \xi^{1/2}$. 
Note that $\xi = o(1)$, $\tau = 1 + o(1)$ and $\xi_{\tau - 1} = o(1)$ by \eqref{eq:xi_xi_tau_minus_1_are_o_1}.

The proof relies on an analogous argument in \cite[Proof of Proposition~3.12]{ADK20}, 
to which we refer in the following.  We begin by introducing some notation. Following \cite[Section 3.1]{ADK20}, we define $\cal V_\alpha \defeq \{x \in [N] \col \alpha_x \geq \alpha \}$ for $\alpha >0$. 
We set $\cal U \deq \cal V_{2 + \xi^{1/4}} \setminus \cal V_{\alpha^*(\nu)}$; we note that the set $\cal U \cup \cal V_{\alpha^*(\nu)}$ was denoted by $\cal V$ in \cite{ADK20}, see \cite[eqs.~(1.9) and (3.3)]{ADK20}. Observe that $\cal V_{\alpha^*(\nu)} \subset \cal V_{2 + \xi^{1/4}}$ by the definition of $\alpha^*$ in 
\eqref{eq:largeDegree}. 
From \eqref{eq:norm_H_minus_EH_minus_wh_H_tau} and $\norm{M^{(\cal V_{\alpha^*(\nu)})}} \leq \norm{M}$ for any matrix $M$, we conclude  that $(H- \E H) ^{(\cal V_{\alpha^*(\nu)})} = (\wh{H}^\tau)^{(\cal V_{\alpha^*(\nu)})} + o(1)$, where $o(1)$  is in the sense of operator norm. 
Therefore, in order to prove \eqref{eq:proof_spectral_gap_aux1}, it suffices to prove that, for any constant $\delta > 0$, the matrix $(\wh{H}^\tau)^{(\cal V_{\alpha^*(\nu)})}$ does not have an eigenvalue larger than $\Lambda(\alpha^*(\nu)) + \delta$.

Suppose, by contradiction, that $\lambda$ is an eigenvalue of $(\wh{H}^\tau)^{(\cal V_{\alpha^*(\nu)})}$ with associated normalized eigenvector $\f w$ such that $\abs{\lambda} >  \Lambda(\alpha^*(\nu)) + \delta$. Let $x \in \cal U$. Since $(\supp \f v^\tau_\sigma(x)) \cap \cal V_{\alpha^*(\nu)} = \emptyset$ by 
the definition of $\f v_\sigma^\tau(x)$ in \eqref{eq:def_v_tau}, \cite[Proposition 3.1 (i)]{ADK20} and 
$\alpha^*(\nu) \geq 2 + \kappa \geq \tau$, we deduce from the definition of $\wh H^\tau$ in \eqref{eq:def_wh_H_tau} that
$\f v_\sigma^{\tau}(x)$ is an eigenvector of $(\wh{H}^\tau)^{(\cal V_{\alpha^*(\nu)})}$ with eigenvalue $\sigma \Lambda(\alpha_x)$. 
As $\Lambda(\alpha_x) \leq \Lambda(\alpha^*(\nu))$ for any $x \in \cal U$, we conclude that $\f w \perp \f v^\tau_\sigma(x)$ for all $x \in \cal U$ and $\sigma = \pm$.  Since $\supp \f w \subset (\cal V_{\alpha^*(\nu)})^c$,  we deduce from \eqref{eq:norm_H_minus_EH_minus_wh_H_tau},  
 \cite[Proposition~3.13]{ADK20}\footnote{We stress that the definition of $H$ in \cite{ADK20} differs from that in the current paper; see \cite[Definition~3.6]{ADK20}.} and \eqref{eq:d_range} that    
\begin{equation} \label{eq:proof_spectral_gap_aux2} 
\abs{\lambda} = \abs{\scalar{\f w}{(\wh{H}^\tau)^{(\cal V_{\alpha^*(\nu)})} \f w}} = \abs{\scalar{\f w}{\wh{H}^\tau \f w}}  \leq 1 + \tau + (1 + o(1)) \sum_{x \in \cal V_\tau \setminus \cal V_{\alpha^*(\nu)}} \alpha_x \scalar{\f 1_x}{\f w}^2  + o(1)\,. 
\end{equation}
Note that $\tau = 1 + o(1)$ and $\Lambda(\alpha^*(\nu)) \geq 2$ imply $2 \tau \leq \Lambda(\alpha^*(\nu)) + o(1)$.

We will now apply Lemma~\ref{lem:weak_delocalization} for each $x \in \cal V_\tau\setminus \cal V_{\alpha^*(\nu)}$. 
Note that $\abs{\lambda} \geq \Lambda(\alpha^*(\nu)) +\delta > 2 + o(1) = 2 \tau + C \xi$. In particular, there is a constant $c \equiv c_\delta \in (0,1)$ such that 
$\frac{2\tau + C \xi}{\abs{\lambda}} \leq 1 - c$. 
Owing to the disjointness of the balls $(B_{2r_\star}^\tau(x))_{x\in \cal V_{\tau}}$ (see \cite[Proposition~3.1 (i)]{ADK20}) 
and the locality of $\wh{H}^\tau$, 
we have $((\wh{H}^\tau)^{(\cal V_{\alpha^*(\nu)})}\f w)|_{B_{2r_\star}^\tau(x)} = (\wh{H}^\tau \f w)|_{B_{2r_\star}^\tau(x)}$ 
for all $x \in \cal V_\tau \setminus \cal V_{\alpha^*(\nu)}$. Hence, for these $x$, \eqref{eq:condition_weak_delocalization} holds. 
If $x \in \cal V_{2 + \xi^{1/4}}\setminus \cal V_{\alpha^*(\nu)}$ then $\f w \perp \f v_\pm^\tau(x)$ as demonstrated above. 
Thus, all conditions of Lemma~\ref{lem:weak_delocalization} are satisfied for each $x \in \cal V_{\tau} \setminus \cal V_{\alpha^*(\nu)}$. 
Therefore, Lemma~\ref{lem:weak_delocalization} 
and the disjointness of the balls $(B_{2r_\star}^\tau(x))_{x\in \cal V_{\tau}}$ (compare with the similar argument in \eqref{eq:norm_w_restricted_V_setminus_W})
yield 
\[ 
\sum_{x \in \cal V_\tau \setminus \cal V_{\alpha^*(\nu)}} \scalar{\f 1_x}{\f w}^2 \leq \bigg(1 - \frac{2 \tau + C \xi}{\abs{\lambda}} \bigg)^{-4} 
\bigg( \frac{2 \tau + C \xi}{\abs{\lambda}} \bigg)^{2r_\star} 
\leq c^{-4} ( 1- c) ^{2r_\star} \ll \frac{d}{\log N}\,,  
\] 
where we used $\frac{2\tau + C \xi}{\abs{\lambda}} \leq 1 - c$ in the second step and $r_\star \asymp \sqrt{\log N}$ as well as \eqref{eq:d_range} in the last step. 
Hence, from Proposition~\ref{prop:graphProperties_every_vertex} \ref{item:upper_bound_alpha_x}, we conclude 
$\sum_{x \in \cal V_\tau \setminus \cal V_{\alpha^*(\nu)}} \alpha_x \scalar{\f 1_x}{\f w}^2 = o(1)$.  
Using this in 
 \eqref{eq:proof_spectral_gap_aux2} implies  
$\abs{\lambda} \leq \Lambda(\alpha^*(\nu)) + o(1)$, in contradiction with the assumption on $\lambda$. 
We therefore conclude \eqref{eq:proof_spectral_gap_aux1}. 
\end{proof}

\section{Tools}

In this appendix we summarize several well known results used throughout this paper.

\subsection{Bennett's inequality}

Define the function $h \colon [0,\infty) \to [0,\infty)$ by
\begin{equation} \label{eq:def_h} 
 h(a) \deq (1 + a) \log (1 + a) - a\,. 
\end{equation}
We use the notation $\cal P_\mu$ to denote a Poisson random variable with parameter $\mu \geq 0$, and $\cal B_{n,p}$ to denote a binomial random variable with parameters $n \in \N^*$ and $p \in [0,1]$.
The following estimate is proved in \cite[Section 2.7]{BLM13}.

\begin{lemma}[Bennett] \label{lem:Bennett}
For $0 \leq \mu \leq n$ and  $a > 0$ we have
\begin{equation*}
\P(\cal B_{n,\mu/n} - \mu \geq a \mu) \leq \ee^{-\mu h(a)}\,, \qquad
\P(\cal B_{n,\mu/n} - \mu \leq - a \mu) \leq \ee^{-\mu a^2/2} \leq \ee^{-\mu h(a)}\,,
\end{equation*}
and $\frac{a^2}{2 (1 + a/3)}  \leq h(a) \leq \frac{a^2}{2}$. By taking $n \to \infty$, the same estimates hold with $\cal B_{n,\mu/n}$ replaced with $\cal P_\mu$.
\end{lemma}

\subsection{Perturbation theory} 

The following lemma contains simple perturbation estimates for approximate eigenvalues and eigenvectors. 
Its proof can be found in \cite[Lemma~4.10]{ADK21}. 

\begin{lemma} \label{lem:perturbation_theory} 
Let $M$ be a real symmetric matrix. Let $\eps$, $\Delta >0$ satisfy $5 \eps \leq \Delta$. 
Suppose that $M$ has a unique eigenvalue, $\lambda$, in $[\wh{\lambda}-\Delta, \wh{\lambda} + \Delta]$ for  
some  $\wh{\lambda} \in \R$. 
Let $\f w$ be a corresponding normalized eigenvector of $M$. 
If there exists a normalized vector $\f v$ such that $\norm{(M- \wh{\lambda}) \f v} \leq \eps$ then, for some $\sigma \in \{\pm\}$,  
\[ \lambda - \wh{\lambda} = \scalar{\f v}{(M- \wh{\lambda}) \f v} + O \bigg( \frac{\eps^2}{\Delta} \bigg)\,, 
\qquad \qquad \norm{ \f w - \sigma\f v} = O \bigg( \frac{\eps}{\Delta} \bigg)\,. \] 
\end{lemma}

\subsection{Eigenvalue interlacing}
For an $N \times N$ Hermitian matrix $M$, we use the convention that $\lambda_i(M) = -\infty$ for $i > N$.

\begin{lemma}[Interlacing for minors] \label{lem:interlacing1}
Let $X \subset Y \subset [N]$ and $M$ be an $N \times N$ Hermitian matrix. Then, for all $i \in [N]$,
\begin{equation*}
\lambda_{i + \abs{Y \setminus X}}(M^{(X)})_+ \leq \lambda_i(M^{(Y)})_+ \leq \lambda_i(M^{(X)})_+\,,
\end{equation*}
where $\lambda_+ \deq \max\{\lambda, 0\}$.
\end{lemma}

\begin{lemma}[Interlacing for rank-one perturbations] \label{lem:interlacing2}
Let $M$ be an $N \times N$ Hermitian matrix and $V$ a rank-one positive semidefinite $N \times N$ matrix. Then, for all $i \in [N]$,
\begin{equation*}
\lambda_{i + 1}(M+V) \leq \lambda_i(M) \leq \lambda_{i}(M+V)\,.
\end{equation*}
\end{lemma}

\subsection{Resolvent of adjacency matrix of $\N^*$} 

Define the infinite tridiagonal matrix 
\begin{equation} \label{eq:def_M_adjacency_matrix_N}
M \defeq \begin{pmatrix} 0 & 1 & 0 &    \\ 1 & 0 & 1 & \ddots   \\ 0 & 1 & 0 &  \ddots \\ 
 & \ddots & \ddots & \ddots  
 \end{pmatrix} \,,
\end{equation} 
which we regard as a bounded operator on $\ell^2(\N^*)$. 
Note that $M$ is the adjacency matrix of $\N^*$ with the graph structure induced by regarding adjacent numbers as neighbours.

\begin{lemma} \label{lem:resolvent_adjacency_matrix_positive_integers}
Let $t > 2$. Then $1- t^{-1} M$ 
is invertible with convergent Neumann series and explicit inverse
\[ (1 - t^{-1} M)^{-1}_{1j} = 
t \bigg( \frac{2}{t + \sqrt{t^2 - 4}} \bigg)^j
\] 
for any $j \in \N^*$. 
\end{lemma} 

\begin{proof} 
By the Schur test, $\norm{M} \leq 2$. Therefore,  $t > 2$ implies the invertibility of 
$1 - t^{-1} M$ in $\ell^2(\N^*)$ and the convergence of the Neumann series representation of its inverse.
In particular, 
$a = (a_j)_{j \in \N^*} \in \ell^2(\N^*)$ with $a_j \defeq (1 - t^{-1} M)^{-1}_{1j}$ for 
any $j \in \N^*$. With $e_1 \defeq (\delta_{1j})_{j \in \N^*}$, the definition of the resolvent yields  
$ ( 1 - t^{-1} M) a = e_1$ and, thus, 
\begin{equation} \label{eq:transfer_matrix_relation} 
 a_1 - t^{-1} a_2 = 1\,, \qquad \qquad \begin{pmatrix} a_{j+1} \\ a_j \end{pmatrix} 
= T(t) \begin{pmatrix} a_j \\ a_{j-1} \end{pmatrix}\,, \qquad 
T(t) \defeq \begin{pmatrix} t & - 1 \\ 1 & 0 \end{pmatrix} 
\end{equation} 
for any $j \in \N$ satisfying $j \geq 2$. 
The transfer matrix $T(t)$ has the two eigenvalues, $\gamma$ and $\gamma^{-1}$, where 
\[ \gamma \defeq \frac{2}{t + \sqrt{t^2 - 4}}\,. \] 
As $t > 2$, we have $\gamma < 1$ and, hence, $a \in \ell^2(\N^*)$ implies 
$a_j = \gamma a_{j-1} = \gamma^{j-1}a_1$ for any $j \geq 2$. 
Together with the first relation in \eqref{eq:transfer_matrix_relation}, we obtain  
 $a_1=\frac{2t}{t + \sqrt{t^2  -4}}$, 
which completes the proof.
\end{proof}

\section{Graph properties -- proofs of Propositions~\ref{prop:graphProperties_every_vertex} and~\ref{prop:graphProperties}} \label{sec:proof_prop_graphProperties}

\begin{proof}[Proof of Proposition \ref{prop:graphProperties_every_vertex}] 
We show that each item holds individually with high probability. By Bennett's inequality (see Lemma \ref{lem:Bennett} below) and a union bound, we conclude that
   \[
 \mathbb{P}\pB{\max_{x\in[N]}|S_1(x)|>10 \log N} \leq N \ee^{-10 \log N} \ll 1\,, 
\] 
where we used $d \leq 3 \log N$. This shows \ref{item:upper_bound_alpha_x}.

To show \ref{item:upper_bound_size_B_r}, we abbreviate $D_x = \abs{S_1(x)}$ and note that by part \ref{item:upper_bound_alpha_x} and a union bound, it suffices to show, for any $x \in [N]$,
\begin{equation} \label{Bi_bound_cond}
\P \pb{\abs{B_i(x)} \geq 2 \max\{D_x, d\} d^{i - 1} \,|\, D_x} \leq N^{-2}
\end{equation}
for all $i \leq \frac{1}{3} \frac{\log N}{\log d}$ and $D_x \leq 10 \log N$. If $D_x \geq d$, then \eqref{Bi_bound_cond} follows from from \cite[eq.~(5.12b)]{ADK19}, where the upper and lower bounds in the condition \cite[eq.~(5.13)]{ADK19} follow from $D_x \leq 10 \log N$ and the upper bound on $i$, as well as $D_x \geq d$ and $d \gg \sqrt{\log N}$, respectively. 

For $D_x \leq d$ the bound \eqref{Bi_bound_cond} follows from the monotonicity property 
\begin{equation} \label{eq:monotonicity_degrees} 
 \P (\abs{B_i(x)} \geq L \,|\, D_x = k ) \leq \P( \abs{B_i(x)} \geq L \,|\, D_x = l ) 
\end{equation} 
for any $k \leq l \leq N-1$, combined with \eqref{Bi_bound_cond} on the event $\{D_x = d\}$.

What remains is the proof of \eqref{eq:monotonicity_degrees}. 
To that end, fix $x \in [N]$ and let $k \leq l \leq N -1$. 
For $y \in [N]\setminus \{x \}$, we denote by $B_i^{(x)}(y)$ the ball in $\mathbb{G}|_{[N] \setminus \{x \}}$ 
of radius $i$ around $y$. 
Moreover, let $(K_1, K_2)$ be disjoint random uniformly chosen subsets of $[N]\setminus \{x \}$ of sizes $k$ and $l - k$, 
respectively. 
Then we estimate
\begin{align*}
\mathbb{P}\big(|B_{i}(x)|\geq L\big|D_{x}=k\big) & =\mathbb{P}\bigg(\absbb{\bigcup_{y\in S_{1}(x)}B_{i-1}^{(x)}(y)}\geq L-1\bigg|D_{x}=k\bigg)\\
 & =\mathbb{P}\bigg(\absbb{\bigcup_{y\in K_{1}}B_{i-1}^{(x)}(y)}\geq L-1\bigg)\\
 & \leq\mathbb{P}\bigg(\absbb{\bigcup_{y\in K_{1}\cup K_{2}}B_{i-1}^{(x)}(y)}\geq L-1\bigg)\\
 & =\mathbb{P}\bigg(\absbb{\bigcup_{y\in S_{1}(x)}B_{i-1}^{(x)}(y)}\geq L-1\bigg|D_{x}=l\bigg)\\
 & =\mathbb{P}(|B_{i}(x)|\geq L|D_{x}=l)\,. 
\end{align*}
This proves \eqref{eq:monotonicity_degrees} and, thus, completes the proof of Proposition~\ref{prop:graphProperties_every_vertex}. 
\end{proof}

\begin{proof}[Proof of Proposition \ref{prop:graphProperties}] 
We show that each item holds individually with high probability. 
For the proof of \ref{item:size_cal_V}, we estimate 
$\P ( \abs{\cal V} \geq t N^{\mu}) \leq t^{-1} N^{-\mu} \E \abs{\cal V} = t^{-1} N^{-\mu } N \P(\alpha_1 \geq \alpha^*)\leq t ^{-1} $ 
for $t>0$ by the definitions of $\cal V$ and $\alpha^*$ in \eqref{eq:def_cal_V} and \eqref{eq:largeDegree}, respectively. 
This proves \ref{item:size_cal_V}. 
Item \ref{item:B_r_tree} is a consequence of \cite[Lemma 5.5]{ADK19} with $k=1$, Proposition~\ref{prop:graphProperties_every_vertex} \ref{item:upper_bound_alpha_x},  $r \leq \big(\frac{1}{5} - \frac{\mu}{4} \big) \frac{\log N}{\log d}$, and $d \geq \sqrt{\log N}$.  
Item \ref{item:B_r_disjoint} follows from \cite[eq.~(9.5)]{ADK21} with the choices $\tau=\alpha^*$ and $n = 2$, $\exp(- d h(\alpha^* -1-\frac{3}{d})) \leq N^{\mu-1+o(1)}$ and $r\leq \big(\frac{1}{3}-\mu\big)\frac{\log N}{\log d}$.

To prove \ref{item:upper_degree_in_balls_around_vertices_in_V}, we set $\eta/2 \defeq r \frac{\log d}{\log N}$
and conclude 
from \ref{item:size_cal_V}, Proposition~\ref{prop:graphProperties_every_vertex} \ref{item:upper_bound_alpha_x} and \ref{item:upper_bound_size_B_r} that 
$\abs{\bigcup_{x \in \cal V} B_r(x)} \leq N^{\mu + \eta/2 + o(1)}$.  
Thus, we obtain \ref{item:upper_degree_in_balls_around_vertices_in_V} by arguing similarly 
 as in \cite[Proof of Proposition~4.4]{ADK21}, see especially \cite[eq.~(9.10)]{ADK21}.
\end{proof}

\section{Eigenfunction correlator and dynamical localization -- proof of Corollary~\ref{cor:dynamical_loc}} \label{sec:dynamical}

Before proving \eqref{eq:projection_of_evolution_with_restricted_energy}, 
we use it to show \eqref{dynamical_localization}. 
If $\dd (x,y) \leq 1$ then \eqref{dynamical_localization} holds trivially. If $\dd (x,y) \geq 2$, \eqref{dynamical_localization} follows 
directly from \eqref{eq:projection_of_evolution_with_restricted_energy} by choosing $r = \dd(x,y) - 1$
 and $F(\lambda) = \ee^{-\ii t \lambda}$. This completes the proof of \eqref{dynamical_localization} assuming 
\eqref{eq:projection_of_evolution_with_restricted_energy}. 

We now turn to the proof of \eqref{eq:projection_of_evolution_with_restricted_energy}. First, we introduce some notation. By  \eqref{eq:u_approx_u_x} from the proof of Theorem~\ref{thm:main}, with high probability, there is a one-to-one correspondence between eigenvalues satisfying \eqref{lambda_condition} and vertices $x$ in some subset of $\cal W$. Under this correspondence, we denote such an eigenvalue and its corresponding eigenvector by $\lambda_x$ and 
$\f w(x)$, respectively.

Let $x \in [N]$. Let $J\subset I$ be an interval, where $I \deq [\Lambda(\alpha^*(\mu))+ \kappa, \sqrt{d}/2]$. Let $F$ be a function as in the statement of Corollary~\ref{cor:dynamical_loc}. Then we obtain   
\[ 
\Pi_J(H) F(H) \f 1_x = \sum_{a \in \cal W \colon \lambda_a \in J} F(\lambda_a) \scalar{\f w(a)}{\f 1_x} \f w(a). 
\] 
By following the proof of Lemma \ref{lem:UXBallY}, we see that if $\mu < \frac{1}{24}$ and $r \leq \frac{1}{6} \frac{\log N}{\log d}$ then, with high probability, 
\begin{equation}\label{eq:u_a_small_on_B_r_b} 
\norm{\f u(a) \vert_{B_r(b)}} \leq N^{-3/8}
\end{equation}
for all $a \neq b \in \cal W$.

We now fix $r \geq 0$ as in the statement. 
Without loss of generality, we can assume that $r \leq 2 \frac{\log N}{\log d}$, because $\diam(\bb G) \leq 2 \frac{\log N}{\log d}$ with high probability.
Suppose first that there is $b \in \cal W$ such that $x \in B_{R}(b)$ with $R \deq \min \{ \frac{r}{2}, \frac{1}{6} \frac{\log N}{\log d}\}$. Note that this $b$ is unique by Proposition~\ref{prop:graphProperties} 
\ref{item:B_r_disjoint} and $R \geq r/12$ as $r \leq 2 \frac{\log N}{\log d}$. 
In this case, \eqref{eq:projection_of_evolution_with_restricted_energy} follows from 
\begin{multline} \label{eq:projection_evolution_restricted_x_inside}  
\norm{(\Pi_J(H) F(H)\f 1_x) |_{B_r(x)^c}}^2 \lesssim \norm{\f w(b) |_{B_r(x)^c}}^2 + \sum_{a \in \cal W\setminus \{b \}} 
\abs{\scalar{\f w(a)}{\f 1_x}}^2 \\ 
  \lesssim \ee^{-c R} + N^{-3/4 + \mu} + N^{2(\eta- \zeta) + \mu} \lesssim \ee^{-c r} . 
\end{multline} 
Here, we used that $\norm{\f w(b) |_{B_r(x)^c}} \leq \norm{\f w(b) |_{B_r(x)^c}} \lesssim \ee^{-cR}$, 
which is a consequence of $B_R(b) \subset B_r(x)$ as well as \eqref{eq:exponential_decay_true_eigenvector} 
if $\alpha_x \leq 3$ or $R \geq 10 \log \log N$ due to \eqref{Lambda_alpha_bound} and \eqref{eq:improved_exponential_decay_true_eigenvector} otherwise. 
Moreover, we employed 
$\abs{\scalar{\f w(a)}{\f 1_x}}^2 \leq \norm{\f w(a)|_{B_R(b)}}^2 \lesssim N^{-3/4} + N^{2(\eta- \zeta)}$ 
due to \eqref{eq:u_approx_u_x} and \eqref{eq:u_a_small_on_B_r_b}. 
The last step in \eqref{eq:projection_evolution_restricted_x_inside} follows from $r/12 \leq R \leq \frac{1}{6} \frac{\log N}{\log d}$ and $ N^{-3/4 + \mu} + N^{2(\eta- \zeta) + \mu} \lesssim \ee^{-c \log N/\log d}$ by appropriate choices of $\eta$ and $\zeta$ if $\mu < \frac{1}{24}$. 

Hence, it remains to consider the case $x \notin B_{R} (\cal W)$. Under this condition, we will 
show that 
\begin{equation} \label{eq:norm_Pi_1_x_x_away_from_cal_W} 
 \norm{\Pi_J(H) \f 1_x} \leq C\ee^{-cr} 
\end{equation} 
if $r \leq 2\frac{\log N}{\log d}$, which will directly imply \eqref{eq:projection_of_evolution_with_restricted_energy} in this missing case. 

For the proof of \eqref{eq:norm_Pi_1_x_x_away_from_cal_W}, we note that $\norm{\Pi_J(H) \f 1_x}^2 
\leq \norm{\Pi_I(H) \f 1_x}^2$ as $J \subset I = [\Lambda(\alpha^*)+ \kappa, \sqrt{d}/2]$ 
and, for $k \in \N$ such that $2k < R$, estimate 
\[ 
\norm{\Pi_I(H) \f 1_x}^2 (\Lambda(\alpha^*) + \kappa)^{2k} \leq \scalar{\f 1_x}{\Pi_I(H) H^{2k} \Pi_I(H) \f 1_x}
\leq \scalar{\f 1_x}{H^{2k} \f 1_x} = \scalar{\f 1_x}{(H^{(\cal W)})^{2k} \f 1_x}, 
\] 
where the last step follows from the locality of $H$, $2k < r/2$ and $x \notin B_{R} (\cal W)$. Moreover, by locality of $H$, we have
\begin{align*}
\scalar{\f 1_x}{(H^{(\cal W)})^{2k} \f 1_x} &= \scalar{\f 1_x}{(H^{(\cal W)} - (\E H)^{(\cal W \cup B_k(x))})^{2k} \f 1_x}
\\
&\leq \norm{H^{(\cal W)} - (\E H)^{(\cal W \cup B_k(x))}}^{2k}
\\
&\leq \pb{\norm{H^{(\cal W)} - (\E H)^{(\cal W)}} + \norm{(\E H)^{(\cal W)} - (\E H)^{(\cal W \cup B_k(x))}}}^{2k}
\\
&\leq (\Lambda(\alpha^*) + o(1))^{2k}
\end{align*}
with high probability,
where we used Proposition~\ref{prop:graphProperties_every_vertex}, Lemma~\ref{lem:bound_A_minus_EA}, \eqref{Lambda_alpha_bound} and \eqref{eq:def_cal_W} in the last step. 
Combining these bounds, choosing $k \asymp r$ and using $r \leq 2\frac{\log N}{\log d}$ implies the existence 
of constants $C>0$ and $c>0$ depending on $\kappa$ such that $\norm{\Pi_I(H) \f 1_x} \leq C \ee^{-cr}$ if $d\asymp \log N$ as $\Lambda(\alpha^*) \asymp 1$ in this case. 
For $d \asymp \log N$, this proves \eqref{eq:norm_Pi_1_x_x_away_from_cal_W} since $\norm{\Pi_J(H) \f 1_x} \leq \norm{\Pi_I(H) \f 1_x}$ as explained above. Otherwise, we use $T$ as in Corollary~\ref{cor:Lambda_alpha_1_2_Lambda_alpha_mu} for $\eps = 1/2$. We choose $\nu <\frac{1}{24}$ satisfying the assumptions made on $\mu$ in this proof up to now, and consider $\cal W$ to be defined 
with respect to this $\nu$ instead of $\mu$, i.e.\ $\cal W = \{ x \in [N] \colon \Lambda(\alpha_x) 
\geq \Lambda(\alpha^*(\nu)) + \kappa/2 \}$ (cf.\ \eqref{eq:def_cal_W}). 
For any constant $\mu \in (0,\nu)$ and $I = [\Lambda(\alpha^*(\mu) + \kappa, \sqrt{d}/2]$, we obtain 
$\norm{\Pi_I(H) \f 1_x}^2 \leq \big(\frac{\Lambda(\alpha^*(\nu)) + o(1)}{\Lambda(\alpha^*(\mu)) + \kappa} \big)^{2k} \leq (1 - \eps)^{2k}$ for some constant $\eps \in (0,1)$, where the last step follows from 
Corollary~\ref{cor:Lambda_alpha_1_2_Lambda_alpha_mu} by assuming $\frac{\log N}{d} \geq T$. 
Since $\norm{\Pi_J(H) \f 1_x} \leq \norm{\Pi_I(H) \f 1_x}$, this completes the proof of 
\eqref{eq:norm_Pi_1_x_x_away_from_cal_W} in the missing regime and, thus, 
the one of \eqref{eq:projection_of_evolution_with_restricted_energy}. 

Next, we show \eqref{ef_correlator}. 
First, we assume that $x,y \notin B_r(\cal W)$ with $r \deq \frac{1}{20} \frac{\log N}{\log d}$. 
Then \eqref{ef_correlator} follows from the Cauchy-Schwarz inequality, \eqref{eq:norm_Pi_1_x_x_away_from_cal_W}  and $\dd(x,y) \leq \diam(\mathbb G) \leq 2 \frac{\log N}{\log d}$. 
Finally, we assume that $x \in B_r(\cal W)$ with $r$ as above. By Proposition~\ref{prop:graphProperties} \ref{item:B_r_disjoint}, there is a unique $b \in \cal W$ such that $\dd(x,b) \leq r$ and $\dd(x,a) \geq 2 r$ 
for all $a \in\cal W \setminus \{ b\}$. 
As in the second step of \eqref{eq:projection_evolution_restricted_x_inside}, we conclude from \eqref{eq:exponential_decay_true_eigenvector} and \eqref{eq:improved_exponential_decay_true_eigenvector} that 
$\abs{\scalar{\f 1_x}{\f w(b)}} \lesssim \ee^{-c \dd(x,b)}$ and $\abs{\scalar{\f w(b)}{\f 1_y}} \lesssim 
\ee^{-c \dd(y,b)}$ for some constant $c>0$.
Hence, arguing similarly as in \eqref{eq:projection_evolution_restricted_x_inside} as well as using 
$\dd(x,b) + \dd (y,b) \geq \dd(x,y)$ and $\abs{\scalar{\f w(a)}{\f 1_y}} \leq 1$ for all $a \in \cal W \setminus \{ b \}$ yields \eqref{ef_correlator}. 
This completes the proofs of \eqref{ef_correlator} and Corollary~\ref{cor:dynamical_loc}.

\section{The expected optimal range for $\mu$} \label{sec:optimal_loc}

In this appendix we give a heuristic argument that yields the expected optimal size of the localized phase, by providing the optimal upper bound on the exponent $\mu$. We also explain how the conclusion of this argument relates to that of \cite{tarzia2022}.

The optimal value for $\mu$ can already be seen from the simple case where $\cal V = \cal W = \{x,y\}$ consists of two vertices, $x$ and $y$, which are are in resonance (i.e.\ with normalized degrees close to each other). In that case, as in Proposition~\ref{prop:EstimateExy}, we obtain with high probability
\begin{equation*}
(H - \lambda(x)) \f u(x) = \epsilon_y(x) \f 1_y \,, \qquad \abs{\epsilon_y(x)} \leq N^{-1/2 + o(1)}\,,
\end{equation*}
where the estimate on $\epsilon_y(x)$ follows as in the proof of Proposition~\ref{prop:EstimateExy}.

We conclude that the tunnelling amplitude between $x$ and $y$ is with high probability bounded by $N^{-1/2 + o(1)}$.
By Mott's criterion, we therefore expect localization at vertex $x$ whenever the typical eigenvalue spacing $N^{-\rho_b(\lambda(x)) + o(1)}$ is much larger than $N^{-1/2 + o(1)}$, i.e. $\rho_b(\lambda(x)) < \frac12$. This condition holds precisely for the $N^\mu$ largest eigenvalues of the semilocalized phase for any $\mu < \frac12$ (see Appendix~\ref{sec:alpha_qual}). This conclusion coincides with that of \cite{tarzia2022}, obtained by a different argument.

The preceding argument has two important shortcomings. First, it only considers a single pair of resonant vertices, instead of the $N^{\mu + o(1)}$ vertices associated with the localized phase. Second, it assumes that the eigenvalue spacing around $\lambda(x)$ is typical, i.e. $N^{-\mu}$, which is true for most eigenvalues in the localized phase but is false for all eigenvalues simultaneously.

To address these shortcomings, we define a random matrix model that heuristically captures the main features of the localized phase. Define the $\cal W \times \cal W$ matrix $M = (M_{xy})_{x,y \in \cal W}$ through
\begin{equation*}
M_{xy} \deq \scalar{\f u(x)}{H \f u(y)}\,.
\end{equation*}
Recall from Proposition \ref{prop:mainRigidity} that $(\f u(x))_{x \in \cal W}$ are approximate eigenvectors of $H$, and that all eigenvalues of $H$ in the localized phase arise by perturbation from the approximate eigenvalues $(\lambda(x))_{x \in \cal W}$. Hence, localization is tantamount to the matrix $M$ being close to diagonal.

We make the simplifying assumption that the vectors $(\f u(x) \col x \in \cal W)$ are orthogonal, i.e.\ we neglect their overlaps for the purpose of this heuristic argument. As above, the proof of Proposition~\ref{prop:EstimateExy} implies that
\begin{equation*}
M_{xy} \approx \delta_{xy} \, \lambda(x) + O(N^{-1/2 + o(1)})
\end{equation*}
with high probability for all $x,y \in \cal W$, since (under the orthogonality assumption)
\begin{align*}
M_{xy} &\approx \delta_{xy}\, \lambda(x)  + (1 - \delta_{xy}) \, \scalar{\f u(x)} {(H - \lambda(y)) \f u(y)}
\\
&= \delta_{xy} \, \lambda(x) + (1 - \delta_{xy}) \, \epsilon_x(y) \scalar{\f u(x)}{\f 1_x}
\\
&\approx \delta_{xy} \, \lambda(x) + (1 - \delta_{xy}) \, \epsilon_x(y)\,,
\end{align*}
where we used that $\scalar{\f u(x)}{\f 1_z} = 0$ for any $z \in \cal W \setminus \{x\}$.

Hence, the random variables $M_{xy}$ are typically of order $N^{-1/2 + o(1)}$. Assuming that they are mean zero, independent up to the symmetry constraint $M_{xy} = M_{yx}$, and have light enough tails, we therefore arrive at the following toy model for $M$:
\begin{equation} \label{def_Mt}
M \approx M(t) \deq D + \sqrt{t} \, W\,, \qquad t = N^{\mu - 1 + o(1)}\,,
\end{equation}
where $D = \diag(\lambda(x) \col x \in \cal W)$ and $W$ is an independent $\cal W \times \cal W$ Wigner matrix, normalized so that its entries have variance $1 / \abs{\cal W}$ (and hence $W$ has typically a norm of order $1$).

The matrix $M(t)$ from \eqref{def_Mt} is a deformed Wigner matrix, the kind of which has been extensively studied in the random matrix theory literature. Assuming that $W$ is a GOE matrix (i.e.\ with a Gaussian law), the law of $M(t)$ is governed by Dyson Brownian motion at time $t$ starting from the diagonal matrix $D$. We refer for instance to \cite{VW1, benaych2021eigenvectors, benigni2020eigenvectors, LS1, von2019non}. 
Although not stated explicitly in these works, the analysis of Dyson Brownian motion, as for instance in \cite{benigni2020eigenvectors}, implies that, with high probability, for any $x \in \cal W$, the eigenvector associated with the eigenvalue $\lambda(t,x)$ of the $\cal W \times \cal W$ matrix $M(t)$ remains localized for
\begin{equation} \label{heuristic_localization}
t \ll \min \hb{\Delta(x) \,, \abs{\cal W} \, \Delta(x)^2}\,,
\end{equation}
where $\Delta(x) \deq \min_{y \neq x} \abs{\lambda(x) - \lambda(y)}$ is the eigenvalue spacing of $D$  around $\lambda(x)$. Indeed, the eigenvalues $\lambda(t,x)$ of $M(t)$ satisfy
\begin{equation*}
\dd \lambda(t,x) = \sqrt{\frac{2}{\abs{\cal W}}} \, \dd B(t,x) + \frac{1}{\abs{\cal W}} \sum_{y \neq x} \frac{\dd t}{\lambda(t,x) - \lambda(t,y)}\,,
\end{equation*}
where $(B(t,x) \col t \geq 0, x \in \cal W)$ is a family of independent standard Brownian motions. Heuristically, we therefore see immediately that if \eqref{heuristic_localization} holds then $\abs{\lambda(t,x) - \lambda(0,x)} \ll \Delta(x)$, i.e.\ the shift in $\lambda(x)$ is negligible compared to $\Delta(x)$. A similar calculation using Dyson Brownian motion for the eigenvectors (see e.g.\ \cite{benigni2020eigenvectors}) shows that they remain localized under the condition \eqref{heuristic_localization}. Another heuristic way of arriving at the condition \eqref{heuristic_localization} is to notice that under it all finite-order corrections in perturbation theory to $\lambda(t,x)$ are $o(\Delta(x))$, and similarly all finite-order corrections in perturbation theory to the eigenvectors are $o(1)$.

We remark that the condition \eqref{heuristic_localization} is strictly weaker than Mott's criterion, which reads $t \ll \Delta(x)^2$. Hence, for the model \eqref{def_Mt}, Mott's criterion is sufficient but not necessary for localization. Heuristically, this discrepancy stems from the special relationship between the mean-field matrix $W$ and the diagonal eigenvector basis of $D$, in which all entries of $W$ are of order $N^{-\mu/2}$, while the norm of $W$ is of order $1$.

The typical eigenvalue spacing of $D$ is $N^{-\mu + o(1)}$, i.e.\ $\Delta(x) = N^{-\mu + o(1)}$ for most $x \in \cal W$. In contrast, we expect the minimal eigenvalue spacing $\min_{x \in \cal W} \Delta(x)$ to be $N^{-2 \mu + o(1)}$, as follows from basic extreme value theory (cf.\ the birthday paradox) and the assumption that the random variables $(\lambda(x) \col x \in \cal W)$ are independent. Using $t = N^{\mu - 1 + o(1)}$ and $\abs{\cal W} = N^{\mu + o(1)}$, from the condition \eqref{heuristic_localization} we therefore expect the following behaviour for the top $N^\mu$ eigenvectors.
\begin{enumerate}[label=(\roman*)]
\item
For $\mu < \frac{1}{2}$ most eigenvectors are localized.
\item
For $\mu < \frac{1}{4}$ all eigenvectors are localized.
\end{enumerate}
For $\frac{1}{4} \leq \mu < \frac{1}{2}$, while most eigenvectors are localized, we expect some eigenvectors to hybridize with others associated with resonant vertices. These hybridized eigenvectors have a divergent localization length \eqref{def_ll}. In fact, a more refined analysis of Dyson Brownian motion in \eqref{def_Mt}, as performed for instance in \cite{benigni2020eigenvectors}, yields a precise picture of the hybridized structure of such eigenvectors. We shall not go into further details in this appendix.

In particular, the optimal bound for $\mu$ in Theorem \ref{thm:main} is expected to be $\frac{1}{4}$.
In our proof we obtain the bound $\frac{1}{24}$ because both our upper bound on the tunnelling amplitude (Proposition \ref{prop:mainRigidity}) is worse than $N^{-1/2 + o(1)}$ and because our lower bound on the eigenvalue spacing (Proposition \ref{prop:MNonClustering}) is worse than $N^{-2 \mu + o(1)}$. Both of these bounds, however, can be somewhat improved by a suitable refinement of our proof.

\section*{Acknowledgements}
JA acknowledges funding from the European Union's Horizon 2020 research and innovation programme under the Marie Sklodowska-Curie grant 
agreement No.\ 895698 as well as funding from the Deutsche Forschungsgemeinschaft (DFG, German Research Foundation) under Germany’s Excellence Strategy - GZ 2047/1, project-id 390685813. 
RD acknowledges funding from the European Research Council (ERC) under the European Union's Horizon 2020 research and innovation programme (grant agreement No.\ 884584). AK acknowledges funding from the European Research Council (ERC) and the Swiss State Secretariat for Education, Research and Innovation (SERI) through the consolidator grant ProbQuant. JA and AK acknowledge funding from the European Research Council (ERC) under the European Union’s Horizon 2020 research and innovation programme (grant agreement No.\ 715539), funding from the Swiss National Science Foundation through the NCCR SwissMAP grant, and support from the US National Science Foundation under Grant No.\ DMS-1928930 during their participation in the program \emph{Universality and Integrability in Random Matrix Theory and Interacting Particle Systems} hosted by the Mathematical Sciences Research Institute in Berkeley, California during the Fall semester of 2021.

{
\bibliographystyle{amsplain-nodash}
\bibliography{bibliography0}
}

\bigskip

\noindent
Johannes Alt, University of Bonn, Institute for Applied Mathematics, \href{mailto:johannes.alt@iam.uni-bonn.de}{johannes.alt@iam.uni-bonn.de}.
\\[0.3em]
Rapha\"el Ducatez, Université Lyon 1, Institut Camille Jordan, \href{mailto:ducatez@math.univ-lyon1.fr}{ducatez@math.univ-lyon1.fr}.
\\[0.3em]
Antti Knowles, University of Geneva, Section of Mathematics, \href{mailto:antti.knowles@unige.ch}{antti.knowles@unige.ch}.

\end{document}

%% file: fig11.pdf_tex
\begingroup%
  \makeatletter%
  \providecommand\color[2][]{%
    \errmessage{(Inkscape) Color is used for the text in Inkscape, but the package 'color.sty' is not loaded}%
    \renewcommand\color[2][]{}%
  }%
  \providecommand\transparent[1]{%
    \errmessage{(Inkscape) Transparency is used (non-zero) for the text in Inkscape, but the package 'transparent.sty' is not loaded}%
    \renewcommand\transparent[1]{}%
  }%
  \providecommand\rotatebox[2]{#2}%
  \newcommand*\fsize{\dimexpr\f@size pt\relax}%
  \newcommand*\lineheight[1]{\fontsize{\fsize}{#1\fsize}\selectfont}%
  \ifx\svgwidth\undefined%
    \setlength{\unitlength}{410.3947034bp}%
    \ifx\svgscale\undefined%
      \relax%
    \else%
      \setlength{\unitlength}{\unitlength * \real{\svgscale}}%
    \fi%
  \else%
    \setlength{\unitlength}{\svgwidth}%
  \fi%
  \global\let\svgwidth\undefined%
  \global\let\svgscale\undefined%
  \makeatother%
  \begin{picture}(1,0.64397777)%
    \lineheight{1}%
    \setlength\tabcolsep{0pt}%
    \put(0,0){\includegraphics[width=\unitlength,page=1]{fig11.pdf}}%
    \put(0.07846567,0.54245844){\color[rgb]{0.01568627,0.01568627,0.01568627}\makebox(0,0)[lt]{\lineheight{1.25}\smash{\begin{tabular}[t]{l}delocalized\end{tabular}}}}%
    \put(0,0){\includegraphics[width=\unitlength,page=2]{fig11.pdf}}%
    \put(0.52451254,0.63306815){\color[rgb]{0.01568627,0.01568627,0.01568627}\makebox(0,0)[lt]{\lineheight{1.25}\smash{\begin{tabular}[t]{l}$b$\end{tabular}}}}%
    \put(0.33760082,0.0097442){\color[rgb]{1,1,1}\makebox(0,0)[lt]{\lineheight{1.25}\smash{\begin{tabular}[t]{l}$-2$\end{tabular}}}}%
    \put(0.49827473,0.0097442){\color[rgb]{1,1,1}\makebox(0,0)[lt]{\lineheight{1.25}\smash{\begin{tabular}[t]{l}$0$\end{tabular}}}}%
    \put(0.64525865,0.0097442){\color[rgb]{1,1,1}\makebox(0,0)[lt]{\lineheight{1.25}\smash{\begin{tabular}[t]{l}$2$\end{tabular}}}}%
    \put(0.69474603,0.51399525){\color[rgb]{0.01568627,0.01568627,0.01568627}\makebox(0,0)[lt]{\lineheight{1.25}\smash{\begin{tabular}[t]{l}$b_*$\end{tabular}}}}%
    \put(0.34256647,-0){\makebox(0,0)[lt]{\lineheight{1.25}\smash{\begin{tabular}[t]{l}-2\end{tabular}}}}%
    \put(0.63848566,-0){\makebox(0,0)[lt]{\lineheight{1.25}\smash{\begin{tabular}[t]{l}2\end{tabular}}}}%
    \put(0.07846564,0.45198811){\color[rgb]{0.01568627,0.01568627,0.01568627}\makebox(0,0)[lt]{\lineheight{1.25}\smash{\begin{tabular}[t]{l}localized\end{tabular}}}}%
    \put(0,0){\includegraphics[width=\unitlength,page=3]{fig11.pdf}}%
    \put(0.07846596,0.36360473){\color[rgb]{0.01568627,0.01568627,0.01568627}\makebox(0,0)[lt]{\lineheight{1.25}\smash{\begin{tabular}[t]{l}semilocalized\end{tabular}}}}%
    \put(0,0){\includegraphics[width=\unitlength,page=4]{fig11.pdf}}%
    \put(0.89368933,0.19113373){\color[rgb]{0.01568627,0.01568627,0.01568627}\makebox(0,0)[lt]{\lineheight{1.25}\smash{\begin{tabular}[t]{l}$1$\end{tabular}}}}%
    \put(0,0){\includegraphics[width=\unitlength,page=5]{fig11.pdf}}%
    \put(0.98129994,0.00788218){\makebox(0,0)[lt]{\lineheight{1.25}\smash{\begin{tabular}[t]{l}$\lambda$\end{tabular}}}}%
    \put(0,0){\includegraphics[width=\unitlength,page=6]{fig11.pdf}}%
  \end{picture}%
\endgroup%

%% file: fig12.pdf_tex
\begingroup%
  \makeatletter%
  \providecommand\color[2][]{%
    \errmessage{(Inkscape) Color is used for the text in Inkscape, but the package 'color.sty' is not loaded}%
    \renewcommand\color[2][]{}%
  }%
  \providecommand\transparent[1]{%
    \errmessage{(Inkscape) Transparency is used (non-zero) for the text in Inkscape, but the package 'transparent.sty' is not loaded}%
    \renewcommand\transparent[1]{}%
  }%
  \providecommand\rotatebox[2]{#2}%
  \newcommand*\fsize{\dimexpr\f@size pt\relax}%
  \newcommand*\lineheight[1]{\fontsize{\fsize}{#1\fsize}\selectfont}%
  \ifx\svgwidth\undefined%
    \setlength{\unitlength}{403.30137103bp}%
    \ifx\svgscale\undefined%
      \relax%
    \else%
      \setlength{\unitlength}{\unitlength * \real{\svgscale}}%
    \fi%
  \else%
    \setlength{\unitlength}{\svgwidth}%
  \fi%
  \global\let\svgwidth\undefined%
  \global\let\svgscale\undefined%
  \makeatother%
  \begin{picture}(1,0.56115775)%
    \lineheight{1}%
    \setlength\tabcolsep{0pt}%
    \put(0,0){\includegraphics[width=\unitlength,page=1]{fig12.pdf}}%
    \put(0.1884084,0.54003199){\makebox(0,0)[lt]{\lineheight{1.25}\smash{\begin{tabular}[t]{l}$H$\end{tabular}}}}%
    \put(0.11931699,0.26287866){\makebox(0,0)[lt]{\lineheight{1.25}\smash{\begin{tabular}[t]{l}$\wh H$\end{tabular}}}}%
    \put(0,0){\includegraphics[width=\unitlength,page=2]{fig12.pdf}}%
    \put(0.40572332,0.31235548){\rotatebox{36.477448}{\makebox(0,0)[lt]{\lineheight{1.25}\smash{\begin{tabular}[t]{l}localization\end{tabular}}}}}%
    \put(0.39319113,0.21240754){\rotatebox{-14.571066}{\makebox(0,0)[lt]{\lineheight{1.25}\smash{\begin{tabular}[t]{l}hybridization\end{tabular}}}}}%
    \put(0.92555755,0.35781887){\makebox(0,0)[lt]{\lineheight{1.25}\smash{\begin{tabular}[t]{l}$\spec(\wh H)$\end{tabular}}}}%
    \put(0.92698395,0.32011621){\makebox(0,0)[lt]{\lineheight{1.25}\smash{\begin{tabular}[t]{l}$\spec(H)$\end{tabular}}}}%
    \put(0.76117657,0.54003199){\makebox(0,0)[lt]{\lineheight{1.25}\smash{\begin{tabular}[t]{l}$H$\end{tabular}}}}%
    \put(0.76117657,0.22182745){\makebox(0,0)[lt]{\lineheight{1.25}\smash{\begin{tabular}[t]{l}$H$\end{tabular}}}}%
    \put(0.92555755,0.0396143){\makebox(0,0)[lt]{\lineheight{1.25}\smash{\begin{tabular}[t]{l}$\spec(\wh H)$\end{tabular}}}}%
    \put(0.92698395,0.00191176){\makebox(0,0)[lt]{\lineheight{1.25}\smash{\begin{tabular}[t]{l}$\spec(H)$\end{tabular}}}}%
    \put(0.35278932,0.1032552){\makebox(0,0)[lt]{\lineheight{1.25}\smash{\begin{tabular}[t]{l}$\spec(\wh H)$\end{tabular}}}}%
  \end{picture}%
\endgroup%

%% file: fig8.pdf_tex
\begingroup%
  \makeatletter%
  \providecommand\color[2][]{%
    \errmessage{(Inkscape) Color is used for the text in Inkscape, but the package 'color.sty' is not loaded}%
    \renewcommand\color[2][]{}%
  }%
  \providecommand\transparent[1]{%
    \errmessage{(Inkscape) Transparency is used (non-zero) for the text in Inkscape, but the package 'transparent.sty' is not loaded}%
    \renewcommand\transparent[1]{}%
  }%
  \providecommand\rotatebox[2]{#2}%
  \newcommand*\fsize{\dimexpr\f@size pt\relax}%
  \newcommand*\lineheight[1]{\fontsize{\fsize}{#1\fsize}\selectfont}%
  \ifx\svgwidth\undefined%
    \setlength{\unitlength}{384.33022029bp}%
    \ifx\svgscale\undefined%
      \relax%
    \else%
      \setlength{\unitlength}{\unitlength * \real{\svgscale}}%
    \fi%
  \else%
    \setlength{\unitlength}{\svgwidth}%
  \fi%
  \global\let\svgwidth\undefined%
  \global\let\svgscale\undefined%
  \makeatother%
  \begin{picture}(1,0.37224542)%
    \lineheight{1}%
    \setlength\tabcolsep{0pt}%
    \put(0,0){\includegraphics[width=\unitlength,page=1]{fig8.pdf}}%
    \put(0.88862493,0.00398422){\makebox(0,0)[lt]{\lineheight{1.25}\smash{\begin{tabular}[t]{l}$\lambda$\end{tabular}}}}%
    \put(0.78698585,0.00398422){\makebox(0,0)[lt]{\lineheight{1.25}\smash{\begin{tabular}[t]{l}$\lambda_{\max}$\end{tabular}}}}%
    \put(0.35680697,0.00398422){\makebox(0,0)[lt]{\lineheight{1.25}\smash{\begin{tabular}[t]{l}$2$\end{tabular}}}}%
    \put(0.0690267,0.03928535){\makebox(0,0)[lt]{\lineheight{1.25}\smash{\begin{tabular}[t]{l}$0$\end{tabular}}}}%
    \put(-0.00157539,0.31785653){\makebox(0,0)[lt]{\lineheight{1.25}\smash{\begin{tabular}[t]{l}$\diam(\bb G)$\end{tabular}}}}%
    \put(0.0690267,0.35622725){\makebox(0,0)[lt]{\lineheight{1.25}\smash{\begin{tabular}[t]{l}$\ell$\end{tabular}}}}%
    \put(0.42664165,0.15286254){\makebox(0,0)[lt]{\lineheight{1.25}\smash{\begin{tabular}[t]{l}$\frac{\lambda}{2 \sqrt{\lambda^2 - 4}}$\end{tabular}}}}%
  \end{picture}%
\endgroup%

%% file: fig7.pdf_tex
\begingroup%
  \makeatletter%
  \providecommand\color[2][]{%
    \errmessage{(Inkscape) Color is used for the text in Inkscape, but the package 'color.sty' is not loaded}%
    \renewcommand\color[2][]{}%
  }%
  \providecommand\transparent[1]{%
    \errmessage{(Inkscape) Transparency is used (non-zero) for the text in Inkscape, but the package 'transparent.sty' is not loaded}%
    \renewcommand\transparent[1]{}%
  }%
  \providecommand\rotatebox[2]{#2}%
  \newcommand*\fsize{\dimexpr\f@size pt\relax}%
  \newcommand*\lineheight[1]{\fontsize{\fsize}{#1\fsize}\selectfont}%
  \ifx\svgwidth\undefined%
    \setlength{\unitlength}{457.64463829bp}%
    \ifx\svgscale\undefined%
      \relax%
    \else%
      \setlength{\unitlength}{\unitlength * \real{\svgscale}}%
    \fi%
  \else%
    \setlength{\unitlength}{\svgwidth}%
  \fi%
  \global\let\svgwidth\undefined%
  \global\let\svgscale\undefined%
  \makeatother%
  \begin{picture}(1,0.24211126)%
    \lineheight{1}%
    \setlength\tabcolsep{0pt}%
    \put(0,0){\includegraphics[width=\unitlength,page=1]{fig7.pdf}}%
    \put(0.87809062,0.18491791){\makebox(0,0)[lt]{\lineheight{1.25}\smash{\begin{tabular}[t]{l}$\{\lambda(x) \col x \in \cal W\}$\end{tabular}}}}%
    \put(0.87809062,0.06103817){\makebox(0,0)[lt]{\lineheight{1.25}\smash{\begin{tabular}[t]{l}$\spec(H)$\end{tabular}}}}%
    \put(0,0){\includegraphics[width=\unitlength,page=2]{fig7.pdf}}%
    \put(0.41195646,0.22865919){\makebox(0,0)[lt]{\lineheight{1.25}\smash{\begin{tabular}[t]{l}$N^{-\eta}$\end{tabular}}}}%
    \put(0.40621074,0.00334595){\makebox(0,0)[lt]{\lineheight{1.25}\smash{\begin{tabular}[t]{l}$2 N^{-\zeta}$\end{tabular}}}}%
    \put(-0.00132301,0.02559326){\makebox(0,0)[lt]{\lineheight{1.25}\smash{\begin{tabular}[t]{l}$\Lambda(\alpha^*) + \kappa$\end{tabular}}}}%
  \end{picture}%
\endgroup%

%% file: fig5.pdf_tex
\begingroup%
  \makeatletter%
  \providecommand\color[2][]{%
    \errmessage{(Inkscape) Color is used for the text in Inkscape, but the package 'color.sty' is not loaded}%
    \renewcommand\color[2][]{}%
  }%
  \providecommand\transparent[1]{%
    \errmessage{(Inkscape) Transparency is used (non-zero) for the text in Inkscape, but the package 'transparent.sty' is not loaded}%
    \renewcommand\transparent[1]{}%
  }%
  \providecommand\rotatebox[2]{#2}%
  \newcommand*\fsize{\dimexpr\f@size pt\relax}%
  \newcommand*\lineheight[1]{\fontsize{\fsize}{#1\fsize}\selectfont}%
  \ifx\svgwidth\undefined%
    \setlength{\unitlength}{370.47501787bp}%
    \ifx\svgscale\undefined%
      \relax%
    \else%
      \setlength{\unitlength}{\unitlength * \real{\svgscale}}%
    \fi%
  \else%
    \setlength{\unitlength}{\svgwidth}%
  \fi%
  \global\let\svgwidth\undefined%
  \global\let\svgscale\undefined%
  \makeatother%
  \begin{picture}(1,0.70221332)%
    \lineheight{1}%
    \setlength\tabcolsep{0pt}%
    \put(0,0){\includegraphics[width=\unitlength,page=1]{fig5.pdf}}%
    \put(0.33800131,0.3832998){\makebox(0,0)[lt]{\lineheight{1.25}\smash{\begin{tabular}[t]{l}$x$\end{tabular}}}}%
    \put(0.36663526,0.48718368){\makebox(0,0)[lt]{\lineheight{1.25}\smash{\begin{tabular}[t]{l}$y$\end{tabular}}}}%
    \put(0.53857553,0.29415518){\makebox(0,0)[lt]{\lineheight{1.25}\smash{\begin{tabular}[t]{l}$z$\end{tabular}}}}%
    \put(0,0){\includegraphics[width=\unitlength,page=2]{fig5.pdf}}%
    \put(0.01618084,0.03208813){\makebox(0,0)[lt]{\lineheight{1.25}\smash{\begin{tabular}[t]{l}$1$\end{tabular}}}}%
    \put(0.04501821,0.00380636){\makebox(0,0)[lt]{\lineheight{1.25}\smash{\begin{tabular}[t]{l}$0$\end{tabular}}}}%
    \put(0,0){\includegraphics[width=\unitlength,page=3]{fig5.pdf}}%
    \put(0.12153202,0.00380636){\makebox(0,0)[lt]{\lineheight{1.25}\smash{\begin{tabular}[t]{l}$T_1$\end{tabular}}}}%
    \put(0,0){\includegraphics[width=\unitlength,page=4]{fig5.pdf}}%
    \put(0.54235802,0.00380636){\makebox(0,0)[lt]{\lineheight{1.25}\smash{\begin{tabular}[t]{l}$T_2$\end{tabular}}}}%
    \put(0,0){\includegraphics[width=\unitlength,page=5]{fig5.pdf}}%
    \put(0.01618084,0.07034492){\makebox(0,0)[lt]{\lineheight{1.25}\smash{\begin{tabular}[t]{l}$2$\end{tabular}}}}%
    \put(0,0){\includegraphics[width=\unitlength,page=6]{fig5.pdf}}%
    \put(0.01618084,0.22337266){\makebox(0,0)[lt]{\lineheight{1.25}\smash{\begin{tabular}[t]{l}$j$\end{tabular}}}}%
    \put(0,0){\includegraphics[width=\unitlength,page=7]{fig5.pdf}}%
    \put(0.58061493,0.00380636){\makebox(0,0)[lt]{\lineheight{1.25}\smash{\begin{tabular}[t]{l}$T_3$\end{tabular}}}}%
    \put(0,0){\includegraphics[width=\unitlength,page=8]{fig5.pdf}}%
    \put(0.69538565,0.00380636){\makebox(0,0)[lt]{\lineheight{1.25}\smash{\begin{tabular}[t]{l}$T_4$\end{tabular}}}}%
    \put(0,0){\includegraphics[width=\unitlength,page=9]{fig5.pdf}}%
    \put(0.81015637,0.00380636){\makebox(0,0)[lt]{\lineheight{1.25}\smash{\begin{tabular}[t]{l}$T_5$\end{tabular}}}}%
    \put(0,0){\includegraphics[width=\unitlength,page=10]{fig5.pdf}}%
    \put(0.84841321,0.00380636){\makebox(0,0)[lt]{\lineheight{1.25}\smash{\begin{tabular}[t]{l}$T_6 = k$\end{tabular}}}}%
    \put(0.96553191,0.00380636){\makebox(0,0)[lt]{\lineheight{1.25}\smash{\begin{tabular}[t]{l}$t$\end{tabular}}}}%
    \put(-0.00163431,0.2638821){\makebox(0,0)[lt]{\lineheight{1.25}\smash{\begin{tabular}[t]{l}$\gamma(t)$\end{tabular}}}}%
  \end{picture}%
\endgroup%

%% file: fig6.pdf_tex
\begingroup%
  \makeatletter%
  \providecommand\color[2][]{%
    \errmessage{(Inkscape) Color is used for the text in Inkscape, but the package 'color.sty' is not loaded}%
    \renewcommand\color[2][]{}%
  }%
  \providecommand\transparent[1]{%
    \errmessage{(Inkscape) Transparency is used (non-zero) for the text in Inkscape, but the package 'transparent.sty' is not loaded}%
    \renewcommand\transparent[1]{}%
  }%
  \providecommand\rotatebox[2]{#2}%
  \newcommand*\fsize{\dimexpr\f@size pt\relax}%
  \newcommand*\lineheight[1]{\fontsize{\fsize}{#1\fsize}\selectfont}%
  \ifx\svgwidth\undefined%
    \setlength{\unitlength}{330.32069373bp}%
    \ifx\svgscale\undefined%
      \relax%
    \else%
      \setlength{\unitlength}{\unitlength * \real{\svgscale}}%
    \fi%
  \else%
    \setlength{\unitlength}{\svgwidth}%
  \fi%
  \global\let\svgwidth\undefined%
  \global\let\svgscale\undefined%
  \makeatother%
  \begin{picture}(1,0.427745)%
    \lineheight{1}%
    \setlength\tabcolsep{0pt}%
    \put(0,0){\includegraphics[width=\unitlength,page=1]{fig6.pdf}}%
    \put(0.94681122,0.30979348){\color[rgb]{0,0.68627451,0.16470588}\makebox(0,0)[lt]{\lineheight{1.25}\smash{\begin{tabular}[t]{l}$\cal V \setminus \{x\}$\end{tabular}}}}%
    \put(0.94681122,0.213647){\color[rgb]{0.29019608,0.45490196,1}\makebox(0,0)[lt]{\lineheight{1.25}\smash{\begin{tabular}[t]{l}$(\cal V \setminus \{x\})^c$\end{tabular}}}}%
  \end{picture}%
\endgroup%

%% file: fig1.pdf_tex
\begingroup%
  \makeatletter%
  \providecommand\color[2][]{%
    \errmessage{(Inkscape) Color is used for the text in Inkscape, but the package 'color.sty' is not loaded}%
    \renewcommand\color[2][]{}%
  }%
  \providecommand\transparent[1]{%
    \errmessage{(Inkscape) Transparency is used (non-zero) for the text in Inkscape, but the package 'transparent.sty' is not loaded}%
    \renewcommand\transparent[1]{}%
  }%
  \providecommand\rotatebox[2]{#2}%
  \newcommand*\fsize{\dimexpr\f@size pt\relax}%
  \newcommand*\lineheight[1]{\fontsize{\fsize}{#1\fsize}\selectfont}%
  \ifx\svgwidth\undefined%
    \setlength{\unitlength}{369.343907bp}%
    \ifx\svgscale\undefined%
      \relax%
    \else%
      \setlength{\unitlength}{\unitlength * \real{\svgscale}}%
    \fi%
  \else%
    \setlength{\unitlength}{\svgwidth}%
  \fi%
  \global\let\svgwidth\undefined%
  \global\let\svgscale\undefined%
  \makeatother%
  \begin{picture}(1,0.48617975)%
    \lineheight{1}%
    \setlength\tabcolsep{0pt}%
    \put(0,0){\includegraphics[width=\unitlength,page=1]{fig1.pdf}}%
    \put(0.70437515,0.38878984){\makebox(0,0)[lt]{\lineheight{1.25}\smash{\begin{tabular}[t]{l}$B_i(b)^c$\end{tabular}}}}%
    \put(0.27853097,0.33886293){\makebox(0,0)[lt]{\lineheight{1.25}\smash{\begin{tabular}[t]{l}$B_i(b)$\end{tabular}}}}%
    \put(0.28296675,0.25838559){\makebox(0,0)[lt]{\lineheight{1.25}\smash{\begin{tabular}[t]{l}$b$\end{tabular}}}}%
  \end{picture}%
\endgroup%

%% file: fig2.pdf_tex
\begingroup%
  \makeatletter%
  \providecommand\color[2][]{%
    \errmessage{(Inkscape) Color is used for the text in Inkscape, but the package 'color.sty' is not loaded}%
    \renewcommand\color[2][]{}%
  }%
  \providecommand\transparent[1]{%
    \errmessage{(Inkscape) Transparency is used (non-zero) for the text in Inkscape, but the package 'transparent.sty' is not loaded}%
    \renewcommand\transparent[1]{}%
  }%
  \providecommand\rotatebox[2]{#2}%
  \newcommand*\fsize{\dimexpr\f@size pt\relax}%
  \newcommand*\lineheight[1]{\fontsize{\fsize}{#1\fsize}\selectfont}%
  \ifx\svgwidth\undefined%
    \setlength{\unitlength}{369.34389849bp}%
    \ifx\svgscale\undefined%
      \relax%
    \else%
      \setlength{\unitlength}{\unitlength * \real{\svgscale}}%
    \fi%
  \else%
    \setlength{\unitlength}{\svgwidth}%
  \fi%
  \global\let\svgwidth\undefined%
  \global\let\svgscale\undefined%
  \makeatother%
  \begin{picture}(1,0.48617976)%
    \lineheight{1}%
    \setlength\tabcolsep{0pt}%
    \put(0,0){\includegraphics[width=\unitlength,page=1]{fig2.pdf}}%
    \put(0.70437519,0.38878986){\makebox(0,0)[lt]{\lineheight{1.25}\smash{\begin{tabular}[t]{l}$B_i(b)^c$\end{tabular}}}}%
    \put(0.27853099,0.33886294){\makebox(0,0)[lt]{\lineheight{1.25}\smash{\begin{tabular}[t]{l}$B_i(b)$\end{tabular}}}}%
    \put(0.28296675,0.25838559){\makebox(0,0)[lt]{\lineheight{1.25}\smash{\begin{tabular}[t]{l}$b$\end{tabular}}}}%
  \end{picture}%
\endgroup%

%% file: fig3.pdf_tex
\begingroup%
  \makeatletter%
  \providecommand\color[2][]{%
    \errmessage{(Inkscape) Color is used for the text in Inkscape, but the package 'color.sty' is not loaded}%
    \renewcommand\color[2][]{}%
  }%
  \providecommand\transparent[1]{%
    \errmessage{(Inkscape) Transparency is used (non-zero) for the text in Inkscape, but the package 'transparent.sty' is not loaded}%
    \renewcommand\transparent[1]{}%
  }%
  \providecommand\rotatebox[2]{#2}%
  \newcommand*\fsize{\dimexpr\f@size pt\relax}%
  \newcommand*\lineheight[1]{\fontsize{\fsize}{#1\fsize}\selectfont}%
  \ifx\svgwidth\undefined%
    \setlength{\unitlength}{369.343907bp}%
    \ifx\svgscale\undefined%
      \relax%
    \else%
      \setlength{\unitlength}{\unitlength * \real{\svgscale}}%
    \fi%
  \else%
    \setlength{\unitlength}{\svgwidth}%
  \fi%
  \global\let\svgwidth\undefined%
  \global\let\svgscale\undefined%
  \makeatother%
  \begin{picture}(1,0.48617975)%
    \lineheight{1}%
    \setlength\tabcolsep{0pt}%
    \put(0,0){\includegraphics[width=\unitlength,page=1]{fig3.pdf}}%
    \put(0.70437515,0.38878984){\makebox(0,0)[lt]{\lineheight{1.25}\smash{\begin{tabular}[t]{l}$B_i(b)^c$\end{tabular}}}}%
    \put(0.27853097,0.33886293){\makebox(0,0)[lt]{\lineheight{1.25}\smash{\begin{tabular}[t]{l}$B_i(b)$\end{tabular}}}}%
    \put(0.28296675,0.25838559){\makebox(0,0)[lt]{\lineheight{1.25}\smash{\begin{tabular}[t]{l}$b$\end{tabular}}}}%
  \end{picture}%
\endgroup%

%% file: fig4.pdf_tex
\begingroup%
  \makeatletter%
  \providecommand\color[2][]{%
    \errmessage{(Inkscape) Color is used for the text in Inkscape, but the package 'color.sty' is not loaded}%
    \renewcommand\color[2][]{}%
  }%
  \providecommand\transparent[1]{%
    \errmessage{(Inkscape) Transparency is used (non-zero) for the text in Inkscape, but the package 'transparent.sty' is not loaded}%
    \renewcommand\transparent[1]{}%
  }%
  \providecommand\rotatebox[2]{#2}%
  \newcommand*\fsize{\dimexpr\f@size pt\relax}%
  \newcommand*\lineheight[1]{\fontsize{\fsize}{#1\fsize}\selectfont}%
  \ifx\svgwidth\undefined%
    \setlength{\unitlength}{195.74191246bp}%
    \ifx\svgscale\undefined%
      \relax%
    \else%
      \setlength{\unitlength}{\unitlength * \real{\svgscale}}%
    \fi%
  \else%
    \setlength{\unitlength}{\svgwidth}%
  \fi%
  \global\let\svgwidth\undefined%
  \global\let\svgscale\undefined%
  \makeatother%
  \begin{picture}(1,0.67844023)%
    \lineheight{1}%
    \setlength\tabcolsep{0pt}%
    \put(0,0){\includegraphics[width=\unitlength,page=1]{fig4.pdf}}%
    \put(0.47489782,0.38144304){\makebox(0,0)[lt]{\lineheight{1.25}\smash{\begin{tabular}[t]{l}$b$\end{tabular}}}}%
    \put(0.79679138,0.54712351){\makebox(0,0)[lt]{\lineheight{1.25}\smash{\begin{tabular}[t]{l}$B_r(b)$\end{tabular}}}}%
  \end{picture}%
\endgroup%

%% file: fig9.pdf_tex
\begingroup%
  \makeatletter%
  \providecommand\color[2][]{%
    \errmessage{(Inkscape) Color is used for the text in Inkscape, but the package 'color.sty' is not loaded}%
    \renewcommand\color[2][]{}%
  }%
  \providecommand\transparent[1]{%
    \errmessage{(Inkscape) Transparency is used (non-zero) for the text in Inkscape, but the package 'transparent.sty' is not loaded}%
    \renewcommand\transparent[1]{}%
  }%
  \providecommand\rotatebox[2]{#2}%
  \newcommand*\fsize{\dimexpr\f@size pt\relax}%
  \newcommand*\lineheight[1]{\fontsize{\fsize}{#1\fsize}\selectfont}%
  \ifx\svgwidth\undefined%
    \setlength{\unitlength}{357.4688241bp}%
    \ifx\svgscale\undefined%
      \relax%
    \else%
      \setlength{\unitlength}{\unitlength * \real{\svgscale}}%
    \fi%
  \else%
    \setlength{\unitlength}{\svgwidth}%
  \fi%
  \global\let\svgwidth\undefined%
  \global\let\svgscale\undefined%
  \makeatother%
  \begin{picture}(1,0.40683014)%
    \lineheight{1}%
    \setlength\tabcolsep{0pt}%
    \put(0,0){\includegraphics[width=\unitlength,page=1]{fig9.pdf}}%
    \put(0.03845128,0.03864001){\makebox(0,0)[lt]{\lineheight{1.25}\smash{\begin{tabular}[t]{l}$2$\end{tabular}}}}%
    \put(0.0199161,0.23771621){\makebox(0,0)[lt]{\lineheight{1.25}\smash{\begin{tabular}[t]{l}$2.1$\end{tabular}}}}%
    \put(0.07513586,0.00461038){\makebox(0,0)[lt]{\lineheight{1.25}\smash{\begin{tabular}[t]{l}$0$\end{tabular}}}}%
    \put(0.20175659,0.00461038){\makebox(0,0)[lt]{\lineheight{1.25}\smash{\begin{tabular}[t]{l}$0.05$\end{tabular}}}}%
    \put(0.3576056,0.00461038){\makebox(0,0)[lt]{\lineheight{1.25}\smash{\begin{tabular}[t]{l}$0.1$\end{tabular}}}}%
    \put(0.4958153,0.00461038){\makebox(0,0)[lt]{\lineheight{1.25}\smash{\begin{tabular}[t]{l}$0.15$\end{tabular}}}}%
    \put(0.65086664,0.00461038){\makebox(0,0)[lt]{\lineheight{1.25}\smash{\begin{tabular}[t]{l}$0.2$\end{tabular}}}}%
    \put(0.79131895,0.00461038){\makebox(0,0)[lt]{\lineheight{1.25}\smash{\begin{tabular}[t]{l}$0.25$\end{tabular}}}}%
    \put(0,0){\includegraphics[width=\unitlength,page=2]{fig9.pdf}}%
    \put(0.88086359,0.00461042){\makebox(0,0)[lt]{\lineheight{1.25}\smash{\begin{tabular}[t]{l}$\mu$\end{tabular}}}}%
    \put(-0.00197952,0.38670287){\makebox(0,0)[lt]{\lineheight{1.25}\smash{\begin{tabular}[t]{l}$\alpha_*(\mu)$\end{tabular}}}}%
    \put(0.72829458,0.07303461){\makebox(0,0)[lt]{\lineheight{1.25}\smash{\begin{tabular}[t]{l}$1 - \frac{b}{b_*}$\end{tabular}}}}%
    \put(0,0){\includegraphics[width=\unitlength,page=3]{fig9.pdf}}%
    \put(0.09647726,0.36266939){\makebox(0,0)[lt]{\lineheight{1.25}\smash{\begin{tabular}[t]{l}$\alpha_{\max}(b)$\end{tabular}}}}%
  \end{picture}%
\endgroup%

%% file: fig10.pdf_tex
\begingroup%
  \makeatletter%
  \providecommand\color[2][]{%
    \errmessage{(Inkscape) Color is used for the text in Inkscape, but the package 'color.sty' is not loaded}%
    \renewcommand\color[2][]{}%
  }%
  \providecommand\transparent[1]{%
    \errmessage{(Inkscape) Transparency is used (non-zero) for the text in Inkscape, but the package 'transparent.sty' is not loaded}%
    \renewcommand\transparent[1]{}%
  }%
  \providecommand\rotatebox[2]{#2}%
  \newcommand*\fsize{\dimexpr\f@size pt\relax}%
  \newcommand*\lineheight[1]{\fontsize{\fsize}{#1\fsize}\selectfont}%
  \ifx\svgwidth\undefined%
    \setlength{\unitlength}{337.10349481bp}%
    \ifx\svgscale\undefined%
      \relax%
    \else%
      \setlength{\unitlength}{\unitlength * \real{\svgscale}}%
    \fi%
  \else%
    \setlength{\unitlength}{\svgwidth}%
  \fi%
  \global\let\svgwidth\undefined%
  \global\let\svgscale\undefined%
  \makeatother%
  \begin{picture}(1,0.36836013)%
    \lineheight{1}%
    \setlength\tabcolsep{0pt}%
    \put(0,0){\includegraphics[width=\unitlength,page=1]{fig10.pdf}}%
    \put(0.9246005,0.00171011){\color[rgb]{0,0,0}\makebox(0,0)[lt]{\lineheight{1.25}\smash{\begin{tabular}[t]{l}$\lambda$\end{tabular}}}}%
    \put(0.03591931,0.29543718){\color[rgb]{0,0,0}\makebox(0,0)[lt]{\lineheight{1.25}\smash{\begin{tabular}[t]{l}1\end{tabular}}}}%
    \put(0.03611879,0.03714559){\color[rgb]{0,0,0}\makebox(0,0)[lt]{\lineheight{1.25}\smash{\begin{tabular}[t]{l}0\end{tabular}}}}%
    \put(0.3779508,0.00107716){\color[rgb]{0,0,0}\makebox(0,0)[lt]{\lineheight{1.25}\smash{\begin{tabular}[t]{l}2\end{tabular}}}}%
    \put(0.69148626,0.00171011){\color[rgb]{0,0,0}\makebox(0,0)[lt]{\lineheight{1.25}\smash{\begin{tabular}[t]{l}$\lambda_{\max}(b)$\end{tabular}}}}%
    \put(-0.0017961,0.35009789){\makebox(0,0)[lt]{\lineheight{1.25}\smash{\begin{tabular}[t]{l}$\rho_b(\lambda)$\end{tabular}}}}%
  \end{picture}%
\endgroup%